# Paraconsistent First-Order Logic $\mathbf{LP}^{\#}_{\omega}$ with restricted modus ponens rule and infinite hierarchy levels of contradiction. Axiomatical system $\mathbf{HST}^{\#}_{\omega}$, as paraconsistent generalization of Hrbacek set theory $\mathbf{HST}$.

Jaykov Foukzon
Israel Institute of Technology
jaykovfoukzon@gmail.com

**Abstract**: In this paper paraconsistent first-order logic $\mathbf{LP}^{\#}_{\omega}$ with restricted modus ponens rule and infinite hierarchy levels of contradiction is proposed. Corresponding paraconsistent set theory $\mathbf{KSth}^{\#}_{\omega}$ is proposed. Axiomatical system $\mathbf{HST}^{\#}_{\omega}$, as inconsistent generalization of Hrbacek set theory $\mathbf{HST}$ is considered.

## Content

**PART I**.

**PART II**.

**PART III**.

**List of designations**

$\mathbf{V}^{Con}$ - consistent universum
$\mathbf{V}^{Inc}$ - inconsistent universum
$\mathbf{U}$ - complete universum $\mathbf{U} \triangleq \mathbf{V}^{Con} \bigcup \mathbf{V}^{Inc}$
$(\bullet = \bullet)$ - relation of the classical consistent equivalence
$(\bullet =_{\mathbf{s}} \bullet)$ - relation of the strong consistent nonclassical equivalence or $\mathbf{s}$-equivalence
$(\bullet =_{\mathbf{s}} \bullet) \upharpoonright \mathbf{V}^{Con} = (\bullet = \bullet)$
$(\bullet \in \bullet)$ - classical consistent membership relation
$(\bullet \in_{\mathbf{s}} \bullet)$ - strong consistent membership relation or $\mathbf{s}$-membership relation

$(\bullet \in_{\mathbf{s}} \ \bullet) \upharpoonright \mathbf{V}^{Con} = (\bullet \in \bullet)$
$(\bullet =_{w} \ \bullet)$ - relation of the weak equivalence or $w$-equivalence
$(\bullet =_{w} \ \bullet) \upharpoonright \mathbf{V}^{Con} = (\bullet = \bullet)$
$(\bullet =_{w_1} \ \bullet)$ - relation of the weak inconsistent equivalence order $1$ or $w_1$-equivalence
$(\bullet =_{w_1} \ \bullet) \leftrightarrow (\bullet =_{w} \ \bullet) \wedge (\bullet \neq_{w} \ \bullet)$
$(\bullet =_{w_n} \ \bullet)$ - relation of the weak inconsistent equivalence order $n$ or $w_n$-equivalence
$(\bullet \in_{w} \ \bullet)$ - weak membership relation or $w$-membership relation
$(\bullet \in_{w_1} \ \bullet)$ weak inconsistent membership relation order 1 or $w_1$-membership relation
$\varnothing_w$-weak empty set

.

# I.Introduction.

The real history of non-Aristotelian logic begins on May 18,1910 when N.A. Vasiliev presented to the Kazan University faculty a lecture "On Partial Judgements, the Triangle of Opposition, the Law of Excluded Fourth" [Vasiliev 1910] to satisfy the requirements for obtaining the title of privat-dozent. In this lecture Vasiliev expounded for the first time the key principles of non-Aristotelian, imaginary, logic. In this work he likewise constructed his "imaginary" logic free of the laws of contradiction and excluded middle in the informal, so-to-speak Aristotelian, manner (although imaginary logic is in essense non- Aristotelian).Thus the birthday of new logic was exactly fixed in the annals of history. Vasiliev's reform of logic was radical, and he did his best to determine whether it was possible for the new logic with new laws and new subject to imply a new logical Universe. Vasiliev began the modern non-classical revolution in logic, but he certainly did not complete it. The founder of paraconsistent logic, N.A. Vasiliev, stated as a characteristic feature of his logic, three kinds of sentence, i.e. "**S** *is* **A**", "**S** *is not* **A**", "**S** *is and is not* **A**". Thus Vasiliev logic rejected the *law of non-contradiction:* $\neg(\mathbf{A} \wedge \neg\mathbf{A})$ and the *law of excluded middle:* $\mathbf{A} \vee \neg\mathbf{A}$. However Vasiliev's logic preserve the *law of excluded fourth:* $\mathbf{A} \vee \neg\mathbf{A} \vee (\mathbf{A} \wedge \neg\mathbf{A})$. Possible formalized versions of Vasiliev's logic with one level of contradiction $\mathbf{LP}_1^{\#}$ was proposed by A.I.Arruda [1]. In this paper we proposed paraconsistent first-order logic $\mathbf{LP}_{\boldsymbol{\omega}}^{\#}$ with infinite hierarchy levels of contradiction. Corresponding paraconsistent set theory $\mathbf{KSth}_{\omega}^{\#}$ is discussed.

The postulates (or their axioms schemata) of Vasiliev-Amida propositional paraconsistent logic $VA_1$ are the following:

The language $\mathcal{L}_1$ of paraconsistent logic $VA_1 \triangleq VA_1[\mathbf{V}]$ has as primitive symbols (i) countable set of a clalassical propositional variables, (ii) countable set $\mathbf{V} = \{\mathbf{P}_i\}_{i \in \mathbb{N}}$ of a non clalassical propositional variables,(iii) the connectives $\neg, , \wedge, \vee, \rightarrow$ and (iv) the parentheses (,).

**I**. **Logical postulates**:

(1) $\mathbf{A} \rightarrow (\mathbf{B} \rightarrow \mathbf{A})$,

(2) $(\mathbf{A} \to \mathbf{B}) \to ((\mathbf{A} \to (\mathbf{B} \to \mathbf{C})) \to (\mathbf{A} \to \mathbf{C}))$,

(3) $\mathbf{A} \to (\mathbf{B} \to \mathbf{A} \wedge \mathbf{B})$,

(4) $\mathbf{A} \wedge \mathbf{B} \to \mathbf{A}$,

(5) $\mathbf{A} \wedge \mathbf{B} \to \mathbf{B}$,

(6) $\mathbf{A} \to (\mathbf{A} \vee \mathbf{B})$,

(7) $\mathbf{B} \to (\mathbf{A} \vee \mathbf{B})$,

(6) $\mathbf{A} \to (\mathbf{A} \vee \mathbf{B})$,

(7) $\mathbf{B} \to (\mathbf{A} \vee \mathbf{B})$,

(8) $(\mathbf{A} \to \mathbf{C}) \to ((\mathbf{B} \to \mathbf{C}) \to (\mathbf{A} \vee \mathbf{B} \to \mathbf{C}))$,

(9) $\mathbf{A} \vee \neg\mathbf{A}$,

(10) $\mathbf{B} \to (\neg\mathbf{B} \to \mathbf{A})$ **if** $\mathbf{B} \notin \mathbf{V}$.

**II.Rules of a conclusion**:

**Anrestricted Modus Ponens rule** $\mathbf{MP}$ : $\mathbf{A}, \mathbf{A} \to \mathbf{B} \vdash \mathbf{B}$.

**Theorem 1.1**.[1]. (1) If $\mathbf{B} \notin \mathbf{V}$, then $\mathbf{B}, \neg\mathbf{B} \vdash \mathbf{A}$; (2) $\neg\neg\mathbf{A} \leftrightarrow \mathbf{A}$ **iff** $\mathbf{A} \notin \mathbf{V}$;
(3) $\neg\neg\mathbf{A} \to \mathbf{A}$.

The postulates (or their axioms schemata) of Vasiliev-Amida propositional paraconsistent logic $VA_2$ are the following:

The language $\mathcal{L}_2$ of paraconsistent logic $VA_2 \triangleq VA_2[\mathbf{V}]$ has as primitive symbols (i) countable set of a clalassical propositional variables, (ii) countable set $\mathbf{V} = \{\mathbf{P}_i\}_{i \in \mathbb{N}}$ of a non clalassical propositional variables,(iii) the connectives $\neg, , \wedge, \vee, \to$ and (iv) the parentheses (,).

**I. Logical postulates**:

(1) $\mathbf{A} \to (\mathbf{B} \to \mathbf{A})$,

(2) $(\mathbf{A} \to \mathbf{B}) \to ((\mathbf{A} \to (\mathbf{B} \to \mathbf{C})) \to (\mathbf{A} \to \mathbf{C}))$,

(3) $\mathbf{A} \to (\mathbf{B} \to \mathbf{A} \wedge \mathbf{B})$,

(4) $\mathbf{A} \wedge \mathbf{B} \to \mathbf{A}$,

(5) $\mathbf{A} \wedge \mathbf{B} \to \mathbf{B}$,

(6) $\mathbf{A} \to (\mathbf{A} \vee \mathbf{B})$,

(7) $\mathbf{B} \to (\mathbf{A} \vee \mathbf{B})$,

(6) $\mathbf{A} \to (\mathbf{A} \vee \mathbf{B})$,

(7) $\mathbf{B} \to (\mathbf{A} \vee \mathbf{B})$,

(8) $(\mathbf{A} \to \mathbf{C}) \to ((\mathbf{B} \to \mathbf{C}) \to (\mathbf{A} \vee \mathbf{B} \to \mathbf{C}))$,

(9) $\mathbf{A} \vee \neg\mathbf{A}$,

(10) $\mathbf{B} \to (\neg\mathbf{B} \to \mathbf{A})$ **if** $\mathbf{B} \notin \mathbf{V}$,

(11) $\mathbf{P}_i \wedge \neg\mathbf{P}_i$ **iff** $\mathbf{P}_i \in \mathbf{V}, i = 1, 2, \ldots$ .

**II.Rules of a conclusion**:

**Anrestricted Modus Ponens rule** $\mathbf{MP}$ : $\mathbf{A}, \mathbf{A} \to \mathbf{B} \vdash \mathbf{B}$.

# II.Paraconsistent Logic with $n$ levels of contradiction $\mathbf{LP}_n^\#$.

Modern formalized versions of Vasiliev's logic with one level of contradiction may be found in Amida [1980], [Puga and Da Costa 1988], Smimov [Smirnov 1987], and [Smimov 1987a, 161-169]. There is also the presentation Smimov given at the International Congress of Logic, Methodology and Philosophy of Science in Uppsala in 1991.

.

# Paraconsistent Logic with one levels of a contradiction $\mathbf{LP}_1^\#$.

Let us consider now Vasiliev-Arruda type paraconsistent logic $\mathbf{LP}_1^\# = \mathbf{LP}_1^\#[\mathbf{V},\Delta]$ with one level of contradiction.

The postulates (or their axioms schemata) of propositional paraconsistent logic $\mathbf{LP}_1^\#$ are the following:

The language $\mathcal{L}_1^\#$ of paraconsistent logic $\mathbf{LP}_1^\# \triangleq \mathbf{LP}_1^\#[\mathbf{V},\Delta]$ has as primitive symbols (i) countable set of a clalassical propositional variables, (ii) countable set $\mathbf{V} = \{\mathbf{P}_i\}_{i\in\mathbb{N}}$ of a non clalassical propositional variables, (iii) the connectives $\neg_w, \neg_\mathbf{s}, \wedge, \vee, \rightarrow$ and (iv) the parentheses (,).

**Remark**.**2**.**1**.We distinguish a weak negation $\neg_w$ and a strong negation $\neg_\mathbf{s}$.

The definition of formula is the usual. We denote the set of the all formulae of $\mathbf{LP}_1^\#[\mathbf{V}_1,\Delta]$ by $\mathcal{F}_1^\#$, where $\mathbf{V}_1 = \mathbf{V}^{[0]} \cup \mathbf{V}^{[1]}$ and $\Delta$ is a given subset of $\mathcal{F}_1^\#$. Here we used the following definitions: $\mathbf{V}^{[0]} \triangleq \mathbf{V}, \mathbf{V}^{[1]} \triangleq \{\alpha^{[1]} | (\alpha \in \mathbf{V})\}, \alpha^{[1]} \triangleq (\alpha \wedge \neg_w \alpha)$. $\mathbf{A},\mathbf{B},\mathbf{C}$, ... will be used as metalanguage variables which indicate formulas of $\mathbf{LP}_1^\#[\mathbf{V},\Delta]$. We assume through that $\mathbf{V}_1 \subset \Delta \subsetneq \mathcal{F}_1^\#$.

**I**. **Logical postulates**:

(1) $\mathbf{A} \rightarrow (\mathbf{B} \rightarrow \mathbf{A})$,

(2) $(\mathbf{A} \rightarrow \mathbf{B}) \rightarrow ((\mathbf{A} \rightarrow (\mathbf{B} \rightarrow \mathbf{C})) \rightarrow (\mathbf{A} \rightarrow \mathbf{C}))$,

(3) $\mathbf{A} \rightarrow (\mathbf{B} \rightarrow \mathbf{A} \wedge \mathbf{B})$,

(4) $\mathbf{A} \wedge \mathbf{B} \rightarrow \mathbf{A}$,

(5) $\mathbf{A} \wedge \mathbf{B} \rightarrow \mathbf{B}$,

(6) $\mathbf{A} \rightarrow (\mathbf{A} \vee \mathbf{B})$,

(7) $\mathbf{B} \rightarrow (\mathbf{A} \vee \mathbf{B})$,

(8) $(\mathbf{A} \rightarrow \mathbf{C}) \rightarrow ((\mathbf{B} \rightarrow \mathbf{C}) \rightarrow (\mathbf{A} \vee \mathbf{B} \rightarrow \mathbf{C}))$,

(9) $\mathbf{P}_i \wedge \neg_w \mathbf{P}_i$ **iff** $\mathbf{P}_i \in \mathbf{V}, i = 1,2,\ldots$

(10) $\mathbf{A} \vee \neg_w \mathbf{A}$ **iff** $\mathbf{A} \notin \mathbf{V}$,

(11) $\mathbf{B} \rightarrow (\neg_w \mathbf{B} \rightarrow \mathbf{A})$ **if** $\mathbf{B} \notin \mathbf{V}_1$,

(12) $\mathbf{A} \vee \neg_w \mathbf{A} \vee (\mathbf{A} \wedge \neg_w \mathbf{A})$ **iff** $\mathbf{A} \in \mathbf{V}_1$,

(13) $\mathbf{A} \vee \neg_\mathbf{s} \mathbf{A}$ **if** $\mathbf{A} \in \mathcal{F}_1^\#$,

(14) $\mathbf{B} \rightarrow (\neg_\mathbf{s} \mathbf{B} \rightarrow \mathbf{A})$ **if** $\mathbf{A},\mathbf{B} \in \mathcal{F}_1^\#$.

**II.Rules of a conclusion**:

**Restricted Modus Ponens rule** $\mathbf{MPR}$ :

$\mathbf{A}, \mathbf{A} \to \mathbf{B} \vdash \mathbf{B}$ **iff** $\mathbf{A} \notin \Delta$.

**Unrestricted Modus Tollens rules**: $\mathbf{P} \to \mathbf{Q}, \neg_w \mathbf{Q} \vdash \neg_w \mathbf{P}; \mathbf{P} \to \mathbf{Q}, \neg_s \mathbf{Q} \vdash \neg_s \mathbf{P}$.

**The rule of a strong contradiction**: $\mathbf{A} \wedge \neg_s \mathbf{A} \vdash \mathbf{B}$.

**III.Quantification**

Corresponding to the propositional paraconsistent relevant logic $\mathbf{LP}_1^{\#}[\mathbf{V}, \Delta]$ we construct the corresponding paraconsistent relevant first-order predicate calculus $\overline{\mathbf{LP}}_1^{\#} = \overline{\mathbf{LP}}_1^{\#}[\tilde{\mathbf{V}}, \tilde{\Delta}]$. The language of the paraconsistent predicate calculus $\overline{\mathbf{LP}}_1^{\#}$, denoted by $\overline{\mathcal{L}}_1^{\#}$, is an extension of the language $\mathcal{L}_1^{\#}$ introduced above, by adding, as usually,for every $m$, denumerable families of $m$-ary predicate symbols $R_1^m, R_2^m, \ldots, R_n^m, \ldots,$ and $m$-ary function symbols $f_1^m, f_2^m, \ldots, f_n^m, \ldots,$ and the universal $\forall$ and existential $\exists$ quantifiers.

We assume throughout that: the language $\overline{\mathcal{L}}_1^{\#}$ contains also

(i) the classical numerals $\bar{0}, \bar{1}, \ldots$;

(ii) countable set $\Gamma$ of the classical consistent set variables $\Gamma = \{x, y, z, \ldots\}$;

(iii) countable set $\breve{\Gamma}$ of the non classical inconsistent set variables $\breve{\Gamma} = \{\breve{x}, \breve{y}, \breve{z}, \ldots\}$;

(iv) countable set $\Theta$ of the classical non-logical constants $\Theta = \{a, b, c\ldots\}$;

(iv) countable set $\breve{\Theta}$ of the non classical non-logical constants $\breve{\Theta} = \{\breve{a}, \breve{b}, \breve{c}\ldots\}$;

**Definition 2.1**. An $\overline{\mathbf{LP}}_1^{\#}$ wff $\Phi$ (well-formed formula $\Phi$) is a $\overline{\mathbf{LP}}_1^{\#}$- sentence iff it hasn't free variables; a wff $\Psi$ is open if it has free variables. We'll use the slang '$k$ – place open wff' to mean a wff with $k$ distinct free variables.

**Definition 2.2**. An $\overline{\mathbf{LP}}_1^{\#}$ wff $\Phi$ is a classical iff it hasn't non classical variables and non classical constants.

**Definition 2.3**. An $\overline{\mathbf{LP}}_1^{\#}$ wff $\Phi$ is a non classical iff it has an non classical variables or non classical constants.We denote the set of the all formulae of $\overline{\mathbf{LP}}_1^{\#}[\tilde{\mathbf{V}}, \tilde{\Delta}]$ by $\overline{\mathcal{F}}_1^{\#}$, where $\tilde{\mathbf{V}} \supset \mathbf{V}_1$ and $\tilde{\Delta} \supset \Delta$ is a given subsets of $\overline{\mathcal{F}}_1^{\#}$.We assume through that $\tilde{\mathbf{V}} \subset \tilde{\Delta} \subsetneq \overline{\mathcal{F}}_1^{\#}$.

The postulates of $\overline{\mathbf{LP}}_1^{\#}[\tilde{\mathbf{V}}, \tilde{\Delta}]$ are those of $\mathbf{LP}_1^{\#}[\tilde{\mathbf{V}}, \tilde{\Delta}]$ (suitably adapted) plus the following:

**(I)** $\dfrac{\alpha \to \beta(x)}{\alpha \to \forall x \beta(x)},$

**(II)** $\forall x \alpha(x) \to \alpha(y),$

**(III)** $\alpha(x) \to \exists x \alpha(x),$

**(IV)** $\dfrac{\alpha(x) \to \beta}{\exists x \alpha(x) \to \beta},$

**(V)** $\forall x(\alpha(x))^{(1)} \to (\forall x \alpha(x))^{(1)},$

**(VI)** $\forall x((\alpha(x))^{(1)} \to (\exists x \alpha(x))^{(1)},$

**(VII)** $\forall x(\alpha(x))^{[1]} \to (\forall x\alpha(x)) \wedge (\forall x\neg_w\alpha(x))$,

**(VIII)** $\forall x((\alpha(x))^{[1]} \to (\exists x\alpha(x)) \wedge (\exists x\neg_w\alpha(x))$,

where we used the following definitions:

$\alpha^{(0)} \triangleq \alpha, \alpha^{(1)} \triangleq \neg_w(\alpha \wedge \neg_w\alpha)$ and

$\alpha^{[0]} \triangleq \alpha, \alpha^{[1]} \triangleq \alpha \wedge \neg_w\alpha$

and where the variables $x$ and $y$ and the formulas $\alpha$ and $\beta$ satisfy the usual definition.

From the calculi $\overline{\mathbf{LP}}_{\mathbf{1}}^{\#}[\tilde{\mathbf{V}}, \tilde{\boldsymbol{\Delta}}]$,one can construct the following predicate calculus with equality. This is done by adding to their languages the binary predicates symbol of strong equality $(\bullet = \bullet)$ or $(\bullet =_{\mathbf{s}} \bullet)$ and weak equality$(\bullet =_w \bullet)$ with suitable modifications in the concept of formula, and by adding the following postulates:

**(IX)** $\forall x(x =_{\mathbf{s}} x)$,

**(X)** $\forall x\forall y\left[(x =_{\mathbf{s}} y)^{[1]} \vdash \mathbf{B}\right]$,

**(XI)** $\forall x\forall y[x =_{\mathbf{s}} y \to (\alpha(x) \leftrightarrow \alpha(y))]$,

**(XII)** $\forall x\forall y\forall z[(x =_{\mathbf{s}} y) \wedge (y =_{\mathbf{s}} z) \to x =_{\mathbf{s}} z]$,

**(XIII)** $\forall y\exists x(y =_w x)$,

**(XIV)** $\forall y\exists x(y =_w x)^{[1]}$,

**(XV)** $\forall x\forall y[x =_w y \to (\alpha(x) \leftrightarrow \alpha(y))]$

**(XVI)** $\forall x\forall y\left[(x =_w y)^{[1]} \leftrightarrow (\alpha^{[1]}(x) \leftrightarrow \alpha^{[1]}(y))\right]$,

**(XVII)** $\forall x\forall y\forall z[(x =_w y) \wedge (y =_w z) \to x =_w z]$,

**(XVIII)** $\forall x\forall y\forall z\left[(x =_w y)^{[1]} \wedge (y =_w z)^{[1]} \to (x =_w z)^{[1]}\right]$,

**(XIX)** $\forall x\forall y\forall z[(x =_w y) \wedge (y =_{\mathbf{s}} z) \to x =_w z]$,

**(XX)** $\forall x\forall y\forall z\left[(x =_w y)^{[1]} \wedge (y =_{\mathbf{s}} z) \to (x =_w z)^{[1]}\right]$,

**(XXI)** $\forall x\forall y\forall z[(x =_{\mathbf{s}} y) \wedge (y =_w z) \to x =_w z]$,

**(XXII)** $\forall x\forall y\forall z\left[(x =_{\mathbf{s}} y) \wedge (y =_w z)^{[1]} \to (x =_w z)^{[1]}\right]$.

**II.Rules of a conclusion**:

**Restricted Modus Ponens rule MPR** :

$\mathbf{A}, \mathbf{A} \to \mathbf{B} \vdash \mathbf{B}$ **iff** $\mathbf{A} \notin \tilde{\boldsymbol{\Delta}}$.

**Unrestricted Modus Tollens rules**: $\mathbf{P} \to \mathbf{Q}, \neg_w\mathbf{Q} \vdash \neg_w\mathbf{P}; \mathbf{P} \to \mathbf{Q}, \neg_{\mathbf{s}}\mathbf{Q} \vdash \neg_{\mathbf{s}}\mathbf{P}$.

**The rule of a strong contradiction**: $\mathbf{A} \wedge \neg_{\mathbf{s}}\mathbf{A} \vdash \mathbf{B}$.

**Definition 2.4**. Classical $\mathbf{V}$-object $\Im^{\mathbf{Cl}} = \Im^{\mathbf{Cl}}[\tilde{\mathbf{V}}, \tilde{\boldsymbol{\Delta}}]$ is the object such that from any classical formula of the form $P(\Im^{\mathbf{Cl}}) \wedge \neg_w P(\Im^{\mathbf{Cl}})$, where $P(\Im^{\mathbf{NCl}}) \notin \tilde{\boldsymbol{\Delta}}$ by using principles as in paraconsistent logical calculas $\overline{\mathbf{LP}}_1^{\#}[\tilde{\mathbf{V}}, \tilde{\boldsymbol{\Delta}}]$ using Restricted Modus Ponens rule, one can deduce any formula i.e., classical object $\Im^{\mathbf{Cl}}$ is the object which hasn't any inconsistent property with respect to a weak negation $\neg_w$.

**Definition 2.5**. Non classical $\mathbf{V}$-object $\Im^{\mathbf{NCl}} = \Im^{\mathbf{NCl}}[\tilde{\mathbf{V}}, \tilde{\boldsymbol{\Delta}}]$ of the 1-degree of inconsistency is the object $\Im^{\mathbf{NCl}}$ such that: from any non classical formula of the form $P(\Im^{\mathbf{NCl}}) \wedge \neg_w P(\Im^{\mathbf{NCl}})$, where $P(\Im^{\mathbf{NCl}}) \notin \tilde{\boldsymbol{\Delta}}$ by using principles as in paraconsistent logical calculas $\overline{\mathbf{LP}}_1^{\#}[\tilde{\mathbf{V}}, \tilde{\boldsymbol{\Delta}}]$ using Restricted Modus Ponens rule one can't deduce any formula whatsoever i.e., non classical object of the $1$-degree of inconsistency is the object $\mathfrak{I}^{\mathbf{NCl}}$ which has at least one inconsistent property of the $1$-degree with respect to a weak negation $\neg_w$.

The simplest example of non classical objects $1$-degree inconsistency is inconsistent numbers $\breve{a}$ such that

$$(\breve{a} =_w \bar{1}) \wedge \neg_w(\breve{a} =_w \bar{1}), \quad (2.1)$$

or

$$\left(\breve{b} =_w \bar{1}\right) \wedge \left(\breve{b} =_w \bar{2}\right). \quad (2.2)$$

**Remark.2.2**. Note that: (i) formula (2.1) meant that $(\breve{a} =_w \bar{1}) \in \tilde{\mathbf{V}}$ and (ii) formula (2.2) meant that $\left(\breve{b} =_w \bar{1}\right) \in \tilde{\boldsymbol{\Delta}}$ and $\left(\breve{b} =_w \bar{2}\right) \in \tilde{\boldsymbol{\Delta}}$.

# Paraconsistent Logic with $n$ levels of contradiction $\mathbf{LP}_n^{\#}$.

Let us consider now paraconsistent logic $\mathbf{LP}_n^{\#} = \mathbf{LP}_n^{\#}[\mathbf{V}, \boldsymbol{\Delta}]$ with $n$ levels of contradiction.

The postulates (or their axioms schemata) of propositional paraconsistent logic $\mathbf{LP}_n^{\#} = \mathbf{LP}_n^{\#}\left[\hat{\mathbf{V}}, \hat{\boldsymbol{\Delta}}\right]$ are the following:

The language $\mathcal{L}_n^{\#}$ of paraconsistent logic $\mathbf{LP}_n^{\#}$ has as primitive symbols (i) countable set of a clalassical propositional variables, (ii) countable set $\mathbf{V} = \{\mathbf{P}_i\}_{i\in\mathbb{N}}$ of a non clalassical propositional variables, (iii) the connectives $\neg_w, \neg_\mathbf{s}, \wedge, \vee, \rightarrow$ and (iv) the parentheses (,).

**Remark 2.3**.We distinguish a weak negation $\neg_w$ and a strong negation $\neg_\mathbf{s}$.

The definition of formula is the usual. We denote the set of the all formulae of $\mathbf{LP}_n^{\#}\left[\hat{\mathbf{V}}, \hat{\boldsymbol{\Delta}}\right]$ by $\mathcal{F}_n^{\#}$ where $\hat{\mathbf{V}}$ and $\hat{\boldsymbol{\Delta}}$ is a given subsets of $\mathcal{F}_n^{\#}$. We assume through that $\hat{\mathbf{V}} \subset \hat{\boldsymbol{\Delta}} \subsetneq \mathcal{F}_n^{\#}$.

$\mathbf{A}, \mathbf{B}, \mathbf{C}$, ... will be used as metalanguage variables which indicate formulas of $\mathbf{LP}_n^{\#}\left[\hat{\mathbf{V}}, \hat{\boldsymbol{\Delta}}\right]$.

**Definition 2.6**. (i) $\alpha^{(k)}$ stands for $\alpha^{(k-1)} \wedge (\alpha^{(k-1)})^{(1)}$, where $\alpha^{(0)} \triangleq \alpha$, $\alpha^{(1)} \triangleq \neg_w(\alpha \wedge \neg_w\alpha), 0 \leq k \leq n$.

(ii) the (finite) $k$-order of the level of a weak consistency ($w$-consistency) is: $\alpha^{(k)}, 0 \leq k \leq n$.

**Definition 2.7**. (i) $\alpha^{[k]}$ stands for $\alpha^{[k-1]} \wedge (\alpha^{[k-1]})^{[1]}$, where $\alpha^{[0]} \triangleq \alpha$, $\alpha^{[1]} \triangleq \alpha \wedge \neg_w\alpha, 0 \leq k \leq n$.

(ii) the (finite) $k$-order of the level of a weak inconsistency ($w$-inconsistency) is: $\alpha^{[n]}, 1 \leq k \leq n$.

**I. Logical postulates**:

(1) $\mathbf{A} \to (\mathbf{B} \to \mathbf{A})$,

(2) $(\mathbf{A} \to \mathbf{B}) \to ((\mathbf{A} \to (\mathbf{B} \to \mathbf{C})) \to (\mathbf{A} \to \mathbf{C}))$,
(3) $\mathbf{A} \to (\mathbf{B} \to \mathbf{A} \wedge \mathbf{B})$,
(4) $\mathbf{A} \wedge \mathbf{B} \to \mathbf{A}$,

(5) $\mathbf{A} \wedge \mathbf{B} \to \mathbf{B}$,
(6) $\mathbf{A} \to (\mathbf{A} \vee \mathbf{B})$,
(7) $\mathbf{B} \to (\mathbf{A} \vee \mathbf{B})$,
(8) $(\mathbf{A} \to \mathbf{C}) \to ((\mathbf{B} \to \mathbf{C}) \to (\mathbf{A} \vee \mathbf{B} \to \mathbf{C}))$,

(9) $\mathbf{P} \wedge \neg_w \mathbf{P}$ **iff** $\mathbf{P} \in \mathbf{V}$,
(10) $\mathbf{P}^{[k]}$ **iff** $\mathbf{P} \in \mathbf{V}$,
(11) $\mathbf{A} \vee \neg_w \mathbf{A}$ **if** $\mathbf{A} \notin \hat{\mathbf{V}} = \bigcup_{\mathbf{k}=\mathbf{0}}^{\mathbf{n}} \mathbf{V}^{[\mathbf{k}]}.$,
(12) $\mathbf{A} \vee \neg_{\mathbf{s}} \mathbf{A}$ **if** $\mathbf{A} \in \mathcal{F}_n^{\#}$,
(13) $\mathbf{B} \to (\neg_{\mathbf{s}} \mathbf{B} \to \mathbf{A})$ **if** $\mathbf{A}, \mathbf{B} \in \mathcal{F}_n^{\#}$,
(14) $(\mathbf{A} \vee \neg_w \mathbf{A}) \vee (\mathbf{A} \wedge \neg_w \mathbf{A}) \vee \mathbf{A}^{[2]} \underbrace{\vee \ldots \vee \mathbf{A}^{[k]} \vee \ldots \vee}_{n-1} \mathbf{A}^{[n]}$ **if** $\mathbf{A} \in \mathcal{F}_n^{\#}$,

(15) $\mathbf{B} \to (\neg_w \mathbf{B} \to \mathbf{A})$ **if** $\mathbf{B} \notin \hat{\mathbf{V}} = \bigcup_{k=0}^{n} \mathbf{V}^{[k]}$.

**II.Rules of a conclusion**:

**Restricted Modus Ponens rule** $\mathbf{MPR}$ :

$\mathbf{A}, \mathbf{A} \to \mathbf{B} \vdash \mathbf{B}$ **iff** $\mathbf{A} \notin \hat{\mathbf{V}}$.
**Unrestricted Modus Tollens rule**: $\mathbf{P} \to \mathbf{Q}, \neg_w \mathbf{Q} \vdash \neg_w \mathbf{P}; \mathbf{P} \to \mathbf{Q}, \neg_{\mathbf{s}} \mathbf{Q} \vdash \neg_{\mathbf{s}} \mathbf{P}$.

**The rule of a strong contradiction**: $\mathbf{A} \wedge \neg_{\mathbf{s}} \mathbf{A} \vdash \mathbf{B}$.

**III.Quantification**

Corresponding to the propositional paraconsistent relevant logic $\mathbf{LP}_n^{\#}[\hat{\mathbf{V}}]$ we construct the corresponding paraconsistent relevant first-order predicate calculus. These new calculus will be denoted by $\overline{\mathbf{LP}}_n^{\#}[\hat{\mathbf{V}}]$.

The postulates of $\overline{\mathbf{LP}}_n^{\#}[\hat{\mathbf{V}}]$ are those of $\mathbf{LP}_n^{\#}[\hat{\mathbf{V}}]$ (suitably adapted)
plus the following:

**(I)** $\dfrac{\alpha \to \beta(x)}{\alpha \to \forall x \beta(x)}$,

**(II)** $\forall x \alpha(x) \to \alpha(y)$,

**(III)** $\alpha(x) \to \exists x \alpha(x)$,

**(IV)** $\dfrac{\alpha(x) \to \beta}{\exists x \alpha(x) \to \beta}$,

**(V)** $\forall x (\alpha(x))^{(k)} \to (\forall x \alpha(x))^{(k)}, k = 1, 2, \ldots, n$,

**(VI)** $\forall x((\alpha(x))^{(k)} \to (\exists x\alpha(x))^{(k)}, k = 1,2\ldots,n,$
**(VII)** $\forall x(\alpha(x))^{[k]} \to (\forall x\alpha(x))^{[k]}, k = 1,2\ldots,n.$

From the calculus $\overline{\mathbf{LP}}_n^{\#}\left[\hat{\mathbf{V}}\right]$,we can construct the following predicate calculus with equality.This is done by adding to their languages the binary predicates symbol of strong equality ($\bullet = \bullet$) or ($\bullet =_{\mathbf{s}} \bullet$) and weak equality($\bullet =_w \bullet$) with suitable modifications in the concept of formula, and by adding the following postulates:
**(IX)** $\forall x(x =_{\mathbf{s}} x),$
**(X)** $\forall x\left[\ (x =_{\mathbf{s}} x)^{[1]} \vdash \mathbf{B}\right],$
**(XI)** $\forall x\forall y[x =_{\mathbf{s}} y \to (\alpha(x) \leftrightarrow \alpha(y))],$
**(XII)** $\forall x\forall y\forall z[(x =_{\mathbf{s}} y) \wedge (y =_{\mathbf{s}} z) \to x =_{\mathbf{s}} z],$

**(XIII)** $\forall y\exists x(x =_w x)^{[k]}, k = 0,1,2,\ldots,n,$

**(XIV)** $\forall x\forall y\left[(x =_w y)^{[k]} \leftrightarrow \forall\alpha(\circ)(\alpha^{[k]}(x) \leftrightarrow \alpha^{[k]}(y))\right], k = 1,2,\ldots,n,$

**(XV)** $\forall x\forall y\forall z\left[(x =_w y)^{[k]} \wedge (y =_w z)^{[k]} \to (x =_w z)^{[k]}\right], k = 0,1,2,\ldots,n,$

**(XVI** $\forall x\forall y\forall z\left[(x =_w y)^{[k]} \wedge (y =_{\mathbf{s}} z) \to (x =_w z)^{[k]}\right], k = 0,1,2,\ldots,n,$

**(XVII)** $\forall x\forall y\forall z\left[(x =_{\mathbf{s}} y) \wedge (y =_w z)^{[k]} \to (x =_w z)^{[k]}\right], k = 0,1,2,\ldots,n,$

**(XVIII)** $\forall y\exists x(y =_w x)^{[k]}, k = 0,1,2,\ldots,n.$

# III. Paraconsistent Logic with infinite hierarchy levels of contradiction $\mathbf{LP}_{\boldsymbol{\omega}}^{\#}$.

The postulates (or their axioms schemata) of propositional paraconsistent logic $\mathbf{LP}_{\boldsymbol{\omega}}^{\#} = \mathbf{LP}_{\boldsymbol{\omega}}^{\#}\left[\hat{\mathbf{V}},\hat{\boldsymbol{\Delta}}\right]$ are the following:

The language $\mathcal{L}_{\boldsymbol{\omega}}^{\#}$ of paraconsistent logic $\mathbf{LP}_{\boldsymbol{\omega}}^{\#}$ has as primitive symbols (i) countable set of a clalassical propositional variables, (ii) countable set $\mathbf{V} = \{\mathbf{P}_i\}_{i\in\mathbb{N}}$ of a non clalassical propositional variables, (iii) the connectives $\neg_w, \neg_{\mathbf{s}}, \wedge, \vee, \to$ and (iv) the parentheses (,).

**Remark.3.1.**We distinguish a weak negation $\neg_w$ and a strong negation $\neg_{\mathbf{s}}$.

The definition of formula is the usual. We denote the set of the all formulae of $\mathbf{LP}_{\boldsymbol{\omega}}^{\#}\left[\hat{\mathbf{V}},\hat{\boldsymbol{\Delta}}\right]$ by $\mathcal{F}_{\boldsymbol{\omega}}^{\#}$ where $\hat{\mathbf{V}}$ and $\hat{\boldsymbol{\Delta}}$ is a given subsets of $\mathcal{F}_{\boldsymbol{\omega}}^{\#}$. We assume through that $\hat{\mathbf{V}} \subset \hat{\boldsymbol{\Delta}} \subsetneq \mathcal{F}_{\boldsymbol{\omega}}^{\#}$.

$\mathbf{A}, \mathbf{B}, \mathbf{C}$, ... will be used as metalanguage variables which indicate formulas of $\mathbf{LP}_{\boldsymbol{\omega}}^{\#}\left[\hat{\mathbf{V}},\hat{\boldsymbol{\Delta}}\right]$.

**Definition 3.1**. (i) $\alpha^{(n)}$ stands for $\alpha^{(n-1)} \wedge (\alpha^{(n-1)})^{(1)}$, where
$\alpha^{(0)} \triangleq \alpha, \alpha^{(1)} \triangleq \neg_w(\alpha \wedge \neg_w\alpha), 1 \leq n < \boldsymbol{\omega}.$

(ii) $\alpha^{(\omega)}$ stands for $\forall n[\alpha^{(n)}]$.
(iii) the finite $n$-order of the level of a weak consistency ($w$-consistency) is:
$\alpha^{(0)} \triangleq \alpha, \alpha^{(n)}, 1 \leq n < \boldsymbol{\omega}$.
(iv) the ifinite $\boldsymbol{\omega}$-order of level of a weak consistency ($w$-consistency) is : $\alpha^{(\omega)}$.
**Definition 3.2**. (i) $\alpha^{[n]}$ stands for $\alpha^{[n-1]} \wedge (\alpha^{[n-1]})^{[0]}$,
where $\alpha^{[0]} \triangleq \alpha \wedge \neg_w \alpha, 1 \leq n < \omega$.
(ii) $\alpha^{[\omega]}$ stands for $\forall n[\alpha^{[n]}]$.
(iii) the finite $n$-order of the level of a weak inconsistency ($w$-inconsistency) is:
$\alpha^{[n]}, 1 \leq n < \boldsymbol{\omega}$.
(iv) the ifinite $\boldsymbol{\omega}$-order of the level of a weak inconsistency ($w$-inconsistency) is: $\alpha^{[\omega]}$.
**I**. **Logical postulates**:

(1) $\mathbf{A} \to (\mathbf{B} \to \mathbf{A})$,
(2) $(\mathbf{A} \to \mathbf{B}) \to ((\mathbf{A} \to (\mathbf{B} \to \mathbf{C})) \to (\mathbf{A} \to \mathbf{C}))$,
(3) $\mathbf{A} \to (\mathbf{B} \to \mathbf{A} \wedge \mathbf{B})$,
(4) $\mathbf{A} \wedge \mathbf{B} \to \mathbf{A}$,
(5) $\mathbf{A} \wedge \mathbf{B} \to \mathbf{B}$,
(6) $\mathbf{A} \to (\mathbf{A} \vee \mathbf{B})$,
(7) $\mathbf{B} \to (\mathbf{A} \vee \mathbf{B})$,
(8) $(\mathbf{A} \to \mathbf{C}) \to ((\mathbf{B} \to \mathbf{C}) \to (\mathbf{A} \vee \mathbf{B} \to \mathbf{C}))$,

(9) $\mathbf{P}_i \wedge \neg_w \mathbf{P}_i$ **iff** $\mathbf{P}_i \in \mathbf{V}, i = 1,2,\ldots,$
(10) $\mathbf{P}_i^{[n]}$ **iff** $\mathbf{P}_i \in \mathbf{V}, i = 1,2,\ldots; 1 \leq n < \boldsymbol{\omega}$,
(11) $\mathbf{A} \vee \neg_w \mathbf{A}$ **if** $\mathbf{A} \notin \hat{\mathbf{V}} = \bigcup_{\mathbf{k} \in \mathbb{N}} \mathbf{V}^{[\mathbf{k}]}$,
(12) $\mathbf{A} \vee \neg_{\mathbf{s}} \mathbf{A}$ **if** $\mathbf{A} \in \mathcal{F}^{\#}_{\boldsymbol{\omega}}$,
(14) $\mathbf{B} \to (\neg_{\mathbf{s}} \mathbf{B} \to \mathbf{A})$ **if** $\mathbf{A}, \mathbf{B} \in \mathcal{F}^{\#}_{\boldsymbol{\omega}}$,
(15) $\mathbf{A} \vee \neg_w \mathbf{A} \vee \mathbf{A}^{[1]} \vee \mathbf{A}^{[2]} \underbrace{\vee \ldots \vee}_{n} \mathbf{A}^{[n]}$ **if** $\mathbf{A} \in \mathcal{F}^{\#}_{\boldsymbol{\omega}}, 1 \leq n < \boldsymbol{\omega}$,
(16) $\mathbf{B} \to (\neg_w \mathbf{B} \to \mathbf{A})$ **if** $\mathbf{B} \notin \hat{\mathbf{V}} = \bigcup_{\mathbf{k} \in \mathbb{N}} \mathbf{V}^{[\mathbf{k}]}$.

**II**.**Rules of a conclusion**:

**Restricted Modus Ponens rule MPR** :

$\mathbf{A}, \mathbf{A} \to \mathbf{B} \vdash \mathbf{B}$ **iff** $\mathbf{A} \notin \hat{\mathbf{V}}$.
**Unrestricted Modus Tollens rule**: $\mathbf{P} \to \mathbf{Q}, \neg_w \mathbf{Q} \vdash \neg_w \mathbf{P}; \mathbf{P} \to \mathbf{Q}, \neg_{\mathbf{s}} \mathbf{Q} \vdash \neg_{\mathbf{s}} \mathbf{P}$.

**The rule of a strong contradiction**: $\mathbf{A} \wedge \neg_{\mathbf{s}} \mathbf{A} \vdash \mathbf{B}$.

**III**.**Quantification**

Corresponding to the propositional paraconsistent relevant logic $\mathbf{LP}^{\#}_{\boldsymbol{\omega}}[\hat{\mathbf{V}}]$ we construct the corresponding paraconsistent relevant first-order predicate calculus. These new calculus will be denoted by $\overline{\mathbf{LP}}^{\#}_{\boldsymbol{\omega}}[\hat{\mathbf{V}}]$.

The postulates of $\overline{\mathbf{LP}}^{\#}_{\boldsymbol{\omega}}[\hat{\mathbf{V}}]$ are those of $\mathbf{LP}^{\#}_{\boldsymbol{\omega}}[\hat{\mathbf{V}}]$ (suitably adapted)

plus the following:

(**I**) $\dfrac{\alpha \to \beta(x)}{\alpha \to \forall x\beta(x)}$,

(**II**) $\forall x\alpha(x) \to \alpha(y)$,

(**III**) $\alpha(x) \to \exists x\alpha(x)$,

(**IV**) $\dfrac{\alpha(x) \to \beta}{\exists x\alpha(x) \to \beta}$,

(**V**) $\forall x(\alpha(x))^{(n)} \to (\forall x\alpha(x))^{(n)}, 1 \leq n < \boldsymbol{\omega}$,

(**VI**) $\forall x((\alpha(x))^{(n)} \to (\exists x\alpha(x))^{(n)}, 1 \leq n < \boldsymbol{\omega}$,

(**VII**) $\forall x(\alpha(x))^{[n]} \to (\forall x\alpha(x))^{[n]}, 1 \leq n < \boldsymbol{\omega}$,. .

From the calculus $\overline{\mathbf{LP}}^{\#}_{\omega}\left[\hat{\mathbf{V}}\right]$,we can construct the following predicate calculus with equality.This is done by adding to their languages the binary predicates symbol of strong equality ($\bullet = \bullet$) or ($\bullet =_{\mathbf{s}} \bullet$) and weak equality($\bullet =_w \bullet$) with suitable modifications in the concept of formula, and by adding the following postulates:

(**IX**) $\forall x(x =_{\mathbf{s}} x)$,

(**X**) $\forall x\left[\ (x =_{\mathbf{s}} x)^{[1]} \vdash \mathbf{B}\right]$,

(**XI**) $\forall x\forall y[x =_{\mathbf{s}} y \to (\alpha(x) \leftrightarrow \alpha(y))]$,

(**XII**) $\forall x\forall y\forall z[(x =_{\mathbf{s}} y) \wedge (y =_{\mathbf{s}} z) \to x =_{\mathbf{s}} z]$,

(**XIII**) $\forall y\exists x(x =_w x)^{[n]}, 0 \leq n < \boldsymbol{\omega}$,

(**XIV**) $\forall x\forall y\left[(x =_w y)^{[n]} \leftrightarrow \forall\alpha(\circ)(\alpha^{[n]}(x) \leftrightarrow \alpha^{[n]}(y))\right], 1 \leq n < \boldsymbol{\omega}$,

(**XV**) $\forall x\forall y\forall z\left[(x =_w y)^{[n]} \wedge (y =_w z)^{[n]} \to (x =_w z)^{[n]}\right], 0 \leq n < \boldsymbol{\omega}$,

(**XVI**) $\forall x\forall y\forall z\left[(x =_w y)^{[n]} \wedge (y =_{\mathbf{s}} z) \to (x =_w z)^{[n]}\right], 0 \leq n < \boldsymbol{\omega}$,

(**XVII**) $\forall x\forall y\forall z\left[(x =_{\mathbf{s}} y) \wedge (y =_w z)^{[n]} \to (x =_w z)^{[n]}\right], 0 \leq n < \boldsymbol{\omega}$,

(**XVIII**) $\forall y\exists x(y =_w x)^{[n]}, 0 \leq n < \boldsymbol{\omega}$.

# IV. Paraconsistent Set Theory $\mathbf{ZFC}^{\#}_{\boldsymbol{\omega}}$.

## IV.1. Paraconsistent set theory $\mathbf{KSth}^{\#}_{\boldsymbol{\omega}}$

Cantor’s "naive" set theory **KSth** was based mainly on two fundamental principles: the postulate of extensionality (if the sets $x$ and $y$ have the same elements, then they are equal), and the postulate of comprehension or separation (every property determines a set, composed of the objects that have this property). The latter postulate, in the standard (first-order) language of set theory, becomes the following schema of formulas:

$$\exists y\forall x(x \in y \leftrightarrow F(x,y)). \qquad (4.1.1)$$

Now, it is enough to replaces the formula $F(x,y)$ in (4.1) by $x \notin x$ to derive Russell’s paradox. That is, the principle of comprehension (4.1) entails an inconsistency. Thus, if one adds (4.1) to classical first-order logic, conceived as the logic of a set-theoretic language, a trivial theory is obtained.

**Remark.4.1.** We distinguish a weakly inconsisten membership relation $(\circ \in_w \circ)$ and a strongly consisten membership relation$(\circ \in_{\mathbf{s}} \circ)$.

**Definition 4.1.** (i) the minimal order of the level of a weak consistency ($w$-consistency) is: $\alpha^{(1)} \triangleq \alpha^{(0)} \wedge \neg_w(\alpha^{(0)} \wedge \neg_w\alpha^{(0)}), \alpha^{(0)} \triangleq \alpha = (x \in_w y)$;

(ii) the minimal order of the level of a weak inconsistency ($w$ -inconsistency) is: $\alpha^{[1]} \triangleq (\alpha^{[0]} \wedge \neg_w\alpha^{[0]}), \alpha^{[0]} \triangleq \alpha = (x \in_w y)$.

**Definition 4.2.**(i) $x \in_{w,(n)} y$ is to stands for $(x \in_w y)^{(n)}$ and is to mean "$x$ is a weakly consistent member of $y$ of the $n$-order (of the $n$-level) of $w$-consistency ".

(ii) $x \in_{w,[n]} y$ is to stands for $(x \in_w y)^{[n]}$ and is to mean "$x$ is a weakly inconsistent member of $y$ of the $n$-order (of the $n$-level) of $w$-inconsistency ".

**Definition 4.1.** An $\overline{\mathbf{LP}}_1^{\#}$ wff $\Phi$ is a $w$-wff iff it does not contain the connective: $\neg_w$.

We now replace the formula (4.1) by formulae

$$\exists y \forall x(x \in_{w,(n)} y \leftrightarrow F(x,y)), \quad n = 0,1,2,\ldots \tag{4.1.2}$$

and

$$\exists y \forall x(x \in_{w,[n]} y \leftrightarrow F(x,y)), \quad n = 0,1,2,\ldots . \tag{4.1.3}$$

**Theorem 4.1.** (1) The collections $\Re_n \triangleq \forall x[(x \in_{w,(n)} \Re_n) \leftrightarrow [\neg_w(x \in_{w,(n)} x)]]$ is contradictory of the $n+1$-order of $w$-inconsistency.

(2) The collections $\Re_n \triangleq \forall x[(x \in_{w,[n]} \Re_n) \leftrightarrow [\neg_w(x \in_{w,[n]} x)]]$ is contradictory of the $n+1$-order of $w$-inconsistency.

**Theorem 4.2.** (1) The collection $\Re_{\boldsymbol{\omega}} \triangleq \forall x \forall n[(x \in_{w,(n)} \Re_{\boldsymbol{\omega}}) \leftrightarrow [\neg_w(x \in_{w,(n)} x)]]$ is contradictory of the $\boldsymbol{\omega}+1$-order of $w$-inconsistency.

(2) The collection $\Re_{\boldsymbol{\omega}} \triangleq \forall x \forall n[(x \in_{w,[n]} \Re_{\boldsymbol{\omega}}) \leftrightarrow [\neg_w(x \in_{w,[n]} x)]]$ is contradictory of the $\boldsymbol{\omega}+1$-order of $w$-inconsistency.

The standard non-classical response to these paradoxes is to find fault with the *logical and deduction* principles involved in the deduction. Most standard approaches to the paradoxes take them to be important lessons in the behaviour of a Boolean negation.

However if you wish to define negation non-classically, there are many options available.You can define negation inferentially, taking $\mathbf{A}$ to mean that if $\mathbf{A}$, then something absurd follows,or it can be defined by way of the equivalence between the truth of $\sim\mathbf{A}$ and the falsity of $\mathbf{A}$, and allowing truth and falsity to have rather more independence from one another than is usually taken to be the case: say, allowing statements to be neither true nor false, or both true and false. The former account takes truth as primary, and defines negation in terms of a rejected proposition and implication.

For example, one can to define a strong negation $\sim_{\mathbf{s}}A$ non-classically [16]:

$$\sim_s \mathbf{A} \triangleq \mathbf{A} \to \forall x \forall y[(x \in_w y) \wedge (x =_s y)]. \tag{4.1.4}$$

**Theorem 4.3**. The collection $\Re_{\sim_s}$ such that $[x \in_w \Re_{\sim_s} \leftrightarrow [\sim_s(x \in_w x)]]$ i.e., $\Re_{\sim_s} \triangleq \hat{x}[\sim_s(x \in_w x)]$ is contradictory.

Proof. Replace $F(x,y)$ in the axiom schema of abstraction (4.2) in the definition of collection by $\sim_s(x \in_w x)$, so that the implicit definition of $\Re_{\sim_s}$ becomes

$$x \in_w \Re_{\sim_s} \leftrightarrow [\sim_s(x \in_w x)]. \tag{4.1.5}$$

Instantiating in (4.5) $x$ by $\Re_{\sim_s}$ then by unrestricted modus pones $\mathbf{MP}$, we obtain:

(1) $\vdash \Re_{\sim_s} \in_w \Re_{\sim_s} \leftrightarrow \sim_s(\Re_{\sim_s} \in_w \Re_{\sim_s})$.

By unrestricted modus pones $\mathbf{MP}$ one obtain the contradiction

(2) $\vdash \Re_{\sim_s} \in_w \Re_{\sim_s} \wedge \sim_s(\Re_{\sim_s} \in_\mathbf{w} \Re_{\sim_s})$.

Thus, if we adds (4.2)-(4.3) to first-order logic $\overline{\mathbf{LP}}^{\#}_{\omega}\left[\hat{\mathbf{V}}, \hat{\mathbf{\Delta}}\right]$, conceived as the logic of a set-theoretic language with suitable adapted $\hat{\mathbf{V}}$ and $\hat{\mathbf{\Delta}}$ a nontrivial paraconsistent set theory $\mathbf{KSth}^{\#}_{\boldsymbol{\omega}}$ is obtained.

# IV.2. Paraconsistent Set Theory $\mathbf{ZFC}^{\#}_{\omega}$. Basic Definitions and Elementary Operations on Inconsistent Sets.

**Remark 4.2.1**. In this subsection, we will be, to distinguish:

(i) a weak implication $A \to_w B$, where $A \to_w B$ abbreviates $A, A \to_w B \nvdash B, B \in \hat{\mathbf{V}}$ and

(ii) a strong implication $A \to_s B$, where $A \to_s B$ abbreviates $A, A \to_s B \vdash B$;

(iii) a weak negation $\neg_w A$, where $\neg_w A$ abbreviates $A \nvdash B, B \in \hat{\mathbf{V}}$ and

(iv) a strong negation $\neg_s A$ where $\neg_s A$ abbreviates $A \vdash B$.

**Designations 4.2.1.**We will be write for short:

$x =_{w_{[n]}} y$ instead $(x =_w y)^{[n]}, n = 1,2,\ldots;$

and we will write for short:

$x \in_{w_{[n]}} y$ instead $(x \in_w y)^{[n]}, n = 1,2,\ldots$

**Remark 4.2.2**. Thus in particular we will be write:

$x =_{w_{[1]}} y$ instead $(x =_w y) \wedge [(x =_w y) \wedge \neg_w(x =_w y)]$, etc.

and we will be write:

$x \in_{w_{[1]}} y$ instead $(x \in_w y) \wedge [(x \in_w y) \wedge \neg_w(x \in_w y)]$, etc.

**Remark 4.2.3.**However we will be often write for short:

$x =_w y$ instead $x =_{w[0]} y$ and $x \in_w y$ instead $x \in_{w[0]} y$.

**Remark 4.2.4**. In this subsection, we will be distinguish:

(**I**) the relations:

(i) consistent ($s$-consistent) equality denoted by $(\bullet =_s \bullet)$ and such that

$$\forall x, y[x =_s y \wedge \neg_s(x =_s y) \vdash B]; \tag{4.2.1}$$

(ii) weak or strongly inconsistent (or $w_\infty$-inconsistent) equality denoted by $(\circ =_{w_\infty} \circ)$ or by $(\circ =_w \circ)$ for short, and such that

$$\exists x, y[(x =_w y) \wedge \neg_w(x =_w y) \nvdash B]; \tag{4.2.2}$$

(iii) $w_{[n]}$-inconsistent) equalities denoted by$(\bullet =_{w[1]} \bullet), \ldots, (\bullet =_{w[n]} \bullet), \ldots, n = 1, 2, \ldots$ and such that

$$\begin{gathered} \forall n \exists x, y[(x =_{w_{[n]}} y) \wedge \neg_w(x =_{w_{[n]}} y) \nvdash B], \\ \forall x, y[x =_{w_{[1]}} y \Rightarrow_s x =_{w_{[0]}} y], \\ \forall n \forall x, y[x \in_{w_{[n+1]}} y \Rightarrow_s x =_{w_{[n]}} y], \end{gathered} \tag{4.2.3}$$

where $x =_{w_{[0]}} y \triangleq x =_w y$;

(**II**) the relations:

(i) consistent (or $s$-consistent) membership relation denoted by $(\bullet \in_s \bullet)$, and such that

$$\forall x, y[x \in_s y \wedge \neg_s(x \in_s y) \vdash B]; \tag{4.2.4}$$

(ii) weak or strongly inconsistent (or $w_\infty$-inconsistent) membership relation denoted by $(\bullet \in_w \bullet)$ and such that

$$\begin{gathered} \exists x, y[(x \in_w y) \wedge \neg_w(x \in_w y) \nvdash B], \\ \forall n \exists x, y[x \in_{w_{[n]}} y \nvdash B] \end{gathered} \tag{4.2.5}$$

(iii) $w_{[n]}$-inconsistent membership relations denoted by $(\bullet \in_{w_{[1]}} \bullet), \ldots,$ $(\bullet \in_{w_{[n]}} \bullet), \ldots,$ $n = 1, 2, \ldots$ and such that

$$\begin{gathered} \forall n \exists x, y[(x \in_{w_{[n]}} y) \wedge \neg_w(x \in_{w_{[n]}} y) \nvdash B], \\ \forall x, y[x \in_{w_{[1]}} y \Rightarrow_s x \in_{w_{[0]}} y], \\ \forall n \forall x, y[x \in_{w_{[n+1]}} y \Rightarrow_s x \in_{w_{[n]}} y], \end{gathered} \tag{4.2.6}$$

where $x \in_{w_{[0]}} y \triangleq x \in_w y$;

**Remark 4.2.5**. Note that: (1) in accordance with (4.2.2) the $w_\infty$-inconsistent equality $(\circ =_w \circ)$ admit the infinite levels of a contradiction;

**Definition 4.2.1**. Let $x$ and $X$ be a sets such that:

(i) the statement $x \in_s X$ holds, then we will be say that $x$ is a strong member (or $s$-member) of a set $X$;

(ii) the statement $x \in_w X$ holds, then we will be say that $x$ is a weak member (or $w$-member) of a set $X$;

**Remark 4.2.5**. We note, that in $\mathbf{ZFC}_{\omega}^{\#}$ valid:
(i) $\exists x,y[(x =_s y) \wedge \neg_s(x =_s y)] \vdash B,$
(ii)
(ii) $\exists x,y[(x =_w y) \wedge \neg_w(x =_w y) \nvdash B], B \in \hat{\mathbf{V}}, n = 1,2,\ldots$
(ii) $\forall x,y[(x =_s y) \wedge \neg_w(x =_s y)] \vdash B,$
(ii) $\exists x,y : (x =_w y) \wedge \neg_n^*(x =_w y) \nvdash B, B \in \hat{\mathbf{V}}, n = 1,2,\ldots$
(iii) $\forall x,y : (x \in_s y) \wedge \neg_w(x \in_s y) \vdash B,$
(iv) $\exists x,y : (x \in_w y) \wedge \neg_n^*(x \in_w y) \nvdash B, B \in \hat{\mathbf{V}}, n = 1,2,\ldots$
**Remark 4.2.4**. (i) $\neg_s A$ abbreviates $A \vdash B$, i.e. $\neg_s$ is a strong negation,
(ii) $A \Rightarrow_s B$ abbreviates $A, A \Rightarrow_s B \vdash B,$ i.e. $\Rightarrow_s$ is a strong implication.
**Designations 4.2.2**. (**I**) We will be write for short:
(i) $x =_{w_{\{0\}}} y$ instead $\left[(x =_w y) \vee \left[(x =_w y) \wedge \neg_w(x =_{w_{[1]}} y)\right]\right],$

(ii) $x =_{w_{\{1\}}} y$ instead $\left[(x =_w y) \vee (x =_w y) \vee [x =_{w_{[1]}} y] \wedge \neg_w(x =_{w_{[2]}} y)\right],$
(iii) $x =_{w_{\{n\}}} y$ instead $\left[(x =_w y) \vee (x =_w y) \vee \ldots \vee [(x =_{w_{[n]}} y)] \wedge \neg_w(x =_{w_{[n+1]}} y)\right],$
$n = 1,2,\ldots$
(iv) $x \in_{w_{\{0\}}} y$ instead $\left[(x \in_w y) \vee \left[(x \in_w y) \wedge \neg_w(x \in_{w_{[1]}} y)\right]\right],$
(v)
(vi)
(iii) $x \in_{w_{\{n\}}} y$ instead $\left[(x \in_s y) \vee (x \in_w y) \vee \ldots \vee \left[(x \in_{w_{[n]}} y) \wedge \neg_w(x \in_{w_{[n+1]}} y)\right]\right],$
$n = 1,2,\ldots$
(iv) $x =_{w_{\{\omega\}}} y$ instead $(x =_s y) \vee (x =_{w_{\{0\}}} y) \vee \bigvee_{0<n<\omega}(x =_{w_{[n]}} y),$
(v) $x \in_{w_{\{\omega\}}} y$ instead $(x \in_s y) \vee (x \in_{w_{\{0\}}} y) \vee \bigvee_{0<n<\omega}(x \in_{w_{[n]}} y).$
(**II**) We will be write for short:
(i) $x =_w^s y$ instead $[(x =_s y) \vee [(x =_w y) \wedge \neg_s[\neg_w(x =_w y)]]]$
(ii) $x =_{w_{\{0\}}}^s y$ instead $\left[(x =_s y) \vee \left[(x =_w y) \wedge \neg_s(x =_{w_{[1]}} y)\right]\right],$
(iii) $x =_{w_{\{1\}}}^s y$ instead $\left[(x =_s y) \vee (x =_w y) \vee [x =_{w_{[1]}} y] \wedge \neg_s(x =_{w_{[2]}} y)\right],$
(iv) $x =_{w_{\{n\}}}^s y$ instead $\left[(x =_s y) \vee (x =_w y) \vee \ldots \vee [(x =_{w_{[n]}} y)] \wedge \neg_s(x =_{w_{[n+1]}} y)\right],$ $n = 1,2,\ldots$
(v) $x \in_{w_{\{n\}}}^s y$ instead $\left[(x \in_s y) \vee (x \in_w y) \vee \ldots \vee \left[(x \in_{w_{[n]}} y) \wedge \neg_s(x \in_{w_{[n+1]}} y)\right]\right],$
$n = 1,2,\ldots$
(vi) $x =_{w_{\{\omega\}}}^s y$ instead $(x =_s y) \vee (x =_{w_{\{0\}}}^s y) \vee \bigvee_{0<n<\omega}(x =_{w_{\{n\}}}^s y),$
(vii) $x \in_{w_{\{\omega\}}}^s y$ instead $(x \in_s y) \vee (x \in_{w_{\{0\}}}^s y) \vee \bigvee_{0<n<\omega}(x \in_{w_{\{n\}}}^s y).$
(**III**)We often will be write for short:
(i) $x =_{w_0} y$ instead $x =_{w_{\{0\}}} y,$
(ii) $x =_{w_1} y$ instead $x =_{w_{\{1\}}} y, n = 1,2,\ldots,$
(iii) $x =_{w_n} y$ instead $x =_{w_{\{n\}}} y, n = 1,2,\ldots,$
(iv) $x =_{w_\omega} y$ instead $x =_{w_{\{\omega\}}} y,$
(v) $x \in_{w_\omega} y$ instead $x \in_{w_{\{\omega\}}} y.$
**Definition 4.2.1**. Let $x$ be an object (set). We shall say that $x$ is a strongly consistent object (s-consistent) or classical object iff: $x =_s x$ and $(x =_w x) \vdash B,$ i.e. $x$ is a strongly consistent object (set) iff $(x =_s x) \wedge \neg_s(x =_w x)$.

**Designations 4.2.3**. We will be write for short: $s\text{-}con(x)$ iff $x$ is $s$-consistent object

(set).

**Definition 4.2.2.**Let $x$ be an object (set). We shall say that:

(i) $x$ is a weakly consistent (w-consistent) object (set) iff

(1) $x =_w x$,(2) $(x =_w x) \nvdash B, B \in \hat{\mathbf{V}}$ (3) $\neg_s(x =_s x)$ and (4) $\neg_s(x =_{w_1} x)$, i.e. $(x =_{w_1} x) \vdash B$;

(ii) $x$ is a weakly inconsistent (w-inconsistent) object (set) iff

(1) $x =_w x$, (2) $\neg_s(x =_s x)$ and (3) $x =_{w_1} x \nvdash B, B \in \hat{\mathbf{V}}$.

**Designations 4.2.4**. We will be write for short:

(i) $w\text{-}con(x)$ or $w_0^*\text{-}con(x)$ iff $x$ is w-consistent object (set).

(ii) $w\text{-}inc(x)$ or $w_0\text{-}inc(x)$ iff $x$ is w-inconsistent object (set).

**Definition 4.2.3.**Let $x$ be an object (set). We shall say that:

(i) $x$ is $w_1$-inconsistent object iff $x =_{w_1} x$ and $x =_{w_1} x \nvdash B, B \in \hat{\mathbf{V}}$.

(ii) $x$ is $w_n$-inconsistent object iff $x =_{w_n} x$ and $x =_{w_n} x \nvdash B, B \in \hat{\mathbf{V}}, n = 1,2,\ldots$

(iii) $x$ is $w_\infty$-inconsistent object iff $\forall n\left[(x =_{w_n} x) \wedge \left(x =_{w_n} x \nvdash B, B \in \hat{\mathbf{V}}\right)\right]$.

**Designations 4.2.5**. We will be write for a short:

(i) $w_n\text{-}inc(x)$ iff $x$ is $w_n$-inconsistent object (set), $n = 1,2,\ldots$

(ii) $w_\infty\text{-}inc(x)$ iff $x$ is $w_\infty$-inconsistent object (set).

**Definition 4.2.4.**Let $x$ be an object (set). We shall say that:

(i) $x$ is a weakly $w_1$-inconsistent object iff $w_1\text{-}inc(x)$ and $\neg_s w_2\text{-}inc(x)$.

(ii) $x$ is a weakly $w_n$-inconsistent object iff $w_n\text{-}inc(x)$ and $\neg_s w_{n+1}\text{-}inc(x), n = 1,2,\ldots$

**Designations 4.2.6**. We will be write for a short:

(i) $w_1^*\text{-}inc(x)$ iff $x$ is a weakly $w_1$-inconsistent object (set).

(ii) $w_n^*\text{-}inc(x)$ iff $x$ is a weakly $w_n$-inconsistent object (set).

**Definition 4.2.5**. Let $x$ and $y$ be any $s$-consistent objects (sets),i.e. $s\text{-}con(x)$ and $s\text{-}con(y)$.

We shall say that objects (sets) $x$ and $y$ are strongly equivalent (s-equivalent) iff $x =_s y$.

**Definition 4.2.6**. Let $x$ and $y$ be an objects (sets) such that $w\text{-}con(x)$ and $w\text{-}con(y)$.

We shall say that objects (sets) $x$ and $y$ are weakly equivalent (w-equivalent) iff $x =_w y$.

**Definition 4.2.7**. Let $x$ and $y$ be an objects (sets) such that $w\text{-}con(x)$ and $w\text{-}con(y)$.

We shall say that objects (sets) $x$ and $y$ are weakly equivalent in consistent sense ($w$-equivalent) iff $x =_w y$ and $\neg_s(x =_{w_1} y)$.

**Definition 4.2.8**. Let $x$ and $y$ be any objects (sets) such that $w_n\text{-}inc(x)$ and $w_n\text{-}inc(y)$, then we shall say that:

(i) $x$ and $y$ are $w_n$- equivalent iff $x =_{w_n} y, x =_{w_n} y \nvdash B, B \in \hat{\mathbf{V}}, n = 0,1,2,\ldots$

(ii) $x$ and $y$ are $w_n$- equivalent in consistent sense ($w_n^*$- equivalent) iff $x =_{w_n} y$ and $\neg_s(x =_{w_{n+1}} y), n = 0,1,2,\ldots$

(iii) $x$ and $y$ are $w_\infty$- equivalent iff $\forall n\left[(x =_{w_n} y) \wedge \left(x =_{w_n} y \nvdash B, B \in \hat{\mathbf{V}}\right)\right]$

**Designations 4.2.7**. We will be write for a short:

(i) $x =_{w_n} y$ iff $x$ and $y$ is a $w_n$- equivalent,(ii) $x =_{w_n^*} y$ iff $x$ and $y$ are $w_n^*$- equivalent,

(iii) $x =_{w_\infty} y$ iff $x$ and $y$ are $w_\infty$- equivalent.

**Definition 4.2.9**. Let $x$ and $y$ be an objects (sets) such that $s\text{-}con(x)$ and $w_n\text{-}inc(y)$.

We shall say that objects (sets) $x$ and $y$ are weakly equivalent (w-equivalent) iff $x =_w y$.

**Definition 4.2.10**. Let $x$ and $y$ be any objects (sets) such that $s\text{-}inc(x)$ and $w_n\text{-}inc(y)$, then we shall say that:

(i) $x$ and $y$ are $w_n$- equivalent iff $x =_{w_n} y$ where $x =_{w_n} y \nvdash B, B \in \hat{\mathbf{V}}, n = 0,1,2,\ldots$

(ii) $x$ and $y$ are $w_n$- equivalent in consistent sense ($w_n^*$- equivalent) iff $x =_{w_n} y$ and $\neg_s(x =_{w_{n+1}} y), n = 0,1,2,\ldots$

(iii) $x$ and $y$ are $w_\infty$- equivalent iff $\forall n\left[(x =_{w_n} y) \wedge \left(x =_{w_n} y \nvdash B, B \in \hat{\mathbf{V}}\right)\right]$

**Designations 4.2.8**. We will be write for a short:

(i) $x =_{w_n} y$ iff $x$ and $y$ is a $w_n$- equivalent,

(ii) $x =_{w_n^*} y$ iff $x$ and $y$ are $w_n^*$- equivalent,

(iii) $x =_{w_\infty} y$ iff $x$ and $y$ are $w_\infty$- equivalent.

**Definition 4.2.11**. Let $x$ and $y$ be any objects (sets), then:

(i) we shall say that $x$ is a strongly consistent member (s-member) of $y$ if $x \in_s y$.

(ii) we shall say that $x$ is a weakly consistent member (wc-member) of $y$ if $x \in_w y$ and $\neg(x \in_w y) \vdash B$, i.e. $x$ is a weakly consistent member of $y$ if $x \in_w y$ and $\neg_s[\neg(x \in_w y)]$

(iii) we shall say that $x$ is a weak $w_1$-inconsistent member ($w_1$-member) of $y$ if $(x \in_w y) \wedge \neg(x \in_w y)$

**Designations 4.2.9**.We will be write for a short:

(i) $x \in_{wc} y$ or $x \in_{w_0^c} y$ iff $x$ is a weak consistent member of $y$

**Definition 4.2.12**. We shall say that:

(i) an formula $\varphi$ of Set Theory $\mathbf{ZFC}_{\boldsymbol{\omega}}^{\#}$ is a a strongly consistent formula iff formula $\varphi$ contains only predicates $x =_s y$ and $x \in_s y$. Sometimes we designate such formula by $\varphi_s$.

(ii)

(iii)

**Designations 4.2.10**. Before introducing any set-theoretic axioms at all, we can introduce some important abbreviations. Let $x$,$y$ and $z$ be any consistent objects (sets)

(i) $x \subseteq_s y$ abbreviates $\forall z(z \in_s x \to z \in_s y)$;

(ii) $x \subset_s y$ abbreviates $x \subseteq_s y \wedge x \neq_s y$;

(iii) $x \notin_s y$ abbreviates $\neg_s(x \in_s y)$;

(iv) $x \neq_s y$ abbreviates $\neg_s(x =_s y)$;

(v) $u =_s \ \cup_s x \triangleq \cup_s (x) \triangleq \forall z[z \in_s u \leftrightarrow (\exists y \in_s x)(z \in_s y)]$;

(vi) $u =_s \ \cap_s x \triangleq \cap_s (x) \triangleq \forall z[z \in_s u \leftrightarrow (\forall y \in_s x)(z \in_s y)]$;

(vii) $\exists x \in_s y\phi_s$ abbreviates $\exists x(x \in_s y \wedge \phi_s)$;

(viii) $\forall x \in_s y\phi_s$ abbreviates $\forall x(x \in_s y \to \phi_s)$;

(ix) $(\exists!_s x)\phi_s(x)$ abbreviates $(\exists x)\phi_s(x) \wedge \forall x\forall y[\phi_s(x) \wedge \phi_s(y) \to x =_s y]$

**Designations 4.2.11**.For any terms $r, s,$ and $t,$ we make the following abbreviations of formulas.

(i) $(\forall x \in_s t)\Phi$ or $(\forall x \in_s t)\Phi_s$ for $\forall x(x \in_s t \Rightarrow_s \Phi)$;

(ii) $(\forall x \in_s t)\Phi_w$ for $\forall x(x \in_s t \Rightarrow_w \Phi)$;

(iii) $(\forall x \in_{w_{[n]}} t)\Phi$ or $(\forall x \in_{w_{[n]}} t)\Phi_s$ for $\forall x(x \in_{w_{[n]}} t \Rightarrow_s \Phi)$;

(iv) $(\forall x \in_{w_{[n]}} t)\Phi_w$ for $\forall x(x \in_{w_{[n]}} t \Rightarrow_w \Phi)$;

(v) $(\forall x \in_{w_{\{n\}}} t)\Phi$ or $(\forall x \in_{w_{\{n\}}} t)\Phi_s$ for $\forall x(x \in_{w_{\{n\}}} t \Rightarrow_s \Phi)$;

(vi) $(\forall x \in_{w_{\{n\}}} t)\Phi_w$ for $\forall x(x \in_{w_{\{n\}}} t \Rightarrow_w \Phi)$.

**Designations 4.2.12**.For any terms $r, s,$ and $t,$ we make the following abbreviations of formulas.

(i) $(\exists x \in_s t)\Phi$ for $\exists x(x \in_s t \wedge \Phi)$;

(ii) $(\exists x \in_w t)\Phi$ for $\exists x(x \in_w t \wedge \Phi)$;

(iii) $(\exists x \in_{w_{[n]}} t)\Phi$ for $\exists x(x \in_{w_{[n]}} t \wedge \Phi)$;
(iv) $(\exists x \in_{w_{\{n\}}} t)\Phi$ for $\exists x(x \in_{w_{\{n\}}} t \wedge \Phi)$;
(v) $(\exists x \in^{s}_{w_{\{n\}}} t)\Phi$ for $\exists x\left(x \in^{s}_{w_{\{n\}}} t \wedge \Phi\right)$.

**Designations 4.2.13**.For any terms $r, s,$ and $t,$ we make the following abbreviations of formulas.
(i) $x \notin_s t$ or $x \notin^s_s t$ for $\neg_s(x \in_s t)$;
(ii) $x \notin^w_s t$ for $\neg_w(x \in_s t)$;
(iii) $x \notin^w_w t$ for $\neg_w(x \in_w t)$;
(iv) $x \notin^s_w t$ for $\neg_s(x \in_w t)$;
(iv) $x \notin^w_{w_{[n]}} t$ for $\neg_w(x \in_{w_{[n]}} t)$;
(iv) $x \notin^s_{w_{[n]}} t$ for $\neg_s(x \in_{w_{[n]}} t)$;
(v) $x \notin^w_{w_{\{n\}}} t$ for $\neg_w(x \in_{w_{\{n\}}} t)$;
(vi) $x \notin^s_{w_{\{n\}}} t$ for $\neg_s(x \in_{w_{\{n\}}} t)$;
(vii) $x \notin^{w,s}_{w_{\{n\}}} t$ for $\neg_w\left(x \in^s_{w_{\{n\}}} t\right)$;
(viii) $x \notin^{s,s}_{w_{\{n\}}} t$ for $\neg_s\left(x \in^s_{w_{\{n\}}} t\right)$.

**Designations 4.2.14**.
(i) The notation $\{x|\Phi(x)\}_{s,s}$ will stand for a set $X[\Phi]$ such that $\forall x[x \in_s X[\Phi] \Leftrightarrow_s \Phi(x)]$.
(ii) The notation $\{x|\Phi(x)\}_{w,s}$ will stand for a set $X[\Phi]$ such that $\forall x[x \in_w X[\Phi] \Leftrightarrow_s \Phi(x)]$.
(iii) The notation $\{x|\Phi(x)\}_{w,w}$ will stand for a set $X[\Phi]$ such that $\forall x[x \in_w X[\Phi] \Leftrightarrow_w \Phi(x)]$.
(iv) The notation $\{x|\Phi(x)\}_{w_{[n]},s}$ will stand for a set $X[\Phi]$ such that $\forall x[x \in_{w_{[n]}} X[\Phi] \Leftrightarrow_s \Phi(x)]$.
(v) The notation $\{x|\Phi(x)\}_{w_{[n]},s}$ will stand for a set $X[\Phi]$ such that $\forall x[x \in_{w_{[n]}} X[\Phi] \Leftrightarrow_s \Phi(x)]$.
(vi) The notation $\{x|\Phi(x)\}_{w_{\{n\}},w}$ will stand for a set $X[\Phi]$ such that $\forall x[x \in_{w_{\{n\}}} X[\Phi] \Leftrightarrow_w \Phi(x)]$.

**Designations 4.2.15**.Whenever we have a finite number of terms $t_1, t_2, \ldots, t_n$ then
(i) the notation $\{t_1, t_2, \ldots, t_n\}_{s,s}$ is used as an abbreviation for the class:
$\{x|x =_s t_1 \vee x =_s t_2 \vee \cdots \vee x =_s t_n\}_{s,s}$;
(ii) the notation $\{t_1, t_2, \ldots, t_n\}_{w,s}$ is used as an abbreviation for the class:
$\{x|x =_w t_1 \vee x =_w t_2 \vee \cdots \vee x =_w t_n\}_{w,s}$;
(iii) the notation $\{t_1, t_2, \ldots, t_n\}_{w,w}$ is used as an abbreviation for the class:
$\{x|x =_w t_1 \vee x =_w t_2 \vee \cdots \vee x =_w t_n\}_{w,w}$;
(iv) the notation $\{t_1, t_2, \ldots, t_n\}_{w_{[n]},s}$ is used as an abbreviation for the class:
$\{x|x =_{w_{[n]}} t_1 \vee x =_{w_{[n]}} t_2 \vee \cdots \vee x =_{w_{[n]}} t_n\}_{w_{[n]},s}$;
(v) the notation $\{t_1, t_2, \ldots, t_n\}_{w_{[n]},s}$ is used as an abbreviation for the class:
$\{x|x =_{w_{[n]}} t_1 \vee x =_{w_{[n]}} t_2 \vee \cdots \vee x =_{w_{[n]}} t_n\}_{w_{[n]},w}$;
(vi) the notation $\{t_1, t_2, \ldots, t_n\}_{w_{\{n\}},s}$ is used as an abbreviation for the class:
$\{x|x =_{w_{\{n\}}} t_1 \vee x =_{w_{\{n\}}} t_2 \vee \cdots \vee x =_{w_{\{n\}}} t_n\}_{w_{\{n\}},s}$;

(vii) the notation $\{t_1,t_2,\ldots,t_n\}_{w_{\{n\}},s}$ is used as an abbreviation for the class: $\{x|x =_{w_{\{n\}}} t_1 \vee x =_{w_{\{n\}}} t_2 \vee \cdots \vee x =_{w_{\{n\}}} t_n\}_{w_{\{n\}},w}$.

**Designations 4.2.16.** We abbreviate the following important sets:

(i) $s,s$-union $t_1 \cup_{s,s} t_2$ or $t_1 \cup_s t_2$ for $\{x|x \in_s t_1 \vee x \in_s t_2\}_{s,s}$;

(ii) $w,s$-union $t_1 \cup_{w,s} t_2$ or $t_1 \cup_{w,s} t_2$ for $\{x|x \in_w t_1 \vee x \in_w t_2\}_{w,s}$;

(iii) $w,w$-union $t_1 \cup_{w,w} t_2$ or $t_1 \cup_{w,w} t_2$ for $\{x|x \in_w t_1 \vee x \in_w t_2\}_{w,w}$;

(iv) $w_{[n]},s$-union $t_1 \cup_{w_{[n]},s} t_2$ or $t_1 \cup_{w_{[n]},s} t_2$ for $\{x|x \in_{w_{[n]}} t_1 \vee x \in_{w_{[n]}} t_2\}_{w_{[n]},s}$;

(v) $w_{[n]},w$-union $t_1 \cup_{w_{[n]},w} t_2$ or $t_1 \cup_{w_{[n]},w} t_2$ for $\{x|x \in_{w_{[n]}} t_1 \vee x \in_{w_{[n]}} t_2\}_{w_{[n]},w}$;

(vi) $w_{\{n\}},s$-union $t_1 \cup_{w_{\{n\}},s} t_2$ or $t_1 \cup_{w_{\{n\}},s} t_2$ for $\{x|x \in_{w_{\{n\}}} t_1 \vee x \in_{w_{\{n\}}} t_2\}_{w_{\{n\}},s}$;

(vii) $w_{\{n\}},w$-union $t_1 \cup_{w_{\{n\}},w} t_2$ or $t_1 \cup_{w_{\{n\}},w} t_2$ for $\{x|x \in_{w_{\{n\}}} t_1 \vee x \in_{w_{\{n\}}} t_2\}_{w_{\{n\}},w}$;

**Designations 4.2.17.** We abbreviate the following important sets:

(i) $s,s$-union $s,s\text{-}\bigcup C$ or $s,s\text{-}\bigcup_{S\in_s C} S$ for $\{x|x \in_s S \text{ for some } S \in_s C\}_{s,s}$;

(ii) $w,s$-union $w,s\text{-}\bigcup C$ or $w,s\text{-}\bigcup_{S\in_w C} S$ for $\{x|x \in_w S \text{ for some } S \in_w C\}_{w,s}$;

(iii) $w,w$-union $w,w\text{-}\bigcup C$ or $w,w\text{-}\bigcup_{S\in_w C} S$ for $\{x|x \in_w S \text{ for some } S \in_w C\}_{w,w}$;

**Designations 4.2.18.** We abbreviate the following important sets:

(i) $s,s$-intersection $t_1 \cap_{s,s} t_2$ or $t_1 \cap_s t_2$ for $\{x|x \in_s t_1 \vee x \in_s t_2\}_{s,s}$;

(ii) $w,s$-intersection $t_1 \cap_{w,s} t_2$ for $\{x|x \in_w t_1 \vee x \in_w t_2\}_{w,s}$;

(iii) $w,w$-intersection $t_1 \cap_{w,w} t_2$ for $\{x|x \in_w t_1 \vee x \in_w t_2\}_{w,w}$;

(iv) $w_{[n]},s$-intersection $t_1 \cap_{w_{[n]},w} t_2$ for $\{x|x \in_{w_{[n]}} t_1 \vee x \in_{w_{[n]}} t_2\}_{w_{[n]},w}$;

(v) $w_{[n]},w$-intersection $t_1 \cap_{w_{[n]},w} t_2$ for $\{x|x \in_{w_{[n]}} t_1 \vee x \in_{w_{[n]}} t_2\}_{w_{[n]},w}$;

(vi) $w_{\{n\}},s$-intersection $t_1 \cap_{w_{\{n\}},w} t_2$ for $\{x|x \in_{w_{\{n\}}} t_1 \vee x \in_{w_{\{n\}}} t_2\}_{w_{\{n\}},w}$;

(vii) $w_{\{n\}},w$-intersection $t_1 \cap_{w_{\{n\}},w} t_2$ for $\{x|x \in_{w_{\{n\}}} t_1 \vee x \in_{w_{\{n\}}} t_2\}_{w_{\{n\}},w}$.

**Designations 4.2.19.** We abbreviate the following important sets:

(i) $s,s$-intersection $s,s\text{-}\bigcap C$ or $s,s\text{-}\bigcap_{S\in_s C} S$ for $\{x|x \in_s S \text{ for all } S \in_s C\}_{s,s}$;

(ii) $w,s$-intersection $w,s\text{-}\bigcap C$ or $w,s\text{-}\bigcap_{S\in_s C} S$ for $\{x|x \in_w S \text{ for all } S \in_w C\}_{w,s}$;

(iii) $w,w$-intersection $w,w\text{-}\bigcap C$ or $w,w\text{-}\bigcap_{S\in_s C} S$ for $\{x|x \in_w S \text{ for all } S \in_w C\}_{w,w}$;

## IV.3. The Axioms of Paraconsistent Set Theory $\mathbf{ZFC}^{\#}_{\omega}$.

### IV.3.1.The Axioms of Extensionality.

(1) Strong axiom of $w$-extensionality

$$\forall u[u \in_w X \Leftrightarrow_s u \in_w Y] \Leftrightarrow_s X =_w Y. \qquad (4.3.1)$$

(2) Weak axiom of $w$-extensionality

$$\forall u[u \in_w X \Leftrightarrow_w u \in_w Y] \Leftrightarrow_w X =_w Y. \qquad (4.3.2)$$

(3) Strong axiom of $w_{[n]}$-extensionality

$$\forall u[u \in_{w_{[n]}} X \Leftrightarrow_s u \in_{w_{[n]}} Y] \Leftrightarrow_s X =_{w_{[n]}} Y. \qquad (4.3.3)$$

(4) Weak axiom of $w_{[n]}$-extensionality

$$\forall u\left[u \in_{w_{[n]}} X \Leftrightarrow_w u \in_{w_{[n]}} Y\right] \Leftrightarrow_w X =_{w_{[n]}} Y. \qquad (4.3.4)$$

(5) Strong axiom of $w_{\{n\}}$-extensionality

$$\forall u[u \in_{w_{\{n\}}} X \Leftrightarrow_s u \in_{w_{\{n\}}} Y] \Leftrightarrow_s X =_{w_{\{n\}}} Y. \qquad (4.3.5)$$

(6) Weak axiom of $w_{\{n\}}$-extensionality

$$\forall u\left[u \in_{w_{\{n\}}} X \Leftrightarrow_w u \in_{w_{\{n\}}} Y\right] \Leftrightarrow_w X =_{w_{\{n\}}} Y. \qquad (4.3.6)$$

# IV.3.2.The Axioms of Empty Set.

(1) Axiom of strongly $w$-empty set

$$\exists x \forall u[u \notin_w^s x]. \qquad (4.3.7)$$

The strongly $w$-empty set, denoted $\varnothing_w^s$.

(2) Axiom of weakly $w$-empty set

$$\exists x \forall u[u \notin_w^w x]. \qquad (4.3.8)$$

The weakly $w$-empty set, denoted $\varnothing_w^w$.

(3) Axiom of weakly $w_{\{0\}}$-empty set

$$\exists x \forall u\left[u \notin_{w_{\{0\}}}^w x\right]. \qquad (4.3.9)$$

The weakly $w_{\{0\}}$-empty set, denoted $\varnothing_{w_{\{0\}}}^w$.

# IV.3.3.The Axioms of Pairing.

(1) Strong axiom of $w,s$-pairing.

$$\forall a \forall b \exists c \forall x[x \in_w c \Leftrightarrow_s (x =_w a) \vee (x =_w b)] \quad (4.3.)$$

and we define the $w,s$-pair $\{a,b\}_{w,s}$ by $\{a,b\}_{w,s} =_w c$.

# IV.3.4.The Axioms of Separation.

(**1**) **Strong Separation Schemes**.

(i) Let $\phi(u,p_1,\ldots,p_k)$ be a formula free from symbols $\notin^s_{w_{[n]}}$ , $\notin^s_{w_{\{n\}}}$ , $n = 1,2,\ldots$. For any $X$ and $p_1,\ldots,p_k$, there exists a set $Y =_w \{u \in_w X|\phi(u,p_1,\ldots,p_k)\}_{w,s}$, i.e.

$$\forall X\forall p\exists Y\forall u[u \in_w Y \Leftrightarrow_s (u \in_w X) \wedge \phi(u,p_1,\ldots,p_k)] \qquad (4.3.4)$$

(ii) Let $\phi(u,p_1,\ldots,p_k)$ be a formula free from symbols $\notin^s_{w_{[n]}}$ , $\notin^s_{w_{\{n\}}}$ , $n = 1,2,\ldots$. For any $X$ and $p_1,\ldots,p_k$, there exists a set $Y =_{w_{[n]}} \{u \in_{w_{[n]}} X|\phi(u,p_1,\ldots,p_k)\}_{w_{[n]},s}$, i.e.

$$\forall X\forall p\exists Y\forall u[u \in_{w_{[n]}} Y \Leftrightarrow_s (u \in_{w_{[n]}} X) \wedge \phi(u,p_1,\ldots,p_k)] \qquad (4.3.)$$

(iii) Let $\phi(u,p_1,\ldots,p_k)$ be a formula free from symbols $\notin^s_{w_{[n]}}$ , $\notin^s_{w_{\{n\}}}$ , $n = 1,2,\ldots$. For any $X$ and $p_1,\ldots,p_k$, there exists a set $Y =_{w_{\{n\}}} \{u \in_{w_{\{n\}}} X|\phi(u,p_1,\ldots,p_k)\}_{w_{\{n\}},s}$, i.e.

$$\forall X\forall p\exists Y\forall u[u \in_{w_{\{n\}}} Y \Leftrightarrow_s (u \in_{w_{\{n\}}} X) \wedge \phi(u,p_1,\ldots,p_k)] \qquad (4.3.)$$

(**2**) **Weak Separation Schemes**.

(i) Let $\phi(u,p_1,\ldots,p_k)$ be a formula. For any $X$ and $p_1,\ldots,p_k$, there exists a set $Y =_w \{u \in_w X|\phi(u,p_1,\ldots,p_k)\}_{w,w}$, i.e.

$$\forall X\forall p\exists Y\forall u[u \in_w Y \Leftrightarrow_w (u \in_w X) \wedge \phi(u,p_1,\ldots,p_k)] \qquad (4.3.)$$

(ii) Let $\phi(u,p_1,\ldots,p_k)$ be a formula. For any $X$ and $p_1,\ldots,p_k$, there exists a set $Y =_{w_{[n]}} \{u \in_{w_{[n]}} X|\phi(u,p_1,\ldots,p_k)\}_{w_{[n]},s}$, i.e.

$$\forall X\forall p\exists Y\forall u[u \in_{w_{[n]}} Y \Leftrightarrow_w (u \in_{w_{[n]}} X) \wedge \phi(u,p_1,\ldots,p_k)] \qquad (4.3.)$$

(iii) Let $\phi(u,p_1,\ldots,p_k)$ be a formula. For any $X$ and $p_1,\ldots,p_k$, there exists a set $Y =_{w_{\{n\}}} \{u \in_{w_{\{n\}}} X|\phi(u,p_1,\ldots,p_k)\}_{w_{\{n\}},s}$, i.e.

$$\forall X\forall p\exists Y\forall u[u \in_{w_{\{n\}}} Y \Leftrightarrow_w (u \in_{w_{\{n\}}} X) \wedge \phi(u,p_1,\ldots,p_k)] \qquad (4.3.)$$

# IV.3.5.The Axioms of Replacement.

(**1**) **Strong Replacement Scheme**.

(i) Let $\phi(x,y,u)$ be a formula free from symbols $\notin^s_{w_{[n]}}$ , $\notin^s_{w_{\{n\}}}$ , then for any $n = 1,2,\ldots$.

$$\begin{aligned}&\forall x\forall y\forall y'[\phi(x,y,u) \wedge \phi(x,y',u) \Rightarrow_s y =_w y'] \Rightarrow_s \\ &\Rightarrow_s \forall s\exists z\forall y[y \in_w z \Leftrightarrow_s \exists x(x \in_w s)\phi(x,y,u)].\end{aligned} \qquad (4.3.)$$

The set $z$ is denoted $\{y|\exists x\phi(x,y,u) \wedge (x \in_w s)\}_{w,s}$.

(ii) Let $\phi(x,y,u)$ be a formula free from symbols $\notin^s_{w_{[n]}}$ , $\notin^s_{w_{\{n\}}}$ , then for any $u = (p_1,\ldots,p_k)$, $\quad n = 1,2,\ldots$.

$$\forall x \forall y \forall y'[\phi(x,y,u) \wedge \phi(x,y',u) \Rightarrow_s y =_{w_{[n]}} y'] \Rightarrow_s$$
$$\Rightarrow_s \forall s \exists z \forall y[y \in_{w_{[n]}} z \Leftrightarrow_s \exists x(x \in_{w_{[n]}} s)\phi(x,y,u)]. \quad (4.3.)$$

The set $z$ is denoted $\{y|\exists x\phi(x,y,u) \wedge (x \in_{w_{[n]}} s)\}_{w_{[n]},s}$.

(iii) Let $\phi(x,y,u)$ be a formula free from symbols $\notin^s_{w_{[n]}}$ , $\notin^s_{w_{\{n\}}}$ ,then for any $u = (p_1,\ldots,p_k)$, $n = 1,2,\ldots$.

$$\forall x \forall y \forall y'[\phi(x,y,u) \wedge \phi(x,y',u) \Rightarrow_s y =_{w_{\{n\}}} y'] \Rightarrow_s$$
$$\Rightarrow_s \forall s \exists z \forall y[y \in_{w_{\{n\}}} z \Leftrightarrow_s \exists x(x \in_{w_{\{n\}}} s)\phi(x,y,u)]. \quad (4.3.)$$

The set $z$ is denoted $\{y|\exists x\phi(x,y,u) \wedge (x \in_{w_{\{n\}}} s)\}_{w_{\{n\}},s}$.

(**2**) **Weak Replacement Scheme**.

(i) Let $\phi(x,y,u)$ be a formula, then for any $u = (p_1,\ldots,p_k)$, $n = 1,2,\ldots$

$$\forall x \forall y \forall y'[\phi(x,y,u) \wedge \phi(x,y',u) \Rightarrow_w y =_w y'] \Rightarrow_w$$
$$\Rightarrow_w \forall s \exists z \forall y[y \in_w z \Leftrightarrow_w \exists x(x \in_w s)\phi(x,y,u)]. \quad (4.3.)$$

The set $z$ is denoted $\{y|\exists x\phi(x,y,u) \wedge (x \in_w s)\}_{w,w}$.

(ii) Let $\phi(x,y,u)$ be a formula, then for any $u = (p_1,\ldots,p_k)$, $n = 1,2,\ldots$

$$\forall x \forall y \forall y'[\phi(x,y,u) \wedge \phi(x,y',u) \Rightarrow_w y =_{w_{[n]}} y'] \Rightarrow_w$$
$$\Rightarrow_w \forall s \exists z \forall y[y \in_{w_{[n]}} z \Leftrightarrow_w \exists x(x \in_{w_{[n]}} s)\phi(x,y,u)]. \quad (4.3.)$$

The set $z$ is denoted $\{y|\exists x\phi(x,y,u) \wedge (x \in_{w_{[n]}} s)\}_{w_{[n]},w}$.

(iii) Let $\phi(x,y,u)$ be a formula,then for any $u = (p_1,\ldots,p_k)$, $n = 1,2,\ldots$

$$\forall x \forall y \forall y'[\phi(x,y,u) \wedge \phi(x,y',u) \Rightarrow_s y =_{w_{\{n\}}} y'] \Rightarrow_w$$
$$\Rightarrow_w \forall s \exists z \forall y[y \in_{w_{\{n\}}} z \Leftrightarrow_s \exists x(x \in_{w_{\{n\}}} s)\phi(x,y,u)]. \quad (4.3.)$$

The set $z$ is denoted $\{y|\exists x\phi(x,y,u) \wedge (x \in_{w_{\{n\}}} s)\}_{w_{\{n\}},w}$.

w

# IV.3.6.The Axioms of Union.

# IV.3.6.(1) Strong Axiom of Union.

(i) Strong $w$-union

$$\forall x \exists y_s \forall t[t \in_w y_s \Leftrightarrow_s \exists u(u \in_w x \wedge t \in_w u)]. \quad (4.3.)$$

The set $y_s$ is denoted $\cup_{w,s}\, x$ or $s,w\text{-}\cup x$.

(ii) Strong $w_{[n]}$-union

$$\forall x \exists y_s \forall t[t \in_{w_{[n]}} y_s \Leftrightarrow_s \exists u(u \in_{w_{[n]}} x \wedge t \in_{w_{[n]}} u)]. \quad (4.3.)$$

The set $y_s$ is denoted $\cup_{w_{[n]},s}\, x$ or $s,w_{[n]}\text{-}\cup x$.

(iii) Strong $w_{\{n\}}$-union

$$\forall x \exists y_s \forall t[t \in_{w_{\{n\}}} y_s \Leftrightarrow_s \exists u(u \in_{w_{\{n\}}} x \wedge t \in_{w_{\{n\}}} u)]. \quad (4.3.)$$

The set $y_s$ is denoted $\cup_{w_{\{n\}},s}\, x$ or $s, w_{\{n\}}$-$\cup x$.

## IV.3.6.(2) Weak Axiom of Union.

(i) Weak $w$-union

$$\forall x \exists y_w \forall t[t \in_w y_w \Leftrightarrow_w \exists u(u \in_w x \wedge t \in_w u)]. \quad (4.3.)$$

The set $y_w$ is denoted $\cup_{w,w}\, x$ or $w, w$-$\cup x$.

(ii) Weak $w_{[n]}$-union

$$\forall x \exists y_w \forall t[t \in_{w_{[n]}} y_w \Leftrightarrow_w \exists u(u \in_{w_{[n]}} x \wedge t \in_{w_{[n]}} u)]. \quad (4.3.)$$

The set $y_w$ is denoted $\cup_{w_{[n]},w}\, x$ or $w, w_{[n]}$-$\cup x$.

(iii) Weak $w_{\{n\}}$-union

$$\forall x \exists y_w \forall t[t \in_{w_{\{n\}}} y_w \Leftrightarrow_w \exists u(u \in_{w_{\{n\}}} x \wedge t \in_{w_{\{n\}}} u)]. \quad (4.3.)$$

The set $y_w$ is denoted $\cup_{w_{\{n\}},w}\, x$ or $w, w_{\{n\}}$-$\cup x$.

## IV.3.7.The Axioms of Power Set.

## IV.3.7.(1) Strong Axioms of Power Set.

(i) Strong axiom of $w$-power set.

$$\forall X \exists Y_s \forall t[t \in_w Y_s \Leftrightarrow_s \forall z(z \in_w t \Rightarrow_s z \in_w X)] \quad (4.3.)$$

For any set $X$, a set $Y_s$ is denoted $P_w^s(X)$.

(ii) Strong axiom of $w_{[n]}$-power set.

$$\forall X \exists Y_s \forall t[t \in_{w_{[n]}} Y_s \Leftrightarrow_s \forall z(z \in_{w_{[n]}} t \Rightarrow_s z \in_{w_{[n]}} X)] \quad (4.3.)$$

For any set $X$, a set $Y_s$ is denoted $P_{w_{[n]}}^s(X)$.

(iii) Strong axiom of $w_{\{0\}}$-power set.

$$\forall X \exists Y_s \forall t[t \in_{w_{\{0\}}} Y_s \Leftrightarrow_s \forall z(z \in_{w_{\{0\}}} t \Rightarrow_s z \in_{w_{\{0\}}} X)] \quad (4.3.)$$

For any set $X$, a set $Y_s$ is denoted $P_{w_{\{0\}}}^s(X)$.

(iv) Strong axiom of $w_{\{n\}}$-power set.

$$\forall X \exists Y_s \forall t[t \in_{w_{\{n\}}} Y_s \Leftrightarrow_s \forall z(z \in_{w_{\{n\}}} t \Rightarrow_s z \in_{w_{\{n\}}} X)] \quad (4.3.)$$

For any set $X$, a set $Y_s$ is denoted $P_{w_{\{n\}}}^s(X)$.

## IV.3.7.(2) Weak Axioms of Power Set.

(i) Weak axiom of $w$-power set.

$$\forall X \exists Y_w \forall t[t \in_w Y_w \Leftrightarrow_w \forall z(z \in_w t \Rightarrow_w z \in_w X)] \quad (4.3.)$$

For any set $X$, a set $Y_w$ denoted $P_w^w(X)$.

(ii) Weak axiom of $w_{[n]}$-power set.

$$\forall X \exists Y_w \forall t[t \in_{w_{[n]}} Y_w \Leftrightarrow_w \forall z(z \in_{w_{[n]}} t \Rightarrow_w z \in_{w_{[n]}} X)] \quad (4.3.)$$

For any set $X,$ a set $Y_w$ is denoted $P^w_{w_{[n]}}(X).$

(iii) Weak axiom of $w_{\{0\}}$-power set.

$$\forall X \exists Y_w \forall t[t \in_{w_{\{0\}}} Y_w \Leftrightarrow_w \forall z(z \in_{w_{\{0\}}} t \Rightarrow_s z \in_{w_{\{0\}}} X)] \quad (4.3.)$$

For any set $X,$ a set $Y_s$ is denoted $P^w_{w_{\{0\}}}(X).$

(iv) Weak axiom of $w_{\{n\}}$-power set.

$$\forall X \exists Y_w \forall t[t \in_{w_{\{n\}}} Y_w \Leftrightarrow_w \forall z(z \in_{w_{\{n\}}} t \Rightarrow_w z \in_{w_{\{n\}}} X)] \quad (4.3.)$$

For any set $X,$ a set $Y_w$ is denoted $P^w_{w_{\{n\}}}(X).$

# IV.3.8.The Axioms of Infinity.

# IV.3.8.(1) Strong Axioms of Infinity.

(i) Strong Axiom of $w$-infinity

$$\begin{gathered}\exists X[\exists y(y \in_w X) \wedge \forall z(z \not\in^s_w y) \wedge \forall y(y \in_w X \Rightarrow_s \\ \Rightarrow_s \exists z(z \in_w X \wedge \forall t(t \in z \Leftrightarrow_s (t \in_w y \vee t =_w y))))]\end{gathered} \quad (4.3.)$$

There is a set $X$ such that $\varnothing^s_w \in_w X$ and whenever $y \in_w X,$ then $y \cup_{w,s} \{y\}_{w,s} \in_w X.$ A set $X$ is denoted $\mathbb{N}^s_w.$

(ii) Strong Axiom of $w_{[n]}$-infinity.

$$\begin{gathered}\exists X\left[\exists y(y \in_{w_{[n]}} X) \wedge \forall z\left(z \not\in^s_{w_{[n]}} y\right) \wedge \forall y(y \in_{w_{[n]}} X \Rightarrow_s \\ \Rightarrow_s \exists z(z \in_{w_{[n]}} X \wedge \forall t(t \in_{w_{[n]}} z \Leftrightarrow_s (t \in_{w_{[n]}} y \vee t =_{w_{[n]}} y))))\right]\end{gathered} \quad (4.3.)$$

There is a set $X$ such that $\varnothing^s_{w_{[n]}} \in_{w_{[n]}} X$ and whenever $y \in_{w_{[n]}} X,$ then $y \cup_{w_{[n]},s} \{y\}_{w_{[n]},s} \in_{w_{[n]}} X.$ A set $X$ is denoted $\mathbb{N}^s_{w_{[n]}}.$

(iii) Strong Axiom of $w_{\{0\}}$-infinity.

$$\begin{gathered}\exists X\left[\exists y(y \in_{w_{\{0\}}} X) \wedge \forall z\left(z \not\in^s_{w_{\{0\}}} y\right) \wedge \forall y(y \in_{w_{\{0\}}} X \Rightarrow_s \\ \Rightarrow_s \exists z(z \in_{w_{\{0\}}} X \wedge \forall t(t \in_{w_{\{0\}}} z \Leftrightarrow_s (t \in_{w_{\{0\}}} y \vee t =_{w_{\{0\}}} y))))\right]\end{gathered} \quad (4.3.)$$

There is a set $X$ such that $\varnothing^s_{w_{\{0\}}} \in_{w_{\{0\}}} X$ and whenever $y \in_{w_{\{0\}}} X,$ then $y \cup_{w_{\{0\}},s} \{y\}_{w_{\{0\}},s} \in_{w_{\{0\}}} X.$ A set $X$ is denoted $\mathbb{N}^s_{w_{\{0\}}}.$

(iv) Strong Axiom of $w_{\{n\}}$-infinity.

$$\begin{gathered}\exists X\left[\exists y(y \in_{w_{\{n\}}} X) \wedge \forall z\left(z \not\in^s_{w_{\{n\}}} y\right) \wedge \forall y(y \in_{w_{\{n\}}} X \Rightarrow_s \\ \Rightarrow_s \exists z(z \in_{w_{\{n\}}} X \wedge \forall t(t \in_{w_{\{n\}}} z \Leftrightarrow_s (t \in_{w_{\{n\}}} y \vee t =_{w_{\{n\}}} y))))\right]\end{gathered} \quad (4.3.)$$

There is a set $X$ such that $\varnothing^s_{w_{\{n\}}} \in_{w_{\{n\}}} X$ and whenever $y \in_{w_{\{n\}}} X,$ then $y \cup_{w_{\{n\}},s} \{y\}_{w_{\{n\}},s} \in_{w_{\{n\}}} X.$ A set $X$ is denoted $\mathbb{N}^s_{w_{\{n\}}}.$

# IV.3.8.(2) Weak Axioms of Infinity.

(i) Weak Axiom of $w$-infinity

$$\exists X[\exists y(y \in_w X) \wedge \forall z(z \notin_w^w y) \wedge \forall y(y \in_w X \Rightarrow_w \\ \Rightarrow_s \exists z(z \in_w X \wedge \forall t(t \in z \Leftrightarrow_w (t \in_w y \vee t =_w y))))] \quad (4.3.)$$

There is a set $X$ such that $\varnothing_w^w \in_w X$ and whenever $y \in_w X,$ then $y \cup_{w,s} \{y\}_{w,w} \in_w X.$
A set $X$ is denoted $\mathbb{N}_w^w$.
(ii) Weak Axiom of $w_{[n]}$-infinity.

$$\exists X\big[\exists y(y \in_{w_{[n]}} X) \wedge \forall z\big(z \notin_{w_{[n]}}^w y\big) \wedge \forall y(y \in_{w_{[n]}} X \Rightarrow_w \\ \Rightarrow_s \exists z(z \in_{w_{[n]}} X \wedge \forall t(t \in_{w_{[n]}} z \Leftrightarrow_w (t \in_{w_{[n]}} y \vee t =_{w_{[n]}} y))))] \quad (4.3.)$$

There is a set $X$ such that $\varnothing_{w_{[n]}}^s \in_{w_{[n]}} X$ and whenever $y \in_{w_{[n]}} X,$ then $y \cup_{w_{[n]},w} \{y\}_{w_{[n]},w} \in_{w_{[n]}} X.$ A set $X$ is denoted $\mathbb{N}_{w_{[n]}}^w$.
(iii) Weak Axiom of $w_{\{0\}}$-infinity.

$$\exists X\big[\exists y(y \in_{w_{\{0\}}} X) \wedge \forall z\big(z \notin_{w_{\{0\}}}^w y\big) \wedge \forall y(y \in_{w_{\{0\}}} X \Rightarrow_w \\ \Rightarrow_s \exists z(z \in_{w_{\{0\}}} X \wedge \forall t(t \in_{w_{\{0\}}} z \Leftrightarrow_w (t \in_{w_{\{0\}}} y \vee t =_{w_{\{0\}}} y))))] \quad (4.3.)$$

There is a set $X$ such that $\varnothing_{w_{\{0\}}}^w \in_{w_{\{0\}}} X$ and whenever $y \in_{w_{\{0\}}} X,$ then

$y \cup_{w_{\{0\}},w} \{y\}_{w_{\{0\}},w} \in_{w_{\{0\}}} X.$ A set $X$ is denoted $\mathbb{N}_{w_{\{0\}}}^w$.
(iv) Weak Axiom of $w_{\{n\}}$-infinity.

$$\exists X\big[\exists y(y \in_{w_{\{n\}}} X) \wedge \forall z\big(z \notin_{w_{\{n\}}}^w y\big) \wedge \forall y(y \in_{w_{\{n\}}} X \Rightarrow_w \\ \Rightarrow_w \exists z(z \in_{w_{\{n\}}} X \wedge \forall t(t \in_{w_{\{n\}}} z \Leftrightarrow_w (t \in_{w_{\{n\}}} y \vee t =_{w_{\{n\}}} y))))] \quad (4.3.)$$

There is a set $X$ such that $\varnothing_{w_{\{n\}}}^w \in_{w_{\{n\}}} X$ and whenever $y \in_{w_{\{n\}}} X,$ then $y \cup_{w_{\{n\}},w} \{y\}_{w_{\{n\}},w} \in_{w_{\{n\}}} X.$ A set $X$ is denoted $\mathbb{N}_{w_{\{n\}}}^w$.

# IV.4. $w$-Inconsistent Relations and Functions

# IV.4.1. $w$-Consistent Relations and Functions

**Definition 4.4.1.** An $w$-consistent ordered pair $(a,b)_{s,w}$ (or $s$-$w$-ordered pair ) is defined to
be

$$(a,b)_{s,w} = \big\{\{a\}_{s,w}, \{a,b\}_{s,w}\big\}_{s,w}. \quad (4.4.1)$$

Similarly we define

$$(a,b,c)_{s,w} = ((a,b),c)_{s,w} = \Big\{\{a\}_{s,w}, \{a,b\}_{s,w}, \Big\{\big\{\{a\}_{s,w}, \{a,b\}_{s,w}\big\}_{s,w}, c\Big\}_{s,w}\Big\}_{s,w}, \quad (4.4.2)$$

etc.
**Definition 4.4.2**. Let $R_w^s$ be an $w$-consistent set. An $w$-consistent set $R_w^s$ is a binary $w$-consistent relation (or $s$-$w$-relation) if all $w$-elements of $R_w^s$ are $w$-consistent ordered pairs, i.e. for $z \in_w^s R_w^s$ there exists $x$ and $y$ such that $z =_w^s (x,y)_{s,w}$. We can also denote $(x,y)_{s,w} \in_w^s R_w^s$ as $xR_w^s y$, and say that $x$ is in $s$-$w$-relation $R_w^s$ with $y$ if $xR_w^s y$ holds.
**Designation 4.4.1**.

# IV.5. $w$-Inconsistent $w$-Equivalences and $w$-Orderings

**Definition 4.5.1**. Let $(\bullet R^s_w \bullet)$ be a binary $w$-consistent $w$-relation in $A$.
(i) $R^s_w$ is $s$-$w$-reflexive (or strongly $w$-reflexive) in $A$ if for all $a \in^s_w A,\ aR^s_w a$.
(ii) $R^s_w$ is $s$-$w$-symmetric (or strongly $w$-symmetric) in $A$ if for all $a,b \in^s_w A$ :
$aR^s_w b \Leftrightarrow_s bR^s_w a$.
(iii) $R^s_w$ is $s$-$w$-antisymmetric (or strongly $w$-antysymmetric) in $A$ if
for all $a,b \in^s_w A : aR^s_w b \wedge bR^s_w a \Rightarrow_s a =^s_w b$.
(iv) $R^s_w$ is $s$-$w$-asymmetric (or strongly $w$-asymmetric) in $A$ if
for all $a,b \in^s_w A : aR^s_w b \Rightarrow_s \neg_s bR^s_w a$,i.e. $aR^s_w b$ and $bR^s_w a$ cannot both be true.
(v) $R^s_w$ is $s$-$w$-transitive (or strongly $w$-transitive) in $A$ if
for all $a,b,c \in^s_w A : aR^s_w b \wedge bR^s_w c \Rightarrow_s aR^s_w c$.

**Definition 5.1.1**.An $w$-consistent (or strong) $w$-ordering $<^s_w$ of $A$ is called $s$-$w$-linear or $s$-$w$-total if any two $w$-elements of $A$ are comparable in the ordering $<^s_w$ ;i.e. for any $a,b \in^s_w A$, either $a <^s_w b, b <^s_w a$, or $a =^s_w b$. The pair $(A, <^s_w)$ is called a $s$-$w$-linearly $w$-ordered set.

**Definition 5.1**.The condition that $X \subset^s_w A$ has a strong $<^s_w$ -least element reads

$$\exists x(x \in^s_w X)[\forall y \in^s_w X(x \leq^s_w y)] \tag{4.5.1}$$

or in the following equivalent form

$$\exists x(x \in^s_w X)[\forall y \in^s_w X \neg_s (y <^s_w x)]. \tag{4.5.2}$$

**Definition 5.1.1**.An $w$-inconsistent (or weak) $w$-ordering $<^w_w$ of $A$ is called

# V.The $w$-inconsistent natural numbers

## V.1.The $w$-consistent natural numbers

In defining the $w$-consistent natural numbers we begin by examining the most fundamental set, the strong $w$-empty set $\varnothing^s_w$.We can very easily create a pattern that is a prime candidate for the definition of the $w$-consistent natural numbers:
$w$-empty set $\varnothing^s_w$ has zero elements in the $w$-consistent sense;
$\{\varnothing^s_w\}$ has one element in the $w$-consistent sense;
$\{\varnothing^s_w, \{\varnothing^s_w\}\}$ has two elements in the $w$-consistent sense,etc.
Revisiting our prime candidate for $w$-consistent natural numbers, we can revise it as:
$0^s_w =^s_w \varnothing^s_w$;
$1^s_w =^s_w \{0^s_w\} =^s_w \varnothing^s_w \cup_{w,s} \{\varnothing^s_w\} =^s_w \{\varnothing^s_w\}$;

We see that each number is de ned based on the number that precedes it. This sequence is anchored by $0^s_w$. As long as $0^s_w$ is defined, then $1^s_w$ can be defined. Once $1^s_w$ is defined, $2^s_w$ can also be, and so on. This brings us to the concept of $w$-consistent induction.

**Definition 5.1.1**. The $w$-consistent $w$-successor (or strong $w$-successor) of a set $x$ is the
set

$$\mathbf{S}^s_w(x) =^s_w x \cup_{w,s} \{x\}_{w,s}. \tag{5.1.1}$$

**Definition 5.1**..A set $\mathbf{I}^s_w$ is called $s$-$w$-inductive (or strongly $w$-inductive) if

(i) $\varnothing^s_w \in^s_w \mathbf{I}^s_w$.

(b) If $x \in^s_w \mathbf{I}^s_w$, then $\mathbf{S}^s_w(x) \in^s_w \mathbf{I}^s_w$.

**Definition 5.1**... The $w$-consistent set of all $w$-consistent natural numbers is defined by

**Definition 5.1**..Let $A_{\setminus^s_w}$ by a set every nonempty $w$-subset $X$ of $A_{\setminus^s_w}$ has $w$-consistent (strong) $w$-complement $A\setminus^s_w X$. Any $w$-consistent $s$-$w$-linear $w$-ordering $\prec^s_w$ of a set $A$ is a $w$-consistent well-ordering if every nonempty $w$-subset $X$ of $A_{\setminus^s_w}$ has a $\prec^s_w$ -least element. The structure $(A_{\setminus^s_w}, \prec^s_w)$ is called $w$-consistent well-ordered set.

## V.2.The $w$-inconsistent natural numbers

# V. The Standard and Non-Standard Models of formal Paraconsistent theories.

## V.1. Generalized Incompleteness Theorems.

Let **Th** be some fixed, but unspecified, paraconsistent, i.e. inconsistent but nontrivial formal theory and in these case we wrote **PTh** or $Pcon(\mathbf{PTh})$ instead **Th**. For later convenience, we assume that the encoding is done in some fixed consistent formal theory **S** and that **PTh** contains **S**. We do not specify **S** — it is usually taken to be a formal system of arithmetic, although a weak set theory is often more convenient. The sense in which **S** is contained in **PTh** is better exemplified than explained: If **S** is a formal system of arithmetic and **PTh** is, say, $\widetilde{\mathbf{ZF}}_n, 1 \leqslant n < \omega$ or $\mathbf{ZFC}^\#$ then **PTh** contains **S** in the sense that there is a well-known embedding, or interpretation, of **S** in **PTh**.Since encoding is to take place in **S**, it will have to have a large supply of constants and closed terms to be used as codes. (E.g. in formal arithmetic, one has $\bar{0},\bar{1},\ldots$ .) **S** will also have certain function symbols to be described shortly.

To each formula $\Phi$, of the language of **PTh** is assigned a closed term, $[\Phi]^c$, called the code of $\Phi$. [N.B. If $\Phi(x)$ is a formula with free variable $x$, then $[\Phi(x)]^c$ is a closed term encoding the formula $\Phi(x)$ with $x$ viewed as a syntactic object and not as a parameter.] Corresponding to the logical connectives and quantifiers are function symbols, $\mathbf{neg}(\cdot)$, $\mathbf{imp}(\cdot)$, etc., such that, for all formulae $\Phi,\Psi$ : **S** $\vdash \mathbf{neg}_{(n)}([\Phi]^c) = [\neg_{(n)}\Phi]^c$,
**S** $\vdash \mathbf{imp}([\Phi]^c,[\Psi]^c) = [\Phi \to \Psi]^c$, etc. Of particular importance is the substitution operator $sub$, represented by the function symbol $\mathbf{sub}(\cdot,\cdot)$. For formulae $\Phi(x)$, terms $t$ with codes $[t]^c$ :

$$\mathbf{S} \vdash \mathbf{sub}([\Phi(x)]^c,[t]^c) = [\Phi(t)]^c. \tag{5.1.1}$$

Iteration of the substitution operator $sub$ allows one to define function symbols $\mathbf{sub}_3$, $\mathbf{sub}_4,\ldots$, such that

$$\mathbf{S} \vdash \mathbf{sub}_n([\Phi(x_1,x_2,\ldots,x_n)]^c,[t_1]^c,[t_2]^c,\ldots,[t_n]^c) = [\Phi(t_1,t_2,\ldots,t_n)]^c. \tag{5.1.2}$$

It well known [17] that one can also encode derivations and have a binary relation $\mathbf{Prov_{Th}}(x,y)$ (read "$x$ proves $y$ " or "$x$ is a proof of $y$") such that for closed $t_1,t_2$ :

**S** $\vdash \mathbf{Prov_{Th}}(t_1,t_2)$ iff $t_1$ is the code of a derivation in **PTh** of the formula with code $t_2$.It follows that

$$\mathbf{PTh} \vdash \Phi \text{ iff } \mathbf{S} \vdash \mathbf{Prov_{PTh}}(t,[\Phi]^c) \tag{5.1.3}$$

for some closed term $t$.

**Definition 5.1.1.**Thus one can define

$$\mathbf{Pr_{PTh}}(y) \leftrightarrow \exists x \mathbf{Prov_{PTh}}(x,y), \tag{5.1.4}$$

and therefore one obtain a predicate asserting provability.

**Remark 5.1.1**. We note that it is not always the case that :

$$\mathbf{PTh} \vdash \Phi \text{ iff } \mathbf{S} \vdash \mathbf{Pr_{PTh}}([\Phi]^c). \tag{5.1.5}$$

It well known [17] that the above encoding can be carried out in such a way that the following important conditions $\mathbf{D}1,\mathbf{D}2$ and $\mathbf{D}3$ are met for all sentences:

$$\mathbf{D}1.\mathbf{PTh} \vdash \Phi \text{ implies } \mathbf{S} \vdash \mathbf{Pr}_{\mathbf{PTh}}([\Phi]^c),$$

$$\mathbf{D}2.\mathbf{S} \vdash \mathbf{Pr}_{\mathbf{PTh}}([\Phi]^c) \to \mathbf{Pr}_{\mathbf{PTh}}([\mathbf{Pr}_{\mathbf{PTh}}([\Phi]^c)]^c), \tag{5.1.6}$$

$$\mathbf{D}3.\mathbf{S} \vdash \mathbf{Pr}_{\mathbf{PTh}}([\Phi]^c) \wedge \mathbf{Pr}_{\mathbf{PTh}}([\Phi \to \Psi]^c) \to \mathbf{Pr}_{\mathbf{PTh}}([\Psi]^c).$$

Generalized Incompleteness Theorems depend on the following.

**Theorem 5.1.1**. (Diagonalization Lemma). Let $\Phi(x)$ in the language of $\mathbf{PTh}$ have only the free variable indicated. Then there is a sentence $\psi$ such that

$$\mathbf{S} \vdash \psi \leftrightarrow \Phi([\psi]^c). \tag{5.1.7}$$

**Proof**. Given $\Phi(x)$, let $\vartheta(x) \leftrightarrow \Phi(\mathbf{sub}(x,x))$ be the diagonalization of $\Phi(x)$. Let $m = [\vartheta(x)]^c$ and $\psi = \vartheta(m)$. Then we claim that $\mathbf{S} \vdash \psi \leftrightarrow \Phi([\psi]^c)$. For $\Phi(x)$ in **S**, we see that

$$\psi \leftrightarrow \vartheta(m) \leftrightarrow \Phi(\mathbf{sub}(m,m)) \leftrightarrow \Phi(\mathbf{sub}([\vartheta(x)]^c,m)) \leftrightarrow \Phi([\vartheta(m)]^c) \leftrightarrow \Phi([\psi]^c). \tag{5.1.8}$$

We apply now (5.1.7) to $\neg_{(n)}\mathbf{Pr}_{\mathbf{Th}}(x)$.

**Theorem 5.1.2**. (Generalized First Incompleteness Theorem).Let (1) $Pcon_{(n)}(\mathbf{Th})$ and (2) $\mathbf{Th} \vdash \phi \leftrightarrow \neg_{(n)}\mathbf{Pr}_{\mathbf{Th}}([\phi]^c)$.

Then (**i**)

$$\mathbf{Th} \nvdash \phi, \tag{5.1.9}$$

(ii) under an additional assumption

$$\mathbf{Th} \nvdash \neg_{(n)}\phi. \tag{5.1.10}$$

**Proof**. (i) Observe $\mathbf{Th} \vdash \phi$ implies $\mathbf{Th} \vdash \mathbf{Pr}_{\mathbf{Th}}([\phi]^c)$ by $\mathbf{D}1$, which implies $\mathbf{Th} \vdash \neg_{(n)}\phi$, contradicting the paraconsistency of $\mathbf{Th}$.

(ii) The additional assumption is a strengthening of the converse to $\mathbf{D}1$, namely $\mathbf{Th} \vdash \mathbf{Pr}_{\mathbf{Th}}([\phi]^c)$ implies $\mathbf{Th} \vdash \phi$. We have $\mathbf{Th} \vdash \neg_{(n)}\phi$, hence $\mathbf{Th} \vdash \neg_{(n)}\neg_{(n)}\mathbf{Pr}_{\mathbf{Th}}([\phi]^c)$ so that $\mathbf{Th} \vdash \mathbf{Pr}_{\mathbf{Th}}([\phi]^c)$ and, by the additional assumption, $\mathbf{Th} \vdash \phi$, again contradicting the paraconsistency of $\mathbf{Th}$.

**Theorem 5.1.3**. (Generalized Second Incompleteness Theorem).

Let $Pcon_{(n)}(\mathbf{Th})$ be $\neg_{(n)}\mathbf{Pr}_{\mathbf{Th}}([\Lambda_n]^c)$, where $\Lambda_n = \mathbf{A} \wedge \neg_{(n)}\mathbf{A}$ is any convenient

$n$-contradictory statement. Then

$$\mathbf{Th} \nvdash Pcon_{(n)}(\mathbf{Th}). \tag{5.1.11}$$

**Proof**. Let $\phi$ be as in the statement of Theorem 5.2.. We show: $\mathbf{S} \vdash \phi \leftrightarrow Pcon_{(n)}(\mathbf{Th})$. Observe that $\mathbf{S} \vdash \phi \to \neg_{(n)}\mathbf{Pr_{Th}}([\phi]^c)$ implies $\mathbf{S} \vdash \phi \to \neg_{(n)}\mathbf{Pr_{Th}}([\Lambda_n]^c)$, since $\mathbf{S} \vdash \phi \to \Lambda_n$ implies $\mathbf{S} \vdash \mathbf{Pr_{Th}}([\phi \to \Lambda_n]^c)$, by $\mathbf{D}1$, which implies $\mathbf{S} \vdash \mathbf{Pr_{Th}}([\Lambda_n]^c) \to \mathbf{Pr_{Th}}([\phi]^c)$, by $\mathbf{D}3$. But $\phi \to \neg_{(n)}\mathbf{Pr_{Th}}([\Lambda_n]^c)$ is just $\phi \to Pcon_{(n)}(\mathbf{Th})$ and we have proven half of the equivalence. Conversely, by $\mathbf{D}2$,$\mathbf{S} \vdash \mathbf{Pr_{Th}}([\phi]^c) \to \mathbf{Pr_{Th}}([\mathbf{Pr_{Th}}([\phi]^c)]^c)$ , which implies $\mathbf{S} \vdash \mathbf{Pr_{Th}}([\phi]^c) \to\to \mathbf{Pr_{Th}}([\neg_{(n)}\phi]^c)$, by $\mathbf{D}1,\mathbf{D}3$, since $\phi \to \neg_{(n)}\mathbf{Pr_{Th}}([\phi]^c)$. This yields $\mathbf{S} \vdash \mathbf{Pr_{Th}}([\phi \wedge \neg_{(n)}\phi]^c)$, by $\mathbf{D}1,\mathbf{D}3$, and logic, which implies $\mathbf{S} \vdash \mathbf{Pr_{Th}}([\phi]^c) \to \mathbf{Pr_{Th}}([\Lambda_n]^c)$ by $\mathbf{D}1,\mathbf{D}3$, and logic. By contraposition, $\mathbf{S} \vdash \neg_{(n)}\mathbf{Pr_{Th}}([\Lambda_n]^c) \to \neg_{(n)}\mathbf{Pr_{Th}}([\phi]^c)$, which is $\mathbf{S} \vdash Pcon_{(n)}(\mathbf{Th}) \to \phi$, by definitions.

**Theorem 5.1.4**. $\mathbf{S} \vdash Pcon_{(n)}(\mathbf{Th}) \to Pcon_{(n)}\big(\mathbf{Th} + \neg_{(n)}Pcon_{(n)}(\mathbf{Th})\big)$.

**Proof**.By the proof of Theorem 5.3, (i) $\mathbf{S} \vdash Pcon_{(n)}(\mathbf{Th}) \to \neg_{(n)}\mathbf{Pr_{Th}}([\phi]^c)$,
(ii) $\mathbf{S} \vdash Pcon_{(n)}(\mathbf{Th}) \leftrightarrow \phi$. Using now $\mathbf{D}2,\mathbf{D}3$, it follows that
$\mathbf{S} \vdash Pcon_{(n)}(\mathbf{Th}) \to \neg_{(n)}\mathbf{Pr_{Th}}([Pcon_{(n)}(\mathbf{Th})]^c)$, so that

$$\mathbf{S} \vdash Pcon_{(n)}(\mathbf{Th}) \to \neg_{(n)}\mathbf{Pr_{Th}}([\neg_{(n)}Pcon_{(n)}(\mathbf{Th}) \to \Lambda_n]^c) \tag{5.1.12}$$

which gives $\mathbf{S} \vdash Pcon_{(n)}(\mathbf{Th}) \to Pcon_{(n)}\big(\mathbf{Th} + \neg_{(n)}Pcon_{(n)}(\mathbf{Th})\big)$.

**Definition 5.2**.Define: (i)

$$\begin{gathered}\mathbf{Prov}^{\Re}_{\mathbf{Th}}(x,y) \leftrightarrow \mathbf{Prov}_{\mathbf{Th}}(x,y) \wedge \\ \wedge \forall z(w \leq x)[\mathbf{Prov}_{\mathbf{Th}}(z,w) \to y \neq neg_{(n)}(w) \wedge w \neq neg_{(n)}(y)]\end{gathered} \tag{5.1.13}$$

(ii)

$$\mathbf{Pr}^{\Re}_{\mathbf{Th}}(y) \leftrightarrow \exists x \mathbf{Prov}^{\Re}_{\mathbf{Th}}(x,y) \tag{5.1.14}$$

and
(iii)

$$Pcon^{\Re}_{(n)}(\mathbf{Th}) \leftrightarrow \mathbf{Pr}^{\Re}_{\mathbf{Th}}([\Lambda_n]^c). \tag{5.1.15}$$

**Theorem 5.1.5**. (Generalized Rossers Theorem).Let (1) $Pcon_{(n)}(\mathbf{Th})$ and (2) $\mathbf{Th} \vdash \phi \leftrightarrow \neg_{(n)}\mathbf{Pr}^{\Re}_{\mathbf{Th}}([\phi]^c)$.
Then
(i)

$$\mathbf{Th} \nvdash \phi, \tag{5.1.16}$$

(ii)

$$\mathbf{Th} \nvdash \neg_{(n)}\phi. \tag{5.1.17}$$

(iii)

$$\mathbf{Th} \vdash Pcon^{\Re}_{(n)}(\mathbf{Th}). \tag{5.1.18}$$

**Proof**.(i) By the paraconsistency of $\mathbf{Th}$, $\mathbf{Prov_{Th}}$ and $\mathbf{Prov}^{\Re}_{\mathbf{Th}}$ binumerate the same relation. Hence $\mathbf{D}1^{\Re}$ holds: $\mathbf{Th} \vdash \phi \Rightarrow \mathbf{Th} \vdash \mathbf{Pr}^{\Re}_{\mathbf{Th}}([\phi]^c)$. Thus, the proof of the first part of the First Incompleteness Theorem yields the result.

(ii) This follows from (iii).

(iii) Follows immediately from the remarks that $\mathbf{Th}$ is paraconsistent and $\mathbf{Th} \vdash \neg_{(n)}\Lambda_n$.

**Theorem 5.1.6**. (Generalized Löb's Theorem). Let be (1) $Pcon_{(n)}(\mathbf{Th})$ and (2) $\phi$ be closed. Then

$$\mathbf{Th} \vdash \mathbf{Pr_{Th}}([\phi]^c) \to \phi \text{ iff } \mathbf{Th} \vdash \phi. \tag{5.1.19}$$

**Proof**. The one direction is obvious. For the other, assume that $\mathbf{Th} \nvdash \phi$. Then $\mathbf{Th} + \neg_{(n)}\phi$ is consistent and we may appeal to the Generalized Second Incompleteness Theorem to conclude that $\mathbf{Th} + \neg_{(n)}\phi$ does not yield $Pcon_{(n)}(\mathbf{Th} + \neg_{(n)}\phi)$, hence not $\neg_{(n)}\mathbf{Pr_{Th}}([\phi \to \Lambda_n]^c)$. Thus $\mathbf{Th} + \neg_{(n)}\phi \nvdash \neg_{(n)}\mathbf{Pr_{Th}}([\phi]^c)$. Contraposition yields $\mathbf{Th} \nvdash \mathbf{Pr_{Th}}([\phi]^c) \to \phi$.

Let be $Pcon_{(n)}(\mathbf{Th})$. Now we focuses our attention on the following schemata:

(**I**) Generalized Local Reflection Principle $\mathbf{Rfn}(\mathbf{Th})$ :

$$\mathbf{Pr_{Th}}([\phi]^c) \to \phi, \ \phi \text{ closed.} \tag{5.1.20}$$

(**II**) Generalized First Uniform Reflection Principle $\mathbf{RFN}(\mathbf{Th})$ :

$$\forall x \mathbf{Pr}_{\mathbf{Th}}([\phi(x)]^c) \to \forall x \phi(x),\ \phi(x) \text{ has only } x \text{ free.} \tag{5.1.21}$$

(**III**) Generalized Second Uniform Reflection Principle $\mathbf{RFN}'(\mathbf{Th})$ :

$$\forall x[\mathbf{Pr}_{\mathbf{Th}}([\phi(x)]^c) \to x\phi(x)],\ \phi(x) \text{ has only } x \text{ free.} \tag{5.1.22}$$

**Theorem 5.1.7**. (Generalized First Incompleteness Theorem).Let be $Pcon_{(n)}(\mathbf{Th})$. Then for some true, unprovable $\phi$

$$\mathbf{Th} \nvdash \mathbf{Pr}_{\mathbf{Th}}([\phi]^c) \to \phi \tag{5.1.23}$$

**Theorem 5.1.8**. (Generalized Second Incompleteness Theorem).Let be $Pcon_{(n)}(\mathbf{Th})$. Then for any refutable $\phi$

$$\mathbf{Th} \nvdash \mathbf{Pr}_{\mathbf{Th}}([\phi]^c) \to \phi \tag{5.1.24}$$

Theorem 5.1.6 simply yields

$$\mathbf{Th} \nvdash \mathbf{Pr}_{\mathbf{Th}}([\phi]^c) \to \phi \text{ iff } \mathbf{Th} \nvdash \phi, \tag{5.1.25}$$

# V.2. The Generalized Compactness Theorem corresponding to paraconsistent first-order logic $\mathbf{LP}^{\#}_{\omega}$.

In order to use the Generalized Compactness Theorem, and in fact, even to state it, we must first develop the logical language to which it applies. In this case, we shall use paraconsistent first-order logic $\mathbf{LP}^{\#}_{\omega}$ with infinite hierarchy levels of contradiction. We will begin by listing the requisite definitions.

**Definition 5.2.1**. A language $\mathcal{L}$ is a not necessarily countable collection of relation symbols $P$, function symbols $G$, and constant symbols $c$.

**Definition 5.2.2**. The inconsistent universe $\mathbf{V}^{inc}$ is inconsistent universal set or any inconsistent set $A_{inc} \subset \mathbf{V}^{inc}$.

**Definition 5.2.3**. An interpretation function $\Im$ is a function such that:
(1) For each $n$-place relation symbol $P$ of $\mathcal{L}, \Im(P) = R$ where $R \subset A_{inc}^n$.
(2) For each $m$-place function symbol $G$ of $\mathcal{L}$, $\Im(G) = F$ where $F : A_{inc}^m \to A_{inc}$
is an $m$-place function on $A_{inc}$.
(3) For each constant symbol $c$, $\Im(c) = x$ for some $x \in A_{inc}$.
**Definition 5.2.4**. An inconsistent model $M_{inc}$ for any paraconsistent theory $Th^{inc}$ with a language $\mathcal{L}$ (paraconsistent $\mathcal{L}$-theory) consists of a universe $A_{inc}$ and an interpretation function $\Im$, which we denote by

$$M^{Th^{inc}} = M_{inc} = \langle A_{inc}, \mathcal{L}, \Im \rangle \tag{5.2.1}$$

or $M^{inc}$.
**Definition 5.2.5**. A term is one of four things:
(i) A variable is a term.
(ii) A constant symbol is a term.
(iii) If $F$ is an $m$-placed function symbol, and $t_1, \ldots, t_m$ are terms, then $F(t_1, \ldots, t_m)$ is a term.
(iv) A string of symbols is a term only if it can be shown to be a term by a finite number of applications of (i)-(iii).
**Remark.5.2.1.**The purpose of (iv) is to ensure that there are no infinite terms.
**Remark.5.2.2.**Now before we continue, we note that there is two two-place relation symbol which always belongs to rst-order logic, even though it does not belong to $\mathcal{L}$. This relation is called the identity relation, and is denoted by $=_s$ and $=_w$ .
**Definition 5.2.6**. An atomic formula of $\mathcal{L}$ is a string of the form:

(i) $t_1 =_s t_2$ where $t_1$ and $t_2$ are terms of $\mathcal{L}$.
(ii) $t_1 =_w t_2$ where $t_1$ and $t_2$ are terms of $\mathcal{L}$.
(iii) $P(t_1, \ldots, t_n)$ where $P$ is an $n$-placed relation and $t_1, \ldots, t_n$ are terms of $\mathcal{L}$.
**Definition 5.2.7**. A formula of $\mathcal{L}$ is defined as follows:
(i) An atomic formula is a formula.
(ii) If $\varphi$ and $\psi$ are formulas, then $\varphi \wedge \psi, \varphi \vee \psi, \neg_s\varphi, \neg_w\varphi, \varphi \Rightarrow \psi, \varphi \to \psi$ are formulas.
(iii) If $v$ is a variable and $\varphi$ is a formula, then $(\forall v)\varphi$ is a formula.
(iv) If a string of symbols can be shown to be a formula by a finite number of applications of (i)-(iii), then it is a formula.
**Definition 5.2.8**. A formula is a sentence if every variable in the formula is bound by the quantifier $\forall$ or $\exists$.
**Definition 5.2.9**. A sentence $\varphi$ is true in a model $M_{inc}$, or alternatively, $M_{inc}$ is a model of $\varphi$, denoted $M_{inc} \models \varphi$, if for every possible sequence of elements in $A_{inc}$, substituting these elements in $A_{inc}$ for the variables present in $\varphi$ yields a true sentence at last at inconsistent sence,i.e.,both $\varphi$ and $\neg_w\varphi$ holds in $M_{inc}$ for some $\varphi$.
**Remark.5.2.2.**Note that: (1) this idea of truth precludes the possibility of both $\varphi$ and $\neg_s\varphi$ holding in $A_{inc}$ but (2) the possibility of both $\varphi$ and $\neg_w\varphi$ holds in $M_{inc}$ for some $\varphi$.
**Definition 5.2.10**. We say that $M_{inc}$ is inconsistent model of a set of sentences $\Sigma$ if $M_{inc}$ is a model of $\varphi$ for all $\varphi \in \Sigma$.
**Definition 5.2.11**. A sentence $\varphi$ is a consequence of a set of sentences $\Sigma$, denoted

$\Sigma \models \varphi,$ if every model $M_{inc}$ of $\Sigma$ is a model of $\varphi$.
**Definition 5.2.12**. sentence $\varphi$ is deducible from $\Sigma,$ expressed $\Sigma \vdash_{\mathbf{RMP}} \varphi,$ if there exists a nite chain of sentences $\psi_0,\ldots,\psi_n$ where $\psi_n$ is $\varphi$ and each previous sentence in the chain either belongs to $\Sigma,$ follows from one of the axioms of $\mathbf{Th}^{inc} + \mathbf{LP}^{\#}_{\omega},$ or can be inferred from previous sentences.
**Definition 5.2.13**. A set of sentences $\Sigma$ is paraconsistent if and only if (i) there does not exist a sentence $\varphi$ such that $\Sigma \vdash_{\mathbf{RMP}} \varphi$ and $\Sigma \vdash_{\mathbf{RMP}} \neg_{\mathbf{s}}\varphi,$(ii) there exist at least one sentence $\varphi$ such that $\Sigma_{\mathbf{RMP}} \nvdash \varphi$.
**Lemma 5.2.1**. Let $\Sigma$ be a set of sentences. If $\Sigma \vdash_{\mathbf{RMP}} \varphi,$ then $\Sigma \models \varphi$.
**Theorem 5.2.1**.(**Generalized Soundness Theorem**). Let $\Sigma$ be a set of formulas. If $\Sigma$

has an inconsistent model $M_{inc},$ then $\Sigma$ is paraconsistent.
**Proof**. Suppose that is not paraconsistent. Then there exists some $\varphi$ such that:
(i) $\Sigma \vdash_{\mathbf{RMP}} \varphi$ and $\Sigma \vdash_{\mathbf{RMP}} \neg_{\mathbf{s}}\varphi,$ or (ii) $\Sigma \vdash_{\mathbf{RMP}} \varphi$ for all $\varphi$.
**1**. From assumption (i) by Lemma 5.2.1, we have $\Sigma \models \varphi$ and $\Sigma \models \neg_{\mathbf{s}}\varphi.$This means that $\Sigma$ cannot have a model $M_{inc}$, because if it did, then $\Sigma \models \varphi$ and $\Sigma \models \neg_{\mathbf{s}}\varphi,$ which is impossible.
**2**. From assumption (i) by Lemma 5.2.1, we have $\Sigma \models \varphi$ for all $\varphi.$This means that $\Sigma$ cannot have a model $M_{inc},$ because if it did, then $\Sigma \vdash_{\mathbf{RMP}} \varphi$ for all $\varphi,$ which is impossible.Thus, we have proved our statement.
**Theorem 5.2.1**.(**Generalized Godel's Completeness Theorem**). Let $\Sigma$ be a set of formulas. If $\Sigma$ is paraconsistent, then it has an inconsistent model $M_{inc}$.
**Proof**. Let $\Sigma$ be an arbitrary paraconsistent set of sentences of some language $\mathcal{L}$. Let $\overline{\mathcal{L}}$ be an expansion of $\mathcal{L}$ created by adding a set of new constant symbols not in $\mathcal{L}$ that has the same cardinality as $\mathcal{L}.$ The firrst step is to add sentences to $\Sigma$ to create a paraconsistent set of sentences $\overline{\Sigma}$ in the language $\overline{\mathcal{L}}.$ It is possible to show by canonical way, that $\overline{\Sigma}$ has a model $M_{inc}$ which is a model for $\overline{\mathcal{L}}.$ Now if we let $\mathfrak{R}_{inc}$ be the reduction of $M_{inc}$ to only involve the original language $\mathcal{L},$ it is possible to show by canonical way, that $\mathfrak{R}_{inc}$ is a model for $\Sigma,$ because the sentences in do not involve any constants which belonged to $\overline{\mathcal{L}},$so the reduction of $M_{inc}$ to $\mathfrak{R}_{inc}$ did not affect its ability to model $\Sigma$.
**Theorem 5.2.2**.(**Generalized Compactness Theorem**). Any paraconsistent set of sentences $\Sigma$ has a model $M_{inc}$ if and only if every finite subset of $\Sigma$ has a model.

## V.3. The Non-Standard Models of Paraconsistent second order arithmetic $Z_2^{\#}$.

This subsection presents the terminology and results necessary to prove the existence of non-standard models of paraconsistent second order arithmetic $Z_2^{\#}$ and paraconsistent Peano arithmetic $PA_{inc} \subset Z_2^{\#}$ [24] and that there are $2^{\aleph_0}$ such countable models.
**Definition 5.3.1**.Any paraconsistent theory $Th^{inc}$ with a language $\mathcal{L}$ is called paraconsistent $\mathcal{L}$-theory.
**Definition 5.3.2**.Let $M_{inc}$ be an inconsistent model $M_{inc} = \langle A_{inc}, \mathcal{L}, \mathfrak{I}\rangle$.
(i) The set $A_{inc}$ of all elements of $M_{inc},$ called the domain of $M_{inc}$ and named by $\mathbf{dom}(M_{inc})$.

(ii) For a constant symbol $c$ the constant element $\Im(c)$ named by $c^{M_{inc}}$.
(iii) For a relation symbol $R$, the relation $\Im(R)$ named by $R^{M_{inc}}$.

(iv) For a function symbol $F$, the function $\Im(F)$ named by $F^{M_{inc}}$.

**Notation5.3.1**. Any inconsistent model $M_{inc} = \langle A_{inc}, \mathcal{L}, \Im \rangle$ we often call as $\mathcal{L}$-model.

**Definition 5.3.3**.Let $Th^{inc}$ be an paraconsistent $\mathcal{L}$-theory and let $\varphi \in \mathcal{L}$. If $M_{inc} \vDash \varphi$ for all $\varphi \in Th^{inc}$,then $M_{inc}$ is inconsistent model for $Th^{inc}$, written $M_{inc} \vDash Th^{inc}$ or $M_{Th^{inc}}$.

**Definition 5.3.4**.The second-order language of paraconsistent second order arithmetic $Z_2^{\#}$ named by $\mathcal{L}_2^{\#}$.

**Definition 5.3.5**. Inconsistent model for paraconsistent second order arithmetic $Z_2^{\#}$ is an inconsistent model $M_{inc} \vDash Z_2^{\#}$ such that $M_{inc} \vDash \varphi$ for all $\varphi \in Z_2^{\#}$.

**Definition 5.3.6**. Any inconsistent model $M_{inc} = \langle A_{inc}, \mathcal{L}, \Im \rangle$ of paraconsistent $\mathcal{L}$-theory $Th^{inc}$ is colled inconsistent $\mathcal{L}$-model.

**Definition 5.3.7**. The signature $S_{M_{inc}}$ of an inconsistent model $M_{inc}$ of paraconsistent $\mathcal{L}$-theory $Th^{inc}$ lists the set of functions, relations and constants of that inconsistent model.

**Definition 5.3.8**.Let $S$ be a signature and let $M_1^{inc}$ and $M_2^{inc}$ be $\mathcal{L}$-models with signature $S = S_{M_1^{inc}} = S_{M_2^{inc}}$.A homomorphism $f : M_1^{inc} \to M_2^{inc}$, is a function $f$ from $\mathbf{dom}(M_1^{inc})$ to $\mathbf{dom}(M_2^{inc})$ such that:

1. For each constant $c$ of $S$, $f\left(c^{M_1^{inc}}\right) = c^{M_2^{inc}}$.
2. For each $n > 0, n$-ary relation symbol $R$ of $S$ and $n$-tuple $\bar{a} \in M_1^{inc}$, if $\bar{a} \in R^{M_1^{inc}}$ then $f(\bar{a}) \in R^{M_2^{inc}}$.
3. For each $n > 0, n$-ary function symbol $F$ of S and $n$-tuple $\bar{a} \in R^{M_1^{inc}}$, $f\left(F^{M_1^{inc}}(\bar{a})\right) = F^{M_2^{inc}}(f(\bar{a}))$.

Where $\bar{a} = (a_0,\ldots,a_{n-1})$ and $f(\bar{a}) = (f(a_0),\ldots,f(a_{n-1}))$.

**Definition 5.3.9**.An embedding of $M_1^{inc}$ into $M_2^{inc}$ is a homomorphism $f : M_1^{inc} \to M_2^{inc}$ which is injective and satisfies:

1. For each $n > 0$, each $n$-ary relation symbol $R$ of $S$ and each $n$-tuple $\bar{a} \in M_1^{inc}$, $\bar{a} \in R^{M_1^{inc}} \leftrightarrow f(\bar{a}) \in R^{M_2^{inc}}$.

Furthermore, $M_1^{inc}$ and $M_2^{inc}$ are isomorphic, written $M_1^{inc} \cong M_2^{inc}$, when there exists a surjective embedding $f : M_1^{inc} \to M_2^{inc}$.

The importent question is recapitulated formally as motivation for classifying the properties of standard and non-standard models of inconsistent arithmetic $Z_2^{\#}$ [24].

**Problem**. Given the standard model $\mathbb{N}_{inc}$ of $Z_2^{\#}$, if $M^{inc} \vDash Z_2^{\#}$ and $M^{inc} \ncong \mathbb{N}_{inc}$, then how many such countable models $M^{inc}$ are there and how do they differ from $\mathbb{N}_{inc}$?

The answer to this question requires the theory of inconsistent non-standard models and will be answered by Theorem 5.3.1 and Theorem 5.3.3 below .

**Definition 5.3.10**. Let $\mathbf{A}(Z_2^{\#})$ be a set of the all axioms of $Z_2^{\#}$. A non-standard model $M^{inc}$ of $Z_2^{\#}$ is an $\mathcal{L}_2^{\#}$-model such that $M^{inc} \vDash_w \phi$, for all $\phi \in \mathbf{Ax}(Z_2^{\#})$, and $M^{inc} \ncong \mathbb{N}_{inc}$, where $\mathbb{N}_{inc}$ is the standard model of $Z_2^{\#}$.

That is to say, an model $M^{inc}$ of $Z_2^{\#}$ or $PA_{inc} \subset Z_2^{\#}$ is non-standard when there does not exist a surjective weakly embedding $f : \mathbb{N}_{inc} \rightarrow_{w^*} M^{inc}$. Unpacking the definitions, this means that for any homomorphism $f : \mathbb{N}_{inc} \rightarrow M^{inc}$, either there exists a constant symbol, relation, or function which is not mapped to (i.e., $f$ is not a bijection), or the condition for an embedding - that $\bar{a} \in R^{\mathbb{N}_{inc}} \Leftrightarrow \bar{a} \in R^{M^{inc}}$ - does not holds. The explicit construction of such a non-standard model $M^{inc}$ of $Z_2^{\#}$ will require a connection between the satisfiability of a theory and some new constant symbol which ensures that $M^{inc} \not\cong \mathbb{N}_{inc}$.

**Theorem 5.3.1**. (Generalized Gödel's Completeness Theorem) Let $Th^{inc}$ be an paraconsistent $\mathcal{L}$-theory and let $\varphi \in \mathcal{L}$, where $\phi$ is an $\mathcal{L}$-sentence. Then $M_{inc} \models_w \phi$ if and only if $Th^{inc} \vdash \phi$.

**Corollary 5.3.1**. $Th^{inc}$ is paraconsistent if and only if $Th^{inc}$ is satisfiable.

**Proof**. Assume to the contrary that there exists a theory $Th^{inc}$ such that $Th^{inc}$ is paraconsistent and $Th^{inc}$ is not satisfiable. Since $Th^{inc}$ is not satisfiable, there does not exist a model $M_{inc}$ of $Th^{inc}$. So, any model $M_{inc}$ of $Th^{inc}$ is a model of $\bot_{\mathbf{s}}$ ($\bot_{\mathbf{s}} \triangleq \alpha \wedge \neg_{\mathbf{s}}\alpha$). Then, $Th^{inc} \models \bot_{\mathbf{s}}$ and so by the Completeness Theorem, $Th^{inc} \vdash_{\mathbf{RMP}} \bot_{\mathbf{s}}$ ; yet this contradicts the assumption that $Th^{inc}$ is paraconsistent. Assume to the contrary that there exists a theory $Th^{inc}$ such that $Th^{inc}$ is satisfiable and $Th^{inc}$ is not paraconsistent; this is an immediate contradiction by the definition of satisfiability. Therefore,$Th^{inc}$ is paraconsistent if and only if $Th^{inc}$ is satisfiable.

**Theorem 5.3.2**.(Generalized Compactness Theorem) $Th^{inc}$ is satisfiable if and only if every finite subset of $Th^{inc}$ is satisfiable.

**Theorem 5.3.3**.There exists inconsistent non-standard models of $\mathbf{PA}_{inc}$.

Proof. We want to prove that there exists a model $M_{inc} \triangleq M_{\mathbf{PA}_{inc}}$ for $\mathbf{PA}_{inc}$ which is not isomorphic to the standard model $\mathbb{N}_{inc}$. Let $\bar{n}_w$ be the value of the weak numeral $n_w$ formed by

$\bar{n}_{w_1} =_{w_1} \underbrace{1_{w_1} +_{w_1} \ldots +_{w_1} 1_{w_1}}_{n_{w_1}\ 1's}$ , and let $c$ be a new constant symbol such that $\Im(c) =_{w_1} c_{w_1}$ where $c_{w_1}$ is an $w_1$-inconsistent object, i.e. $w_1\text{-}inc(c_{w_1})$. Then we set:

$$Th^{inc}_{k_{w_1}} = \{\mathbf{Ax}(\mathbf{PA}_{inc})\} \cup \{\neg_w(c_{w_1} =_{w_1-inc} \bar{n}_{w_1} | \bar{n}_{w_1-inc} <_{w_1-inc} k_{w_1})\} \tag{5.3.1}$$

be a set of axioms in the language $\mathcal{L}_2^{\#} \cup \{c_{w_1}\}$, where $n,k \in_{w\text{-}con} \mathbb{N}_{inc}$. For a given $k$, give the interpretation $c_{w_1}^{\mathbb{N}_{inc}} =_{w_1-inc} \bar{k}_{w_1}^{\mathbb{N}_{inc}}$ Then, since $\mathbf{PA}_{inc}$ is paraconsistent,$Th^{inc}_{k_{w_1-inc}}$ is paraconsistent, and thus satisfiable by Corollary 5.3.1, for each $k \in_{w\text{-}con} \mathbb{N}_{inc}$. Therefore, the standard model of $\mathbf{PA}_{inc}$ is a model for $Th^{inc}_{k_{w_1}}$; that is to say, $\mathbb{N}_{inc} \models Th^{inc}_{k_{w_1}}$. Since

$$Th^{inc}_{k^{(1)}_{w_1}} \subseteq Th^{inc}_{k^{(2)}_{w_1}} \subseteq \ldots \subseteq Th^{inc}_{k^{(i)}_{w_1}} \subseteq \ldots \tag{5.3.2}$$

and each $Th^{inc}_{k^{(i)}_{w_1}}$ is satisfiable for $i \in \mathbb{N}$. We set now

$$Th^{\#}_{inc} = \bigcup_{i\in\mathbb{N}} Th^{inc}_{k^{(i)}_{w_1}}. \qquad (5.3.3)$$

is satisfiable by the Generalized Compactness Theorem. So, there exists an $\mathcal{L}_2^{\#} \cup \{c_{w_1}\}$-model $M_{\mathbf{PA}_{inc}}$ such that $M_{\mathbf{PA}_{inc}} \models Th^{\#}_{inc}$ and thus $M_{\mathbf{PA}_{inc}} \models \mathbf{Ax}(\mathbf{PA}_{inc})$. Assume now to the contrary that $M_{\mathbf{PA}_{inc}} \cong \mathbb{N}_{inc}$, then there exists a surjective embedding $f : \mathbb{N}_{inc} \to M_{\mathbf{PA}_{inc}}$ and so $f(n) =_{w^*} n$, for all $n \in \mathbb{N}_{inc}$. But since $c_{w_1} \neq_{w^*} n$, for all $n \in \mathbb{N}_{inc}$, there does not exist an image in $M_{\mathbf{PA}_{inc}}$ of $c_{w_1}$ under $f$, which contradicts that $f$ is a $w^*$-surjective embedding. Therefore, $M_{\mathbf{PA}_{inc}}$ is a model for Peano arithmetic $\mathbf{PA}_{inc}$.

V.3

# VI.Paralogical Nonstandard Analysis.

## VI.1. The inconsistent ultrafilter.

### VI.1.1. The consistent ultrafilter.

We remind some classical definitions.

**Definition 6.1.1.** Let $\wp^{con}$ be an infinite classical set $\wp^{con} \in \mathbf{V}^{con}$. Any consistent filter $\mathcal{F}^{con}$ that is a family of subsets of $\wp^{con}$ satisfying the following properties:

(i) $\wp^{con} \in_s \mathcal{F}^{con}, \varnothing_s \notin_s \mathcal{F}^{con}$.

(ii) $A_1,\ldots,A_n \in_s \mathcal{F}^{con} \Rightarrow_s A_1 \cap_s \ldots \cap_s A_n \in_s \mathcal{F}^{con}$.

(iii) $A \in_s \mathcal{F}^{con}$ and $A \subset_s B \Rightarrow_s B \in_s \mathcal{F}^{con}$.

**Definition 6.1.2.** A filter $\mathcal{F}^{con}$ on $\wp^{con}$ is called free if it contains no finite set.

**Definition 6.1.3.** A filter $\mathcal{F}^{con}$ is called an consistent ultrafilter over $\wp^{con}$ if for all $E \subseteq_s \wp^{con}$ either $E \in_s \mathcal{F}^{con}$ or $\wp^{con}\backslash_s E \in_s \mathcal{F}^{con}$, i.e.,

$$E \in_s \mathcal{F}^{con} \Leftrightarrow_s \wp^{con}\backslash_s E \notin^s_s \mathcal{F}^{con}, \qquad (6.1.1)$$

where we abrraviate: $\wp^{con}\backslash_s E \notin^s_s \mathcal{F}^{con} \triangleq \neg_s(\wp^{con}\backslash_s E \in_s \mathcal{F}^{con})$.

**Remark.6.1.1.** Notice that from the nontriviality condition [Definition 6.1.1(i)] it follows that
if $\mathcal{F}^{con}$ is an ultrafilter on $\wp^{con}$ and $E \subseteq_s \wp^{con}$, then exactly one of the sets $E$ and $\wp^{con}\backslash_s E$
belongs $\mathcal{F}^{con}$.

### VI.1.2. The inconsistent $w$-ultrafilter.

**Definition 6.1.4.** Let $\wp^w$ be a infinite weakly inconsistent ($w$-inconsistent ) set $\wp^w \in \mathbf{V}^w$.

Any weakly inconsistent filter ($w$-filter) $\mathcal{F}^w$ that is $w_{[1]}$-inconsistent family of $w$-subsets of

$\wp^w$ satisfying the following properties:
(i) $\wp^w \in_w \mathcal{F}^w, \varnothing_w^s \notin_w^s \mathcal{F}^w$.
(ii.a) $A_1,\ldots,A_n \in_w \mathcal{F}^w \Rightarrow_s A_1 \cap_{w,s} \ldots \cap_{w,s} A_n \in_w \mathcal{F}^w$.
(ii.b) $A_1,\ldots,A_n \in_{w_{[1]}} \mathcal{F}^w \Rightarrow_s A_1 \cap_{w,s} \ldots \cap_{w,s} A_n \in_{w_{[1]}} \mathcal{F}^w$.
(iii.a) $A \in_w \mathcal{F}^w$ and $A \subset_w B \Rightarrow_s B \in_w \mathcal{F}^w$.
(iii.b) $A \in_{w_{[1]}} \mathcal{F}^w$ and $A \subset_w B \Rightarrow_s B \in_{w_{[1]}} \mathcal{F}^w$.

**Definition 6.1.5**.A $w$-filter $\mathcal{F}^w$ on $\wp^w$ is called free if it contains no finite set.
**Definition 6.1.6**.A $w$-filter $\mathcal{F}^w$ is called an $w$-ultrafilter over $\wp^w$ if for all $E \subseteq_w \wp^w$ either $E \in_w \mathcal{F}^w$ or $\wp^w \backslash_s E \in_w \mathcal{F}^w$, i.e.,

$$E \in_w \mathcal{F}^w \Leftrightarrow_s \wp^w \backslash_s E \notin_w^s \mathcal{F}^w, \qquad (6.1.2)$$

**Remark.6.1.2**.Notice that from the nontriviality condition [Definition 6.1.4(i)] it follows that
if $\mathcal{F}^w$ is an $w$-ultrafilter on $\wp^w$ and $E \subseteq_w \wp^w$,then exactly one of the sets $E$ and $\wp^w \backslash_s E$ $w$-belongs $\mathcal{F}^w$.

# VI.1.3. The weakly consistent $w_0$-ultrafilter.

**Definition 6.1.7**.Let $\wp^{w_0}$ be a infinite weakly consistent ($w_0$-consistent) set $\wp^{w_0} \in \mathbf{V}^{w_0}$.
Any weakly consistent filter ($w_0$-filter) $\mathcal{F}^{w_0}$ that is a family of $w_0$-subsets of $\wp^{w_0}$ satisfying
the following properties:
(i) $\wp^{w_0} \in \mathcal{F}^{w_0}, \varnothing_{w_0}^s \notin_{w_0}^s \mathcal{F}^{w_0}$.
(ii.a) $A_1,\ldots,A_n \in_{w_0} \mathcal{F}^{w_0} \Rightarrow_s A_1 \cap_{w_0} \ldots \cap_{w_0} A_n \in_{w_0} \mathcal{F}^{w_0}$,
(ii.b) $(A_1,\ldots,A_n \in_{w_0} \mathcal{F}^{w_0}) \wedge (A_1,\ldots,A_n \notin_{w_0}^w \mathcal{F}^{w_0}) \Rightarrow_s$
$(A_1 \cap_{w_0} \ldots \cap_{w_0} A_n \in_{w_0} \mathcal{F}^{w_0}) \wedge (A_1 \cap_{w_0} \ldots \cap_{w_0} A_n \notin_{w_0}^w \mathcal{F}^{w_0})$.
(iii.a) $A \in_{w_0} \mathcal{F}^{w_0}$ and $A \subset_{w_0} B \Rightarrow_s B \in_{w_0} \mathcal{F}^{w_0}$,
(iii.b) $(A \in_{w_0} \mathcal{F}^{w_0}) \wedge (A \in_{w_0}^w \mathcal{F}^{w_0})$ and $A \subset_{w_0} B \Rightarrow_s (B \in_{w_0} \mathcal{F}^{w_0}) \wedge_s (B \notin_{w_0}^w \mathcal{F}^{w_0})$.
**Definition 6.1.8**.Any $w_0$-filter $\mathcal{F}^{w_0}$ on $\wp^{w_0}$ is called free if it contains no finite set.
**Definition 6.1.9**.Weakly consistent filter $\mathcal{F}^{w_0}$ is called a $w_0$-ultrafilter over $\wp^{w_0}$ if for all $E \subseteq_{w_0} \wp^{w_0}$ either $E \in_{w_0} \mathcal{F}^{w_0}$ or $\wp^{w_0} \backslash_{w_0} E \in_{w_0} \mathcal{F}^{w_0}$.

$$E \in_{w_0} \mathcal{F}^{w_0} \Leftrightarrow_s \wp^{w_0} \backslash_s E \notin_w^s \mathcal{F}^{w_0}, \qquad (6.1.3)$$

**Remark.6.1.3**.Notice that from the nontriviality condition [Definition 6.1.7(i)] it follows that
if $\mathcal{F}^{w_0}$ is an ultrafilter on $\wp^{w_0}$ and $E \subseteq_{w_0} \wp^{w_0}$,then $E \in_{w_0} \mathcal{F}^{w_0}$ or $\wp^{w_0} \backslash E_{w_0} \in_{w_0} \mathcal{F}^{w_0}$.
We can now construct an $w$-inconsistent and $w_0$-consistent nonstandard extensions.

# VI.1.4 The $w$-inconsistent nonstandard extension.

**Definition 6.1.10**.Let be $\mathcal{F}^w$ a free $w$-ultrafilter on $\wp^w$ and introduce a strong $w$-equivalence relation $f^w \sim_{\mathcal{F}^w}^s g^w$ on $w$-sequences $f^w \in_w \mathbb{R}_w^{\wp^w}$ by

$$f^w \sim_{\mathcal{F}^w}^s g^w \Leftrightarrow_s \left\{ v \in_w \wp^w \mid f^w(v) =_w g^w(v) \right\}_w \in_w \mathcal{F}^w. \qquad (6.1.4)$$

**Remark.6.1.4**.Note that for any $f^w, g^w, h^w \in_w \mathbb{R}^{\wp^w}$

$$\left(f^w \sim_{\mathcal{F}^w} g^w\right) \wedge \left(g^w \sim_{\mathcal{F}^w} h^w\right) \Rightarrow_s f^w \sim_{\mathcal{F}^w} h^w. \qquad (6.1.5)$$

**Definition 6.1.11.** $\mathbb{R}_w^{\wp^w}$ "divided" out by the $w$-equivalence relation $\sim_{\mathcal{F}^{inc}}$ on clases $f_{\mathcal{F}^w}^w$ by
formula

$$f^w \in_w f_{\mathcal{F}^w}^w,$$
$$\forall g^w \left\{ g^w \in_w f_{\mathcal{F}^w}^w \Leftrightarrow_s f^w \sim_{\mathcal{F}^w} g^w \right\}, \quad (6.1.6)$$

gives us the inconsistent nonstandard extension ${}^{\#_w}\mathbb{R}$, the inconsistent hyperreals; in symbols,

$${}^{\#_w}\mathbb{R}_w \triangleq \mathbb{R}_w^{\wp^w} \big/ \sim_{\mathcal{F}^w}, \quad (6.1.7)$$

which mean a natural $w$-embedding:

$$f^w \mapsto_w f_{\mathcal{F}^w}^w \quad (6.1.8)$$

If $f^w \in_w \mathbb{R}_w^{\wp^{inc}}$, we denote its image in ${}^{\#_w}\mathbb{R}_w$ by $f_{\mathcal{F}^w}^w$, and, of course,every element in ${}^{\#_w}\mathbb{R}_w$
is of the form $f_{\mathcal{F}^{inc}}^w$, for some $f^w : \wp^w \to_w \mathbb{R}_w$.

**Remark.6.1.5.** Note that for any $f^w, g^w \in_w \mathbb{R}_w^{\wp^w}$

$$f^w \sim_{\mathcal{F}^w} g^w \Rightarrow_s f_{\mathcal{F}^w}^w =_w g_{\mathcal{F}^w}^w. \quad (6.1.9)$$

For any $w$-inconsistent real number $r_w \in_w \mathbb{R}_w$, such that $r_w =_w r, r \in_s \mathbb{R}$, let $\mathbf{r}^w$ denote the
constant $w$-function with value $r_w$ in $\mathbb{R}_w$, i.e., $\mathbf{r}^w(v) =_w r_w$, for all $v \in_w \wp^w$. We then have a
natural $w$-embedding:

$$\#_w : \mathbb{R}_w \hookrightarrow_w {}^{\#_w}\mathbb{R}_w \quad (6.1.10)$$

by setting ${}^{\#_w}r_w =_w \mathbf{r}_{\mathcal{F}^w}^w$, for all $r_w \in_w \mathbb{R}_w$. We must now lift the structure of $\mathbb{R}_w$ to the $w$-inconsistent hyperreals ($w$-hyperreals) ${}^{\#_w}\mathbb{R}_w$

**Remark.6.1.6.** Notice that as an algebraic $w$-inconsistent structure, $\mathbb{R}_w$ is a $w$-complete $w$-odered field,i.e., a $w$-structure of the form

$$\{\mathbb{R}_w, +_w, \times_w, 0_w, 1_w\}, \quad (6.1.11)$$

where $\mathbb{R}_w =_w \mathbb{R}$ is the set of elements of the structure, $+_w$ and $\times_w$ are the binary operations of addition and multiplication, $<_w$ is the ordering relation, and $0_w =_w 0 \in_s \mathbb{R}$ and $1_w =_w 1 \in_s \mathbb{R}$ are two distinguished elements of the domain. And it is complete in the sense that every nonempty set $w$-bounded from above has a $w$-least $w$-upper bound.

(**I**) The $\#_w$-embedding of (6.1.10) sends $0_w$ to ${}^{\#_w}0 =_w \mathbf{0}_{\mathcal{F}^w}^w \triangleq \mathbf{0}_w$ and $1$ to ${}^{\#_w}1 =_w \mathbf{1}_{\mathcal{F}^w}^w \triangleq 1_w$. We must lift the operations and relations of $\mathbb{R}_w$ to ${}^{\#_w}\mathbb{R}_w$. We get the clue from (6.1.9), which tells us when:

(i) two elements $f_{\mathcal{F}^w}$ and $g_{\mathcal{F}^w}$, of ${}^{\#_w}\mathbb{R}_w$ are $w$-equal:

$$f_{\mathcal{F}^w}^w =_w g_{\mathcal{F}^w}^w \Leftrightarrow_s \left\{ v \in_w \wp^w \mid f^w(v) =_w g^w(v) \right\} \in_w \mathcal{F}^w, \quad (6.1.12)$$

(ii) two elements $f_{\mathcal{F}^w}$ and $g_{\mathcal{F}^w}$, of ${}^{\#_w}\mathbb{R}_w$ are not $w$-equal in strong consistent sense:

$$f_{\mathcal{F}^w}^w \neq_w^s g_{\mathcal{F}^w}^w \Leftrightarrow_s \left\{ v \in_w \wp^w \mid f^w(v) =_w g^w(v) \right\} \notin_w^s \mathcal{F}^w, \quad (6.1.13)$$

(iii) two elements $f_{\mathcal{F}^w}$ and $g_{\mathcal{F}^w}$, of ${}^{\#_w}\mathbb{R}_w$ are not $w$-equal in a weak sense:

$$f^{w}_{\mathcal{F}^w} \neq^{w}_{w} g^{w}_{\mathcal{F}^w} \Leftrightarrow_s \left\{v \in_w \wp^w | f^w(v) =_w g^w(v)\right\} \notin^{w}_{w} \mathcal{F}^w, \qquad (6.1.14)$$

(iv) two elements $f_{\mathcal{F}^w}$ and $g_{\mathcal{F}^w}$, of ${}^{\#_w}\mathbb{R}_w$ are are $w$-equal and are not $w$-equal in a weak $w$-inconsistent sense:

$$\left(f^{w}_{\mathcal{F}^w} =_w g^{w}_{\mathcal{F}^w}\right) \wedge \left(f^{w}_{\mathcal{F}^w} \neq^{w}_{w} g^{w}_{\mathcal{F}^w}\right) \Leftrightarrow_s f^{w}_{\mathcal{F}^w} =_{w_{[1]}} g^{w}_{\mathcal{F}^w} \Leftrightarrow_s$$
$$\Leftrightarrow_s \left[\left\{v \in_w \wp^w | f^w(v) =_w g^w(v)\right\} \in_w \mathcal{F}^w\right] \wedge \qquad (6.1.15)$$
$$\wedge\left[\left\{v \in_w \wp^w | f^w(v) =_w g^w(v)\right\} \notin^{w}_{w} \mathcal{F}^w\right].$$

In a similar way we extend $<_w$ to ${}^{\#_w}\mathbb{R}$ by setting for arbitrary $f^{w}_{\mathcal{F}^w}$,and $g^{w}_{\mathcal{F}^w}$,in ${}^{\#_w}\mathbb{R}_w$:

$$f^{w_0}_{\mathcal{F}^w} <_w g^{w}_{\mathcal{F}^w} \Leftrightarrow_s \{v \in_w \wp^w | f^w(v) <_w g^w(v)\} \in_w \mathcal{F}^w,$$
$$f^{w_0}_{\mathcal{F}^w} \nless^{s}_{w} g^{w}_{\mathcal{F}^w} \Leftrightarrow_s \{v \in_w \wp^w | f^w(v) <_w g^w(v)\} \notin^{s}_{w} \mathcal{F}^w,$$
$$f^{w_0}_{\mathcal{F}^w} \nless^{w}_{w} g^{w}_{\mathcal{F}^w} \Leftrightarrow_s \{v \in_w \wp^w | f^w(v) <_w g^w(v)\} \notin^{w}_{w} \mathcal{F}^w,$$
$$\left(f^{w_0}_{\mathcal{F}^w} <_w g^{w}_{\mathcal{F}^w}\right) \wedge \left(f^{w_0}_{\mathcal{F}^w} \nless^{s}_{w} g^{w}_{\mathcal{F}^w}\right) \Leftrightarrow_s f^{w_0}_{\mathcal{F}^w} <_{w_{[1]}} g^{w}_{\mathcal{F}^w} \Leftrightarrow_s \qquad (6.1.16)$$
$$\Leftrightarrow_s [\{v \in_w \wp^w | f^w(v) <_w g^w(v)\} \in_w \mathcal{F}^w] \wedge$$
$$\wedge[\{v \in_w \wp^w | f^w(v) <_w g^w(v)\} \notin^{w}_{w} \mathcal{F}^w].$$

(**II**) With this definition of $<_w$ in ${}^{\#_w}\mathbb{R}_w$ we easily show that the extended domain ${}^{\#_w}\mathbb{R}_w$ is $w$-linearly $w$-ordered $w$-inconsistent field. As an example we verify $w$-transitivity of $<_w$ in ${}^{\#_w}\mathbb{R}_w$. Let $f^{w}_{\mathcal{F}^w} <_w g^{w}_{\mathcal{F}^w}$, and $g^{w}_{\mathcal{F}_w} <_w h^{w}_{\mathcal{F}^w}$,i.e.,

$$D^{w}_{1} =_w \left\{v \in_w \wp^w | f^w(v) <_w g^w(v)\right\} \in_w \mathcal{F}^w,$$
$$D^{w}_{2} =_w \{v \in_w \wp^w | g^w(v) <_w h^w(v)\} \in_w \mathcal{F}^w \qquad (6.1.17)$$

By the finite intersection property,[see Definition 6.1.4.(ii)] $D^{w}_{1} \cap_w D^{w}_{2} \in_w \mathcal{F}^w$. If $v \in_w D^{w}_{1} \cap_w D^{w}_{2}$, then $f^w(v) <_w g^w(v)$ and $g^w(v) <_w h^w(v)$; hence by transitivity of $<_w$ in $\mathbb{R}_w$,

$$\left[f^w(v) <_w g^w(v)\right] \wedge [g^w(v) <_{w_0} h^w(v)] \quad \Rightarrow_s \quad f^w(v) <_{w_0} h^w(v). \qquad (6.1.18)$$

Thus

$$D^{w}_{1} \cap_w D^{w}_{2} \subseteq_{w_0} \{v \in_w \wp^w | f^w(v) <_w h^w(v)\} \qquad (6.1.19)$$

The closure property [Definition 6.1.4(iii)] then tells us that $f^{w}_{\mathcal{F}^w} <_w h^{w}_{\mathcal{F}^w}$. Similarly one can to prove that given any $f^{w}_{\mathcal{F}^w}, g^{w}_{\mathcal{F}^w} \in_w {}^{\#_w}\mathbb{R}_w$, then either $f^{w}_{\mathcal{F}^w} <_w g^{w}_{\mathcal{F}^w}$,or $g^{w}_{\mathcal{F}^w} <_w f^{w}_{\mathcal{F}^w}$, or $f^{w}_{\mathcal{F}^w} =_w g^{w}_{\mathcal{F}^w}$.

**Remark.6.1.7.**The $w$-relation $<_w$ on ${}^{\#_w}\mathbb{R}_w$ introduced in (6.1.16) extends the relation $<_w$ on $\mathbb{R}_w$, i.e., given any $r_1, r_2 \in_w \mathbb{R}_w$ we see that $r_1 <_w r_2$ in $\mathbb{R}_w$ iff ${}^{\#_w}r_1 <_w {}^{\#_w}r_2$ in ${}^{\#_w}\mathbb{R}_w$.
We now have $w$-inconsistent $w$-linear order on ${}^{\#_w}\mathbb{R}_w$ and can verify that ${}^{\#_w}\mathbb{R}_w$ contains $w$-inconsistent infinitesimals ($w$-infinitesimals) and weakly consistent infinite numbers ($w$-infinite numbers). A (positive) $w$-infinitesimal $\delta_w$ in ${}^{\#_w}\mathbb{R}_w$ is an $w$-element $\delta_w \in_w {}^{\#_w}\mathbb{R}_w$ such that ${}^{\#_w}0_w <_w \delta_w <_w {}^{\#_w}r_w$ for all $r_w\,_w > 0_w$ in $\mathbb{R}_w$.
Notice that $w$-nfinitesimals exist. Let $\mathcal{F}^w$ a free $w$-ultrafilter on $\wp^w =_w \mathbb{N}_w$

and let $f_1^w(n_w) =_w n_w^{-1}$ and $f_2^w(n_w) =_w n_w^{-2}$ for $n_w \in \mathbb{N}_w$. Then $\delta_{1,w} =_w f_{1,\mathcal{F}^w}^w$ and $\delta_{2,w} =_w f_{2,\mathcal{F}^w}^w$,is a positive $w$-infinitesimals and $\delta_{2,w} <_w \delta_{1,w}$.
In the same way $g_1^w(n_{w_0}) =_w n_w$ and $g_2^w(n_w) =_w n_w^2$ introduce a weakly consistent infinite numbers, $\omega_{1,w} =_w g_{1,\mathcal{F}^w}^w$, and $\omega_{2,w} =_w g_{2,\mathcal{F}^w}^w$, and we have that $\omega_{1,w} <_w \omega_{2,w}$ in ${}^{\#_w}\mathbb{R}_w$.
(**III**) It remains to extend the operations $+_w$ and $\times_w$ to ${}^{\#_w}\mathbb{R}_w$. Looking back to (6.1.12) and (6.1.16) we have nothing to do but to set

$$f_{\mathcal{F}^w}^w +_w g_{\mathcal{F}^w}^w =_w h_{\mathcal{F}^w}^w \Leftrightarrow_s$$
$$\Leftrightarrow_s \left\{ v \in_w \wp^w | f^w(v) +_w g^w(v)_w =_w h^w \right\} \in_w \mathcal{F}^w, \tag{6.1.19}$$

and

$$f_{\mathcal{F}^w}^w \times_w g_{\mathcal{F}^w}^w =_w h_{\mathcal{F}^w}^w \Leftrightarrow_w$$
$$\Leftrightarrow_w \left\{ v \in_w \wp^w | f^w(v) \times_w g^w(v) =_w h^w(v) \right\} \in_w \mathcal{F}^w. \tag{6.1.20}$$

With these definitions one can to proves easily that ${}^{\#_w}\mathbb{R}_w$ is an $w$-inconsistent extension of
$\mathbb{R}_w$. And these definitions introduce an $w$-inconsistent algebra on the $w$-infinitesimals and
on the $w$-infinitely large numbers. One may wish to verify easily that if $f_{\mathcal{F}^w}^w <_w g_{\mathcal{F}^w}^w$, and ${}^{\#_w}0_w <_w h_{\mathcal{F}^w}^w$,then

$$f_{\mathcal{F}^w}^w \times_w h_{\mathcal{F}^w}^w <_w g_{\mathcal{F}^w}^w \times_w h_{\mathcal{F}^w}^w. \tag{6.1.21}$$

One should also notice that for the $w_0$-infinitesimals $\delta_{1,w}$ and $\delta_{2,w}$ and the $w_0$-infinite $\omega_{1,w_0}$
and $\omega_{2,w_0}$ introduced above, we have, e.g., $\omega_{2,w} =_w \omega_{1,w}^2$ , $\delta_{1,w_0} \times_{w_0} \omega_{1,w_0} =_{w_0} \mathbf{1}_{w_0}$, is infinitesimal, and 6 w' is infinite.Thus the $w_0$-infinitely small and the $w_0$-infinitely large have
a decent wekly consistent arithmetic.
The way we extended the particular operators $+_{w_0}$ and $\times_{w_0}$ and the particular relation $=_{w_0}$ from $\mathbb{R}_{w_0}$ to ${}^{\#_{w_0}}\mathbb{R}_{w_0}$ can be used to extend any function and relation on $\mathbb{R}_{w_0}$ to ${}^{\#_{w_0}}\mathbb{R}_{w_0}$. Let $F$ be an $n$-ary $w_0$-function on $\mathbb{R}_{w_0}$, i.e.,

$$F^{w_0} : \underbrace{\mathbb{R}_{w_0} \times_{w_0} \cdots \times_{w_0} \mathbb{R}_{w_0}}_{n \text{ times}} \to_{w_0} \mathbb{R}_{w_0}. \tag{6.1.22}$$

Then we introduce the extended $w_0$-function ${}^{\#_{w_0}}F$ by the $w_0$-equivalence

$${}^{\#_{w_0}}F^{w_0}\left(f_{\mathcal{F}^{w_0}}^{1,w_0}, \ldots, f_{\mathcal{F}^{w_0}}^{n,w_0}\right) =_{w_0} g_{\mathcal{F}^{w_0}}^{w_0} \Leftrightarrow_s$$
$$\left\{ v \in_{w_0} \wp^{w_0} | F^{w_0}\left(f^{1,w_0}(v), \ldots, f^{n,w_0}(v) =_{w_0} g^{w_0}(v)\right) \right\} \in_{w_0} \mathcal{F}^{w_0} \tag{6.1.23}$$

The reader may want to verify that ${}^{\#_{w_0}}F$ is a $w_0$-function and that ${}^{\#_{w_0}}F$ really extends $F$, i.e., ${}^{\#_{w_0}}F({}^{\#_{w_0}}r_{1,w_0}, \ldots, {}^{\#_{w_0}}r_{n,w_0}) =_{w_0} {}^{\#_{w_0}}r_{w_0}$ iff $F^{w_0}(r_{1,w_0}, \ldots, r_{n,w_0}) = r_{w_0}$. In the same way we extend any $n$-ary $w_0$-relation $S^{w_0}$ on $\mathbb{R}_{w_0}$ to a $w_0$-relation ${}^{\#_{w_0}}S^{w_0}$ on ${}^{\#_{w_0}}\mathbb{R}_{w_0}$.
Note that since a $w_0$-subset $E \subseteq_{w_0} \mathbb{R}_{w_0}$ corresponds to an unary $w_0$-relation, we have an
$w_0$-extension ${}^{\#_{w_0}}E$ characterized by the condition

$$f^{w_0}_{\mathcal{F}^{w_0}} \in_{w_0} {}^{\#_{w_0}}E \Leftrightarrow_s \{v \in_{w_0} \wp^{w_0} | f(v) \in E\} \in_{w_0} \mathcal{F}^{w_0}. \quad (6.1.24)$$

Thus if $E = (0_{w_0}, 1_{w_0}]$, then ${}^{\#_{w_0}}E$ as a subset of ${}^{\#_{w_0}}\mathbb{R}_{w_0}$ will have every positive $w_0$-infinitesimal as an $w_0$-element, but not ${}^{\#_{w_0}}0_{w_0}$, a fact which can be read off immediately from condition (6.1.24).But first a few elementary observations on the $w_0$-extension of subsets of $\mathbb{R}_{w_0}$ : ${}^{\#_{w_0}}\varnothing_{w_0}$ is the $w_0$-empty set in ${}^{\#_{w_0}}\mathbb{R}_{w_0}$. If $E \subseteq_{w_0} \mathbb{R}_{w_0}$, then ${}^{\#_{w_0}}r_{w_0} \in_{w_0} {}^{\#_{w_0}}E$ for all $r_{w_0} \in_{w_0} E$, but in general(see the example $E = (0_{w_0}, 1_{w_0}]$ above)${}^{\#_{w_0}}E$ will contain elements not of the form ${}^{\#_{w_0}}r_{w_0}$ for any $r_{w_0} \in_{w_0} E$. Furthermore ${}^{\#_{w_0}}$ is a Boolean homomorphism in the sense that ${}^{\#_{w_0}}(E_1 \cup_{w_0} E_2) =_{w_0} ({}^{\#_{w_0}}E_1) \cup_{w_0} ({}^{\#_{w_0}}E_2)$ and ${}^{\#_w}(E_1 \cap_w E_2) =_{w_0} ({}^{\#_w}E_1) \cap_w ({}^{\#_w}E_2)$ for arbitrary sets $E_1, E_2 \subseteq_w \mathbb{R}_w$. Finally, we note that ${}^{\#_w}E_1 =_w {}^{\#_w}E_2$ iff $E_1 =_w E_2$, and ${}^{\#_w}r_w \in_w {}^{\#_w}E$ iff $r_w \in_w E$.

Before proceeding we need to discuss the important concept of standard part. By virtue of (6.1.23) the absolute-value function $|\cdot|_w$ on $\mathbb{R}_w$ has an extension to ${}^{\#_w}\mathbb{R}_w$ that we will denote by ${}^{\#_w}|\cdot|_w$.

**Definition 6.1.12.**An $w$-element $x \in_w {}^{\#_w}\mathbb{R}_w$is called $w$-finite if ${}^{\#_w}|x|_w <_w {}^{\#_w}r_w$ for some $0_w <_w r_w$.

As we shall see, every finite $x \in_w {}^{\#_w}\mathbb{R}_w$ is $w$-infinitely close to some $r_w \in_w \mathbb{R}_w$ in the sense that ${}^{\#_w}\left|x -_w {}^{\#_w}r_w\right|_w$ is either ${}^{\#_w}0_w$ or $w_0$-positively $w_0$-infinitesimal in ${}^{\#_{w_0}}\mathbb{R}_w$.

**Definition 6.1.13.**This $w$-unique $r_w$ is called the $w$-standard part of $x$ and is denoted by $st_w(x)$ or ${}^{\circ_w}x$.

The proof of existence of the standard part is simple. Let $x \in_w {}^{\#_w}\mathbb{R}_w$ be finite. Let $D_1$ be the set of $r_w \in_w \mathbb{R}_w$ such that ${}^{\#_w}r_w <_w x$ and $D_2$ the set of $r'_w \in_w \mathbb{R}_w$ such that $x <_w {}^{\#_w}r_w$.The pair $\{D_1, D_2\}$ forms a Dedekind cut in $\mathbb{R}_w$, hence determines a unique $\bar{r}_w \in_w \mathbb{R}_w$. A simple argument shows that $st_w(x) =_w \bar{r}_w$.

# VI.2.1 The $w_0$-consistent nonstandard extension.

**Definition 6.2.1.**Let be $\mathcal{F}^{w_0}$ a free $w_0$-consistent ultrafilter on $\wp^{w_0}$ and introduce an $w_0$-equivalence relation $f^{w_0} \sim^{w_0}_{\mathcal{F}^{w_0}} g^{w_0}$ on $w_0$-sequences $f^{w_0} \in_{w_0} \mathbb{R}_{w_0}^{\wp^{w_0}}$ by

$$f^{w_0} \sim_{\mathcal{F}^{w_0}} g^{w_0} \Leftrightarrow_s \left\{v \in_{w_0} \wp^{w_0} | f(v) =_{w_0} g(v)\right\} \in_{w_0} \mathcal{F}^{w_0}. \quad (6.2.1)$$

**Definition 6.2.2.**$\mathbb{R}_{w_0}^{\wp^{w_0}}$ divided out by the $w_0$-equivalence relation $\sim_{\mathcal{F}^{w_0}}$ gives us the $w_0$-consistent nonstandard extension ${}^{\#_{w_0}}\mathbb{R}_{w_0}$, the hyperreals; in symbols,

$${}^{\#_{w_0}}\mathbb{R}_{w_0} \triangleq \mathbb{R}_{w_0}^{\wp^{w_0}} \Big/ \sim_{\mathcal{F}^{w_0}} . \quad (6.2.2)$$

**Remark.6.2.3.**Note that for any $f^{w_0}, g^{w_0}, h^{w_0} \in_{w_0} \mathbb{R}_{w_0}^{\wp^{w_0}}$ it follows

$$\left(f^{w_0} \sim_{\mathcal{F}^{w_0}} g^{w_0}\right) \wedge \left(g^{w_0} \sim_{\mathcal{F}^{w_0}} h^{w_0}\right) \Rightarrow_s f^{w_0} \sim_{\mathcal{F}^{w_0}} h^{w_0}. \quad (6.2.3)$$

**Remark.6.2.4.**If $f^{w_0} \in_{w_0} \mathbb{R}_{w_0}^{\wp^{w_0}}$, we denote its image in ${}^{\#_{w_0}}\mathbb{R}_{w_0}$ by $f^{w_0}_{\mathcal{F}^{w_0}}$, and, of course,every $w_0$-element in ${}^{\#_{w_0}}\mathbb{R}_{w_0}$ is of the form $f^{w_0}_{\mathcal{F}^{w_0}}$, for some $f^{w_0} : \wp^{w_0} \to \mathbb{R}_{w_0}$.

**Remark.6.2.5.**Note that for any $f^{w_0}, g^{w_0}, h^{w_0}$ it follows by definitions

$$f^{w_0} \sim^{w_0}_{\mathcal{F}^{w_0}} g^{w_0} \Rightarrow_s f^{w_0}_{\mathcal{F}^{w_0}} =_{w_0} g^{w_0}_{\mathcal{F}^{w_0}},$$
$$f^{w_0} \sim^{w_0}_{\mathcal{F}^{w_0}} g^{w_0} \wedge g^{w_0} \sim h^{w_0}_{\mathcal{F}^{w_0}} \Rightarrow_s f^{w_0} \sim h^{w_0}_{\mathcal{F}^{w_0}}. \quad (6.2.4)$$

For any real number $r_{w_0} \in_{w_0} \mathbb{R}_{w_0}$ let $\mathbf{r}^{w_0}$ denote the constant $w_0$-function with value $r_{w_0}$ in $\mathbb{R}_{w_0}$, i.e.,$\mathbf{r}^{w_0}(v) =_{w_0} r_{w_0}$, for all $v \in_{w_0} \wp^{w_0}$. We then have a natural $w_0$-embedding:

$$\#_{w_0} : \mathbb{R}_{w_0} \to_{w_0} {}^{\#_{w_0}}\mathbb{R}_{w_0} \quad (6.2.5)$$

by setting ${}^{\#_{w_0}}r_{w_0} =_{w_0} \mathbf{r}^{w_0}_{\mathcal{F}^{w_0}}$, for all $r_{w_0} \in_{w_0} \mathbb{R}_{w_0}$.

(**I**) The $\#_{w_0}$-embedding of (6.2.5) sends $0$ to ${}^{\#_{w_0}}0 =_{w_0} \mathbf{0}^{w_0}_{\mathcal{F}^{w_0}} \triangleq 0_{w_0}$ and $1$ to ${}^{\#_{w_0}}1 =_{w_0} \mathbf{1}^{w_0}_{\mathcal{F}^{w_0}} \triangleq 1_{w_0}$.We must lift the operations and relations of $\mathbb{R}$ to ${}^{\#_{w_0}}\mathbb{R}_{w_0}$. We get the clue from (6.2.1), which tells us when two elements $f_{\mathcal{F}^{w_0}}$ and $g_{\mathcal{F}^{w_0}}$, of ${}^{\#_w}\mathbb{R}_{w_0}$ are $w_0$-equal:

$$f^{w_0}_{\mathcal{F}^{w_0}} =_{w_0} g^{w_0}_{\mathcal{F}^{w_0}} \Leftrightarrow_s \left\{v \in_{w_0} \wp^{w_0} | f^{w_0}(v) =_{w_0} g^{w_0}(v)\right\} \in_{w_0} \mathcal{F}^{w_0}. \quad (6.2.6)$$

In a similar way we extend $<_{w_0}$ to ${}^{\#_{w_0}}\mathbb{R}_{w_0}$ by setting for arbitrary $f^{w_0}_{\mathcal{F}^{w_0}}$,and $g^{w_0}_{\mathcal{F}^{w_0}}$,in ${}^{\#_{w_0}}\mathbb{R}_{w_0}$

$$f^{w_0}_{\mathcal{F}^{w_0}} <_{w_0} g^{w_0}_{\mathcal{F}^{w_0}} \Leftrightarrow_s \{v \in_{w_0} \wp^{w_0} | f^{w_0}(v) <_{w_0} g^{w_0}(v)\} \in_{w_0} \mathcal{F}^{w_0},$$
$$f^{w_0}_{\mathcal{F}^{w_0}} \nless^s_{w_0} g^{w_0}_{\mathcal{F}^{w_0}} \Leftrightarrow_s \{v \in_{w_0} \wp^{w_0} | f^{w_0}(v) <_{w_0} g^{w_0}(v)\} \notin^s_{w_0} \mathcal{F}^{w_0},$$
$$f^{w_0}_{\mathcal{F}^{w_0}} \nless^w_{w_0} g^{w_0}_{\mathcal{F}^{w_0}} \Leftrightarrow_s \{v \in_{w_0} \wp^{w_0} | f^{w_0}(v) <_{w_0} g^{w_0}(v)\} \notin^w_{w_0} \mathcal{F}^{w_0},$$
$$\left(f^{w_0}_{\mathcal{F}^{w_0}} <_{w_0} g^{w_0}_{\mathcal{F}^{w_0}}\right) \wedge \left(f^{w_0}_{\mathcal{F}^{w_0}} \nless^w_{w_0} g^{w_0}_{\mathcal{F}^{w_0}}\right) \Leftrightarrow_s \quad (6.2.7)$$
$$\Leftrightarrow_s (\{v \in_{w_0} \wp^{w_0} | f^{w_0}(v) <_{w_0} g^{w_0}(v)\} \in_{w_0} \mathcal{F}^{w_0}) \wedge$$
$$\wedge(_s\{v \in_{w_0} \wp^{w_0} | f^{w_0}(v) <_{w_0} g^{w_0}(v)\} \notin^w_{w_0} \mathcal{F}^{w_0}).$$

(**II**) With this definition of $<_{w_0}$ in ${}^{\#_{w_0}}\mathbb{R}_{w_0}$ we easily show that the extended domain ${}^{\#_{w_0}}\mathbb{R}_{w_0}$ is linearly $w_0$-ordered $w_0$-inconsistent field. As an example we verify $w_0$-transitivity of $<_{w_0}$ in ${}^{\#_{w_0}}\mathbb{R}$. Let $f^{w_0}_{\mathcal{F}^{w_0}} <_{w_0} g^{w_0}_{\mathcal{F}^{w_0}}$, and $g^{w_0}_{\mathcal{F}_{w_0}} <_{w_0} h^{w_0}_{\mathcal{F}^{w_0}}$,i.e.,

$$D_1^{w_0} =_{w_0} \left\{v \in_{w_0} \wp^{w_0} | f^{w_0}(v) <_{w_0} g^{w_0}(v)\right\} \in_{w_0} \mathcal{F}^{w_0},$$
$$D_2^{w_0} =_{w_0} \{v \in_{w_0} \wp^{w_0} | g^{w_0}(v) <_{w_0} h^{w_0}(v)\} \in_{w_0} \mathcal{F}^{w_0} \quad (6.2.8)$$

By the finite intersection property,[see Definition 6.2.7.(ii)] $D_1^{w_0} \cap_{w_0} D_2^{w_0} \in_{w_0} \mathcal{F}^{w_0}$. If $v \in_{w_0} D_1^{w_0} \cap_{w_0} D_2^{w_0}$, then $f^{w_0}(v) <_{w_0} g^{w_0}(v)$ and $g^{w_0}(v) <_{w_0} h^{w_0}(v)$; hence by transitivity of $<_{w_0}$ in $\mathbb{R}_{w_0}$,

$$\left[f^{w_0}(v) <_{w_0} g^{w}(v)\right] \wedge [g^{w_0}(v) <_{w_0} h^{w_0}(v)] \quad \Rightarrow_s \quad f^{w_0}(v) <_{w_0} h^{w_0}(v). \quad (6.2.9)$$

Thus

$$D_1^{w_0} \cap_{w_0} D_2^{w_0} \subseteq_{w_0} \{v \in_{w_0} \wp^{w_0} | f^{w_0}(v) <_{w_0} h^{w_0}(v)\} \quad (6.2.10)$$

The closure property [Definition 6.1.7(iii)] then tells us that $f^{w_0}_{\mathcal{F}^{w_0}} <_{w_0} h^{w_0}_{\mathcal{F}^{w_0}}$. Similarly one can to prove that given any $f^{w_0}_{\mathcal{F}^{w_0}}, g^{w_0}_{\mathcal{F}^{w_0}} \in_{w_0} \#_{w_0}\mathbb{R}$, then either $f^{w_0}_{\mathcal{F}^{w_0}} <_{w_0} g^{w_0}_{\mathcal{F}^{w_0}}$,or $g^{w_0}_{\mathcal{F}^{w_0}} <_{w_0} f^{w_0}_{\mathcal{F}^{w_0}}$, or $f^{w_0}_{\mathcal{F}^{w_0}} =_{w_0} g^{w_0}_{\mathcal{F}^{w_0}}$.

**Remark.6.2.6.**The $w_0$-relation $<_{w_0}$ on ${}^{\#_{w_0}}\mathbb{R}_{w_0}$ introduced in (6. 2.5) extends the relation $<_{w_0}$ on $\mathbb{R}_{w_0}$, i.e., given any $r_1, r_2 \in_{w_0} \mathbb{R}_{w_0}$ we see that $r_1 <_{w_0} r_2$ in $\mathbb{R}_{w_0}$ iff ${}^{\#_{w_0}}r_1 <_{w_0} {}^{\#_{w_0}}r_2$ in ${}^{\#_{w_0}}\mathbb{R}_{w_0}$.

We now have a weakly consistent linear order on ${}^{\#_{w_0}}\mathbb{R}$ and can verify that ${}^{\#_{w_0}}\mathbb{R}_{w_0}$ contains
weakly consistent infinitesimals ($w_0$-infinitesimals) and weakly consistent infinite numbers
($w_0$-infinite numbers). A (positive) $w_0$-infinitesimal $\delta_{w_0}$ in ${}^{\#_{w_0}}\mathbb{R}_{w_0}$ is an
$w_0$-element $\delta_{w_0} \in_{w_0} {}^{\#_{w_0}}\mathbb{R}$ such that $0_{w_0} <_{w_0} \delta_{w_0} <_{w_0} {}^{\#_{w_0}}r$ for all $r > 0$ in $\mathbb{R}$.
Notice that $w_0$-nfinitesimals exist. Let $\mathcal{F}^{w_0}$ a free $w_0$-consistent ultrafilter on $\wp^{w_0} =_{w_0} \mathbb{N}_{w_0}$
and let $f_1^{w_0}(n_{w_0}) =_{w_0} n_{w_0}^{-1}$ and $f_2^{w_0}(n_{w_0}) =_{w_0} n_{w_0}^{-2}$ for $n \in \mathbb{N}_{w_0}$. Then $\delta_{1,w_0} =_{w_0} f_{1,\mathcal{F}^{w_0}}^{w_0}$ and
$\delta_{2,w_0} =_{w_0} f_{2,\mathcal{F}^{w_0}}^{w_0}$,is a positive $w_0$-infinitesimals and $\delta_{2,w_0} <_{w_0} \delta_{1,w_0}$.
In the same way $g_1^{w_0}(n_{w_0}) =_{w_0} n_{w_0}$ and $g_2^{w_0}(n_{w_0}) =_{w_0} n_{w_0}^2$ introduce a weakly consistent
infinite numbers, $\omega_{1,w_0} =_{w_0} g_{1,\mathcal{F}^{w_0}}^{w_0}$, and $\omega_{2,w_0} =_{w_0} g_{2,\mathcal{F}^{w_0}}^{w_0}$, and we have that $\omega_{1,w_0} <_{w_0} \omega_{2,w_0}$
in ${}^{\#_{w_0}}\mathbb{R}_{w_0}$.

(**III**) It remains to extend the operations $+_{w_0}$ and $\times_{w_0}$ to ${}^{\#_{w_0}}\mathbb{R}_{w_0}$. Looking back to (6.2.1) and (6.2.2) we have nothing to do but to set

$$f_{\mathcal{F}^{w_0}}^{w_0} +_{w_0} g_{\mathcal{F}^{w_0}}^{w_0} =_{w_0} h_{\mathcal{F}^{w_0}}^{w_0} \Leftrightarrow_{w_0}$$
$$\Leftrightarrow_{w_0} \left\{ v \in_{w_0} \wp^{w_0} | f^{w_0}(v) +_{w_0} g^{w_0}(v)_{w_0} =_{w_0} h^{w_0} \right\} \in_{w_0} \mathcal{F}^{w_0}, \tag{6.2.11}$$

and

$$f_{\mathcal{F}^{w_0}}^{w_0} \times_{w_0} g_{\mathcal{F}^{w_0}}^{w_0} =_{w_0} h_{\mathcal{F}^{w_0}}^{w_0} \Leftrightarrow_{w_0}$$
$$\Leftrightarrow_{w_0} \left\{ v \in_{w_0} \wp^{w_0} | f^{w_0}(v) \times_{w_0} g^{w_0}(v) =_{w_0} h^{w_0}(v) \right\} \in_{w_0} \mathcal{F}^{w_0}. \tag{6.2.12}$$

With these definitions one can to proves easily that ${}^{\#_{w_0}}\mathbb{R}$ is an $w_0$-consistent extension of
$\mathbb{R}$. And these definitions introduce an $w_0$-consistent algebra on the
$w_0$-infinitesimals and on the $w_0$-infinitely large numbers. One may wish to verify easily that
if $f_{\mathcal{F}^{w_0}}^{w_0} <_{w_0} g_{\mathcal{F}^{w_0}}^{w_0}$, and ${}^{\#_{w_0}}0 <_{w_0} h_{\mathcal{F}^{w_0}}^{w_0}$,then

$$f_{\mathcal{F}^{w_0}}^{w_0} \times_{w_0} h_{\mathcal{F}^{w_0}}^{w_0} <_{w_0} g_{\mathcal{F}^{w_0}}^{w_0} \times_{w_0} h_{\mathcal{F}^{w_0}}^{w_0}. \tag{6.2.13}$$

One should also notice that for the $w_0$-infinitesimals $\delta_{1,w_0}$ and $\delta_{2,w_0}$ and the $w_0$-infinite $\omega_{1,w_0}$
and $\omega_{2,w_0}$ introduced above, we have, e.g., $\omega_{2,w_0} =_{w_0} \omega_{1,w_0}^2$ , $\delta_{1,w_0} \times_{w_0} \omega_{1,w_0} =_{w_0} \mathbf{1}_{w_0}$, is
infinitesimal, and 6 w' is infinite.Thus the $w_0$-infinitely small and the $w_0$-infinitely large have
a decent wekly consistent arithmetic.
The way we extended the particular operators $+_{w_0}$ and $\times_{w_0}$ and the particular
relation $=_{w_0}$ from $\mathbb{R}_{w_0}$ to ${}^{\#_{w_0}}\mathbb{R}_{w_0}$ can be used to extend any function and relation
on $\mathbb{R}_{w_0}$ to ${}^{\#_{w_0}}\mathbb{R}_{w_0}$. Let $F$ be an $n$-ary $w_0$-function on $\mathbb{R}_{w_0}$, i.e.,

$$F^{w_0} : \underbrace{\mathbb{R}_{w_0} \times_{w_0} \cdots \times_{w_0} \mathbb{R}_{w_0}}_{n \text{ times}} \to_{w_0} \mathbb{R}_{w_0}. \tag{6.2.14}$$

Then we introduce the extended $w_0$-function ${}^{\#_{w_0}}F$ by the $w_0$-equivalence

$$
\begin{gathered}
{}^{\#_{w_0}}F^{w_0}\left(f^{1,w_0}_{\mathcal{F}^{w_0}},\ldots,f^{n,w_0}_{\mathcal{F}^{w_0}}\right) =_{w_0} g^{w_0}_{\mathcal{F}^{w_0}} \Leftrightarrow_s \\
\left\{v \in_{w_0} \wp^{w_0} | F^{w_0}\left(f^{1,w_0}(v),\ldots,f^{n,w_0}(v) =_{w_0} g^{w_0}(v)\right)\right\} \in_{w_0} \mathcal{F}^{w_0}
\end{gathered} \tag{6.2.15}
$$

The reader may want to verify that ${}^{\#_{w_0}}F$ is a $w_0$-function and that ${}^{\#_{w_0}}F$ really extends $F$, i.e., ${}^{\#_{w_0}}F({}^{\#_{w_0}}r_{1,w_0},,\ldots,{}^{\#_{w_0}}r_{n,w_0}) =_{w_0} {}^{\#_{w_0}}r_{w_0}$ iff $F^{w_0}(r_{1,w_0},\ldots,r_{n,w_0}) = r_{w_0}$. In the same way we extend any $n$-ary $w_0$-relation $S^{w_0}$ on $\mathbb{R}_{w_0}$ to a $w_0$-relation ${}^{\#_{w_0}}S^{w_0}$ on ${}^{\#_{w_0}}\mathbb{R}_{w_0}$. Note thatsince a $w_0$-subset $E \subseteq_{w_0} \mathbb{R}_{w_0}$ corresponds to an unary $w_0$-relation, we have an $w_0$-extension ${}^{\#_{w_0}}E$ characterized by the condition

$$
f^{w_0}_{\mathcal{F}^{w_0}} \in_{w_0} {}^{\#_{w_0}}E \Leftrightarrow_s \{v \in_{w_0} \wp^{w_0} | f(v) \in E\} \in_{w_0} \mathcal{F}^{w_0}. \tag{6.2.16}
$$

Thus if $E =_{w_0} (0_{w_0}, 1_{w_0}]_{w_0}$, then ${}^{\#_{w_0}}E$ as a subset of ${}^{\#_{w_0}}\mathbb{R}_{w_0}$ will have every positive $w_0$-infinitesimal as an $w_0$-element, but not ${}^{\#_{w_0}}0_{w_0}$, a fact which can be read off immediately from condition (6.2.16).But first a few elementary observations on the $w_0$-extension of subsets of $\mathbb{R}_{w_0}$ : ${}^{\#_{w_0}}\varnothing_{w_0}$ is the $w_0$-empty set in ${}^{\#_{w_0}}\mathbb{R}_{w_0}$. If $E \subseteq_{w_0} \mathbb{R}_{w_0}$, then ${}^{\#_{w_0}}r_{w_0} \in_{w_0} {}^{\#_{w_0}}E$ for all $r_{w_0} \in_{w_0} E$, but in general (see the example $E = (0_{w_0}, 1_{w_0}]_{w_0}$ above) ${}^{\#_{w_0}}E$ will contain elements not of the form ${}^{\#_{w_0}}r_{w_0}$ for any $r_{w_0} \in_{w_0} E$. Furthermore ${}^{\#_{w_0}}$ is a Boolean homomorphism in the sense that ${}^{\#_{w_0}}(E_1 \cup_{w_0} E_2) =_{w_0} ({}^{\#_{w_0}}E_1) \cup_{w_0} ({}^{\#_{w_0}}E_2)$ and ${}^{\#_{w_0}}(E_1 \cap_{w_0} E_2) =_{w_0} ({}^{\#_{w_0}}E_1) \cap_{w_0} ({}^{\#_{w_0}}E_2)$ for arbitrary sets $E_1, E_2 \subseteq_{w_0} \mathbb{R}_{w_0}$. Finally, we note that ${}^{\#_{w_0}}E_1 =_{w_0} {}^{\#_{w_0}}E_2$ iff $E_1 =_{w_0} E_2$, and ${}^{\#_{w_0}}r_{w_0} \in_{w_0} {}^{\#_{w_0}}E$ iff $r_{w_0} \in_{w_0} E$.

Before proceeding we need to discuss the important concept of standard part. By virtue of (6.2.15) the absolute-value function $|\cdot|_{w_0}$ on $\mathbb{R}_{w_0}$ has an extension to ${}^{\#_{w_0}}\mathbb{R}_{w_0}$ that we will denote by ${}^{\#_{w_0}}|\cdot|_{w_0}$.

**Definition 6.2.3**.An $w_0$-element $x \in_{w_0} {}^{\#_{w_0}}\mathbb{R}_{w_0}$is called finite if ${}^{\#_{w_0}}|x|_{w_0} <_{w_0} {}^{\#_{w_0}}r_{w_0}$ for some $0_{w_0} <_{w_0} r_{w_0}$.

As we shall see, every finite $x \in_{w_0} {}^{\#_{w_0}}\mathbb{R}_{w_0}$ is $w_0$-infinitely close to some $r_{w_0} \in_{w_0} \mathbb{R}_{w_0}$ in the sense that ${}^{\#_{w_0}}\left|x -_{w_0} {}^{\#_{w_0}}r_{w_0}\right|_{w_0}$ is either ${}^{\#_{w_0}}0_{w_0}$ or $w_0$-positively $w_0$-infinitesimal in ${}^{\#_{w_0}}\mathbb{R}_{w_0}$.

**Definition 6.2.4**.This $w_0$-unique $r_{w_0}$ is called the $w_0$-standard part of $x$ and is denoted by $st_{w_0}(x)$ or ${}^{\circ_{w_0}}x$.

The proof of existence of the standard part is simple. Let $x \in_{w_0} {}^{\#_{w_0}}\mathbb{R}_{w_0}$ be finite. Let $D_1$ be the set of $r_{w_0} \in_{w_0} \mathbb{R}_{w_0}$ such that ${}^{\#_{w_0}}r_{w_0} <_{w_0} x$ and $D_2$ the set of $r'_{w_0} \in_{w_0} \mathbb{R}_{w_0}$ such that $x <_{w_0} {}^{\#_{w_0}}r_{w_0}$.The pair $\{D_1, D_2\}$ forms a Dedekind cut in $\mathbb{R}_{w_0}$, hence determines a unique $\bar{r}_{w_0} \in_{w_0} \mathbb{R}_{w_0}$. A simple argument shows that $st_{w_0}(x) =_{w_0} \bar{r}_{w_0}$.

# VI.2.2.The classical transfer principle.

We now remind the construction of the nonstandard extension. Let $\mathcal{F}$ be a free ultrafilter on $\wp$ and introduce an equivalence relation on sequences in $\mathbb{R}^{\wp}$ as

$$f \sim_{\mathcal{F}} g \Leftrightarrow_s \left\{v \in \wp \mid f(v) = g(v)\right\} \in_s \mathcal{F}. \tag{6.2.17}$$

$\mathbb{R}^{\wp}$ divided out by the equivalence relation $\sim_{\mathcal{F}}$ gives us the nonstandard extension ${}^*\mathbb{R}$, the hyperreals: ${}^*\mathbb{R} \triangleq \mathbb{R}^{\wp}/\mathcal{F}$. Two elements $f_{\mathcal{F}}, g_{\mathcal{F}} \in {}^*\mathbb{R}$ are equal:

$$f_{\mathcal{F}} =_s g_{\mathcal{F}} \Leftrightarrow_s \{v \in \wp | f(v) = g(v)\} \in_s \mathcal{F}. \tag{6.2.18}$$

In a similar way we extend $<$ to ${}^*\mathbb{R}$ by setting for arbitrary $f_{\mathcal{F}}, g_{\mathcal{F}} \in {}^*\mathbb{R}$

$$f_{\mathcal{F}} <_s g_{\mathcal{F}} \Leftrightarrow_s \{v \in \wp | f(v) < g(v)\} \in_s \mathcal{F}. \tag{6.2.19}$$

It remains to extend the operations $+, \times$ to ${}^*\mathbb{R}$ by

$$\begin{gathered} f_{\mathcal{F}} +_s g_{\mathcal{F}} =_s h_{\mathcal{F}} \Leftrightarrow_s \{v \in \wp | f(v) + g(v) = h(v)\} \in_s \mathcal{F}, \\ f_{\mathcal{F}} \times_s g_{\mathcal{F}} =_s h_{\mathcal{F}} \Leftrightarrow_s \{v \in \wp | f(v) \times g(v) = h(v)\} \in_s \mathcal{F} \end{gathered} \tag{6.2.20}$$

With these definitions we can prove easily that ${}^*\mathbb{R}$ is an ordered field extension of $\mathbb{R}$.

Let $F$ be an n-ary function on $\mathbb{R}$.We introduce the extended function ${}^*F$ by the equivalence

$${}^*F\left(f^1_{\mathcal{F}}, \ldots, f^n_{\mathcal{F}}\right) = g_{\mathcal{F}} \Leftrightarrow_s \left\{v \in_s \wp | F\left(f^1(v), \ldots, f^n(v) = g(v)\right)\right\} \in_s \mathcal{F} \tag{6.2.21}$$

Note thatsince a $w_0$-subset $E \subseteq_{w_0} \mathbb{R}_{w_0}$ corresponds to an unary $w_0$-relation, we have

an $w_0$-extension ${}^{\#_{w_0}}E$ characterized by the condition

$$f_{\mathcal{F}} \in_s {}^*E \Leftrightarrow_s \left\{v \in \wp \mid f(v) \in E\right\} \in_s \mathcal{F}. \tag{6.2.22}$$

We consider now the <u>standard consistent reals</u> as a structure

$$\{\mathbb{R}, +_s, \times_s, =_s, <_s, |\cdot|, 0, 1\}, \tag{6.2.23}$$

# The properties of ordered fields in classical consistent case

Any consistent ordered field $F$ is a field together with a total ordering of its elements that is compatible with the field operations. The basic example of an ordered field is the field of real numbers, and every Dedekind-complete ordered field is isomorphic to the reals $\mathbb{R}$.

**Definition 6.2.5**. A field is a nonempty set $F$ containing at least 2 elements alongside the two binary operations of addition, $f_+ : F \times_s F \to F$ such that $f_+(x,y) = x +_s y$ and multiplication $f_\times(x,y) = x \times y$ that satisfy all of the axioms below.

**I.Basic Properties of Equality**

1. $\forall x \in_s F[x =_s x]$.

2.$\forall x, y \in_s F[x =_s y \Rightarrow y =_s x]$.

3.For any function $f(x_1, \ldots, x_n) : F \times_s \ldots \times_s F \to F$, if $x_1 =_s y_1, \ldots, x_n =_s y_n$ then $f(x_1, \ldots, x_n) =_s f(y_1, \ldots, y_n)$.

**II.Axioms for Addition**

*Field Axiom for Addition* 1. The operation of addition is closed, that is $\forall x, \forall y(x +_s y \in_s F)$.

*Field Axiom for Addition* 2. The operation of addition is commutative, that is $\forall x \forall y(x +_s y =_s y +_s x)$ (Commutativity of addition).

*Field Axiom for Addition* 3. The operation of addition is associative, that is $\forall x \forall y \forall z[x +_s (y +_s z) =_s (x +_s y) +_s z]$ (Associativity of addition).

*Field Axiom for Addition* 4. The operation of addition has the additive identity element of

$0_s$ such that $\forall x(x +_s 0_s =_s x)$ (Existence of an additive identity).

**III.Axioms for Multiplication**

*Field Axiom for Multiplication* 1. The operation of multiplication is closed, that is $\forall x \forall y(x \times_s y \in_s F)$.

*Field Axiom for Multiplication* 2. The operation of multiplication is commutative, that is $\forall x \forall y\ (x \times_s y = y \times_s x)$ (Commutativity of multiplication).

*Field Axiom for Multiplication* 3 The operation of multiplication is associative, that is $\forall x \forall y \forall z\ [x \times_s (y \times_s z) =_s (x \times_s y) \times_s z]$ (Associativity of multiplication).

*Field Axiom for Multiplication* 4 The operation of multiplication has the multiplicative identity element of $1_s$ such that

$\forall x\ (1_s \times_s x =_s x)$ (Existence of an multiplicative identity).

*Field Axiom for Multiplication* 5 The operation of multiplication has the multiplicative inverse element of $1_s/x$ such that

$\forall x\ (x \times_s 1_s/x =_s 1_s)$ (Existence of a multiplicative inverse).

**IV.Field Axiom for Distributivity**

The operation of multiplication is distributive over addition, that is $\forall x \forall y \forall z[x \times_s (y +_s z) =_s x \times_s y +_s x \times_s z]$ (Distributive law).

**V.Order Axioms**

1. Either $x =_s y$ or $x <_s y$ or $y <_s x$ (Trichotemy)
2. $x <_s y$ if and only if $x +_s z <_s y +_s z$ (Addition Law)
3. If $z_s > 0_s$, then $x \times_s z <_s y \times_s z$ if and only if $x <_s y$. If $c <_s 0_s$, then $a \times_s c <_s b \times_s c$ if and only if $b <_s a$ (Multiplication Law)
4. If $x <_s y$ and $y <_s z$, then $x <_s z$ (Transitivity)

# The upper and lower bounds in classical consistent case

**Definition 6.2.5.** If $A \subset \mathbb{R}$ is a set of real numbers, then:

(i) $a \in_s \mathbb{R}$ is an upper bound for $A$ if $x \leq_s a$ for all $x \in_s A$,

and we shall denote this relation by $(\cdot U_s \cdot)$, so $aU_sA$ meant that $a$ is an upper bound of $A$;

(ii) $b$ is the least upper bound or supremum ($s$-$\sup(A)$) for $A$ if $b$ is an upper bound, and moreover $b \leq_s a$ whenever $a$ is any upper bound for $A$, and we shall denote this relation by $(\cdot L_sU_s \cdot)$, so $bL_sU_sA$ meant that $b$ is least upper bound of $A$.

One similarly defines lower bound and greatest lower bound or infinum $(\inf(A))$ for $A$ by replacing $\leq_s$ by $_s \geq$.

**Definition 6.2.6.** If $A \subset \mathbb{R}$ is a set of real numbers, then:

(i) $a \in_s \mathbb{R}$ is an lower bound for $A$ if $a \leq_s x$ for all $x \in_s A$,

and we shall denote this relation by $(\cdot L_s \cdot)$, so $aL_sA$ meant that $a$ is an lower bound of $A$;

(ii) $b$ is the gratest lower bound or infinum ($s$-$\inf(A)$) for $A$ if $b$ is an lower bound, and moreover $a \leq_s b$ whenever $a$ is any lower bound for $A$, and we shall denote this relation by $(\cdot G_sL_s \cdot)$, so $bG_sL_sA$ meant that $b$ is gratest lower bound of $A$.

**Remark.6.2.7.**The following second-order sentence expresses the least upper bound property:

$$(\forall A \subset_s \mathbb{R})([(\exists w)(w \in_s A) \wedge (\exists z)(\forall u)(u \in_s A \Rightarrow u \leq_s z)] \ \Rightarrow\searrow$$
$$\Rightarrow (\exists x)(\forall y)([(\forall w)(w \in_s A \Rightarrow w \leq_s x) \wedge\searrow \quad (6.2.23)$$
$$\wedge(\forall u)(u \in_s A \to u \leq_s y)] \ \Rightarrow x \leq_s y)).$$

The structure $\mathbb{R}$ has an associated consistent simple language $\mathcal{L}_s(\mathbb{R})$ that can be used to describe the kind of properties of $\mathbb{R}$ that are preserved under the $*$-embedding:

$$* : \mathbb{R} \hookrightarrow {}^*\mathbb{R}. \quad (6.2.23)$$

The elementary formulas of $\mathcal{L}(\mathbb{R})$ are expressions of the form
(i) $t_1 + t_2 =_s t_3$,(ii) $t_1 \times t_2 =_s t_3$,(iii) $|t_1| =_s t_2$,(iv) $t_1 =_s t_2$,(v) $t_1 <_s t_2$,(vi) $t_1 \in_s X$,
where $t_1, t_2, t_3$ are either the constants $0$ or $1$ or a variable for an arbitrary number $r \in_s \mathbb{R}$, and $X$ is a variable for a subset $A \subseteq_s \mathbb{R}$.
From the elementary formulas we generate the class of all formulas or expressions of $\mathcal{L}_s(\mathbb{R})$ using the propositional connectives:$\wedge, \vee, \neg_s, \Rightarrow_s$ ,
and the number quantifiers: $\forall x(x \in_s \mathbb{R}), \exists x(x \in_s \mathbb{R})$ by the rules:
(vii) If $\Phi$ and $\Psi$ are formulas of $\mathcal{L}_s(\mathbb{R})$, then $\Phi \wedge \Psi, \Phi \vee \Psi, \Phi \Rightarrow_s \Psi, \neg_s\Phi$,
are formulas of $\mathcal{L}_s(\mathbb{R})$,and the consistent number quantifiers $\forall x(x \in_s \mathbb{R}), \exists x(x \in_s \mathbb{R})$ are
formulas of $\mathcal{L}_s(\mathbb{R})$.
(viii) If $\Phi$ is a formula of $\mathcal{L}_s(\mathbb{R})$ and $x$ is a consistent number variable, then $\forall x\Phi, \exists x\Phi$ are formulas of $\mathcal{L}_s(\mathbb{R})$.
The language $\mathcal{L}_s(\mathbb{R})$ is basically a first-order consistent language; i.e., we allow number quantification but not set quantification.

We give a few examples: in the language $\mathcal{L}_s(\mathbb{R})$ we can write down conditions which express that $<_s$ is a strongly consistent linear ordering:
(1) $s$-transitive $\forall x\forall y\forall z[(x <_s y) \wedge (y <_s z) \Rightarrow_s x <_s z]$
(2) $s$-irreflexive $\forall x[\neg_s(x <_s x)]$
(3) $s$-linear $\forall x\forall y[(x <_s y) \vee (x =_s y) \vee (y <_s x)]$
A formula $\Phi$ of $\mathcal{L}_s(\mathbb{R})$ is in general of the form

$$\Phi = \Phi(X_1, \ldots, X_m, x_1, \ldots, x_q), \quad (6.2.24)$$

where $x_1, \ldots, x_q$ are the free consistent number variables of $\Phi$,i.e., variables not bound by a quantifier $\forall, \exists$ and $X_1, \ldots, X_m$ are the (free) consistent set variables of $\Phi$. Every formula in $\mathcal{L}(\mathbb{R})$ has an standard interpretation in the structure $\mathbb{R}$;e.g., let $\Phi(X)$ be the formula

$$\Phi(X) \equiv \forall y\big[y \in_s X \Rightarrow_s \exists z\big[\big(z \ {}_s > 0\big) \wedge \forall y_1[|y -_s y_1| <_s z \Rightarrow_s y_1 \in_s X]\big]\big] \quad (6.2.25)$$

and let $A \subseteq_s \mathbb{R}$,then $\Phi(A)$ expresses the fact that $A$ is open in $\mathbb{R}$.

# The classical Łȯs Theorem

Remind the following theorem.

**Theorem 6.2.1**. (**Łoś Theorem**) Let $\Phi(X_1, \ldots, X_m, x_1, \ldots, x_q)$ be a formula of $\mathcal{L}(\mathbb{R})$. Then for any $A_1,\ldots,A_m \subseteq_s \mathbb{R}$ and $f^1_{\mathcal{F}},\ldots,f^q_{\mathcal{F}} \in_s {}^*\mathbb{R}$

$$\Phi\left(A_1,\ldots,A_m,f^1_{\mathcal{F}},\ldots,f^q_{\mathcal{F}}\right) \Leftrightarrow_s$$
$$\Leftrightarrow_s \left\{v \in_s \wp|\Phi\left(A_1,\ldots,A_m,f^1(v),\ldots,f^q(v)\right)\right\} \in_s \mathcal{F}. \tag{6.2.26}$$

**Proof**.The proof is by induction on the number of logical symbols in $\Phi$. If $\Phi$ has no logical symbols, it is an elementary formula of the form (i)-(vi),and (6.2.26) then reduces to one of (6.2.18),(6.2.19),(6.2.20),(6.2.21), or (6.2.22). If $\Phi$ contains logical symbols, then $\Phi$ is of the form $\Phi \wedge \Psi, \Phi \vee \Psi, \Phi \Rightarrow_s \Psi, \neg_s\Phi, \forall x(x \in_s \mathbb{R}), \exists x(x \in_s \mathbb{R})$. The verification of (6.2.26) is, by induction, in each case reduced to an elementary property of the consistent ultrafilter $\mathcal{F}$.For example, if $\Phi = \Phi_1 \wedge \Phi_2$, (6.2.26) follows from the finite intersection property of the consistent ultrafilter $\mathcal{F}$. The case $\Phi = \neg_s\Phi_1$ uses in an essential way that $\mathcal{F}$ is an consistent ultrafilter namely, that

$$\wp\backslash_s\Im \in_s \mathcal{F} \Leftrightarrow_s \Im \notin^s_s \mathcal{F}. \tag{6.2.27}$$

Quantifiers offer no special difficulties,For example, if $\Phi = \exists x\Phi_1$ and let $\Phi$ have one free variable; we shall prove

$$\Phi(f_{\mathcal{F}}) \text{ iff } \{v \in \wp|\Phi(f(v))\} \in_s \mathcal{F}, \tag{6.2.28}$$

where $\Phi(f_{\mathcal{F}})$ is of the form $\exists x\Phi_1(x,f_{\mathcal{F}})$. Now $\Phi(f_{\mathcal{F}})$ is true in ${}^*\mathbb{R}$ iff there is some $g_{\mathcal{F}} \in_s {}^*\mathbb{R}$ such that $\Phi_1(g_{\mathcal{F}},f_{\mathcal{F}})$ is true in ${}^*\mathbb{R}$. By the induction hypothesis this means that

$$\{v \in \wp|\Phi_1(g(v),f(v))\} \in_s \mathcal{F}. \tag{6.2.29}$$

But if $\Phi_1(g(v),f(v))$ is true in $\mathbb{R}$, then $\exists x\Phi_1(x,f(v))$ is also true in $\mathbb{R}$,i.e.,

$$\{v \in \wp|\Phi_1(g(v),f(v))\} \subseteq \{v \in \wp|\exists x\Phi_1(x,f(v))\}. \tag{6.2.30}$$

From (6.2.30) and the property (3) of consistent filters it follows that

$$\{v \in \wp|\Phi(f(v))\} \in_s \mathcal{F}. \tag{6.2.21}$$

In order to prove the converse, assume that: $\Re_\Phi = \{v \in \wp|\Phi(f(v))\} \in_s \mathcal{F}$. For each $v$ such that $v \in \Re_\Phi$ we choose some $\alpha_v \in_s \mathbb{R}$ such that $\Phi_1(\alpha_v,f(v))$ is true in $\mathbb{R}$. Let $g \in \mathbb{R}^\wp$ be a function $g : \wp \to \mathbb{R}$ such that $g(v) =_s \alpha_v$ for all $v \in \Re_\Phi$ and $g(v) =_s \beta$ otherwise, where $\beta$ is some arbitrary $s$-element of $\mathbb{R}$. Then we have:

$$\{v \in \wp|\Phi_1(g(v),f(v))\} \in_s \mathcal{F}. \tag{6.2.22}$$

Hence by the induction hypothesis we have $\Phi_1(g_{\mathcal{F}},f_{\mathcal{F}})$ is true in ${}^*\mathbb{R}$, i.e., we have $\Phi(f_{\mathcal{F}}) = \exists x\Phi_1(x,f_{\mathcal{F}})$ is true in ${}^*\mathbb{R}$.

The theorem of Łoś has the consistent transfer principle as an immediate corollary.

**Theorem 6.2.2**.(**CONSISTENT TRANSFER PRINCIPLE**).Let $\Phi(X_1,\ldots,X_m,x_1,\ldots,x_q)$ be a formula of $\mathcal{L}(\mathbb{R})$. Then for any $A_1,\ldots,A_m \subseteq_s \mathbb{R}$ and $r_1,\ldots,r_n \in_s \mathbb{R}$, $\Phi(A_1,\ldots,A_m,r_1,\ldots,r_n)$ holds in $\mathbb{R}$ iff ${}^{*}\Phi = \Phi({}^{*}A_1,\ldots,{}^{*}A_m,{}^{*}r_1,\ldots,{}^{*}r_n)$ holds in ${}^{*}\mathbb{R}$,i.e.,

$$\Phi(A_1,\ldots,A_m,r_1,\ldots,r_n) \Leftrightarrow_s \Phi({}^{*}A_1,\ldots,{}^{*}A_m,{}^{*}r_1,\ldots,{}^{*}r_n). \quad (6.2.23)$$

**Proof.**From (6.2.16) we get at once

$$\Phi({}^{*}A_1,\ldots,{}^{*}A_m,{}^{*}r_1,\ldots,{}^{*}r_n) \Leftrightarrow_s \{v \in \wp | \Phi(A_1,\ldots,A_m,r_1,\ldots,r_n)\} \in_s \mathcal{F}^{con} \quad (6.2.24)$$

But the set $\{v \in \wp | \Phi(A_1,\ldots,A_m,r_1,\ldots,r_n)\}$ is equal to $\wp \in \mathcal{F}^{con}$ if $\Phi$ is true of $A_1,\ldots,A_m,r_1,\ldots,r_n$ in $\mathbb{R}$, and is equal to $\varnothing \notin_s \wp$ if $\Phi$ is not true of $A_1,\ldots,A_m,r_1,\ldots,r_n$ in $\mathbb{R}$.

Thus $\Phi(A_1,\ldots,A_m,r_1,\ldots,r_n)$ holds in $\mathbb{R}$ iff $\Phi({}^{*}A_1,\ldots,{}^{*}A_m,{}^{*}r_1,\ldots,{}^{*}r_n)$ holds in ${}^{*}\mathbb{R}$.

# VI.4.The Generalized Łoś Theorem

We will consider the standard $w$-inconsistent reals as an $w$-inconsistent algebraic structure $\mathfrak{R}_w$.As $w$-inconsistent algebraic structure, $\mathfrak{R}_w$ is a $w$-complete $w$-ordered field, i.e., $w$-inconsistent structure of the form

$$\mathfrak{R}_w \triangleq \{\mathbb{R}_w, +_w, \times_w, =_w, <_w, 0_w, 1_w\}, \quad (6.4.1)$$

where $\mathbb{R}_w =_w \{r | r =_w r' \in_s \mathbb{R}\}_w$ is the set of $w$-elements of the inconsistent structure,$+_w$ and $\times_w$ are the binary operations of addition and multiplication, $<_w$ is the $w$-ordering relation, and $0_w$ and $1_w$ are two distinguished elements of the domain such that $\neg_s(0_w =_w 1_w)$ but note that

$$\neg_w(0_w =_w 1_w) \wedge (0_w =_w 1_w) \nVdash A, \quad (6.4.2)$$

i.e. sentence $\neg_w(0_w =_w 1_w) \wedge (0_w =_w 1_w)$ holds in $\mathfrak{R}_w$.

And it is complete in the sense that every $w$-nonempty set $w$-bounded from above has a $w$-least $w$-upper bound. We consider now the standard inconsistent $w$-reals as $w$-inconsistent structure

$$\{\mathbb{R}_w, +_w, \times_w, =_w, <_w, |\bullet|_w, 0_w, 1_w\}, \quad (6.4.3)$$

where, in addition to the information in (6.4.1), we have added the absolute value $|\bullet|_w$ that defines the metric on $\mathbb{R}_w$. Of course,$|\bullet|_w$ is definable in terms of the other entities in (6.4.1), but it makes things a bit easier to include it explicitly in the specification.

The structure $\mathbb{R}_w$ has an associated simple language $\mathcal{L}_w = \mathcal{L}_w(\mathbb{R}_w)$ that can be used to describe the kind of properties of $\mathbb{R}_w$ that are preserved under the $\#_w$-embedding:

$$\#_w : \mathbb{R}_w \hookrightarrow {}^{\#_w}\mathbb{R}_w. \quad (6.4.4)$$

The elementary formulas of $\mathcal{L}_w(\mathbb{R}_w)$ are expressions of the form:

(i) $t_1 +_w t_2 =_w t_3$,(ii) $t_1 \times_w t_2 =_w t_3$,(iii) $|t_1|_w =_w t_2$,(iv) $t_1 =_w t_2$,(v) $t_1 <_w t_2$,(vi) $t_1 \in_w X$, where $t_1,t_2,t_3$ are either the constants $0_w$ or $1_w$ or a variable for an arbitrary number $r \in_w \mathbb{R}_w$, and $X$ is a variable for a $w$-subset $A \subseteq_w \mathbb{R}_w$.

From the elementary formulas we generate the class of all formulas or expressions of $\mathcal{L}_w(\mathbb{R}_w)$ using the propositional connectives:$\wedge, \vee, \neg_s, \Rightarrow_s, \neg_w, \Rightarrow_w$, and the inconsistent number quantifiers: $\forall x(x \in_w \mathbb{R}_w), \exists x(x \in_w \mathbb{R}_w)$ by the rules:

(vii) If $\Phi$ and $\Psi$ are formulas of $\mathcal{L}_w(\mathbb{R}_w)$, then

$\Phi \wedge \Psi, \Phi \vee \Psi, \Phi \Rightarrow_s \Psi, \neg_s\Phi, \Phi \Rightarrow_w \Psi, \neg_w\Phi$

are formulas of $\mathcal{L}_w(\mathbb{R}_w)$,and the consistent or inconsistent number quantifiers $\forall x(x \in_w \mathbb{R}_w), \exists x(x \in_w \mathbb{R}_w)$ are formulas of $\mathcal{L}_w(\mathbb{R}_w)$.

(viii) If $\Phi$ is a formula of $\mathcal{L}_w(\mathbb{R}_w)$ and $x$ is a consistent or inconsistent number variable, then $\forall x\Phi, \exists x\Phi$ are formulas of $\mathcal{L}_w(\mathbb{R}_w)$.

The language $\mathcal{L}_w(\mathbb{R}_w)$ is basically a first-order inconsistent language; i.e., we allow number quantification but not set quantification.

We give a few examples: in the language $\mathcal{L}_w(\mathbb{R}_w)$ we can write down conditions which express that $<_w$ is a $w$-inconsistent linear $w$-ordering:

(1) $w$-transitive $\forall x \forall y \forall z[(x <_w y) \wedge (y <_w z) \Rightarrow_s x <_w z]$,

(2) $w_{[1]}$-transitive $\forall x \forall y \forall z[(x <_{w_{[1]}} y) \wedge (y <_{w_{[1]}} z) \Rightarrow_s x <_{w_{[1]}} z]$,

(3) $w$-reflexive $\forall x[\neg_s(x <_w x)]$,

(4) $\forall x \forall y[(x <_w y) \Leftrightarrow_s \neg_s(y <_w x)]$,

(5) $w_{[1]}$-linear $\forall x \forall y[(x <_w y) \vee (x =_w y) \vee (y <_w x)(x <_{w_{[1]}} y) \vee (x =_{w_{[1]}} y) \vee (y <_{w_{[1]}} x)]$

A formula $\Phi$ of $\mathcal{L}_w(\mathbb{R}_w)$ is in general of the form

$$\Phi = \Phi(X_1,\ldots,X_m,x_1,\ldots,x_q), \quad (6.4.5)$$

where $x_1,\ldots,x_q$ are the free consistent and inconsistent number variables of $\Phi$, i.e., variables not bound by a quantifier $\forall,\exists$ and $X_1,\ldots,X_m$ are the (free) consistent and inconsistent set variables of $\Phi$. Every formula in $\mathcal{L}_w(\mathbb{R}_w)$ has an standard interpretation in

the structure $\mathbb{R}_w$; e.g., (i) let $\Phi(X)$ be the formula

$$\Phi(X) \equiv \forall y\big[y \in_w X \Rightarrow_s \exists z\big[\big(z\ _w > 0_w\big) \wedge \forall y_1[|y -_w y_1|_w <_w z \Rightarrow_s y_1 \in_w X]\big]\big] \quad (6.4.6)$$

and let $A \subseteq_w \mathbb{R}_w$, then $\Phi(A)$ expresses the fact that $A$ is open in $\mathbb{R}_w$;

(ii) let $\Phi_1(X)$ be the formula

$$\Phi(X) \equiv \forall y\big[y \in_w X \Rightarrow_s \exists z\big[\big(z\ _w > 0_w\big) \wedge \forall y_1[|y -_w y_1|_w <_w z \Rightarrow_s y_1 \in_w X]\big]\big] \quad (6.4.7)$$

and let $A \subseteq_w \mathbb{R}_w$, then $\Phi(A)$ expresses the fact that $A$ is open in $\mathbb{R}_w$;

**Remark.6.4.1**. Note that $w$-inconsistent algebraic structure $\mathfrak{R}_w$ mentioned above is a $w$-complete $w$-ordered field.

## VI.4.2.The properties of $w$-inconsistent $w$-ordered field $F_w$

**Definition 6.4.1.**A $w$-inconsistent field is a nonempty $w$-inconsistent set $F_w$ containing at least 2 elements along side the two binary operations of $w$-addition, $f_{+_w} : F_w \times_w F_w \to F_w$ such that $f_{+_w}(x,y) =_w x +_w y$ and $w$-multiplication $f_{\times_w}(x,y) =_w x \times_w y$ that satisfy all of the axioms below.

**I.Basic properties of $w$-inconsistent $w$-equality**

1. $\forall x \in_w F_w[(x =_w x) \vee (x =_{w_{[1]}} x)]$.

2.

3.$\forall x,y \in_w F_w[x =_w y \Rightarrow_s y =_w x]$.

4.

5.For any function $f(x_1,\ldots,x_n) : F_w \times_w \ldots \times_w F_w \to F_w$, if $x_1 =_w y_1,\ldots,x_n =_w y_n$ then $f(x_1,\ldots,x_n) =_w f(y_1,\ldots,y_n)$.

6.

**II.Axioms for $w$-addition**

*Field axiom for w-addition* 1. The operation of $w$-addition is closed, that is $\forall x \forall y(x +_w y \in_w F_w)$.

*Field axiom for w-addition* 2. The operation of $w$-addition is $w$-commutative, that is $\forall x \forall y(x +_w y =_w y +_w x)$ ($w$-commutativity of $w$-addition).

*Field axiom for w-addition* 3. The operation of addition is associative, that is $\forall x \forall y \forall z[x +_w (y +_w z) =_w (x +_w y) +_w z]$ (Associativity of addition).

*Field axiom for w-addition* 4. The operation of $w$-addition has the $w$-additive $w$-identity $w$-element of $0_w$ such that

$\forall x(x +_w 0_w =_w x)$ (Existence of an $w$-additive $w$-identity).

**III.Axioms for $w$-multiplication**

*Field axiom for w-multiplication* 1. The operation of $w$-multiplication is closed, that is $\forall x \forall y(x \times_w y \in_w F_w)$.

*Field axiom for w-multiplication* 2. The operation of $w$-multiplication is $w$-commutative, that is

$\forall x \forall y\ (x \times_w y = y \times_w x)$ ($w$-commutativity of $w$-multiplication).

*Field axiom for w-multiplication* 3 The operation of multiplication is associative, that is $\forall x \forall y \forall z\ [x \times_w (y \times_w z) =_s (x \times_s y) \times_s z]$ ($w$-associativity of $w$-multiplication).

*Field axiom for w-multiplication* 4 The operation of $w$-multiplication has the $w$-multiplicative

$w$-identity element of $1_w$ such that

$\forall x\ (1_w \times_w x =_w x)$ (Existence of an $w$-multiplicative $w$-identity).

*Field axiom for w-multiplication* 5 The operation of multiplication has the $w$-multiplicative

$w$-inverse $w$-element of $1_w/x$ such that

$\forall x(x \neq_s 0_w)\ [(x \times_w 1_w/x =_s 1_s)]$ (Existence of a multiplicative inverse).

**IV.Field axiom for $w$-distributivity**

The operation of $w$-multiplication is $w$-distributive over $w$-addition, that is

$\forall x \forall y \forall z[x \times_s (y +_s z) =_s x \times_s y +_s x \times_s z]$ (Distributive law).

**V.Order Axioms**

1. Either $x =_w y$ or $x <_w y$ or $y <_w x$ or $x =_{w_{[1]}} y$ or ($w$-trichotemy)
2. $x <_w y$ if and only if $x +_w z <_w y +_w z$ ($w$-addition law)
3. If $z_w > 0_w$, then $x \times_w z <_w y \times_w z$ if and only if $x <_w y$.

If $c <_w 0_w$, then $x \times_w c <_w y \times_w c$ if and only if $y <_s x$ (Multiplication Law)

4.If $x <_w y$ and $y <_w z$, then $x <_w z$ ($w$-transitivity)

# VI.4.3.The $w$-upper and $w$-lower bounds in $w$-inconsistent case.

**Definition 6.4.2.**If $A \subset_w \mathbb{R}_w$ is a $w$-set of $w$-inconsistent real numbers, then:

(i) $a \in_w \mathbb{R}$ is an strong $w$-upper bound for $A$ if

$$\forall x(x \in_w A)\big[\ x \leq_w a\ \big], \tag{6.4.8}$$

and we shall denote this relation by $(\cdot SU_w \cdot)$, so $aSU_wA$ meant that $a$ is an strong $w$-upper bound of $A$;

(ii) $b$ is the least strong $w$-upper bound or strong $w$-supremum ($w_S\text{-sup}(A)$) for $A$ if $b$ is an strong $w$-upper bound, and moreover

$$b \leq_w a \tag{6.4.9}$$

whenever $a$ is any strong $w$-upper bound for $A$, and we shall denote this relation by $(\cdot L_w SU_w \cdot)$, so $bL_wSU_wA$ meant that $b$ is least strong $w$-upper bound of $A$.

**Remark.6.4.2**. One similarly defines strong $w$-lower bound and greatest strong $w$-lower bound or strong $w$-infinum ($w_S$-inf($A$)) for $A$ by replacing $\leq_w$ by $_w \geq$.

**Definition 6.4.3.**If $A \subset_w \mathbb{R}_w$ is a $w$-set of $w$-inconsistent real numbers, then:

(i) $a \in_w \mathbb{R}_w$ is an strong $w$-lower bound for $A$ if

$$\forall x(x \in_w A)[a \leq_w x] \tag{6.4.10}$$

for all $x \in_s A$,

and we shall denote this relation by $(\bullet SL_w \bullet)$, so $aSL_wA$ meant that $a$ is an strong $w$-lower bound of $A$;

(ii) $b$ is the gratest strong $w$-lower bound or strong $w$-infinum ($w_S$-inf($A$)) for $A$ if $b$ is an strong $w$-lower bound, and moreover

$$a \leq_w b \tag{6.4.11}$$

whenever $a$ is any strong $w$-lower bound for $A$, and we shall denote this relation by $(\bullet G_wSL_w \bullet)$, so $bG_wSL_wA$ meant that $b$ is gratest strong $w$-lower bound of $A$.

**Remark.6.4.3.**We rewrite now the inequality (6.4.8) in the following equivalent form

$$\forall x(x \in_w A)\left[\neg_s\left( a <_w x\right)\right]. \tag{6.4.12}$$

From the statement (6.4.12) by using logical postulate $\neg_s A \Rightarrow_s \neg_w A$ we obtain

$$\forall x(x \in_w A)\left[\neg_w\left( a <_w x\right)\right]. \tag{6.4.13}$$

Note that by using (6.4.13) one obtains more weakened conditions then required above in Definition 6.4.2-6.4.3.

**Definition 6.4.4.**If $A \subset_w \mathbb{R}_w$ is a $w$-set of $w$-inconsistent real numbers, then:

(i) $a \in_w \mathbb{R}$ is an weak $w$-upper bound for $A$ if

$$\forall x(x \in_w A)\left[\neg_w\left( a <_w x\right)\right], \tag{6.4.12}$$

and we shall denote this relation by $(\bullet WU_w \bullet)$, so $aWU_wA$ meant that $a$ is an weak $w$-upper bound of $A$;

(ii) $b$ is the least weak $w$-upper bound or weak $w$-supremum ($w_W$-sup($A$)) for $A$ if $b$ is an weak $w$-upper bound, and moreover

$$\neg_w\left( a <_w b\right) \tag{6.4.13}$$

whenever $a$ is any weak $w$-upper bound for $A$, and we shall denote this relation by $(\bullet L_wWU_w \bullet)$, so $bL_wWU_wA$ meant that $b$ is least weak $w$-upper bound of $A$.

**Remark.6.4.4**. One similarly defines weak $w$-lower bound and greatest weak $w$-lower bound or weak $w$-infinum ($w_W$-inf($A$)) for $A$ by replacing $<_w$ by $_w >$.

**Definition 6.4.5.**If $A \subset_w \mathbb{R}_w$ is a $w$-set of $w$-inconsistent real numbers, then:

(i) $a \in_w \mathbb{R}_w$ is an strong $w$-lower bound for $A$ if

$$\forall x(x \in_w A)[a \leq_w x] \tag{6.4.14}$$

for all $x \in_s A$,

and we shall denote this relation by $(\bullet SL_w \bullet)$, so $aSL_wA$ meant that $a$ is an strong $w$-lower bound of $A$;

(ii) $b$ is the gratest strong $w$-lower bound or strong $w$-infinum ($w_S\text{-}\inf(A)$) for $A$ if $b$ is an strong $w$-lower bound, and moreover

$$a \leq_w b \qquad (6.4.11)$$

# VI.4.4.$w$-complete $w$-inconsistent $w$-ordered field.

# VI.3.3.The properties of $w$-inconsistent naturals $\mathbb{N}_w$

**Remark.6.3.1**. The $w$-inconsistent structure $\Re_w$ has an consistent substructure $\Re_w^s \subset_w^s \Re_w$

$$\{\mathbb{R}_w^s, +_w^s, \times_w^s, =_w^s, <_w^s, |\bullet|_w^s, 0_w^s, 1_w^s\}, \qquad (6.3.)$$

denoted below by $\mathbb{R}_w^s$ or by $\mathbb{R}_w^{con}$. The structure $\mathbb{R}_w^s$ has an associated simple language $\mathcal{L}_w^s = \mathcal{L}_w^s(\mathbb{R}_w^s)$ that can be used to describe the kind of properties of $\mathbb{R}_w^s$ that are preserved under the $\#_w$-embedding:

$$\#_w : \mathbb{R}_w^s \hookrightarrow {}^{\#_w}\mathbb{R}_w^s \subset_w^s {}^{\#_w}\mathbb{R}_w. \qquad (6.3.)$$

The elementary formulas of $\mathcal{L}_w^s(\mathbb{R}_w^s)$ are expressions of the form:

(i) $t_1 +_w^s t_2 =_w^s t_3$,(ii) $t_1 \times_w^s t_2 =_w^s t_3$,(iii) $|t_1|_w^s =_w^s t_2$,(iv) $t_1 =_w^s t_2$,(v) $t_1 <_w^s t_2$,(vi) $t_1 \in_w^s X$, where $t_1, t_2, t_3$ are either the constants $0_w^s$ or $1_w^s$ or a variable for an arbitrary number $r \in_w^s \mathbb{R}_w^s$, and $X$ is a variable for a $w$-consistent $w$-subset $A \subseteq_w^s \mathbb{R}_w^s$.From the elementary formulas we generate the class of all formulas or expressions of $\mathcal{L}_w(\mathbb{R}_w)$ using the propositional connectives:$\wedge, \vee, \neg_s, \Rightarrow_s, \neg_w, \Rightarrow_w$ ,and the $w$-consistent number quantifiers:

$\forall x(x \in_w^s \mathbb{R}_w), \exists x(x \in_w^s \mathbb{R}_w^s)$ by the rules:

(vii) If $\Phi$ and $\Psi$ are formulas of $\mathcal{L}_w^s(\mathbb{R}_w^s)$, then

$\Phi \wedge \Psi, \Phi \vee \Psi, \Phi \Rightarrow_s \Psi, \neg_s\Phi, \Phi \Rightarrow_w \Psi, \neg_w\Phi$

are formulas of $\mathcal{L}_w^s(\mathbb{R}_w^s)$,and the $w$-consistent number quantifiers $\forall x(x \in_w^s \mathbb{R}_w), \exists x(x \in_w^s \mathbb{R}_w^s)$ are formulas of $\mathcal{L}_w^s(\mathbb{R}_w^s)$.

(viii) If $\Phi$ is a formula of $\mathcal{L}_w^s(\mathbb{R}_w^s)$ and $x$ is a $w$-consistent number variable, then $\forall x\Phi, \exists x\Phi$ are formulas of $\mathcal{L}_w^s(\mathbb{R}_w^s)$.

The language $\mathcal{L}_w^s(\mathbb{R}_w^s)$ is basically a first-order $w$-consistent language; i.e., we allow number quantification but not set quantification.

We give a few examples: in the language $\mathcal{L}_w^s(\mathbb{R}_w^s)$ we can write down conditions which express that $<_w^s$ is a $w$-consistent linear $w$-ordering:

We give a few examples: in the language $\mathcal{L}_w^s(\mathbb{R}_w^s)$ we can write down conditions which express that $<_w^s$ is a strongly consistent linear ordering:

(1) $s$-$w$-transitive $\forall x\forall y\forall z[(x <_w^s y) \wedge (y <_w^s z) \Rightarrow_s x <_w^s z]$
(2) $s$-$w$-irreflexive $\forall x[\neg_s(x <_w^s x)]$
(3) $s$-$w$-linear $\forall x\forall y[(x <_w^s y) \vee (x =_w^s y) \vee (y <_w^s x)]$

A formula $\Phi$ of $\mathcal{L}_w^s(\mathbb{R}_w^s)$ is in general of the form

$$\Phi = \Phi(X_1,\ldots,X_m,x_1,\ldots,x_q), \tag{6.3.6}$$

where $x_1,\ldots,x_q$ are the free consistent number variables of $\Phi$, i.e., variables not bound by a quantifier $\forall,\exists$ and $X_1,\ldots,X_m$ are the (free) consistent set variables of $\Phi$. Every formula in $\mathcal{L}(\mathbb{R}_w^s)$ has an standard interpretation in the structure $\mathbb{R}_w^s$; e.g., let

$\Phi(X)$ be the formula

$$\Phi(X) \equiv \forall y\Big[y \in_w^s X \Rightarrow_s \exists z\Big[\Big(z \,{}_w^s\!> 0_w^s\Big) \wedge \forall y_1[|y -_w^s y_1|_w^s <_w^s z \Rightarrow_s y_1 \in_w^s X]\Big]\Big] \tag{6.3.7}$$

and let $A \subseteq_w^s \mathbb{R}_w^s$, then $\Phi(A)$ expresses the fact that $A$ is open in $\mathbb{R}_w^s$;

# The properties of $w$-consistent $w$-ordered fields

**Definition 6.2.**. A $w$-consistent field is a nonempty $w$-consistent set $F_w^s$ containing at least

2 elements alongside the two binary operations of addition, $f_{+_w^s} : F \times_w^s F \to F$ such that $f_{+_w^s}(x,y) =_w^s x +_w^s y$ and multiplication $f_{\times_w^s}(x,y) =_w^s x \times_w^s y$ that satisfy all of the axioms below.

**I.Basic Properties of $w$-Consistent Equality**

1. $\forall x \in_w^s F[x =_w^s x]$.
2.$\forall x,y \in_w^s F[x =_w^s y \Rightarrow_s y =_w^s x]$.
3.For any $w$-consistent function $f(x_1,\ldots,x_n) : F \times_w^s \ldots \times_s F \to F,$ if $x_1 =_w^s y_1,\ldots,x_n =_w^s y_n$ then $f(x_1,\ldots,x_n) =_w^s f(y_1,\ldots,y_n)$.

,

**Theorem 6.3.1**. (**Generalized Łoś Theorem**) Let $\Phi_{m,q} = \Phi(X_1,\ldots,X_m,x_1,\ldots,x_q)$ be a formula of $\mathcal{L}_w(\mathbb{R}_w)$.

(**I**) Assume that $\Phi_{m,q}$ is not of the form $\Psi \wedge \neg_w\Psi$. Then for any $A_1,\ldots,A_m \subseteq_w \mathbb{R}$ and $f^1_{\mathcal{F}^{con}},\ldots,f^q_{\mathcal{F}^{con}} \in_w {}^{\#_w}\mathbb{R}$ :

$$\Phi\Big(A_1,\dots,A_m,f^1_{\mathcal{F}^w},\dots,f^q_{\mathcal{F}^w}\Big) \Leftrightarrow_s$$
$$\Leftrightarrow_s \Big\{v \in_w \wp^w | \Phi\Big(A_1,\dots,A_m,f^1(v),\dots,f^q(v)\Big)\Big\} \in_w \mathcal{F}^w. \quad (6.3.6)$$

(**II**) Assume that $\Phi_{m,q}$ is of the form $\Psi \wedge \neg_w \Psi$. Then for any $A_1,\dots,A_m \subseteq_w \mathbb{R}_w$ and $f^1_{\mathcal{F}^{con}},\dots,f^q_{\mathcal{F}^{con}} \in_w {}^{\#_w}\mathbb{R}_w$ :

$$\Phi_{[1]}\Big(A_1,\dots,A_m,f^1_{\mathcal{F}^w},\dots,f^q_{\mathcal{F}^w}\Big) \Leftrightarrow_s$$
$$\Leftrightarrow_s \Big\{v \in_w \wp^w | \Phi\Big(A_1,\dots,A_m,f^1(v),\dots,f^q(v)\Big)\Big\} \in_{w_{[1]}} \mathcal{F}^w. \quad (6.3.7)$$

**Proof.**(**I**)The proof is by induction on the number of logical symbols in $\Phi$. If $\Phi$ has no logical symbols, it is an elementary formula of the form (i)-(vi),and (6.3.6) then reduces to one of (), (), ( ) ,( ) , or (). If $\Phi$ contains logical symbols, then $\Phi$ is of the form $\Phi \wedge \Psi, \Phi \vee \Psi, \Phi \Rightarrow_s \Psi, \neg_s\Phi, \Phi \Rightarrow_w \Psi, \neg_w\Phi, \forall x(x \in_w \mathbb{R}_w), \exists x(x \in_w \mathbb{R}_w)$. The verification of (6.3.6) is, by induction, in each case reduced to an elementary property of the inconsistent $w$-ultrafilter $\mathcal{F}^w$.For example, if $\Phi = \Phi_1 \wedge \Phi_2$, (6.3.6) follows from the finite $w$-intersection property of the inconsistent $w$-ultrafilter $\mathcal{F}^w$. The case $\Phi = \neg_s\Phi_1$ uses in an essential way that $\mathcal{F}^w$ is an inconsistent $w$-ultrafilter namely, that

$$\wp^w \backslash_{w,s} \mathfrak{I} \in_w \mathcal{F}^w \Leftrightarrow_s \mathfrak{I} \notin^s_w \mathcal{F}^w \quad (6.3.7)$$

and

$$\wp^w \backslash_{w,s} \mathfrak{I} \in_{w_{[1]}} \mathcal{F}^w \Leftrightarrow_s \mathfrak{I} \notin^s_{w_{[1]}} \mathcal{F}^w \Leftrightarrow_s \mathfrak{I} \notin^s_w \mathcal{F}^w, \quad (6.3.8)$$

we remind that $\neg_s(a \wedge \neg_w a) \Leftrightarrow_s \neg_s a \vee a$.
The case $\Psi = \neg_w\Phi$ immediately from definition

$$\neg_w\Phi\Big(A_1,\dots,A_m,f^1_{\mathcal{F}^w},\dots,f^q_{\mathcal{F}^w}\Big) \Leftrightarrow_s$$
$$\Leftrightarrow_s \Big\{v \in_w \wp^w | \Phi\Big(A_1,\dots,A_m,f^1(v),\dots,f^q(v)\Big)\Big\} \notin^w_w \mathcal{F}^w. \quad (6.3.9)$$

Quantifiers offer no special difficulties.For example, if $\Phi = \exists x\Phi_1$ and let $\Phi$ have one free variable; we shall prove

$$\Phi(f_{\mathcal{F}^w}) \Leftrightarrow_s \{v \in_w \wp^w | \Phi(f(v))\} \in_w \mathcal{F}^w, \quad (6.3.10)$$

where $\Phi(f_{\mathcal{F}^w})$ is of the form $\exists x\Phi_1(x,f_{\mathcal{F}^w})$. Now $\Phi(f_{\mathcal{F}^w})$ is true in ${}^{\#_w}\mathbb{R}_w$ iff there is some $g_{\mathcal{F}^w} \in_w {}^{\#_w}\mathbb{R}_w$ such that $\Phi_1(g_{\mathcal{F}^w},f_{\mathcal{F}^w})$ is true in ${}^{\#_w}\mathbb{R}_w$. By the induction hypothesis this means that

$$\{v \in_w \wp | \Phi_1(g(v),f(v))\} \in_w \mathcal{F}^w. \quad (6.3.11)$$

But if $\Phi_1(g(v),f(v))$ is true in $\mathbb{R}_w$, then $\exists x\Phi_1(x,f(v))$ is also true in $\mathbb{R}_w$,i.e.,

$$\{v \in_w \wp^w | \Phi_1(g(v),f(v))\} \subseteq_w \{v \in_w \wp^w | \exists x\Phi_1(x,f(v))\}. \quad (6.3.12)$$

From (6.3.12) and the property (3) of inconsistent $w$-filters it follows that

$$\{v \in_w \wp^w | \Phi(f(v))\} \in_w \mathcal{F}^w. \quad (6.3.13)$$

In order to prove the converse, assume that: $\Re_\Phi =_w \{v \in \wp | \Phi(f(v))\} \in_w \mathcal{F}^w$. For each $v$ such that $v \in_w \Re_\Phi$ we choose some $\alpha_v \in_w \mathbb{R}$ such that $\Phi_1(\alpha_v, f(v))$ is true in $\mathbb{R}_w$. Let $g \in_w \mathbb{R}_w^{\wp^w}$ be a $w$-function $g : \wp^w \to \mathbb{R}_w$ such that $g(v) =_w \alpha_v$ for all $v \in_w \Re_\Phi$ and $g(v) =_w \beta$ otherwise, where $\beta$ is some arbitrary $w$-element of $\mathbb{R}$. Then we have:

$$\{v \in \wp^w | \Phi_1(g(v), f(v))\} \in_w \mathcal{F}^w. \quad (6.3.14)$$

Hence by the induction hypothesis we have $\Phi_1(g_{\mathcal{F}^{con}}, f_{\mathcal{F}^{con}})$ is true in ${}^{\#_w}\mathbb{R}_w$, i.e., we have $\Phi(f_{\mathcal{F}^w}) = \exists x \Phi_1(x, f_{\mathcal{F}^w})$ is true in ${}^{\#_w}\mathbb{R}_w$. Thus $\Phi(A_1, \ldots, A_m, r_1, \ldots, r_n)$ holds in $\mathbb{R}_w$ iff $\Phi({}^{\#_w}A_1, \ldots, {}^{\#_w}A_m, {}^{\#_w}r_1, \ldots, {}^{\#_w}r_n)$ holds in ${}^{\#_w}\mathbb{R}_w$.

(**II**) The case $\Psi = \neg_w \Phi$ immediately from definition

$$\begin{gathered}
\Psi\left(A_1, \ldots, A_m, f^1_{\mathcal{F}^w}, \ldots, f^q_{\mathcal{F}^w}\right) \Leftrightarrow_s \\
\Phi\left(A_1, \ldots, A_m, f^1_{\mathcal{F}^w}, \ldots, f^q_{\mathcal{F}^w}\right) \wedge \neg_w \Phi\left(A_1, \ldots, A_m, f^1_{\mathcal{F}^w}, \ldots, f^q_{\mathcal{F}^w}\right) \Leftrightarrow_s \\
\Leftrightarrow_s \left[\left\{v \in_w \wp^w | \Phi\left(A_1, \ldots, A_m, f^1(v), \ldots, f^q(v)\right)\right\} \in_w \mathcal{F}^w\right] \wedge \\
\wedge \left[\left\{v \in_w \wp^w | \Phi\left(A_1, \ldots, A_m, f^1(v), \ldots, f^q(v)\right)\right\} \notin^w_w \mathcal{F}^w\right].
\end{gathered} \quad (6.3.15)$$

Quantifiers offer no special difficulties. For example, if $\Psi = \exists x \Phi_1$, where $\Phi_1 = \Phi \wedge \neg_w \Phi$ and let $\Phi$ have one free variable; we shall prove

$$\begin{gathered}
\Psi(f_{\mathcal{F}^w}) \Leftrightarrow_s \\
[\{v \in_w \wp^w | \Phi(f(v))\} \in_w \mathcal{F}^w] \wedge [\{v \in_w \wp^w | \Phi(f(v))\} \notin^w_w \mathcal{F}^w],
\end{gathered} \quad (6.3.16)$$

where $\Psi(f_{\mathcal{F}^w})$ is of the form $\exists x \Phi_1(x, f_{\mathcal{F}^w})$. Now $\Psi(f_{\mathcal{F}^w})$ is true in ${}^{\#_w}\mathbb{R}_w$ iff there is some $g_{\mathcal{F}^w} \in_w {}^{\#_w}\mathbb{R}_w$ such that $\Psi(g_{\mathcal{F}^w}, f_{\mathcal{F}^w})$ is true in ${}^{\#_w}\mathbb{R}_w$. By the induction hypothesis this means that

$$\left[\{v \in_w \wp | \Phi(g(v), f(v))\} \in_w \mathcal{F}^w\right] \wedge \left[\{v \in_w \wp | \Phi(g(v), f(v))\} \notin^w_w \mathcal{F}^w\right]. \quad (6.3.17)$$

and therefore

$$\left[\{v \in_w \wp | \Phi(g(v), f(v))\} \in_w \mathcal{F}^w\right] \wedge \left[\{v \in_w \wp | \neg_w \Phi(g(v), f(v))\} \in_w \mathcal{F}^w\right]. \quad (6.3.18)$$

Remind that (6.3.18) means that $\Phi_1(g(v), f(v))$ is true in $\mathbb{R}_w$, i.e., both $\Phi(g(v), f(v))$ and $\neg_w \Phi(g(v), f(v))$ is true in $\mathbb{R}_w$

But if both $\Phi(g(v), f(v))$ and $\neg_w \Phi(g(v), f(v))$ is true in $\mathbb{R}_w$, then both $\exists x \Phi(x, f(v))$ and $\exists x[\neg_w \Phi(x, f(v))]$ is also true in $\mathbb{R}_w$, i.e.,

$$\{v \in_w \wp^w | \Phi(g(v), f(v))\} \subseteq_w \{v \in_w \wp^w | \exists x \Phi(x, f(v))\}. \quad (6.3.19)$$

and

$$\{v \in_w \wp^w | \neg_w \Phi(g(v), f(v))\} \subseteq_w \{v \in_w \wp^w | \exists x[\neg_w \Phi(x, f(v))]\}. \quad (6.3.20)$$

From (6.3.19)-(6.3.20) and the property (3) of inconsistent $w$-filters it follows that

$$\{v \in_w \wp^w | \Phi(f(v))\} \in_w \mathcal{F}^w \quad (6.3.13)$$

and

$$\{v \in_w \wp^w | \neg_w \Phi(f(v))\} \in_w \mathcal{F}^w \quad (6.3.13)$$

In order to prove the converse, assume that: $\Re_\Phi =_w \{v \in \wp | \Phi(f(v))\} \in_w \mathcal{F}^w$. For each $v$

such that $\nu \in_w \Re_\Phi$ we choose some $\alpha_\nu \in_w \mathbb{R}$ such that $\Phi_1(\alpha_\nu, f(\nu))$ is true in $\mathbb{R}_w$. Let $g \in_w \mathbb{R}_w^{\wp^w}$ be a $w$-function $g : \wp^w \to \mathbb{R}_w$ such that $g(\nu) =_w \alpha_\nu$ for all $\nu \in_w \Re_\Phi$ and $g(\nu) =_w \beta$ otherwise,where $\beta$ is some arbitrary $w$-element of $\mathbb{R}$. Then we have:

$$\{\nu \in \wp^w | \Phi_1(g(\nu), f(\nu))\} \in_w \mathcal{F}^w. \qquad (6.3.14)$$

# VI.3.2.The Generalized Transfer Principle

**Theorem 6.3.1**.(Transfer principle).Let $\Phi(X_1,\ldots,X_m,x_1,\ldots,x_q)$ be a formula of

$\mathcal{L}(\mathbb{R}_w)$. Then for any $A_1,\ldots,A_m \subseteq_w \mathbb{R}$ and $r_1,\ldots,r_n \in_w \mathbb{R}_w$, $\Phi(A_1,\ldots,A_m,r_1,\ldots,r_n)$ holds in $\mathbb{R}_w$

iff ${}^*\Phi = \Phi({}^*A_1,\ldots,{}^*A_m,{}^*r_1,\ldots,{}^*r_n)$ holds in ${}^*\mathbb{R}$,i.e.,

$$\Phi(A_1,\ldots,A_m,r_1,\ldots,r_n) \Leftrightarrow_s \Phi({}^*A_1,\ldots,{}^*A_m,{}^*r_1,\ldots,{}^*r_n). \qquad (6.2.23)$$

**Proof**.From (6.2.16) we get at once

$$\Phi({}^*A_1,\ldots,{}^*A_m,{}^*r_1,\ldots,{}^*r_n) \Leftrightarrow_s \{\nu \in \wp | \Phi(A_1,\ldots,A_m,r_1,\ldots,r_n)\} \in_s \mathcal{F}^{con} \qquad (6.2.24)$$

But the set $\{\nu \in \wp | \Phi(A_1,\ldots,A_m,r_1,\ldots,r_n)\}$ is equal to $\wp \in \mathcal{F}^{con}$ if $\Phi$ is true of

We consider now the $w_0$-consistent $w_0$-reals $\mathbb{R}_{w_0}$ as a structure

$$\{\mathbb{R}_{w_0}, +_{w_0}, \times_{w_0}, <_{w_0}, |\cdot|_{w_0}, 0_{w_0}, 1_{w_0}\}, \qquad (6.2.17)$$

The structure $\mathbb{R}_{w_0}$ has an associated simple language $\mathcal{L}(\mathbb{R}_{w_0})$ that can be used to describe

the kind of properties of $\mathbb{R}$ that are preserved in wekly consistent sense

under the embedding:

$$\#_{w_0} : \mathbb{R}_{w_0} \hookrightarrow {}^{\#_{w_0}}\mathbb{R}. \qquad (6.2.18)$$

**Theorem 6.2.2**. (**Generalized Łoś Theorem**) Let $\Phi(X_1,\ldots,X_m,x_1,\ldots,x_q)$ be a formula of

$\mathcal{L}(\mathbb{R}_{w_0})$. Then for any $A_1,\ldots,A_m \subseteq_s \mathbb{R}_{w_0}$ and $f^1_{\mathcal{F}^{w_0}},\ldots,f^q_{\mathcal{F}^{w_0}} \in_{w_0} {}^{\#_{w_0}}\mathbb{R}$

$$\Phi\left(A_1,\ldots,A_m,f^1_{\mathcal{F}^{w_0}},\ldots,f^q_{\mathcal{F}^{w_0}}\right) \Leftrightarrow_s$$

$$\Leftrightarrow_s \left\{\nu \in_{w_0} \wp^{w_0} | A_1,\ldots,A_m,f^1(\nu),\ldots,f^q(\nu)\right\} \in_{w_0} \mathcal{F}^{w_0}. \qquad (6.2.16)$$

$$\Phi\left(A_1,\ldots,A_m,f^1_{\mathcal{F}^{w_0}},\ldots,f^q_{\mathcal{F}^{w_0}}\right) \wedge \neg_w \Phi\left(A_1,\ldots,A_m,f^1_{\mathcal{F}^{w_0}},\ldots,f^q_{\mathcal{F}^{w_0}}\right) \Leftrightarrow_s$$

**Remark**.**6**.**2**.**2**.

# VI.3.The $\#_w$-transfer.

## VI.3.1.The $\#_w$-embedding

We consider the inconsistent $w$-reals as a structure

$$\{\mathbb{R}_w, +_w, \times_w, <_w, |\cdot|, 0_w, 1_w\}, \qquad (6.2.)$$

(**II**)The $\#_w$-embedding of (6.1.10) sends $0_w$ to ${}^{\#_w}0_w =_w \mathbf{0}^w_{\mathcal{F}^{inc}} \triangleq \mathbf{0}_w$ and $1_w$ to ${}^{\#_w}1 =_w \mathbf{1}^w_{\mathcal{F}^{inc}} \triangleq \mathbf{1}_w$. We must lift the operations and relations of $\mathbb{R}$ to ${}^{\#_w}\mathbb{R}$.

**Definition 6**.**3**.**1**.We get the clue from (6.1.6), which tells us when two elements $f^w_{\mathcal{F}^{inc}}$ and $g^w_{\mathcal{F}^{inc}}$, of ${}^{\#_w}\mathbb{R}$ are weakly $w$-equal in a weak paraconsistent sense iff :

$$f^w_{\mathcal{F}^{inc}} =_w g^w_{\mathcal{F}^{inc}} \Leftrightarrow_w f^w \sim_{\mathcal{F}^w} g^w \qquad (6.3.1)$$

e.g.,

$$f^w_{\mathcal{F}^{inc}} =_w g^w_{\mathcal{F}^{inc}} \Leftrightarrow_w \left\{v \in_w \wp^w | f^w(v) =_w g^w(v)\right\} \in_w \mathcal{F}^{inc}. \qquad (6.3.2)$$

**Remark**.**6**.**3**.**1**.

**Definition 6**.**3**.**2**.Two elements $f^w_{\mathcal{F}^{inc}}$ and $g^w_{\mathcal{F}^{inc}}$, of ${}^{\#_w}\mathbb{R}$ are $w$-equivalent in a weak inconsistent sense :

In a similar way we extend $<$ to ${}^{\#_w}\mathbb{R}$ by setting for arbitrary $f^w_{\mathcal{F}^{inc}}$,and $g^w_{\mathcal{F}^{inc}}$,in ${}^{\#_w}\mathbb{R}$:

$$f^w_{\mathcal{F}^{inc}} <_w g^w_{\mathcal{F}^{inc}} \Leftrightarrow \left\{v \in_w \wp^{inc} | f^w(v) <_w g^w(v)\right\} \in_w \mathcal{F}^{inc}. \qquad (6.3.2)$$

With this definition of $<_w$ in ${}^{\#_w}\mathbb{R}$ we easily show that the extended domain ${}^{\#_w}\mathbb{R}$ is linearly $w$-ordered $w$-inconsistent field. As an example we verify $w$-transitivity of $<_w$ in ${}^{\#_w}\mathbb{R}$. Let $f^w_{\mathcal{F}^{inc}} <_w g^w_{\mathcal{F}^{inc}}$, and $g^w_{\mathcal{F}^{inc}} <_w h^w_{\mathcal{F}^{inc}}$,i.e.,

$$\begin{aligned} D^w_1 &=_w \{v \in_w \wp^{inc} | f^w(v) <_w g^w(v)\} \in_w \mathcal{F}^{inc}, \\ D^w_2 &=_w \{v \in_w \wp^{inc} | g^w(v) <_w h^w(v)\} \in_w \mathcal{F}^{inc} \end{aligned} \qquad (6.3.3)$$

By the finite intersection property (ii),[see Definition 6.1.1.(ii)] $\Rightarrow_w D^w_1 \cap_w D^w_2 \in_w \mathcal{F}^{inc}$. If $v \in_w D_1 \cap_w D_2$, then $f(v) <_w g(v)$ and $g^w(v) <_w h^w(v)$; hence by transitivity of $<$ in $\mathbb{R}$,

$$\left[f^w(v) <_w g^w(v)\right] \wedge [g^w(v) <_w h^w(v)] \quad \Rightarrow_{w_0} \quad f^w(v) <_w h^w(v). \qquad (6.3.4)$$

Thus

$$\Rightarrow_{w_0} D_1 \cap_w D_2 \subseteq_w \{v \in_w \wp^{inc} | f^w(v) <_w h^w(v)\} \qquad (6.3.5)$$

The closure property (3) then tells us that:

$$\Rightarrow_{w_0} f^w_{\mathcal{F}^{inc}} <_w h^w_{\mathcal{F}^{inc}} \quad (6.3.6)$$

## VI.4.The $\#_{w_n}$ transfer and $\#_{w_n}$-embedding

## VI.5. The Extendent Paralogical Universe.

### VI.5.1. The inconsistent superstructures over universal set.

**Definition 6.5.1**. The superstructure over inconsistent set, or inconsistent universe $\mathbf{S}^{inc}$, denoted by $\mathbf{V}^{w}(\mathbf{S}^{inc}),\mathbf{V}^{w_0}(\mathbf{S}^{inc}),\mathbf{V}^{w_1}(\mathbf{S}^{inc}),$etc. is defined by the following canonical recursion:

$$\begin{gathered}\mathbf{V}_1^{w}(\mathbf{S}^{inc}) =_w \mathbf{S}^{inc},\\ \mathbf{V}_{n+1}^{w}(\mathbf{S}^{inc}) =_w \mathbf{V}_n^{w}(\mathbf{S}^{inc}) \cup_w \{X|X \subseteq_w \mathbf{V}_n^{w}(\mathbf{S}^{inc})\}_w,\\ \mathbf{V}^{w}(\mathbf{S}^{inc}) =_w w\text{-}\bigcup_n \mathbf{V}_n(\mathbf{S}^{inc}).\end{gathered} \tag{6.5.1}$$

$$\begin{gathered}\mathbf{V}_1^{w_0}(\mathbf{S}^{inc}) =_{w_0} \mathbf{S}^{inc},\\ \mathbf{V}_{n+1}^{w_0}(\mathbf{S}^{inc}) =_{w_0} \mathbf{V}_n^{w_0}(\mathbf{S}^{inc}) \cup_{w_0} \{X|X \subseteq_{w_0} \mathbf{V}_n^{w_0}(\mathbf{S}^{inc})\}_{w_0},\\ \mathbf{V}^{w_0}(\mathbf{S}^{inc}) =_w w_0\text{-}\bigcup_n \mathbf{V}_n^{w_0}(\mathbf{S}^{inc}).\end{gathered} \tag{6.5.2}$$

$$\begin{gathered}\mathbf{V}_1^{w_m}(\mathbf{S}^{inc}) =_{wm} \mathbf{S}^{inc},\\ \mathbf{V}_{n+1}^{w_m}(\mathbf{S}^{inc}) =_{w_0} \mathbf{V}_n^{w_m}(\mathbf{S}^{inc}) \cup_{w_0} \{X|X \subseteq_{w_m} \mathbf{V}_n^{w_m}(\mathbf{S}^{inc})\}_{w_0},\\ \mathbf{V}^{w_m}(\mathbf{S}^{inc}) =_w w_m\text{-}\bigcup_n \mathbf{V}_n^{w_m}(\mathbf{S}^{inc}),\\ m = 1,2,\ldots\ .\end{gathered} \tag{6.5.3}$$

The extended inconsistent nonstandard universe of paraconsistent nonstandard analysis will be obtained by postulating: the extensions ${}^{\#_w}\mathbb{R} \supset \mathbb{R}, {}^{\#_{w_0}}\mathbb{R} \supset \mathbb{R}, {}^{\#_{w_n}}\mathbb{R} \supset \mathbb{R},$ and postulating the embeddings

$$\begin{gathered}{}^{\#_w} : \mathbf{V}^{w}(\mathbf{S}^{inc} \supset \mathbb{R}) \hookrightarrow \mathbf{V}^{w}({}^{\#_w}\mathbf{S}^{inc}),\\ {}^{\#_w} : \mathbf{V}^{w_0}(\mathbf{S}^{inc} \supset \mathbb{R}) \hookrightarrow \mathbf{V}^{w_0}({}^{\#_{w_0}}\mathbf{S}^{inc}),\\ {}^{\#_w} : \mathbf{V}^{w_n}(\mathbf{S}^{inc} \supset \mathbb{R}) \hookrightarrow \mathbf{V}^{w_n}({}^{\#_{w_n}}\mathbf{S}^{inc}).\end{gathered} \tag{6.5.4}$$

We shall now extend the construction of the inconsistent ultrafilter to demonstrate that will have properties similar to the embedding ${}^{\#_w} : \mathbb{R} \hookrightarrow {}^{\#_w}\mathbb{R}$ constructed in Subsections 6.1-6.4.

**Remark.6.5.1.**First of all we assume the following principle.**EXTENSION PRINCIPLE**.

(i) ${}^{\#_w}\mathbb{R}$ is a proper $w$-inconsistent extension of $\mathbb{R}$ and ${}^{\#_w}r =_w r$ for all $r \in \mathbb{R}$,

(ii) ${}^{\#_{w_0}}\mathbb{R}$ is a proper $w$-consistent extension of $\mathbb{R}$ and ${}^{\#_{w_0}}r =_{w_0} r$ for all $r \in \mathbb{R}$,

(iii) ${}^{\#_{w_n}}\mathbb{R}$ is a proper $w_n$-consistent extension of $\mathbb{R}$ and ${}^{\#_{w_n}}r =_{w_n} r$ for all $r \in \mathbb{R}, n \in \mathbb{N}$.

## VI.5.2. The Bounded Paralogical Ultrapowers.

## VI.5.2.1. The Bounded Consistent Ultrapowers.

Remind the following definitions.

**Definition 6.5.2**. A sequence $\mathbf{A} = \langle A_v \rangle_{v \in \wp}$ of elements of $\mathbf{V}(\mathbf{S}^{con}) = \mathbf{V}(\mathbf{S}^{con} \supsetneqq \mathbb{R})$ is bounded if there is a fixed $n \geq 1$ such that each $A_v \in \mathbf{V}_n(\mathbf{S}^{con} \supsetneqq \mathbb{R})$.

**Remark.6.5.2.**

**Definition 6.5.3.**Two bounded sequences $\mathbf{A}$ and $\mathbf{B}$ are equivalent with respect to the free
consistent ultrafilter $\mathcal{F}$, in symbols $\mathbf{A} \sim_{\mathcal{F}} \mathbf{B}$, iff

$$\{v \in \wp | A_v = B_v\} \in \mathcal{F}. \tag{6.5.5}$$

We let $\mathbf{A}_{\mathcal{F}}$ denote the equivalence class of $A$ and define the bounded ultrapower by

$$\mathbf{V}(\mathbf{S}^{con})^{\wp}/\mathcal{F} = \left\{\mathbf{A}_{\mathcal{F}} | \mathbf{A} \text{ is a bounded } \mathbf{V}(\mathbf{S}^{con})\text{-sequence}\right\}. \tag{6.5.6}$$

**Definition 6.5.4.**We define the membership relation $\in_{\mathcal{F}}$ in the ultrapower by

$$\mathbf{A}_{\mathcal{F}} \in_{\mathcal{F}} \mathbf{B}_{\mathcal{F}} \text{ iff } \{v \in \wp | A_v \in B_v\} \in \mathcal{F}. \tag{6.5.7}$$

There is a natural proper embedding

$$i : \mathbf{V}(\mathbf{S}^{con}) \hookrightarrow \mathbf{V}(\mathbf{S}^{con})^{\wp}/\mathcal{F} \tag{6.5.8}$$

namely let $i(\mathbf{A}) = \langle A \rangle_{\mathcal{F}}$, the equivalence class corresponding to the constant sequence $\mathbf{A} = \langle A \rangle$.

## VI.5.2.2. The Bounded Paralogical $w$-Ultrapowers.

**Definition 6.5.5**. (i) A $w$-sequence $\mathbf{A} \triangleq \langle A_v \rangle_{v \in_w \wp^w}$ of $w$-elements of $\mathbf{V}^w(\mathbf{S}^{inc}{}_w \supset \mathbb{R}_w)$ is $w$-bounded if there is a fixed $n \geq 1$ such that each $A_v \in_w \mathbf{V}_n^w(\mathbf{S}^{inc}{}_w \supset \mathbb{R}_w)$.

**Remark.6.5.3.**

**Definition 6.5.6.**Let $\mathcal{F}^w$ be a free $w$-ultrafilter on $\wp_w$. Two $w$-bounded $w$-sequences $\mathbf{A}_w$ and $\mathbf{B}_w$ are $w$-equivalent with respect to the free inconsistent ultrafilter $\mathcal{F}^w$, in symbols

$$\mathbf{A}_w \sim_{\mathcal{F}^w}^w \mathbf{B}_w \text{ iff } \{i \in_w \wp_w | A_v =_w B_v\} \in_w \mathcal{F}^w. \tag{6.5.9}$$

**Definition 6.5.7.**We let $\mathbf{A}_{\mathcal{F}^w}^w$ denote the $w$-equivalence class of $\mathbf{A}$ and define the $w$-bounded $w$-ultrapower by

$$\mathbf{V}^w(\mathbf{S}^{inc})^{\wp_w}/\mathcal{F}^w =_w \left\{\mathbf{A}_{\mathcal{F}^w}^w | \mathbf{A}_w \text{ is a } w\text{-bounded } \mathbf{V}^w(\mathbf{S}^{inc})\text{-sequence}\right\}_w \tag{6.5.10}$$

There is a natural proper embedding

$$i_w : \mathbf{V}(\mathbf{S}^{inc}) \hookrightarrow \mathbf{V}^w(\mathbf{S}^{inc})^{\wp_w}/\mathcal{F}^w \tag{6.5.11}$$

namely let $i_w(\mathbf{A}_w) =_w \langle A \rangle_{\mathcal{F}^w}$, the equivalence class corresponding to the constant $w$-sequence $\mathbf{A}_w =_w \langle A \rangle_w$.

**Definition 6.5.8.**We define the $w$-membership relation $\in_{\mathcal{F}^w}$ in the $w$-ultrapower by

$$\mathbf{A}_{\mathcal{F}^w}^w \in_{\mathcal{F}^w} \mathbf{B}_{\mathcal{F}^w}^w \text{ iff } \{i \in_w \wp_w | A_v \in_w B_v\} \in_w \mathcal{F}^w, \tag{6.5.12}$$

## VI.5.2.3. The Bounded Paralogical $w_0$-Ultrapowers.

**Definition 6.5.9.**A $w_0$-sequence $\mathbf{A}_{w_0} \triangleq \langle A_v \rangle_{v \in_{w_0} \wp^{w_0}}$ of $w_0$-elements of $\mathbf{V}^{w_0}(\mathbf{S}^{inc}{}_{w_0} \supset \mathbb{R}_{w_0})$ is $w_0$-bounded if there is a fixed $n \geq 1$ such that each $A_v \in_{w_0} \mathbf{V}_n^{w_0}(\mathbf{S}^{inc}{}_{w_0} \supset \mathbb{R}_{w_0})$.

**Remark**.6.5.4.

**Definition 6.5.10**.Let $\mathcal{F}^{w_0}$ be a free $w_0$-ultrafilter on $\wp_{w_0}$.Two $w_0$-bounded $w_0$-sequences $\mathbf{A}_{w_0}$ and $\mathbf{B}_{w_0}$ are $w_0$-equivalent with respect to the free $w_0$-consistent ultrafilter $\mathcal{F}^{w_0}$, in symbols

$$\mathbf{A}_{w_0} \sim^{w_0}_{\mathcal{F}^{w_0}} \mathbf{B}_{w_0} \text{ iff } \{i \in_{w_0} \wp_{w_0} | A_v =_{w_0} B_v\} \in_{w_0} \mathcal{F}^{w_0}. \qquad (6.5.13)$$

**Definition 6.5.11**.We let $A^{w_0}_{\mathcal{F}^{w_0}}$ denote the $w_0$-equivalence class of $A$ and define the $w_0$-bounded $w_0$-ultrapower by

$$\mathbf{V}^{w_0}(\mathbf{S}^{inc})^{\wp_{w_0}}/\mathcal{F}^{w_0} =_{w_0} \left\{A^{w_0}_{\mathcal{F}^{w_0}} | \mathbf{A}_{w_0} \text{ is a } w\text{-bounded } \mathbf{V}^{w_0}(\mathbf{S}^{inc})\text{-sequence}\right\} \qquad (6.5.14)$$

There is a natural proper $w_0$-embedding

$$i_{w_0} : \mathbf{V}^{w_0}(\mathbf{S}^{inc}) \hookrightarrow \mathbf{V}^{w_0}(\mathbf{S}^{inc})^{\wp_{w_0}}/\mathcal{F}^{w_0}, \qquad (6.5.15)$$

namely let $i_{w_0}(\mathbf{A}_{w_0}) = \langle A \rangle_{\mathcal{F}^{w_0}}$,the $w_0$-equivalence class corresponding to the constant $w_0$-sequence $\mathbf{A}_{w_0} \triangleq \langle A \rangle_{w_0}$.

**Definition 6.5.12**.We define the $w_0$-membership relation $\in_{\mathcal{F}^{w_0}}$ in the $w_0$-ultrapower by

$$\mathbf{A}^{w_0}_{\mathcal{F}^{w_0}} \in_{\mathcal{F}^{w_0}} \mathbf{B}^{w_0}_{\mathcal{F}^{w_0}} \text{ iff } \{v \in_{w_0} \wp_{w_0} | A_v \in_{w_0} B_v\} \in_{w_0} \mathcal{F}^{w_0}. \qquad (6.5.16)$$

# VI.5.2.4. The Bounded Paralogical $w_n$-Ultrapowers.

**Definition 6.5.13**.A $w_n$-sequence $\mathbf{A} \triangleq \langle A_v \rangle_{v \in_{w_n} \wp^{w_n}}$ of $w_n$-elements of $\mathbf{V}^{w_n}(\mathbf{S}^{inc}{}_w \supset \mathbb{R}_w)$ is $w_n$-bounded if there is a fixed $n \geq 1$ such that each $A_v \in_{w_n} \mathbf{V}^{w_n}_n(\mathbf{S}^{inc})$.

**Remark**.6.5.5.

**Definition 6.5.14**.Let $\mathcal{F}^{w_n}$ be a free $w_n$-ultrafilter on $\wp_{w_n}$.Two $w_n$-bounded $w_n$-sequences $\mathbf{A}$ and $\mathbf{B}$ are $w_n$-equivalent with respect to the free inconsistent ultrafilter $\mathcal{F}^{w_n}$, in symbols $\mathbf{A} \sim^{w_n}_{\mathcal{F}^{w_n}} \mathbf{B}$, iff $\{i \in_{w_n} \wp_{w_n} | A_v =_{w_n} B_v\} \in_{w_n} \mathcal{F}^{w_n}$.

**Definition 6.5.15**.We let $\mathbf{A}^{w_n}_{\mathcal{F}^{w_n}}$ denote the $w_n$-equivalence class of $\mathbf{A}$ and define the $w_n$-bounded $w_n$-ultrapower by

$$\mathbf{V}^{w_n}(\mathbf{S}^{inc})^{\wp_{w_n}}/\mathcal{F}^{w_n} =_{w_n} \left\{\mathbf{A}^{w_n}_{\mathcal{F}^{w_n}} | \mathbf{A} \text{ is a } w_n\text{-bounded } \mathbf{V}^{w_n}(\mathbf{S}^{inc})\text{-sequence}\right\}. \qquad (6.5.17)$$

There is a natural proper $w_n$-embedding

$$i_{w_n} : \mathbf{V}^{w_n}(\mathbf{S}^{inc}) \hookrightarrow \mathbf{V}^{w_n}(\mathbf{S}^{inc})^{\wp_{w_n}}/\mathcal{F}^{w_n}, \qquad (6.5.18)$$

$n \in \mathbb{N}$ namely let $i_{w_n}(\mathbf{A}) =_{w_n} \langle A \rangle_{\mathcal{F}^{w_n}}$,the $w_n$-equivalence class corresponding to the constant $w_n$-sequence $\mathbf{A}_{w_n} \triangleq \langle A \rangle_{w_n}$.

**Definition 6.5.16**.we define the $w_n$-membership relation $\in_{\mathcal{F}^{w_n}}$ in the $w_n$-ultrapower by

$$\mathbf{A}^{w_n}_{\mathcal{F}^{w_n}} \in_{\mathcal{F}^{w_n}} \mathbf{B}^{w_n}_{\mathcal{F}^{w_n}} \text{ iff } \{v \in_{w_n} \wp_{w_n} | A_v \in_{w_0} B_v\} \in_{w_n} \mathcal{F}^{w_n}. \qquad (6.5.19)$$

# VI.5.3.The embedding $\mathbf{V}^w(\mathbf{S}^{inc})^{\wp_w}/\mathcal{F}^w$ into $\mathbf{V}^w({}^{\#_w}\mathbf{S}^{inc})$,etc.

# VI.5.3.1.Classical embedding $\mathbf{V}(\mathbf{S}^{con})^{\wp}/\mathcal{F}^{con}$ into $\mathbf{V}({}^*\mathbf{S}^{con})$

Let us consider now the classical embedding $\mathbf{V}(\mathbf{S}^{con})^{\wp}/\mathcal{F}$ into $\mathbf{V}({}^*\mathbf{S}^{con})$.Remind that ${}^*\mathbf{S}^{con} \supsetneqq {}^*\mathbb{R}$ is the bounded ultrapower $\mathbf{V}(\mathbf{S}^{con})^{\wp}/\mathcal{F}$.

**Remark**.6.5.6.Note that in classical case the bounded ultrapower $\mathbf{V}(\mathbf{S}^{con})^{\wp}/\mathcal{F}$ alwais will

not be the same as the full superstructure $\mathbf{V}({}^*\mathbf{S}^{con}) \supsetneqq ({}^*\mathbb{R})$

Remind the construction of canonical embedding

$$j : \mathbf{V}(\mathbf{S}^{con})^{\wp}/\mathcal{F} \hookrightarrow \mathbf{V}({}^{*}\mathbf{S}^{con}) \tag{6.5.20}$$

such that: (i) $j$ is the identity on ${}^{*}\mathbf{S}^{con} \supsetneqq {}^{*}\mathbb{R}$ and (ii) if $\mathbf{A}_{\mathcal{F}} \notin {}^{*}\mathbf{S}^{con}$ then

$$j(\mathbf{A}_{\mathcal{F}}) = \{j(\mathbf{B}_{\mathcal{F}})|\mathbf{B}_{\mathcal{F}} \in \mathbf{A}_{\mathcal{F}}\}. \tag{6.5.21}$$

This means that the relation $\in_{\mathcal{F}}$ in the ultrapower is mapped into the ordinary membership
relation in $\mathbf{V}({}^{*}\mathbf{S}^{con})$. The embedding $j$ is constructed in stages. Let

$$\mathbf{V}_k(\mathbf{S}^{con})^{\wp}/\mathcal{F} = \left\{\mathbf{A}_{\mathcal{F}}|\mathbf{A} \text{ is a sequence from } \mathbf{V}_k(\mathbf{S}^{con})\right\} \tag{6.5.22}$$

Then the bounded ultrapower is the union of the chain

$${}^{*}\mathbb{R} \subset {}^{*}\mathbf{S}^{con} = \mathbf{V}_1(\mathbf{S}^{con})^{\wp}/\mathcal{F} \subseteq \ldots \mathbf{V}_k(\mathbf{S}^{con})^{\wp}/\mathcal{F} \subseteq \ldots \tag{6.5.23}$$

and we can define $j$ by induction. For $k = 1$, the embedding $j$ must be the identity. If $\mathbf{A}_{\mathcal{F}} \in \mathbf{V}_{k+1}(\mathbf{S}^{con})^{\wp}/\mathcal{F}$ and $\mathbf{A}_{\mathcal{F}} \notin_{\mathcal{F}} {}^{*}\mathbf{S}^{con}$ we simply set $j(\mathbf{A}_{\mathcal{F}}) = \left\{j(\mathbf{B}_{\mathcal{F}})|\mathbf{B}_{\mathcal{F}} \in_{\mathcal{F}} \mathbf{A}_{\mathcal{F}}\right\}$ This
makes sense: if $\mathbf{B}_{\mathcal{F}} \in_{\mathcal{F}} \mathbf{A}_{\mathcal{F}}$ it follows from (6.5.7 ) that $\{v \in \wp|B_v \in \mathbf{V}_k(\mathbf{S}^{con})\} \in \mathcal{F}$, i.e., $\mathbf{B}_{\mathcal{F}} \in \mathbf{V}_k(\mathbf{S}^{con})^{\wp}/\mathcal{F}$, which means that $j(\mathbf{B}_{\mathcal{F}})$ is defined at a previous stage of the inductive
construction. Combining $i$ and $j$ we get a model of the extended nonstandard univers

$$\begin{array}{ccc} \mathbf{V}(\mathbf{S}^{con})^{\wp}/\mathcal{F} & & \mathbf{V}({}^{*}\mathbf{S}^{con}) \\ \nwarrow & & \nearrow \\ i & & * \\ & \mathbf{V}(\mathbf{S}^{con}) & \end{array} \tag{6.5.24}$$

where ${}^{\#_w}\mathbf{A} =_w j_w(i_w(\mathbf{A}))$, for any $\mathbf{A} \in_w \mathbf{V}^w(\mathbf{S}^{con})$. Here $\mathbf{V}^w(\mathbf{S}^{con} \supsetneqq \mathbb{R})$ and $\mathbf{V}^w({}^{\#_w}\mathbf{S}^{inc})$ are connected by a classical transfer principle.

**Theorem 6.5.1**.(**TRANSFER PRINCIPLE**) Let $A_1,\ldots,A_n \in \mathbf{V}(\mathbf{S}^{con})$. Any $\mathcal{L}(\mathbf{V}(\mathbf{S}^{con}))$ statement $\Phi$ that is true of $A_1,\ldots,A_n$ in $\mathbf{V}(\mathbf{S}^{con})$ is true of ${}^{*}A_1,\ldots,{}^{*}A_n$ in $\mathbf{V}({}^{*}\mathbf{S}^{con})$ .

**Proof**. In the ultrapower model there are three structures involved, $\mathbf{V}(\mathbf{S}^{con})$, $\mathbf{V}(\mathbf{S}^{con})^{\wp}/\mathcal{F}$ and $\mathbf{V}({}^{*}\mathbf{S}^{con})$. Given any $\mathcal{L}(\mathbf{V}(\mathbf{S}^{con}))$ formula $\Phi(X,Y)$ (see Remark.6.5.7), we have explained
how to interpret it in the three structures. Notice that Łòs theorem, 1.1.3, immediately extends to the bounded ultrapower $\mathbf{V}(\mathbf{S}^{con})^{\wp}/\mathcal{F}$ by exactly the same proof; i.e., for any $\mathbf{A}_{\mathcal{F}},\mathbf{B}_{\mathcal{F}} \in \mathbf{V}(\mathbf{S}^{con})^{\wp}/\mathcal{F}$ we have

$$(\mathbf{i})\ \Phi(\mathbf{A}_{\mathcal{F}},\mathbf{B}_{\mathcal{F}}) \text{ iff } \{v \in \wp|\Phi(A_v,B_v)\} \in \mathcal{F},$$

from which transfer follows between $\mathbf{V}(\mathbf{S}^{con})$ and $\mathbf{V}(\mathbf{S}^{con})^{\wp}/\mathcal{F}$ exactly as in 1.1.4.
But Principle 1.2.4 asserts transfer between $\mathbf{V}(\mathbf{S}^{con})$ and $\mathbf{V}({}^{*}\mathbf{S}^{con})$. And in order to prove
this we need to replace the equivalence (i) by

$$(\mathbf{ii})\ \Phi(j(\mathbf{A}_{\mathcal{F}}),j(\mathbf{B}_{\mathcal{F}})) \text{ iff } \{v \in \wp|\Phi(A_v,B_v)\} \in \mathcal{F}.$$

But this is a rather immediate extension which follows from the fact that every element of,
say, $j(\mathbf{A}_{\mathcal{F}})$ in $\mathbf{V}({}^{*}\mathbf{S}^{con})$ is of the form $j(\mathbf{A}'_{\mathcal{F}})$ for some $\mathbf{A}'_{\mathcal{F}} \in \mathbf{V}(\mathbf{S}^{con})^{\wp}/\mathcal{F}$; see the construction of the $j$-map above. And once we have (1 1) the Transfer Principle 1.2.4 follows by the same argument as in 1.1.4.

**Remark**.6.5.7.The structure $\mathbf{S}^{con} \supset \mathbb{R}$ has an associated elementary language $\mathcal{L}(\mathbf{S}^{con})$, which we used to give the necessary precision to the transfer principle. We need a similar formal tool to state the extended transfer principle.The language $\mathcal{L}(\mathbf{V}(\mathbf{S}^{con}))$ will be an extension of the language $\mathcal{L}(\mathbf{S}^{con}) = \mathcal{L}(\mathbf{S}^{con} \supset \mathbb{R})$.

# VI.5.3.2.The $w$-embedding $\mathbf{V}^{w}(\mathbf{S}^{inc})^{\wp_w}/\mathcal{F}^{w}$ into $\mathbf{V}^{w}({}^{\#_w}\mathbf{S}^{inc})$

Let ${}^{\#_w}\mathbf{S}^{inc}$ be the ($w$-bounded) $w$-ultrapower $\mathbf{V}^{w}(\mathbf{S}^{inc})^{\wp_w}/\mathcal{F}^{w}$.

**Remark**.6.5.8.Note that in contrast with a classical case $\mathbf{V}^{w}(\mathbf{S}^{inc})^{\wp_w}/\mathcal{F}^{w}$ not alwais will not be the same as the full $w$-superstructure $\mathbf{V}^{w}({}^{\#_w}\mathbf{S}^{inc})$.

We shall now construct an $w$-embedding

$$j_w : \mathbf{V}^{w}(\mathbf{S}^{inc})^{\wp_w}/\mathcal{F}^{w} \hookrightarrow \mathbf{V}^{w}({}^{\#_w}\mathbf{S}^{inc}) \quad (6.5.25)$$

such that: (i) $j_w$ is the $w$-identity on ${}^{\#_w}\mathbf{S}^{inc}$ and (ii) if $A^{w}_{\mathcal{F}^w} \notin^{w}_{\mathcal{F}^w} {}^{\#_w}\mathbf{S}^{inc}$, then

$$j_w(\mathbf{A}^{w}_{\mathcal{F}^w}) =_w \left\{ j_w(\mathbf{B}^{w}_{\mathcal{F}^w}) | \mathbf{B}^{w}_{\mathcal{F}^w} \in_{\mathcal{F}^w} \mathbf{A}^{w}_{\mathcal{F}^w} \right\}_w . \quad (6.5.26)$$

This means that the relation $\in_{\mathcal{F}^w}$ in the $w$-ultrapower is mapped into the ordinary $w$-membership relation in $\mathbf{V}^{w}({}^{\#_w}\mathbf{S}^{inc})$.The $w$-embedding $j_w$ is constructed in stages. Let

$$\mathbf{V}^{w}_{k}(\mathbf{S}^{inc})^{\wp_w}/\mathcal{F}^{w} =_w \left\{ \mathbf{A}^{w}_{\mathcal{F}^w} | \ \mathbf{A} \text{ is a } w\text{-sequence from } \mathbf{V}^{w}_{k}(\mathbf{S}^{inc}) \right\}_w . \quad (6.5.27)$$

Then the bounded $w$-ultrapower is the $w$-union of the $w$-chain

$${}^{\#_w}\mathbf{S}^{inc} =_w \mathbf{V}^{w}_{1}(\mathbf{S}^{inc})^{\wp_w}/\mathcal{F}^{w} \subseteq_w \ldots \subseteq_w \mathbf{V}^{w}_{k}(\mathbf{S}^{inc})^{\wp_w}/\mathcal{F}^{w} \subseteq_w \ldots, \quad (6.5.28)$$

and we can define $j_w$ by induction. For $k = 1$, $j_w$ must be the $w$-identity. If $\mathbf{A}^{w}_{\mathcal{F}^w} \in_w \mathbf{V}^{w}_{k+1}(\mathbf{S}^{inc})^{\wp_w}/\mathcal{F}^{w}$ and $\mathbf{A}^{w}_{\mathcal{F}^w} \notin^{w}_{\mathcal{F}^w} {}^{\#_w}\mathbf{S}^{inc}$ we simply set $j_w(\mathbf{A}^{w}_{\mathcal{F}^w}) =_w \left\{ j_w(\mathbf{B}^{w}_{\mathcal{F}^w}) | \mathbf{B}^{w}_{\mathcal{F}^w} \in_{\mathcal{F}^w} \mathbf{A}^{w}_{\mathcal{F}^w} \right\}_w$ This makes sense: if $\mathbf{B}^{w}_{\mathcal{F}^w} \in_{\mathcal{F}^w} \mathbf{A}^{w}_{\mathcal{F}^w}$it follows from (6.5.12 ) that $\{v \in_w \wp_w | B_v \in_w \mathbf{V}^{w}_{k}(\mathbf{S}^{inc})\}_w \in_w \mathcal{F}^{w}$,i.e.,$B^{w}_{\mathcal{F}^w} \in_w \mathbf{V}^{w}_{k}(\mathbf{S}^{inc})^{\wp_w}/\mathcal{F}^{w}$,which means that $j_w(\mathbf{B}^{w}_{\mathcal{F}^w})$ is defined at a previous stage of the inductive construction.

Combining $i_w$ and $j_w$ we get a model of the extended $w$-inconsistent nonstandard universe

$$\begin{array}{ccc} \mathbf{V}^{w}(\mathbf{S}^{inc})^{\wp_w}/\mathcal{F}^{w} & & \mathbf{V}^{w}({}^{\#_w}\mathbf{S}^{inc}) \\ & & {}_{w}\!\supseteq \\ \nwarrow & & \swarrow\nearrow \\ i_w & & \#_w \\ & \mathbf{V}^{w}(\mathbf{S}^{inc}) & \end{array} \quad (6.5.29)$$

where ${}^{\#_w}\mathbf{A} =_w j_w(i_w(A))$, for any $A \in_w \mathbf{V}^{w}(\mathbf{S}^{inc})$.Here $\mathbf{V}^{w}(\mathbf{S}^{inc} \supset \mathbb{R})$ and $\mathbf{V}^{w}({}^{\#_w}\mathbf{S}^{inc})$ are connected by $w$-inconsistent transfer principle.

# VI.5.3.3.The embedding $\mathbf{V}^{w_0}(\mathbf{S}^{inc})^{\wp_{w_0}}/\mathcal{F}^{w_0}$ into $\mathbf{V}^{w_0}({}^{\#_{w_0}}\mathbf{S}^{inc})$

Let ${}^{\#_{w_0}}\mathbf{S}^{inc}$ be the ($w_0$-bounded) $w_0$-ultrapower $\mathbf{V}^{w_0}(\mathbf{S}^{inc})^{\wp_{w_0}}/\mathcal{F}^{w_0}$.

**Remark**.6.5.9.Note that in contrast with a classical case $\mathbf{V}^{w_0}(\mathbf{S}^{inc})^{\wp_{w_0}}/\mathcal{F}^{w_0}$ not alwais will not be the same as the full $w_0$-superstructure $\mathbf{V}^{w_0}({}^{\#_{w_0}}\mathbf{S}^{inc})$.

We shall now construct an $w_0$-embedding

$$j_{w_0} : \mathbf{V}^{w_0}(\mathbf{S}^{inc})^{\wp_{w_0}}/\mathcal{F}^{w_0} \hookrightarrow \mathbf{V}^{w_0}({}^{\#_{w_0}}\mathbf{S}^{inc}) \quad (6.5.30)$$

such that: (i) $j_{w_0}$ is the $w_0$-identity on ${}^{\#_{w_0}}\mathbf{S}^{inc}$ and (ii) if $\mathbf{A}^{w_0}_{\mathcal{F}^{w_0}} \notin^{w}_{\mathcal{F}^{w_0}} {}^{\#_{w_0}}\mathbf{S}^{inc}$, then

$$j_{w_0}(\mathbf{A}^{w_0}_{\mathcal{F}^{w_0}}) =_{w_0} \left\{ j_{w_0}(\mathbf{B}^{w_0}_{\mathcal{F}^{w_0}}) | \mathbf{B}^{w_0}_{\mathcal{F}^{w_0}} \in_{\mathcal{F}^{w_0}} \mathbf{A}^{w_0}_{\mathcal{F}^{w_0}} \right\}_{w_0}. \quad (6.5.31)$$

This means that the relation $\in_{\mathcal{F}^{w_0}}$ in the $w_0$-ultrapower is mapped into the ordinary $w_0$-membership relation in $\mathbf{V}^{w_0}({}^{\#_{w_0}}\mathbf{S}^{inc})$. The $w_0$-embedding $j_{w_0}$ is constructed in stages. Let

$$\mathbf{V}^{w_0}_k(\mathbf{S}^{inc})^{\wp_{w_0}}/\mathcal{F}^{w_0} =_{w_0} \left\{ \mathbf{A}^{w_0}_{\mathcal{F}^{w_0}} | \ \mathbf{A} \text{ is a } w_0\text{-sequence from } \mathbf{V}^{w_0}_k(\mathbf{S}^{inc}) \right\}_{w_0}. \quad (6.5.32)$$

Then the bounded $w_0$-ultrapower is the $w_0$-union of the $w_0$-chain

$${}^{\#_{w_0}}\mathbf{S}^{inc} =_{w_0} \mathbf{V}^{w_0}_1(\mathbf{S}^{inc})^{\wp_{w_0}}/\mathcal{F}^{w_0} \subseteq_{w_0} \ldots \subseteq_{w_0} \mathbf{V}^{w_0}_k(\mathbf{S}^{inc})^{\wp_{w_0}}/\mathcal{F}^{w_0} \subseteq_{w_0} \ldots, \quad (6.5.33)$$

and we can define $j_{w_0}$ by induction. For $k = 1$, $j_{w_0}$ must be the $w_0$-identity. If $\mathbf{A}^{w_0}_{\mathcal{F}^{w_0}} \in_{w_0} \mathbf{V}^{w_0}_{k+1}(\mathbf{S}^{inc})^{\wp_{w_0}}/\mathcal{F}^{w_0}$ and $\mathbf{A}^{w_0}_{\mathcal{F}^{w_0}} \notin^{w}_{\mathcal{F}^{w_0}} {}^{\#_{w_0}}\mathbf{S}^{inc}$ we simply set $j_w(\mathbf{A}^{w_0}_{\mathcal{F}^{w_0}}) =_{w_0} \left\{ j_w(\mathbf{B}^{w_0}_{\mathcal{F}^{w_0}}) | \mathbf{B}^{w_0}_{\mathcal{F}^{w_0}} \in_{\mathcal{F}^{w}} \mathbf{A}^{w_0}_{\mathcal{F}^{w_0}} \right\}_{w}$ This makes sense: if $\mathbf{B}^{w_0}_{\mathcal{F}^{w_0}} \in_{\mathcal{F}^{w_0}} \mathbf{A}^{w_0}_{\mathcal{F}^{w_0}}$ it follows from (6.5.16) that $\{v \in_{w_0} \wp_{w_0} | B_v \in_{w_0} \mathbf{V}^{w_0}_k(\mathbf{S}^{inc})\}_{w_0} \in_{w_0} \mathcal{F}^{w_0}$, i.e., $\mathbf{B}^{w_0}_{\mathcal{F}^{w_0}} \in_{w_0} \mathbf{V}^{w_0}_k(\mathbf{S}^{inc})^{\wp_{w_0}}/\mathcal{F}^{w_0}$, which means that $j_{w_0}(B^{w_0}_{\mathcal{F}^{w_0}})$ is defined at a previous stage of the inductive construction.

Combining $i_{w_0}$ and $j_{w_0}$ we get a model of the extended $w_0$-consistent nonstandard universe

$$\begin{array}{ccc} \mathbf{V}^{w_0}(\mathbf{S}^{inc})^{\wp_{w_0}}/\mathcal{F}^{w_0} & & \mathbf{V}^{w_0}({}^{\#_{w_0}}\mathbf{S}^{inc}) \\ & & {}_{w_0}\supseteq \\ \nwarrow & & \swarrow \nearrow \\ i_{w_0} & & \#_{w_0} \\ & \mathbf{V}^{w_0}(\mathbf{S}^{inc}) & \end{array} \quad (6.5.34)$$

where ${}^{\#_{w_0}}\mathbf{A} =_{w_0} j_{w_0}(i_{w_0}(\mathbf{A}))$, for any $\mathbf{A} \in_{w_0} \mathbf{V}^{w_0}(\mathbf{S}^{inc})$. Here $\mathbf{V}^{w_0}(\mathbf{S}^{inc})$ and $\mathbf{V}^{w_0}({}^{\#_{w_0}}\mathbf{S}^{inc})$ are connected by $w_0$-consistent transfer principle.

# VI.5.3.4.The embedding $\mathbf{V}^{w_n}(\mathbf{S}^{inc})^{\wp_{w_n}}/\mathcal{F}^{w_n}$ into $\mathbf{V}^{w_n}({}^{\#_{w_n}}\mathbf{S}^{inc})$.

Let ${}^{\#_{w_n}}\mathbf{S}^{inc}$ be the ($w_n$-bounded) $w_n$-ultrapower $\mathbf{V}^{w_n}(\mathbf{S}^{inc})^{\wp_{w_n}}/\mathcal{F}^{w_n}$.

**Remark.6.5.10.**Note that in contrast with a classical case $\mathbf{V}^{w_n}(\mathbf{S}^{inc})^{\wp_{w_n}}/\mathcal{F}^{w_n}$ not alwais will not be the same as the full $w_n$-superstructure $\mathbf{V}^{w_n}({}^{\#_{w_n}}\mathbf{S}^{inc})$.

We shall now construct an $w_n$-embedding

$$j_{w_n} : \mathbf{V}^{w_n}(\mathbf{S}^{inc})^{\wp_{w_n}}/\mathcal{F}^{w_n} \hookrightarrow \mathbf{V}^{w_n}({}^{\#_{w_n}}\mathbf{S}^{inc}) \quad (6.5.35)$$

such that: (i) $j_{w_n}$ is the $w_n$-identity on ${}^{\#_{w_n}}\mathbf{S}^{inc}$ and (ii) if $\mathbf{A}^{w_n}_{\mathcal{F}^{w_n}} \notin^{w}_{\mathcal{F}^{w_n}} {}^{\#_{w_n}}\mathbf{S}^{inc}$, then

$$j_{w_n}(\mathbf{A}^{w_n}_{\mathcal{F}^{w_n}}) =_{w_n} \left\{ j_{w_n}(\mathbf{B}^{w_n}_{\mathcal{F}^{w_n}}) | \mathbf{B}^{w_n}_{\mathcal{F}^{w_n}} \in_{\mathcal{F}^{w_n}} \mathbf{A}^{w_n}_{\mathcal{F}^{w_n}} \right\}_{w_n}. \quad (6.5.36)$$

This means that the relation $\in_{\mathcal{F}^{w_n}}$ in the $w_n$-ultrapower is mapped into the ordinary $w_n$-membership relation in $\mathbf{V}^{w_n}({}^{\#_{w_n}}\mathbf{S}^{inc})$. The $w_n$-embedding $j_{w_n}$ is constructed in stages. Let

$$\mathbf{V}^{w_n}_k(\mathbf{S}^{inc})^{\wp_{w_n}}/\mathcal{F}^{w_n} =_{w_n} \left\{ \mathbf{A}^{w_n}_{\mathcal{F}^{w_n}} | \ \mathbf{A} \text{ is a } w_n\text{-sequence from } \mathbf{V}^{w_n}_k(\mathbf{S}^{inc}) \right\}_{w_n}. \quad (6.5.37)$$

Then the bounded $w_n$-ultrapower is the $w_n$-union of the $w_n$-chain

$$^{\#_{w_n}}\mathbf{S}^{inc} =_{w_n} \mathbf{V}_1^{w_n}(\mathbf{S}^{inc})^{\wp w_n}/\mathcal{F}^{w_n} \subseteq_w \ldots \subseteq_w \mathbf{V}_k^{w_n}(\mathbf{S}^{inc})^{\wp w_n}/\mathcal{F}^{w_n} \subseteq_{w_n} \ldots, \tag{6.5.38}$$

and we can define $j_{w_n}$ by induction. For $k = 1$, $j_{w_n}$ must be the $w_n$-identity. If $\mathbf{A}^{w_n}_{\mathcal{F}^{w_n}} \in_{w_0} \mathbf{V}^{w_0}_{k+1}(\mathbf{S}^{inc})^{\wp w_0}/\mathcal{F}^{w_0}$ and $\mathbf{A}^{w_n}_{\mathcal{F}^{w_n}} \notin^{w}_{\mathcal{F}^{w_0}} {}^{\#_{w_0}}\mathbf{S}^{inc}$ we simply set $j_w(\mathbf{A}^{w_0}_{\mathcal{F}^{w_0}}) =_{w_0} \left\{ j_w(\mathbf{B}^{w_n}_{\mathcal{F}^{w_n}}) | \mathbf{B}^{w_n}_{\mathcal{F}^{w_n}} \in_{\mathcal{F}^w} \mathbf{A}^{w_n}_{\mathcal{F}^{w_n}} \right\}_w$ This makes sense: if $\mathbf{B}^{w_n}_{\mathcal{F}^{w_n}} \in_{\mathcal{F}^w} \mathbf{A}^{w_n}_{\mathcal{F}^{w_n}}$ it follows from (6.5.19) that $\{v \in_{w_n} \wp_{w_n} | B_v \in_{w_0} \mathbf{V}_k^{w_n}(\mathbf{S}^{inc})\}_{w_n} \in_{w_n} \mathcal{F}^{w_n}$, i.e., $\mathbf{B}^{w_n}_{\mathcal{F}^{w_n}} \in_w \mathbf{V}_k^{w_n}(\mathbf{S}^{inc})^{\wp w_n}/\mathcal{F}^{w_n}$, which means that $j_{w_n}(B^{w_n}_{\mathcal{F}^{w_n}})$ is defined at a previous stage of the inductive construction.

Combining $i_{w_n}$ and $j_{w_n}$ we get a model of the extended $w_0$-consistent nonstandard universe

$$\begin{array}{ccc} \mathbf{V}^{w_n}(\mathbf{S}^{inc})^{\wp w_n}/\mathcal{F}^{w_n} & & \mathbf{V}^{w_n}(^{\#_{w_n}}\mathbf{S}^{inc}) \\ & & {}_{w_0}\supseteq \\ \nwarrow & & \swarrow\nearrow \\ i_{w_n} & & \#_{w_n} \\ & \mathbf{V}^{w_n}(\mathbf{S}^{inc}) & \end{array} \tag{6.5. }$$

where $^{\#_{w_n}}\mathbf{A} =_{w_n} j_{w_n}(i_{w_n}(\mathbf{A}))$, for any $\mathbf{A} \in_{w_n} \mathbf{V}^{w_n}(\mathbf{S}^{inc})$. Here $\mathbf{V}^{w_n}(\mathbf{S}^{inc})$ and $\mathbf{V}^{w_n}(^{\#_{w_n}}\mathbf{S}^{inc})$ are connected by $w_n$-inconsistent transfer principle.

# VI.6.The Paralogical Transfer Principle

## VI.6.1.The restricted inconsistent language

The structure $\mathbf{S}^{inc} \supset \mathbb{R}_w$ has an associated elementary language $\mathcal{L}(\mathbf{S}^{inc} \supset \mathbb{R}_w)$, which we used to give the necessary precision to the transfer principle. We need a similar formal tool to state the extended transfer principle.The language $\mathcal{L}(\mathbf{V}^w(\mathbf{S}^{inc}))$ will be an extension of the language $\mathcal{L}(\mathbf{S}^{inc}) = \mathcal{L}(\mathbf{S}^{inc} \supset \mathbb{R}_w)$. We add to our stock of elementary formulas [see (i)-(vi) in Section 1.11 expressions of the form

$$X =_s Y, X =_w Y, X =_{w_0} Y, X =_{w_n} Y, \ldots \tag{6.6.1}$$

and

$$X \in_s Y, X \in_w Y, X \in_{w_0} Y, X \in_{w_n} Y, \ldots \tag{6.6.2}$$

We keep the logical symbols of $\mathcal{L}(\mathbf{S}^{inc})$, but in addition to the number quantifiers we add bounded set quantifiers

$$\begin{gathered} \forall X(X \in_s Y), \forall X(X \in_w Y), \forall X(X \in_{w_0} Y), \forall X(X \in_{w_n} Y), \ldots \\ \exists X(X \in_s Y), \exists X(X \in_w Y), \exists X(X \in_{w_0} Y), \exists X(X \in_{w_n} Y), \ldots \end{gathered} \tag{6.6.3}$$

Formulas $\Phi$ of $\mathcal{L}(\mathbf{V}^w(\mathbf{S}^{inc}))$ are then constructed in exactly the same way as formulas of $\mathcal{L}(\mathbf{S}^{inc})$. A formula $\Phi$ of of $\mathcal{L}(\mathbf{V}^w(\mathbf{S}^{inc}))$ can be interpreted in a natural way in any of the structures $\mathbf{V}^w(\mathbf{S}^{inc})$, $\mathbf{V}^w(\mathbf{S}^{inc})^{\wp w}/\mathcal{F}^w$, and $\mathbf{V}^w(^{\#_w}\mathbf{S}^{inc})$; note that in $\mathbf{V}^w(\mathbf{S}^{inc})$ and $\mathbf{V}^w(^{\#_w}\mathbf{S}^{inc})$ we have the standard interpretation of the $\in_w$ symbol, in $\mathbf{V}^w(\mathbf{S}^{inc})^{\wp w}/\mathcal{F}^w$ we use $\in_{\mathcal{F}^w}$ as introduced in (6.5.8) to interpret $w$-membership. Given any formula $\Phi(X_1, \ldots, X_n)$ with $X_1, \ldots, X_n$ as the only free set parameters, and given sets $A_1, \ldots, A_n \in_w \mathbf{V}^w(\mathbf{S}^{inc})$, we mean by $\Phi(A_1, \ldots, A_n)$ the statement about $\mathbf{V}^w(\mathbf{S}^{inc})$ obtained by giving the variables $X_1, \ldots, X_n$ the values $A_1, \ldots, A_n$, respectively. In a similar way we interpret $\Phi(^{\#_w}A_1, \ldots, {}^{\#_w}A_n)$ as a condition about $\mathbf{V}^w(^{\#_w}\mathbf{S}^{inc})$ obtained by giving each $X_k$ the value $^{\#_w}A_k =_w j_w(i_w(A_k))$.

**Theorem 6.6.1.(TRANSFER PRINCIPLE)** (i) Let $A_1, \ldots, A_n \in_w \mathbf{V}^w(\mathbf{S}^{inc})$. Any

$\mathcal{L}(\mathbf{V}^{w}(\mathbf{S}^{inc}))$

statement $\Phi$ that is true of $A_1,\ldots,A_n$ in $\mathbf{V}^{w}(\mathbf{S}^{inc})$ is true of ${}^{\#_w}A_1,\ldots,{}^{\#_w}A_n$ in $\mathbf{V}^{w}({}^{\#_w}\mathbf{S}^{inc})$ .

(ii) Let $A_1,\ldots,A_n \in_{w_0} \mathbf{V}^{w_0}(\mathbf{S}^{inc})$. Any $\mathcal{L}(\mathbf{V}^{w_0}(\mathbf{S}^{inc}))$ statement $\Phi$ that is true of $A_1,\ldots,A_n$ in $\mathbf{V}^{w_0}(\mathbf{S}^{inc})$ is true of ${}^{\#_{w_0}}A_1,\ldots,{}^{\#_{w_0}}A_n$ in $\mathbf{V}^{w_0}({}^{\#_{w_0}}\mathbf{S}^{inc})$.

(iii) Let $A_1,\ldots,A_n \in_{w_n} \mathbf{V}^{w_n}(\mathbf{S}^{inc})$. Any $\mathcal{L}(\mathbf{V}^{w_n}(\mathbf{S}^{inc}))$ statement $\Phi$ that is true of $A_1,\ldots,A_n$ in $\mathbf{V}^{w_n}(\mathbf{S}^{inc})$ is true of ${}^{\#_{w_n}}A_1,\ldots,{}^{\#_{w_n}}A_n$ in $\mathbf{V}^{w_n}({}^{\#_{w_n}}\mathbf{S}^{inc})$.

**Remark 6.6.1.** Let $A \subseteq {}^{*}\mathbb{R}$ then $A \in V_2({}^{*}\mathbb{R})$. The canonical embedding

$* : V(\mathbb{R}) \overset{*}{\hookrightarrow} V({}^{*}\mathbb{R})$

maps $V_2(\mathbb{R})$ to a set ${}^{*}V_2(\mathbb{R}) \subsetneq V({}^{*}\mathbb{R})$. Will $A$ belong to this set: ${}^{*}V_2(\mathbb{R})$ ? It well known that is not necessarily except if $A \in {}^{*}B$ for some $B \in V_2(\mathbb{R})$, i.e. then $A \in {}^{*}V_2(\mathbb{R})$.

We thus want to prove that

$$\forall A\left[\left(A \in {}^{*}B\right) \wedge \left(A \subseteq {}^{*}\mathbb{R}\right) \Rightarrow A \in {}^{*}V_2(\mathbb{R})\right]. \qquad (6.6.4)$$

As it stands, (6.6.4) is not an $\mathcal{L}(\mathbf{V}(\mathbb{R}))$ formula. However, it is equivalent to

$$\forall A\left(A \in {}^{*}B\right)\left[\forall r \in A\left(r \in {}^{*}\mathbb{R}\right) \Rightarrow A \in {}^{*}V_2(\mathbb{R})\right]. \qquad (6.6.5)$$

This is genuine $\mathcal{L}(\mathbf{V}^{w}(\mathbf{S}^{inc}))$; i.e., we have only bounded set quantifiers. Now (6.6.5) is a condition $\Phi({}^{*}B, {}^{*}\mathbb{R}, {}^{*}V_2(\mathbb{R}))$, which by transfer is true in $V({}^{*}\mathbb{R})$ iff the corresponding $\Phi(B,\mathbb{R},V_2(\mathbb{R}))$ is true in $V(\mathbb{R})$. But the latter condition is trivially true. Thus we have shown that if a subset of ${}^{*}\mathbb{R}$ is an element of some ${}^{*}B$ in $V({}^{*}\mathbb{R})$, then it is already an element of the $*$-image of $V_2(\mathbb{R})$.

Remind the following definition.

**Definition 6.6.1.** Let $A \in V({}^{*}\mathbb{R})$, then

(i) $A$ is called $*$-standard if $A = {}^{*}B$ for some $B \in V(\mathbb{R})$,

(ii) $A$ is called $*$-internal if $A \in {}^{*}B$ for some $B \in V(\mathbb{R})$, and

(iii) $A$ is called $*$-external if $A$ is not $*$-internal.

**Remark6.6.2.** It well known that every $*$-standard set is $*$-internal and that every element of an $*$-internal set is $*$-internal.

**Definition 6.6.2.** Let $A \in_w V^w({}^{\#_w}\mathbb{R}_w)$, then:

(i) $A$ is called $w$-standard if $A =_w {}^{\#_w}B$ for some $B \in_w V^w(\mathbb{R}_w)$,

(ii) $A$ is called $w$-internal if $A \in_w {}^{\#_w}B$ for some $B \in_w V^w(\mathbb{R}_w)$, and

(iii) $A$ is called $w$-external if $A$ is not $w$-internal.

**Definition 6.6.3.** Let $A \in_w V^w({}^{\#_w}\mathbb{R}_w)$, then:

(i) $A$ is called weakly $w$-standard or $w_{[1]}$-standard if $A =_{w_{[1]}} {}^{\#_w}B$ for some $B \in_w V^w(\mathbb{R}_w)$,

(ii) $A$ is called weakly $w$-internal or $w_{[1]}$-internal if $A \in_{w_{[1]}} {}^{\#_w}B$ for some $B \in_w V^w(\mathbb{R}_w)$, and

(iii) $A$ is called $w_{[1]}$-external if $A$ is not $w_{[1]}$-internal.

**Definition 6.6.4.** Let $A \in_{w_{\{0\}}} V^{w_{\{0\}}}\left({}^{\#_{w_{\{0\}}}}\mathbb{R}_{w_{\{0\}}}\right)$, then:

(i) $A$ is called $w_{\{0\}}$-standard if $A =_{w_{\{0\}}} {}^{\#_{w_{\{0\}}}}B$ for some $B \in_{w_{\{0\}}} V^{w_{\{0\}}}(\mathbb{R}_{w_{\{0\}}})$,

(ii) $A$ is called $w_{\{0\}}$-internal if $A \in_{w_{\{0\}}} {}^{\#_{w_{\{0\}}}}B$ for some $B \in_{w_{\{0\}}} V^{w_{\{0\}}}(\mathbb{R}_{w_{\{0\}}})$, and

(iii) $A$ is called $w_{\{0\}}$-external if $A$ is not $w_{\{0\}}$-internal.

**Definition 6.6.5.** Let $A \in_{w_{\{0\}}} V^{w_{\{0\}}}\left({}^{\#_{w_{\{0\}}}}\mathbb{R}_{w_{\{0\}}}\right)$, then:

(i) $A$ is called weakly $w_{\{0\}}$-standard or $w^{w}_{\{0\}}$-standard if

$\left(A =_{w_{\{0\}}} {}^{\#_{w_{\{0\}}}}B\right) \wedge \left(A \neq^{w}_{w_{\{0\}}} {}^{\#_{w_{\{0\}}}}B\right)$ for some $B \in_{w_{\{0\}}} V^{w_{\{0\}}}(\mathbb{R}_{w_{\{0\}}})$,

(ii) $A$ is called weakly $w_{\{0\}}$-internal or $w^{w}_{\{0\}}$-internal if

$\left(A \in_{w_{\{0\}}} {}^{\#_{w_{\{0\}}}}B\right) \wedge \left(A \notin^{w}_{w_{\{0\}}} {}^{\#_{w_{\{0\}}}}B\right)$ for some $B \in_{w_{\{0\}}} V^{w_{\{0\}}}(\mathbb{R}_{w_{\{0\}}})$, and

(iii) $A$ is called $w^{w}_{\{0\}}$-external if $A$ is not $w^{w}_{\{0\}}$-internal.

**Definition 6.6.6.** Let $A \in_{w_{\{n\}}} V^{w_{\{n\}}}\left({}^{\#_{w_{\{n\}}}}\mathbb{R}_{w_{\{n\}}}\right), n = 1,2,\ldots,$ then:

(i) $A$ is called weakly $w_{\{n\}}$-standard or $w^{w}_{\{n\}}$-standard if

$\left(A =_{w_{\{n\}}} {}^{\#_{w_{\{0\}}}}B\right) \wedge \left(A \neq^{w}_{w_{\{n\}}} {}^{\#_{w_{\{0\}}}}B\right)$ for some $B \in_{w_{\{n\}}} V^{w_{\{n\}}}(\mathbb{R}_{w_{\{n\}}})$,

(ii) $A$ is called $w_{\{n\}}$-internal if $A \in_{w_{\{n\}}} {}^{\#_{w_{\{0\}}}}B$ for some $B \in_{w_{\{n\}}} V^{w_{\{n\}}}(\mathbb{R}_{w_{\{n\}}})$, and

(iii) $A$ is called $w_{\{n\}}$-external if $A$ is not $w_{\{n\}}$-internal.

**Definition 6.6.7.** Assume now that for any $A \in_{w} V^{w}({}^{\#_{w}}\mathbb{R}_{w})$:

(i) $A$ is $w_{[1]}$-standard or (ii) $A$ is $w_{[1]}$-internal,

then superstructure $V^{w}({}^{\#_{w}}\mathbb{R}_{w})$ is called purely $w_{[1]}$-internal and we abbreviate $V^{w}_{[1]\text{-}Int}({}^{\#_{w}}\mathbb{R}_{w})$.

**Definition 6.6.8.** Assume now that for any $A \in_{w_{\{0\}}} V^{w_{\{0\}}}\left({}^{\#_{w_{\{0\}}}}\mathbb{R}_{w}\right)$ :

(i) $A$ is $w^{w}_{\{0\}}$-standard or (ii) $A$ is $w^{w}_{\{0\}}$-internal,

then superstructure $V^{w_{\{0\}}}\left({}^{\#_{w_{\{0\}}}}\mathbb{R}_{w_{\{0\}}}\right)$ is called purely $w^{w}_{\{0\}}$-internal and we abbreviate $V^{w_{\{0\}}}_{w^{w}_{\{0\}}\text{-}Int}\left({}^{\#_{w_{\{0\}}}}\mathbb{R}_{w_{\{0\}}}\right)$.

**Definition 6.6.9.** Assume now that for any $A \in_{w_{\{n\}}} V^{w_{\{n\}}}\left({}^{\#_{w_{\{n\}}}}\mathbb{R}_{w_{\{n\}}}\right)$ :

(i) $A$ is $w^{w}_{\{n\}}$-standard or (ii) $A$ is $w^{w}_{\{n\}}$-internal,

then superstructure $V^{w_{\{n\}}}\left({}^{\#_{w_{\{n\}}}}\mathbb{R}_{w_{\{n\}}}\right)$ is called purely $w^{w}_{\{n\}}$-internal and we abbreviate $V^{w_{\{n\}}}_{w^{w}_{\{n\}}\text{-}Int}\left({}^{\#_{w_{\{n\}}}}\mathbb{R}_{w_{\{n\}}}\right)$.

**Remark.6.6.3.** We remind now the details of the description the $*$-internal sets in the consistent model. Let $A$ be $*$-internal; thus $A \in {}^{*}V_{k+1}(\mathbb{R})$ for some $k \geq 1$. This means that $A$ will be of the form $A = j(A_{\mathcal{F}})$, for some $A_{\mathcal{F}}$. By the By the construction of $j$, one then gets

$$A \in {}^{*}V_{k+1}(\mathbb{R}) \text{ iff } j(A_{\mathcal{F}}) \in j(i(V_{k+1}(\mathbb{R}))),$$
$$\text{iff } A_{\mathcal{F}} \in i(V_{k+1}(\mathbb{R})),$$

where $i$ is the embedding of $V(\mathbb{R})$ into the ultrapower. The definition of $\in_{\mathcal{F}}$ then gives

$$A \in {}^{*}V_{k+1}(\mathbb{R}) \text{ iff } \{v \in \wp | A_{v} \in V_{k+1}(\mathbb{R})\} \in \mathcal{F},$$

where $\langle A_{v}\rangle_{v\in\wp}$ is the bounded sequence defining $A_{\mathcal{F}}$. Thus the $*$-internal sets are precisely the objects we obtain by starting with an arbitrary bounded sequence $\langle A\rangle_{v\in\wp}$ and the standard objects are obtained by starting from a constant sequence $\langle A\rangle_{v\in\wp}$.

**Remark.6.6.4.** Because of their importance we will describe in detail the $w$-internal sets in the models. Let $A$ be $w$-internal; thus $A \in_{w} {}^{\#_{w}}V^{w}_{k+1}\left(\mathbf{S}^{inc}_{w}{}_{w} \supset \mathbb{R}_{w}\right)$ for some $k \geq 1$. This means that $A$ will be of the form $A =_{w} j_{w}(A^{w}_{\mathcal{F}^{w}})$, for some $A^{w}_{\mathcal{F}^{w}}$. By the construction of $j_{w}$, we then get

$$A \in_w {}^{\#_w}V^w_{k+1}\left(\mathbf{S}^{inc}_{w}{}_{\ w} \supset \mathbb{R}_w\right) \text{ iff } j_w(A^w_{\mathcal{F}^w}) \in_w j_w\left(i_w\left(V^w_{k+1}\left(\mathbf{S}^{inc}_{w}{}_{\ w} \supset \mathbb{R}_w\right)\right)\right),$$
$$\text{iff } A^w_{\mathcal{F}^w} \in_{\mathcal{F}^w} i_w\left(V^w_{k+1}\left(\mathbf{S}^{inc}_{w}{}_{\ w} \supset \mathbb{R}_w\right)\right),$$

where $i_w$ is the $w$-embedding of $V^w(\mathbb{R}_w)$ into the $w$-ultrapower. The definition of $\in_{\mathcal{F}^w}$ then

gives

$$A \in_w {}^{\#_w}V^w_{k+1}\left(\mathbf{S}^{inc}_{w}{}_{\ w} \supset \mathbb{R}_w\right) \text{ iff } \{v \in_w \wp_w|\} \in_w \mathcal{F}^w,$$

where $\langle A_v \rangle_{v \in_w \wp}$ is the bounded $w$-sequence defining $A_{\mathcal{F}}$. Thus the $w$-internal sets are precisely the objects we obtain by starting with an arbitrary bounded sequence $\langle A \rangle_{v \in_w \wp}$ and the standard objects are obtained by starting from a constant sequence $\langle A \rangle_{v \in_w \wp}$.

**Remark.6.6.5.**

Remind the following Theorem.

**Theorem 6.6.2.**(i) Every nonempty $*$-internal subset of ${}^*\mathbb{N}$ has a least element.

(ii) Every nonempty $*$-internal subset of ${}^*\mathbb{R}$ with an upper bound has a $<$-least upper bound.

**Proof.**We prove (i), so let $A \subseteq {}^*\mathbb{N}$ be internal. Then $A \in {}^*V(\mathbb{R});$ see (6.6.5). We can express the fact that an internal subset of ${}^*\mathbb{N}$ has a least element by the condition

$$\Phi \triangleq$$
$$\forall X\left(X \in {}^*V_2(\mathbb{R})\right)\left\{\left[(X \neq \varnothing) \wedge \left(X \subseteq {}^*\mathbb{N}\right)\right] \Rightarrow X \text{ has a } < \text{-least element}\right\}, \quad (6.6.6)$$

where the condition: $X$ has a $<$-least element, in detail is

$$\exists x(x \in X)\forall y[y \in X \Rightarrow \neg(y < x)]. \quad (6.6.7)$$

Finally we have a condition

$$\Phi \triangleq \forall X\left(X \in {}^*V_2(\mathbb{R})\right)$$
$$\left\{\left[(X \neq \varnothing) \wedge \left(X \subseteq {}^*\mathbb{N}\right)\right] \Rightarrow \exists x(x \in X)\forall y[y \in X \Rightarrow \neg(y < x)]\right\}. \quad (6.6.8)$$

We thus have a condition $\Phi(\mathbb{N}, V_2(\mathbb{R}))$ such that $\Phi(\mathbb{N}, V_2(\mathbb{R}))$ is true in $V(\mathbb{R})$. By $*$-transfer condition $\Phi({}^*\mathbb{N}, {}^*V_2(\mathbb{R}))$ is true in $V({}^*\mathbb{R})$ proving (i).

**Remark.6.6.6.**It follows from Theorem 6.6.2(i) that:

(1) ${}^*\mathbb{N}\backslash\mathbb{N}$ is $*$-external since there is no $<$-least element in ${}^*\mathbb{N}\backslash\mathbb{N}$ : if $x \in {}^*\mathbb{N}\backslash\mathbb{N}$ then also $x - 1 \in {}^*\mathbb{N}\backslash\mathbb{N}$.

(2) We also see that $\mathbb{N}$ is external; thus $\mathbb{N} \in V_2({}^*\mathbb{R})\backslash{}^*V_2(\mathbb{R})$

(3) From Theorem 6.6.2(ii) it follows that $\mathbb{R}$ as a subset of ${}^*\mathbb{R}$ is $*$-external.

(4) Note that Theorem 6.6.2 is valid only for $*$-internal sets; the positive infinitesimals in

${}^*\mathbb{R}$ is bounded but has no least upper bound.

**Theorem 6.6.3.**(i) If $A$ is $*$-internal and $\mathbb{N} \subseteq A,$ then $A$ contains some infinite natural number, i.e., an element of ${}^*\mathbb{N}\backslash\mathbb{N}$.

(ii) If $A$ is internal and every infinite $n \in {}^*\mathbb{N}$ belongs to $A,$ then $A$ contains some standard

$n \in \mathbb{N}$.

(iii) If an internal set $A$ contains every positive infinitesimal, then $A$ contains some

positive

standard real $r \in \mathbb{R}$.

(iv) If an internal set $A$ contains every standard positive real, then $A$ contains some positive infinitesimal.

**Remark 6.6.7.** Let $A \subseteq_w {}^{\#_w}\mathbb{R}_w$ then $A \in_w V_2^w({}^{\#_w}\mathbb{R}_w)$. The canonical embedding $\#_w : V^w(\mathbb{R}_w) \overset{\#_w}{\hookrightarrow} V({}^{\#_w}\mathbb{R}_w)$ maps $V_2^w(\mathbb{R}_w)$ to a set ${}^{\#_w}V_2^w(\mathbb{R}_w) \subsetneq_w V^w({}^{\#_w}\mathbb{R}_w)$. Will $A$ $w$-belong to this set: ${}^{\#_w}V_2^w(\mathbb{R}_w)$ ? That is not necessarily except if $A \in {}^{\#_w}B$ for some $B \in_w V_2^w(\mathbb{R}_w)$, i.e. then $A \in {}^{\#_w}V_2(\mathbb{R})$. We thus want to prove that

$$\forall A\Big[\Big(A \in_w {}^{\#_w}B\Big) \wedge \Big(A \subseteq_w {}^{\#_w}\mathbb{R}_w\Big) \Rightarrow_s A \in_w {}^{\#_w}V_2^w(\mathbb{R}_w)\Big]. \tag{6.6.9}$$

As it stands, (6.6.4) is not an $\mathcal{L}(\mathbf{V}^w(\mathbb{R}_w))$ formula. However, it is equivalent to

$$\forall A\Big(A \in_w {}^{\#_w}B\Big)\Big[\forall r \in_w A\Big(r \in_w {}^{\#_w}\mathbb{R}_w\Big) \Rightarrow A \in_w {}^{\#_w}V_2(\mathbb{R}_w)\Big]. \tag{6.6.10}$$

This is genuine $\mathcal{L}(\mathbf{V}^w(\mathbb{R}_w))$; i.e., we have only bounded set quantifiers. Now (6.6.5) is a condition $\Phi({}^{\#_w}B, {}^{\#_w}\mathbb{R}, {}^{\#_w}V_2^w(\mathbb{R}_w))$, which by transfer is true in $V({}^{\#_w}\mathbb{R}_w)$ iff the corresponding $\Phi(B, \mathbb{R}_w, V_2^w(\mathbb{R}_w))$ is true in $V^w(\mathbb{R}_w)$. But the latter condition is trivially true. Thus we have shown that if a subset of ${}^{\#_w}\mathbb{R}_w$ is an element of some ${}^{\#_w}B$ in $V^w({}^{\#_w}\mathbb{R}_w)$, then it is already an element of the $w$-image of $V_2^w(\mathbb{R}_w)$.

The following Theorem very similar to Theorem 6.6.2.

**Theorem 6.6.4.** (i) Every nonempty $w$-internal subset of ${}^{\#_w}\mathbb{N}_w$ has a $<_w$ -least element.

(ii) Every nonempty $w$-internal subset of ${}^{\#_w}\mathbb{R}_w$ with an $<_w$ -upper bound has a $<_w$ -least upper bound.

**Proof.** We prove (i), so let $A \subseteq_w {}^{\#_w}\mathbb{N}_w$ be $w$-internal. Then $A \in {}^{\#_w}V^w(\mathbb{R}_w)$; see (6.6.5). We

can express the fact that an internal subset of ${}^{\#_w}\mathbb{N}_w$ has a least element by the condition

$$\Phi \triangleq \forall X\Big(X \in_w {}^{\#_w}V_2(\mathbb{R}_w)\Big) \searrow \Big\{\Big[(X \neq_w \varnothing_w^w) \wedge \Big(X \subseteq_w {}^{\#_w}\mathbb{N}_w\Big)\Big] \Rightarrow_s X \text{ has a } <_w \text{-least } w\text{-element}\Big\}, \tag{6.6.11}$$

where the condition: $X$ has a $<_w$ -least $w$-element, in detail is

$$\exists x(x \in_w X)\forall y[y \in_w X \Rightarrow_s \neg_w(y <_w x)]. \tag{6.6.12}$$

Finally we have a condition

$$\Phi \triangleq \forall X\Big(X \in_w {}^{\#_w}V_2^w(\mathbb{R}_w)\Big) \searrow \Big\{\Big[(X \neq_w \varnothing) \wedge \Big(X \subseteq_w {}^{\#_w}\mathbb{N}_w\Big)\Big] \Rightarrow_s \exists x(x \in_w X)\forall y[y \in_w X \Rightarrow_s \neg_w(y <_w x)]\Big\}. \tag{6.6.13}$$

We thus have a condition $\Phi(\mathbb{N}_w, V_2(\mathbb{R}_w))$ such that $\Phi(\mathbb{N}_w, V_2^w(\mathbb{R}_w))$ is true in $V^w(\mathbb{R}_w)$. By $\#_w$-transfer condition $\Phi({}^{\#_w}\mathbb{N}_w, {}^{\#_w}V_2^w(\mathbb{R}_w))$ is true in $V^w({}^{\#_w}\mathbb{R}_w)$ proving (i).

# VII. Set theory $\mathbf{HST}^{\#}_{\omega}$.

## VII.1.Axiomatical system $\mathbf{HST}^{\#}_{\omega}$, as inconsistent generalization of Hrbacek set theory **HST**.

In this chapter we introduces $\mathbf{HST}^{\#}_{\omega}$, inconsistent generalization of Hrbacek set theory **HST** and describes the basic structure of the $\mathbf{HST}^{\#}_{\omega}$ set universe. Syntactically, $\mathbf{HST}^{\#}_{\omega}$ is a theory in the $\mathbf{st_s}$-$\in_{\mathbf{s}}$ -$\mathbf{st}_w$-$\in_w$ -language, which contains: (1) a binary consistent predicate of strong or consistent membership $\in_{\mathbf{s}}$ and consistent unary predicate of strong or consistent standardness $\mathbf{st_s}$ (and strong or consistent equality $=_{\mathbf{s}}$ of course) as the consistent primary notions and (2) a binary inconsistent predicate of weak or inconsistent membership $\in_w$ and inconsistent unary predicate of weak or inconsistent standardness $\mathbf{st}_w$ (and weak or inconsistent equality $=_w$ of course) as the inconsistent primary notions. Formula $x \in_w y$ reads: $x$ weakly belongs to $y$, or $x$ is an weak element of $y$, with the usual set theoretic understanding of inconsistent membership. The formula $\mathbf{st}_w x$ reads: $x$ is a weakly standard, its meaning will be explained below. A $\mathbf{st_s}$-$\in_{\mathbf{s}}$ -$\mathbf{st}_w$-$\in_w$ -formula is a formula of the $\mathbf{st_s}$-$\in_{\mathbf{s}}$ -$\mathbf{st}_w$-$\in_w$ -language. An $\in_w$ -formula is a formula of the $\in_w$ -language having $\in_w$ as the only atomic predicate. Thus an $\in_w$ -formula is a $\mathbf{st}_w$-$\in_w$ -formula in which the standardness predicate does not occur. $\in_w$ -formulas are also called weak internal formulas, in opposition to weak external formulas, i.e., those $\mathbf{st}_w$-$\in_w$ -formulas containing $\mathbf{st}_w$.

## VII.2. The universe of $\mathbf{HST}^{\#}_{\omega}$

Inconsistent set theory $\mathbf{HST}^{\#}_{\omega}$ deals with eight major types of sets: (i) strongly external or $\mathbf{s}$-external,(ii) strongly internal or $\mathbf{s}$-internal, (iii) strongly standard or $\mathbf{s}$-standard, (iv) strongly well-founded or $\mathbf{s}$-well-founded,(v) weakly external or $w$-external,(vi) weakly internal or $w$-internal, (vii) weakly well-founded or $w$-well-founded.

First of all, strongly standard sets are those consistent sets $x$ which satisfy $\mathbf{st_s}x$ and weakly standard sets are those inconsistent sets $x$ which satisfy $\mathbf{st}_w x$.Strongly internal sets are those consistent sets $y$ which satisfy $\mathbf{int_s}y$, where $\mathbf{int_s}y$ is the formula $\exists^{\mathbf{st_s}}x(y \in_{\mathbf{s}} x) \equiv \exists x[\mathbf{st_s}x \wedge (y \in_{\mathbf{s}} x)]$ (saying: $y$ strongly belongs to a strongly standard set), weakly internal sets are those inconsistent sets $y$ which satisfy $\mathbf{int}_w y$, where $\mathbf{int}_w y$ is the formula $\exists \mathbf{st}_w x(y \in_w x)$ (saying: $y$ weakly belongs to a weakly standard set). Thus,

(**i**) $\mathbf{S_s} = \{x : \mathbf{st_s}x\}_{\mathbf{s}}$ is the class of all consistent standard sets,

(**ii**) $\mathbf{I_s} = \{y : \mathbf{int_s}y\}_{\mathbf{s}} = \{y : \exists^{\mathbf{st_s}}x(y \in_{\mathbf{s}} x)\}_{\mathbf{s}}$ is the class of all consistent internal sets,

(**iii**) $\mathbf{S}_w = \{x : \mathbf{st}_w x\}_w$ is the class of all inconsistent standard sets,

(**iv**) $\mathbf{I}_w = \{y : \mathbf{int}_w y\}_w = \{y : \exists^{\mathbf{st}_w}x(y \in_w x)\}_w$ is the class of all inconsistent internal sets,

(**v**) $\mathbf{S}^{\#} = \mathbf{S}_s \cup_s \mathbf{S}_w = \{x : \mathbf{st}_s x\}_s \cup_s \{x : \mathbf{st}_w x\}_w$ is the class of all consistent and inconsistent standard sets,

(**vi**) $\mathbf{I}^{\#} = \mathbf{I}_s \cup_s \mathbf{I}_w = \{y : \mathbf{int}_s y\}_s \cup_s \{y : \mathbf{int}_w y\}_w$ is the class of all consistent and inconsistent internal sets.

The class $\mathbf{I}_s$ is the source of some typical objects of consistent "nonstandard" mathematics like consistent hyperintegers and consistent hyperreals, the class $\mathbf{I}_w$ is the source of some typical objects of inconsistent "nonstandard" mathematics like inconsistent hyperintegers and inconsistent hyperreals [],

**Blanket agreement 1.1**. Thus, internal sets are precisely all sets which are elements of consistent or inconsistent standard sets. This understanding of the notion of internality and the associated notions like $\mathbf{I}^{\#}, \exists^{\mathbf{st}_s}, \exists^{\mathbf{st}_w}, \exists^{\mathbf{st}^{\#}} \equiv \exists^{\mathbf{st}_s} \vee \exists^{\mathbf{st}_w}, \forall^{\mathbf{st}_s}, \forall^{\mathbf{st}^{\#}} \equiv \forall^{\mathbf{st}_s} \wedge \forall^{\mathbf{st}_w}$is default throughout this paper. All exceptions (e.g., when $\mathbf{IST}^{\#}_{\omega}$ is considered) will be explicitly indicated.

External sets consistent and inconsistent, are simply all sets in the nonstandard universe of $\mathbf{HST}^{\#}_{\omega}$. We shall use $\mathbf{H}^{\#}_{\omega}$ to denote the class of all consistent and inconsistent external sets. Thus, $\mathbf{H}^{\#}_{\omega}$ is the "universe of discourse", the universe of all sets considered by the theory, including the class $\mathbf{WF}^{\#}_{\omega}$ of all well-founded sets. $\mathbf{WF}^{\#}_{\omega}$ will satisfy all axioms of $ZFC^{\#}_{\omega}$. The class $\mathbf{S}^{\#}$ of all standard sets {determined by the predicate st, as above) will be shown to be $\in_s$ -$\in_w$ -isomorphic to $\mathbf{WF}^{\#}_{\omega}$. In a sense, $\mathbf{S}^{\#}$ is an "isomorphic expansion" of $\mathbf{WF}^{\#}_{\omega}$ into $\mathbf{H}^{\#}_{\omega}$. Given that $\mathbf{S}^{\#}$ is not transitive, $\mathbf{I}^{\#}$ arises naturally as the class of all elements of sets in $\mathbf{S}^{\#}$. It is viewed as an elementary extension of $\mathbf{S}^{\#}$ {in $\in_s$ -$\in_w$ -language), and thereby also of $\mathbf{WF}^{\#}_{\omega}$. Finally, $\mathbf{H}^{\#}_{\omega}$ is a comprehensive universe in which all these classes coexist in a reasonable common set theoretic structure, with $\in_s$ -$\in_w$ having the natural meaning in all mentioned universes.

# VI.3. The axioms of the external inconsistent universe.

This group includes the $ZFC^{\#}_{\omega}$ Extensionality, Pair, Union, Infinity axioms and the schemata of Separation and Collection (therefore also Replacement, which is a consequence of Collection, as usual) for all $\mathbf{st}_s$-$\in_s$ $\mathbf{st}_w$-$\in_w$ -formulas or for all $\mathbf{st}_{\#}$-$\in_{\#}$ -formulas for short.

# VI.4. Axioms for standard and internal sets

**Notation 4.1**. (1).Let quantifiers $\exists^{\mathbf{st}_s}, \forall^{\mathbf{st}_s}, \exists^{\mathbf{st}_w}$ and $\forall^{\mathbf{st}_w}$ be shortcuts meaning: there

exists a

strongly standard..., for all strongly standard,there exists a weakly standard..., for all weakly standard, ..., formally:

(i) $\exists^{\mathbf{st}_s}x\Phi(x)$ means $\exists x[\mathbf{st}_s x \wedge \Phi(x)]$, (ii) $\forall^{\mathbf{st}_s}x\Phi(x)$ means $\forall x[\mathbf{st}_s x \Rightarrow \Phi(x)]$,

(iii) $\exists^{\mathbf{st}_w}x\Phi(x)$ means $\exists x[\mathbf{st}_w x \wedge \Phi(x)]$, (iv) $\forall^{\mathbf{st}_w}x\Phi(x)$ means $\forall x[\mathbf{st}_w x \Rightarrow \Phi(x)]$.

Quantifiers $\exists^{\mathbf{int}}$ and $\forall^{\mathbf{int}}$ (meaning there exists an internal ... , for all internal ...) are introduced similarly. If g, is an E-formula then g,st, the relativization of g> to S, is the formula obtained by restriction of all quantifiers in g> to the class S, so that all occurrences of 3 x ... are changed to 3stx ... while all occurrences of V x ... are changed to ystx .... In other words, g,st says that g> is true in S. Rela-tivization g,int, which displays the truth of an e-formula g> in the universe 0, is defined similarly: the quantifiers 3, V change to 3int, yint. The following axioms specify the behaviour of standard and internal sets.

**Notation 4.2**.For all $\mathbf{st}_s$-$\in_s$ $\mathbf{st}_w$-$\in_w$ -formulas or for all $\mathbf{st}_\#$-$\in_\#$ - formulas for short.

$ZFC^{\mathbf{st}_\#}_\omega$ : The collection of all formulas of the form g,st, where g> is an e- statement which is an axiom of $ZFC^\#_\omega$. In other words, it is postulated that the universe $\mathbf{S}^\#$ is a $ZFC^\#_\omega$ universe. (Note that the $ZFC^\#_\omega$ axioms are assumed to be formulated as certain closed $\in_\#$ - formulas in this definition.) This is enough to prove the following statement:

**Lemma 4.1**. (1) $\mathbf{S}_s \subseteq \mathbf{I}_s$, (2) $\mathbf{S}_w \subseteq \mathbf{I}_w$.

Proof.(1) See [18] Lemma 1.1.3.

(2) Let $x \in_w \mathbf{S}_w$. The formula $\exists y(x \in_w y)$ is a theorem of $ZFC^\#_\omega$, therefore $[\exists y(x \in_w y)]^{\mathbf{st}_w}$ that is the formula $\exists^{\mathbf{st}_w}y(x \in_w y)$,is true. In other words, $x$ is an element of a standard set, which means $x \in_w \mathbf{I}_w$.

1.**Strong or Consistent Transfer** (**s-Transfer**): $\Phi^{\mathbf{int}_s} \Leftrightarrow \Phi^{\mathbf{st}_s}$, where $\Phi$ is an arbitrary closed $\in_s$ -formula containing only consistent standard sets as parameters.

To be more exact, Consistent Transfer is the collection of all statements of the form $\forall^{\mathbf{st}_s}x_1\ldots\forall^{\mathbf{st}_s}x_n[\Phi^{\mathbf{int}_s}(x_1,\ldots,x_n) \Leftrightarrow \Phi^{\mathbf{st}_s}(x_1,\ldots,x_n)]$

2.**Strong Consistent Transitivity of** $\mathbf{I}_s$ : $\forall^{\mathbf{int}_s}x\forall y(y \in_s x \Rightarrow \mathbf{int}_s y)$.

3.**Consistent Regularity over** $\mathbf{I}_s$ : For any non empty consistent set $X$ there exists $x \in_s X$ such that $x \cap_s X \subseteq_s \mathbf{I}_s$. (The full Regularity of $ZFC$ requires $x \cap_s X = \emptyset_s$.)

4.**Consistent Standardization**: $\forall X\exists^{\mathbf{st}_w}y(X \cap_s \mathbf{S}_s = \cap_s \mathbf{S}_s)$. (Such consistent standard set

$Y$, unique by Consistent Transfer and Consistent Extensionality, is sometimes denoted by ${}^{\mathbf{S}_s}X$.)

5.**Weak Transfer** ($w$-**Transfer**): $\Phi^{\mathbf{int}_w} \Rightarrow \Phi^{\mathbf{st}_w}$, where $\Phi$ is an arbitrary closed $\in_s$ -$\in_w$ -formula containing only consistent and inconsistent standard sets as parameters.

To be more exact, Weak Transfer is the collection of all statements of the form $\forall^{\mathbf{st}_s}x_1\ldots\forall^{\mathbf{st}_s}x_n\forall^{\mathbf{st}_w}y_1\ldots\forall^{\mathbf{st}_w}y_m[\Phi^{\mathbf{int}_s}(x_1,\ldots,x_n;y_1,\ldots,y_m) \Rightarrow \Phi^{\mathbf{st}_w}(x_1,\ldots,x_n;y_1,\ldots,y_m)]$

6.**Weak Transitivity of** $\mathbf{I}_w$ : $\forall^{\mathbf{int}_s}x\forall y(y \in_w x \Rightarrow \mathbf{int}_w y)$.

7.**Weak Regularity over** $\mathbf{I}_w$ : For any non empty consistent set $X$ there exists $x \in_w X$ such that $x \cap_w X \subseteq_w \mathbf{I}_w$. (The full Regularity of $ZFC$ requires $x \cap_w X = \emptyset_w$.)

8.**Strictly Weak Regularity** (**Strictly** $w$-**Regularity**): For any non empty inconsistent set $X$ there exists $x \in_w X$ such that $x \cap_w X =_w \emptyset_w \wedge \neg(x \cap_w X =_w \emptyset_w)$.

9.**Weak Standardization** ($w$-**Standardization**): $\forall X\exists^{\mathbf{st}_w}Y(X \cap_w \mathbf{S}_w =_w Y \cap_w \mathbf{S}_w)$.

9.**Weak Standardization**: $\forall X\exists^{\mathbf{st}_w}y(X \cap_w \mathbf{S}_w = \cap_w \mathbf{S}_w)$. (Such consistent standard set $Y$, unique by Consistent Transfer and Consistent Extensionality, is sometimes denoted by

$^{\mathbf{s}_w}X.$)

Such inconsistent standard set $Y$, $w$-unique by $w$-Transfer and weak Extensionality, is sometimes denoted by $^{\mathbf{s}_w}X$.

**Remark 4.1**. (**i**) $w$-Transfer can be considered as saying that: $\mathbf{I}_w$, the universe of all inconsistent internal sets, is an elementary extension of $\mathbf{S}_w$ in the $\in_{\mathbf{s}}$ -$\in_w$ -language. It follows, by $(ZFC^{\#}_{\omega})^{\mathbf{st}_w}$, that the class $\mathbf{I}_w$ of all inconsistent internal sets satifies $ZFC^{\#}_{\omega}$ (in the $\in_{\mathbf{s}}$ -$\in_w$ -language ), in fact, we can replace $(ZFC^{\#}_{\omega})^{\mathbf{st}_w}$ by $(ZFC^{\#}_{\omega})^{\mathbf{int}_w}$, with relativization to $\mathbf{I}_w$, in the list of $\mathbf{HST}^{\#}_{\boldsymbol{\omega}}$ axioms. See also Theorem 1.3.9 below.

(**ii**) $w$-Transitivity of $\mathbf{I}_w$ postulates that: inconsistent internal sets to form the basement of the $\in_{\mathbf{s}}$ -$\in_w$ -structure of the universe $\mathbf{H}^{\#}_{\boldsymbol{\omega}}$. This axiom is very important since it implies that some set operations in $\mathbf{I}_w$ retain their sense in the whole universe $\mathbf{H}^{\#}_{\boldsymbol{\omega}}$.

(**iii**) $w$-Regularity over $\mathbf{I}_w$ organizes the $\mathbf{HST}^{\#}_{\boldsymbol{\omega}}$ set universe $\mathbf{H}^{\#}_{\boldsymbol{\omega}}$ in general case as a sort of hierarchy over the internal universe $\mathbf{I}_w$, in the same way as the $w$-Regularity axiom organizes the universe in the von Neumann $w$-hierarchy over the $w$-empty set $\varnothing_w$ in $ZFC^{\#}_{\omega}$.

(**iv**) Strictly $w$-Regularity organizes the $\mathbf{HST}^{\#}_{\boldsymbol{\omega}}$ set universe $\mathbf{H}^{\#}_{\boldsymbol{\omega}}$ in the von Neumann $w$-hierarchy over the $w$-empty set $\varnothing_w$, but in a strictly inconsistent sense only.

(**v**) $w$-Standardization postulates that $\mathbf{H}^{\#}_{\boldsymbol{\omega}}$ does not contain collections of standard sets other than those of the form $S \cap_w \mathbf{S}_w$ for inconsistent standard set $S$.

**Remark 4.2**. It well known that the $ZFC$ Regularity fails in $\mathbf{H} = \mathbf{H}_{\mathbf{s}}$ : the set of all nonstandard $\mathbf{I}_{\mathbf{s}}$-natural numbers does not contain an $\in_{\mathbf{s}}$ -minimal element, (see for example [18], Exercise 1.2. 15(3)). In contrast with a classical case, $ZFC^{\#}_{\omega}$ $w$-Regularity valid in $\mathbf{H}^{\#}_{\boldsymbol{\omega}}$, but in a strictly inconsistent sense only. For example the set of all nonstandard $\mathbf{I}_w$-natural contain an inconsistent $\in_w$ -minimal element, see [22]-[23].

# VII.5. Well-founded inconsistent sets.

Now we can introduce the last principal class: well-founded inconsistent sets. Recall the following notions from general inconsistent set theory.

**Definition 5.1**. (i) A binary weak relation $\prec_w$ on inconsistent set or inconsistent class $X$ is a strictly well-founded if any nonempty set $Y \subseteq_w X$ contains consistent $\prec_w$ -minimal $w$-element $x^* \in_w Y$, that is there exists $x \in_w Y$ such that no $y \in_w Y$ satisfies $y \prec_w x$.

(ii) A binary weak relation $\prec_w$ on inconsistent set or inconsistent class $X$ is weakly well-founded (or $w$-well-founded) if:

(1) $\prec_w$ is not a strictly well-founded and

(2) any nonempty set $Y \subseteq_w X$ contains a $\prec_w$ -minimal $w$-element $x^* \in_w Y$, that is there exists $y \in_w Y$ satisfies: $y \prec_w x \wedge x \prec_w y$, i.e. $y \prec_w x \wedge \neg(y \prec_w x)$.

(iii) Inconsistent set or inconsistent class $X$ is $w$-transitive if any $x \in_w X$ satisfies $x \subseteq_w X$, i.e., weak elements of weak elements of $X$ are weak elements of $X$.

(iv) Inconsistent set or inconsistent class $X$ is $w$-complete if we have $y \in_w X$ whenever $y \subseteq_w x \in_w X$, that is a weak subsets of weak elements of $X$ are weak elements of $X$.

(v) Inconsistent set $x$ is a strictly well-founded if there is a $w$-transitive set $X$ such that $x \subseteq_w X$ and the restriction $\in_w \restriction X$ is a strictly well-founded weak relation.

(vi) Inconsistent set $x$ is $w$-well-founded if there is a $w$-transitive set $X$ such that $x \subseteq_w X$ and the restriction $\in_w \restriction X$ is a $w$-well-founded weak relation.

**Remark 5.1**. It is known that all sets are well-founded in $ZFC$ by the Regularity axiom.

This is not the case in **HST** : the set $^*\mathbb{N}$ of all $\mathbf{I_s}$-natural numbers is ill-founded [18].

**Remark 5.2**. In contrast with a classical case, all inconsistent sets are $w$-well-founded in $\mathbf{HST}^{\#}_{\boldsymbol{\omega}}$ by the Strictly $w$-Regularity axiom. For example, the set $^{\#}\mathbb{N} = {}^*\mathbb{N}_{\mathbf{inc}}$ of all $\mathbf{I}_w$-natural numbers is $w$-well-founded by the Strictly Weak Regularity axiom.

**Definition 5.2**.($\mathbf{HST}^{\#}_{\boldsymbol{\omega}}$). (i) Let $\mathbf{s\text{-}wf}_w x$ mean that $x$ is a strictly well-founded. We put $\mathbf{s\text{-}WF}_w =_w \{x : \mathbf{s\text{-}wf}_w x\}_w$, the class of all strictly well-founded inconsistent sets and (ii) let $w\text{-}\mathbf{wf}_w x$ mean that $x$ is a $w$-well-founded. We put $w\text{-}\mathbf{WF}_w =_w \{x : w\text{-}\mathbf{wf}_w x\}_w$, the class of all $w$-well-founded inconsistent sets.

**Notation 5.1**.We introduce quantifiers $\exists^{\mathbf{s\text{-}wf}_w}, \forall^{\mathbf{s\text{-}wf}_w}, \exists^{w\text{-}\mathbf{wf}_w}$ and $\forall^{w\text{-}\mathbf{wf}_w}$ (meaning: there is a well-founded ... , for any well-founded ... ) and the relativization (1) $\Phi^{\mathbf{s\text{-}wf}_w}$ to $\mathbf{s\text{-}WF}_w$, (2) $\Phi^{w\text{-}\mathbf{wf}_w}$ to $w\text{-}\mathbf{WF}_w$ similarly to $\exists^{\mathbf{st_s}}, \exists^{\mathbf{st_s}}, \Phi^{\mathbf{st_s}}, \exists^{\mathbf{st}_w}, \exists^{\mathbf{st}_w}, \Phi^{\mathbf{st}_w}$ in §VII.1.3. In other words, $\Phi^{\mathbf{s\text{-}wf}_w}$ says that gj is true in WIF. The main property of the classes $\mathbf{s\text{-}WF}_w$ and $w\text{-}\mathbf{WF}_w$ in $\mathbf{HST}^{\#}_{\boldsymbol{\omega}}$ is that it admits a definable $\in_w$ -isomorphism $\mathbf{w} \mapsto {}^{\#}\mathbf{w}$ onto the class S of all standard sets.

# VIII.Inconsistent number systems

.

## VIII.1.The Graham Priest argument.

**Remark 8.1.1.**In contrast with ordinary classical (unrestricted) modus ponens rule

$$A, A \to B \vdash_{\mathbf{UMP}} B \tag{8.1.1}$$

the restricted modus ponens rule reads

$$A, A \to B \vdash_{\mathbf{RMP}} B \text{ if and only if } A \notin \Delta_1 \text{and } A \to B \notin \Delta_2, \tag{8.1.2}$$

where $\Delta_1, \Delta_2 \subsetneqq \overline{\mathcal{F}}^{\#}_{\omega}$.Thus it is not in general true by using paraconsistent first order logic $\overline{LP}^{\#}_{\omega}$ that if $A \to B$ holds and $A$ holds then $B$ holds.

**Remark 8.1.2**. In adition in logic $\overline{LP}^{\#}_{\omega}$ we distingvish a strong negation $\neg_s A$ and a weak negation $\neg_w A$. A strong negation that is ordinary classical negation, i.e. $\neg_s A$ holds if and only if $A \vdash B, B \in \overline{\mathcal{F}}^{\#}_{\omega}$.A weak negation that is nonclassical negation, i.e. $A \wedge \neg_w A$ might holds.

**Remark.8.1.3**. In particular, $((A \wedge \neg_w A) \to B) \wedge (A \wedge \neg_w A)$ might hold while $B$ does not.

**Remark.8.1.4.**In particular is that this permit a solution to the following problem raised by Graham Priest [5],[6]. Ordinarily one wants postulates such as the Cancellation Law

$$\forall x \left[ \neg(x = 0) \to \forall y \forall z (x \times y = x \times z \to y = z) \right] \tag{8.1.3}$$

to hold when moving from the classical theory of the rings to nonclassical theory of the inconsistent rings.But canonical inconsistent fields [i.e. inconsistent fields based on canonical inconsistent logic with unique negation $\neg$] have both

$$\neg(x_{\mathbf{Inc}} = 0) \text{ and } x_{\mathbf{Inc}} = 0 \tag{8.1.4}$$

for some $x_{\mathbf{Inc}}$ [for examle $(x_{\mathbf{Inc}} = 1) \wedge (x_{\mathbf{Inc}} = 0)$] and therefore

$$\forall y \forall z(x_{\mathbf{Inc}} \times y = x_{\mathbf{Inc}} \times z = 0) \quad (8.1.5)$$

holding.Yet one does not want to detach every $y = z$ or the theory is trivial. Yet one also does not want to forbid detachment for those $x$ which are classically not identical with zero.

**Remark 8.1.5.**Obviously Priest paradox arises from the statement:

$$x \text{ classically not identical with zero.} \quad (8.1.6)$$

The statement (8.1.6) completely does not well defined by using canonical inconsistent logic with unique negation $\neg$. Note that if the statement $\neg(x = 0)$ treated classically i.e. under definition of the strong negation $\neg_s$ (see Remark.8.1.1), this meant impossibilyty $x = 0$, i.e. the statement $x = 0$ is not holds classically and assuming that both $\neg(x = 0)$ and $x = 0$ holding we conclude that the statement $\neg(x = 0)$ is not holds classically in contrary with Priest assumption (8.1.6).Cancellation Law (8.1.3) breaks down.

**Remark 8.1.6.**In order to avoid the difficultness mentioned above we apply the logics $\overline{LP}^{\#}_{n}, n \geq 1$and postulate the Cancellation Law in the following form

$$\forall x\left[\neg_s(x = 0) \to \forall y \forall z(x \times y = x \times z \to y = z)\right], \quad (8.1.7)$$

and

$$\forall x\left[\neg_w(x = 0) \to \forall y, \forall z(x \times y = x \times z \to y = z)\right]. \quad (8.1.8)$$

We set now instead (2.4) that

$$\neg_w(x = 0) \text{ and } x = 0. \quad (8.1.9)$$

From (8.1.7) one to detach $y = z$ only for such $x$ which are classically not identical with zero as it should be. However if we set $(x_{\mathbf{Inc}} = 1) \wedge (x_{\mathbf{Inc}} = 0) \in \Delta_1$ and $(((x_{\mathbf{Inc}} = 1) \wedge (x_{\mathbf{Inc}} = 0)) \to 0 = 1) \in \Delta_1$ then

$$((x_{\mathbf{Inc}} = 1) \wedge (x_{\mathbf{Inc}} = 0)) \to 0 = 1 \nvdash_{\mathbf{RMP}} 0 = 1 \quad (8.1.10)$$

as it should be.

# VIII.2.Inconsistent logic with restricted modus ponens rule based on RM3-assignment.

The classical example of the inconsistent logic with restricted modus ponens rule has been proposed by C.E. Mortensen, see ref.[5-7].

Let $\mathcal{L}$ be an canonical language consisting of simple terms (names), one for each real number; function symbols $+, \times, -, \div$; atomic predicates $=, <, \in$; variables $x, y, z, \ldots$ and operators $\neg, \wedge, \forall$.

Remind that any RM3-assignment [5],[6] is a function $I$ assigning to the wff's of $\mathcal{L}$, or the appropriate sublanguage of $\mathcal{L}$ under investigation at the time, values from the set $\{\mathbf{T}, \mathbf{N}, \mathbf{F}\}$ in accordance with the following definition:

(i) for any atomic wff with terms $t_1, t_2$, we have $I(t_1 = t_2), I(t_1 < t_2)$ and $I(t_1 \in t_2)$ all

belong to $\{\mathbf{T}, \mathbf{N}, \mathbf{F}\}$, (read 'true, neuter, false');

(ii) $I(\neg A)$ and $I(A \wedge B)$ are given by the RM3-matrices:

$$\begin{array}{ccccc} \wedge & \mathbf{T} & \mathbf{N} & \mathbf{F} & \neg \\ {}^*\mathbf{T} & \mathbf{T} & \mathbf{N} & \mathbf{F} & \mathbf{F} \\ {}^*\mathbf{N} & \mathbf{N} & \mathbf{N} & \mathbf{F} & \mathbf{N} \\ \mathbf{F} & \mathbf{F} & \mathbf{F} & \mathbf{F} & \mathbf{T} \end{array} \quad (8.2.1)$$

(iii) $I((x)A) = \min\{y$ : for some term $t, I(A(t|x)) = y\}$, where min is relative to the ordering: false < neuter < true. A sentence $A$ holds in an assignment $I$ iff $I(A) \in \{\mathbf{T}, \mathbf{N}\}$.

Let us consider now the classical standard model of the natural numbers, equipped with names for the natural numbers. In view of the Extendability Lemma [8],[9], the set of sentences holding therein can be extended by adding any collection of sentences of the form $\neg(n = n)$ and evaluating in an RM3-assignment. Note that the contradiction does not spread to other sentences of the form $\neg(m = m)$. Similarly, collections of sentences of the form $n = m$ for distinct $n, m,$ may be added with the same result. This raises the following question [5]. If we add, for example, $0 = 2$ to the standard model of the natural numbers, then, in virtue of the substitutivity of identity and the fact that $\neg(0 = 2)$ also holds, have we not imported the further sentence $\neg(0 = 0)$? The answer is no, and it illustrates the generality of the Extendability Lemma. The rule of substitutivity of identity $(SI)$ in the form if $t_1 = t_2$ holds, then $Ft_1$ holds iff $Ft_2$ holds (all terms $t_1, t_2,$ with $t_2$ replacing $t_1$ in $F_1$ at least one place) does not always hold in RM3-assignments. What is the case, if the sentences holding in an RM3-assignment include those holding in the standard model of the natural numbers, is that $(t_1 = t_2 \wedge Ft_1) \to Ft_2$ holds, since it holds in the standard model.

**Remark 8.2.1.** But it is not in general true that if $A \to B$ holds and $A$ holds then $B$ holds. In particular, $((A \wedge \neg A) \to B) \wedge (A \wedge \neg A)$ might hold while $B$ does not, i.e.

$$(A \wedge \neg A) \to B \nvdash_{RM3} B, \quad (8.2.2)$$

where by $\vdash_{RM3}$ we denote the rule of conclusion corresponding to RM3-assignments. However, this leads to no loss of information from classical arithmetic, since we do have that if $(A \to B) \wedge A$ holds, and if moreover $(A \to B) \wedge A$ holds back in the standard model for arithmetic, then $B$ holds (trivial).

**Remark 8.2.2.** (i) A special case of interest is this: if $t_1 = t_2 \wedge Ft_1$ holds and if moreover $\neg(t_1 = t_2)$ and $\neg Ft_1$ both do not hold, then $Ft_2$ holds. (Reason: for then $t_1 = t_2 \wedge Ft_1$ holds back in the classical complete subtheory, wherein $Ft_2$ could be detached.)

(ii) Thus the rule $SI$ does not hold in all RM3-assignments.

# VIII.3. The da Costa type paraconsistent logic $\mathbf{C}_\infty^\#$ with infinite levels of a contradictions and restricted modus ponens rule can to save da Costa Set Theory from a trivality.

.

It well known that canonical da Costa's paraconsistent logics is invalid in order to obtain non trivial paraconsistent set theory,see [8]. In order to resolve this tension we consider the da Costa type paraconsistent logics $\mathbf{C}_{\infty}^{\#}$ with infinite levels of a contradictions and restricted modus ponens rule mentioned above.

We remind that da Costa paraconsistent set theory is a paraconsistent set theory whose underlying logic is one of da Costa's paraconsistent logics $C_n^{=}, 1 \leq n \leq \omega$.

**Definition 8.3.1.**The postulates of $C_{\omega}^{=}$ are those of the positive intuitionistic first-order logic with equality, plus:

(1) $\neg\neg A \Rightarrow A$,

(2) $A \vee \neg A$.

(3) unrestricted modus ponens rule : $A, A \Rightarrow B \vdash_{\mathbf{UMP}} B$.

**Definition 8.3.2.**The postulates of $C_n^{=}, 1 \leq n \leq \omega$, are those of $C_{\omega}^{=}$, plus:

(1) $B^{(n)} \wedge (A \Rightarrow B) \wedge (A \Rightarrow \neg B) \Rightarrow \neg A$,

(2) $A^{(n)} \wedge B^{(n)} \Rightarrow (A \Rightarrow B)^{(n)} \wedge (A \wedge B)^{(n)} \wedge (A \vee B)^{(n)}$,

(3) $(\forall x)(A(x))^{(n)} \Rightarrow ((\forall x)A(x))^{(n)} \wedge ((\exists x)A(x))^{(n)}$,

where $A^{(n)}$ defined as follows: $A^1 = A^0 = \neg(A \wedge \neg A), A^{n+1} = (A^n)^0, A^{(n)} = A^1 \wedge \ldots \wedge A^n$,

(4) unrestricted modus ponens rule :

$A, A \Rightarrow B \vdash_{\mathbf{UMP}} B$.

In each $C_n^{=}, 1 \leq n \leq \omega$, $\neg^* A$ is defined as $\neg A \wedge A^{(n)}$,and it is proved that satisfies all the properties of the classical negation. Then classical logic can be obtained inside these systems; consequently, they are finitely trivializable. For, from any formula of the form $A \wedge \neg A \wedge A^{(n)}$ one can deduce any formula whatsoever. Nonetheless, $C_{\omega}^{=}$ is not finitely trivializable. Moreover, each system terns in the hierarchy $C_1^{=}, C_2^{=}, \ldots, C_n^{=}, \ldots, C_{\omega}^{=}$ is strictly stronger than the following ones. Thus, we may construct a hierarchy of da Costa's paraconsistent set theories in which, at least intuitively, it seems that each system may admit more nonclassical sets than the preceding ones.

**Definition 8.3.3.**Let $\mathcal{L}_{\omega} = \mathcal{L}_{\omega}(C_{\omega}^{=})$ be the formal language corresponding to logic $C_{\omega}^{=}$ and let $W_{\omega} = W_{\omega}(\mathcal{L}_{\omega})$ be the set of the all wff's of $\mathcal{L}_{\omega}$The postulates of $C_{\omega}^{\#}$ are those of the postulates of $C_{\omega}^{=}$ but with restricted modus ponens rule : $A, A \Rightarrow B \vdash_{\mathbf{RMP}_{\omega}} B$ if and only if $(A,B) \in \Delta_{\omega} \times \Delta_{\omega}$ [instead unrestricted modus ponens rule (3)], where $\Delta_{\omega} \subsetneq W_{\omega}$.

**Definition 8.3.4.**Let $\mathcal{L}_n = \mathcal{L}_n(C_n^{=})$ be the formal language corresponding to logic $C_n^{=}$ and let $W_n = W_n(\mathcal{L}_n)$ be the set of the all wff's of $\mathcal{L}_{\omega}$The postulates of $C_n^{\#}, 1 \leq n \leq \omega$, are those of $C_n^{=}, 1 \leq n \leq \omega$,but with restricted modus ponens rule : $A, A \Rightarrow B \vdash_{\mathbf{RMP}_n} B$ if and only if $(A,B) \in \Delta_n \times \Delta_n$instead unrestricted modus ponens rule (4), where $\Delta_n \subsetneq W_n(C_n^{=})$

**Definition 8.3.5.**The postulates of $C_{\infty}^{\#}$,are $\bigcup_{n=1}^{\omega} C_n^{\#}$ with restricted modus ponens rule :

$A, A \Rightarrow B \vdash_{\mathbf{RMP}_{\infty}} B$ if and only if $(A,B) \in \bigcup_{n=1}^{\omega} \Delta_n \times \Delta_n$.

We remind that da Costa paraconsistent set theories $NF_n^C$ are constructed very similarly

to $NF$. The main postulates of $NF_{\omega}^C$ are the following [8]:

**I.Extensionality**

$$\forall \alpha \forall \beta \forall x[x \in \alpha \Leftrightarrow x \in \beta \Rightarrow \alpha = \beta]. \tag{8.3.1}$$

**II.Abstraction**

$$\exists\alpha\forall x[x \in \alpha \Leftrightarrow F(x)], \qquad (8.3.2)$$

where $\alpha$ does not occur free in $F(x)$ and $F(x)$ is stratified or it does not contain any formula of the form $A \Rightarrow B$.

## VIII.4. The paraconsistent set theory $NF^{\#}_{\infty}$based on logic $C^{\#}_{\infty}$ with restricted modus ponens rule and infinite levels of a contradiction.

**Definition 8.4.1**.The main postulates of $NF^{\#}_{\omega}$ are those of the postulates of $NF^{C}_{\omega}$ but with logic of type $C^{\#}_{\omega}$ instead logic $C^{=}_{\omega}$.

**Definition 8.4.2**.The main postulates of $NF^{\#}_{n}, 1 \leq n \leq \omega$, are those of the postulates of $NF^{C}_{n}$ but with logic of type $C^{\#}_{n}, 1 \leq n \leq \omega$,instead logic $C^{=}_{n}$.

**Definition 8.4.3**.The main postulates of $NF^{\#}_{\infty}$,are those of the postulates of $NF^{C}_{\omega}$ but with logic of type $C^{\#}_{\infty}$.

Da Costa's paraconsistent set theories of type $NF^{C}_{\omega}$ and $NF^{C}_{n}$, $1 \leq n \leq \omega$. has been studying A.I. Arruda [8].A.I. Arruda has been proved that da Costa's formulation of the axiom schema of abstraction (1.2.2) for the systems $NF_n$, $1 \leq n < \omega$, leads to the trivialization of the systems (see [8]).

**Remark 8.4.1**.Note that in $NF^{C}_{\omega}$, the restrictions regarding the use of non-stratified formulas obstruct a direct proof of the paradox of Curry. Russell's set $\mathfrak{R}$, defined as $\hat{x}\neg(x \in x)$, exists as well as many other non-classical sets. The paradox of Russell in the form $\mathfrak{R} \in \mathfrak{R} \wedge \neg(\mathfrak{R} \in \mathfrak{R})$ is derivable but apparently, it causes no ham to the system.Due to its weakness, the primitive negation of $NF^{C}_{\omega}, \neg$, is almost useless for set-theoretical purposes. Thus, let us define

$$\sim A \text{ for } A \Rightarrow \forall x\forall y[x \in y \wedge x = y]. \qquad (8.4.1)$$

The universal set $\mathbf{V}$ is defined as $\hat{x}(x = x)$, the empty set $\emptyset$ as $\hat{x} \sim (x = x)$,and the complement of a set $\alpha$, $\bar{\alpha}$, as $\hat{x}\sim(x \in a)$.

**Theorem 5.1**.[8]. In $NF^{C}_{\omega},\sim$ is a minimal intuitionistic negation.

**Corollary 1**. $\vdash A \Rightarrow (\sim A \Rightarrow \sim B)$, $\vdash (A \Rightarrow B) \Rightarrow (\sim B \Rightarrow \sim A)$.

**Corollary 2**. All the theorems of $NF$ whose proofs depend only on the laws of the minimal intuitionistic first-order logic with equality and on the postulates of extensionality and abstraction of $NF$ are valid in $NF^{C}_{\omega}$.

**Theorem 5.2**.[8].(Cantor's Theorem) $NF^{C}_{\omega} \vdash \sim(\alpha \leq P(\alpha))$.

**Corollary**.[8].(Cantor's Paradox) $NF^{C}_{\omega} \vdash (\mathbf{V} \leq P(\mathbf{V})) \wedge \sim(\mathbf{V} \leq P(\mathbf{V}))$.

**Remark 5.2**.Note that Cantor's paradox does not trivialize $NF^{C}_{\omega}$, since from $A$ and $\neg A$ we cannot obtain any formula $B$ whatsoever. For instance, apparently, we cannot obtain
any formula of the form $\neg B$, where $B$ is a nonatomic formula.

**Theorem 5.3**.[8].(i) $NF^{C}_{\omega} \vdash \forall\alpha\forall\beta[(\alpha = \beta) \wedge \sim(\alpha = \beta)]$,(ii) $NF^{C}_{\omega} \vdash [(\alpha \in \beta) \wedge \sim(\alpha \in \beta)]$, (iii) $NF^{C}_{\omega} \vdash [(\alpha \in \alpha) \wedge \sim(\alpha \in \alpha)]$.

**Proof**. By the corollaries of theorems 5.1 and 5.2, we obtain

$$x = x \Rightarrow \delta,$$
$$\delta \Leftrightarrow \forall\alpha\forall\beta[(\alpha \in \beta) \wedge (\alpha = \beta)] \quad (5.24)$$

Thus, as $x = x,$ then $\forall\alpha\forall\beta(\alpha = \beta)$. By the same corollaries we also obtain $\forall\alpha\forall\beta[\sim(\alpha = \beta)]$.The proof of part (ii) is similar to that of part (i). Part (iii) is an immediate consequence of part (ii).

**Remark 5.3**. We introduce now the logic of type $C_\omega^\#$ with $\Delta_\omega = \Delta_\omega^\delta$ such that $(x = x,\delta) \notin \Delta_\omega^\delta$.Thus in $NF_\omega^\#$ with $\Delta_\omega = \Delta_\omega^\delta$ Theorem 1.2.3 no longer holds.

**Remark 5.4**. Note that: (i) By Theorem 1.2.3, it could seem that $NF_\omega^C$ is trivial. Nonetheless, apparently this is not the case.

(ii) However, though it is nontrivial, $NF_\omega^C$ is without interest, for not only are every two sets identical, but also every set belongs and does not belong to itself.

**Remark 5.5**.In order to avoid the results mentioned in Theorem 5.3, one could think of introducing more restrictions in da Costa's formulation of the axiom schema of abstraction when $F(x)$ is non-stratified. Nonetheless, we be live that this is a worthless effort. For:

(i)The only non-stratified formula used in the proof of Cantor's Theorem (which is fundamental in the proof of Theorem 5.3) is a non-stratified formula of the form $\alpha \in \beta$. Then, the new restrictions must avoid those nonstratified atomic formulas of the form $\alpha \in \beta$ which determine a set.

(ii) A new proof of Theorem 5.3 may be obtained in the following way: in $NF$ the formula $y = \{x\}$ cannot determine a relaticm because $\langle x,y\rangle = \langle x,y\rangle \wedge (y = \{x\})$ is non-stratified. But, such a formula does not contain any subformula of the form $A \Rightarrow B$; then, in $NF_\omega^C$ it determines a relation $S$ such that $S \in 1-1$, see [7],pp.12. With such a relation we prove that $\alpha \leq P(\alpha)$. In $NF_\omega^C$ we also prove that $\sim(\alpha \leq P(\alpha))$. Then, these new restrictions must also avoid that those non-stratified formulas whose atomic sub formulas are of the form $\alpha = \beta$ determine a set.

(iii) From the above remarks (i) and (ii) we conclude that, in order to avoid the counterintuive results mentioned in Theorem 5.3, the axiom schema of abstraction in $NF_\omega^C$ should be formulated as in $NF$.

**Remark 5.6**.Due to the paradoxes obtained in $NF_\omega^C$, we conclude that in these system the axiom schema of abstraction should be formulated as in $NF$. Thus, if we want these theories to be paraconsistent set theories, we need to postulate directly the existence of contradictory sets. .Apparently, we may postulate the existence of Russell's set without any problem. Nonetheless, due to the two above considerations about the non-stratified formulas that lead to the proof of the paradox of identity, we believe that, besides Russell's set, very few other non-classical sets may exist in $NF_n,\ 1 \leq n \leq \omega$.

**Definition 5.4**.[8].Let us denote by $\mathbf{DC}_n$ any da Costa set theory based on the respective $C_n^=$, where Russell's class is a set. Thus, in $\mathbf{DC}_n, 1 < n < \omega,$ the defined negation $\neg^* A \Leftrightarrow \neg A \wedge A^{(n)}$ is a classical negation; and in $\mathbf{DC}_\omega$ the defined negation $\sim A \Leftrightarrow A \Rightarrow \forall x\forall y[x \in y \wedge x = y]$ is a minimal intuitionistic negation.

**Theorem 5.4**.[8].Let $\mathfrak{R}$ be Russell's set. In $\mathbf{DC}_n, 1 < n < \omega, \cup\cup\mathfrak{R}$ is the universal set.

**Definition 5.5**.[8].Let $\mathbf{DC}_{\omega}^{\mathbf{V}}$ be a $\mathbf{DC}_{\omega}$ with universal set $\mathbf{V}$ defined as $\hat{x}(x = x)$.
Let us define $x \neq \mathbf{V}$ for $\exists y[\sim(y \in x)]$. We introduce now the postulate:
**P**1.$\forall x[(x = \mathbf{V}) \vee (x \neq \mathbf{V})]$.
**Theorem 1.2.5**.[8]. In $\mathbf{DC}_{\omega}^{\mathbf{V}} + \mathbf{P}1$ is derivable $\cup\cup\Re = \mathbf{V}$.

# VIII.5. The paraconsistent set theory $ZF_{\infty}^{\#}$ based on logic $C_{\infty}^{\#}$ with restricted modus ponens rule and infinite levels of a contradiction.

In this section da Costa's set theories of type $ZF$ we denoted by $ZF_n, 1 \leq n \leq \omega$.

A.I. Arruda has been proved that da Costa's set theories of type $ZF$ incompatible with the existence of Russell's set $\Re$ [8].

Le us consider the set theories $ZF_n, 1 \leq n \leq \omega$, in which the axioms of pairing and inion are postulated in general, and in which we also postulate the existence of the espty set amI of RusselltsS set. Moreover, let us suppose that there is no universal set, i.e.,

$S_n$. $\forall x \exists y[\neg^*(y \in x)]$, in $ZF_n, 1 \leq n < \omega$;

$S_{\omega}$. $\forall x \exists y\ [\sim(y \in x)]$, in $ZF_{\omega}$.

**Theorem 6.1**.[8].The set theories $ZF_n, 1 \leq n < \omega$ plus $S_n$ are trivial.

**Proof**. By $S_n$ there exists $y$ such that $\neg^*(y \in \cup\cup\Re)$. By part II of Lemma 4.3 [7], and part I of Lemma 4.4 [8], we obtain $\forall x(x \in \cup\cup\Re)$. Consequently,

$$y \in \cup\cup\Re \wedge \neg^*(y \in \cup\cup\Re), \qquad (8.5.1)$$

and this formula trivializes the system.

**Theorem 6.2**.[8]. The paradox of identity is derivable in $ZF_{\omega}$ plus $S_{\omega}$.

**Proof**. By $S_{\omega}$ there exists $y$ such that $\sim(y \in \cup\cup\Re)$. Using part II of Lemma 4.3 [8], we obtain $y \in \cup\cup\Re$. Consequently, by the definition of ~ it follows $\forall x \forall y(x \in y \wedge x = y)$ and therefore, the paradox of identity, $\forall x \forall y(x = y)$,follows.

**Theorem 6.3**.[8].. The systems $ZF_n$ ,$1 \leq n < \omega$ with Russell's set and the axiom schema of separation postulate for all sets are trivial.

**Proof**.If the axiom schema of separation is postulated for all sets then there exists a subset $\alpha$ of $\Re$ such that (1) $\forall x\left[x \in \alpha \Leftrightarrow (x \in \Re) \wedge (x \in x)^{(n)}\right]$. From (1) we obtain (2) $\alpha \in \alpha \Leftrightarrow \neg(\alpha \in \alpha) \wedge (\alpha \in \alpha)^{(n)}$. Consequently, we have $(\alpha \in \alpha) \wedge \neg^*(\alpha \in \alpha)$, and this formula trivializes the system.

**Definition 6.1**.We introduce now paraconsistent logic of type $C_{\infty}^{\#}, = \bigcup_{n=1}^{\omega} C_n^{\#}$ with restricted modus ponens rule : $A, A \Rightarrow B \vdash_{\mathbf{RMP}_{\infty}} B$ if and only if $(A,B) \in \Delta_{\infty} = \bigcup_{n=1}^{\omega} \Delta_n$, and we choose a set $\Delta_{\infty}$ such that

$$(\alpha \in \alpha) \wedge \neg^*(\alpha \in \alpha) \vdash_{\mathbf{RMP}_{\infty}} \varnothing, \qquad (6.2)$$

i.e.,from $(\alpha \in \alpha) \wedge \neg^*(\alpha \in \alpha)$ we cannot obtain any formula whatsoever.

**Definition 6.2**.(i)The main postulates of $ZF_n^{\#}$ are those of the postulates of $ZF_n^{\#}$ but with logic of type $C_n^{\#}$ instead logic $C_n^{=}$.(ii) The main postulates of $ZF_{\infty}^{\#}$ are $\bigcup_{n=1}^{\omega} ZF_n^{\#}$ with logic of type $C_{\infty}^{\#}$ mentioned above in Definition 6.1

. **PART III**.

# I.Introduction

## I.1.Carleson's theorem and generalizations in dimention $N = 1$.

L.Carleson's celebrated theorem of 1965 [25] asserts the pointwise convergence of the partial Fourier sums of square integrable functions. The Fourier transform has a formulation on each of the Euclidean groups $\mathbb{R}, \mathbb{Z}$ and $\mathbf{T}$.Carleson's original proof worked on $\mathbf{T}$.bFefferman's proof translates very easily to $\mathbb{R}$ . Máté [26] extended Carleson's proof to $\mathbb{Z}$.Each of the statements of the theorem can be stated in terms of a maximal Fourier multiplier theorem [27]. Inequalities for such operators can be transferred between these three Euclidean groups, and was done P. Auscher and M.J. Carro [28]. But Carleson's original proof and another proofs very long and very complicated. We give a very short and very "simple" proof of this fact. Our proof uses PNSA technique only, developed in part I, and does not uses complicated technical formations unavoidable by the using of purely standard approach to the present problems. In contradiction to Carleson's method, which is based on profound properties of trigonometric series, the proposed approach is quite general and allows to research a wide class of analogous problems for the general orthogonal series. Let us suppose that there are general orthogonal series in space $\mathscr{L}_2(\Omega)$ $\Omega \subseteq \mathbb{R}^d, d = 1, 2\ldots$

$$\sum_{n=0}^{\infty} c_n f_n(x), \{c_n\}_{n=0}^{n=\infty} \in l_2, f_n \in \mathscr{L}_2(\Omega), n \in \mathbb{N}.$$
$$\int_{\Omega} f_i(x) \cdot f_j(x) d^N x = \delta_{ij}. \tag{1.1.1}$$

We shall say that a sequence $\{f_n\}_{n=0}^{n=\infty}$ or series (1.1.1) admit $LC$-property if series (1.1.1) converges a.e. It is well known that a general orthogonal series does not admit $LC$-property [29-30].

**Definition 1.1.1**. We shall say that for orthogonal series (1.1.1) $LC$-property holds iff series (1.1.1) converges a.e. on a set $\Omega$.

A problem corresponding to $LC$-property is still open for many orthogonal series, as example for the series by Jakoby's polynomial. In the present work we shall obtain

a

general sufficient condition guaranteeing the $LC$-property for series (1.1).

**Definition 1.1.2**. We shall say that orthogonal series (1.1.1) in a space $\mathbf{L}_2(\Omega)$ is a

strongly paraorthogonal series, iff the following condition is satisfies

$$
{}^{\#_w}\int_{{}^{\#_w}\Omega}\left[{}^{\#_w}f_i(x)\cdot{}^{\#_w}f_j(x)\right]d^N x =_{w^*} ({}^{\#_w}\delta_{ij}),\quad i,j \in_w {}^{\#_w}\mathbb{N}, \tag{1.1.2}
$$

$$
{}^{\#_w}f_{\mathbf{i}}(x) \in_w {}^{\#_w}L_2({}^{\#_w}\Omega), i \in_w {}^{\#_w}\mathbb{N}.
$$

Here

$$
{}^{\#_w}\delta_{ij} =_{w^*} 1_{w^*} \Leftrightarrow i =_{w^*} j; {}^{\#_w}\delta_{ij} =_{w^*} 0_{w^*} \Leftrightarrow i \neq_{w^*} j \quad \text{and} \quad {}^{\#_w}\delta_{ij} =_w 1_{w^*} \Leftrightarrow i =_w j; {}^{\#}\delta_{ij} =_w 0_{w^*} \Leftrightarrow i \neq_w j. \tag{1.1.3}
$$

# I.2.Carleson's theorem and generalizations in dimentions $N \geq 2$.

Carleson's results are trivially transferred on $N$ -harmonic Fourier series, for the case of convergence by cubes, but in the case of arbitral convergence Carleson methods does not works and,in general,the problem for $N$ -harmonic Fourier series is still open.

Particularly,this problem is open for the case of the spherical sum $E_M[f(x)], x \in \mathbb{R}^N$ :

$$
E_M[f(x)] = (2\pi)^N \sum_{\|n\|^2 \leq M} f_n \cdot \exp(inx),\quad n \in \mathbb{Z}^N,\quad \|n\| = \sqrt{\textstyle\sum_{i=1}^N n_i^2}. \tag{1.2.1}
$$

In 1971 R. Cooke proved Cantor-Lebesque theorem in two dimentions [30]: if

$$
\lim_{k\to\infty} \sum_{|n|^2 = k} c_n \exp(inx) = 0 \tag{1.2.2}
$$

a.e. on $\mathbf{T}^2$, then

$$
\lim_{k\to\infty} \sum_{|n|^2 = k} |c_n|^2 = 0. \tag{1.2.3}
$$

# I.3.The uniqueness problem of the trigonometric expansion

## in dimention $N = 1$. Cantor-Lebesque theorem in dimention $N = 1$.

The uniqueness problem of the trigonometric expansion in dimention $N = 1$ can be stated as follows. Suppose the series

$$\frac{a_0}{2} + \sum_{n=1}^{\infty}(a_n \cos nx + b_n \sin nx) \qquad (1.3.1.)$$

converges to zero for every $x \in [-\pi, \pi]$, does it follow that $a_n = b_n = 0$ for all $n$ ? The answer is not obvious, but was found to be affirmative by Cantor in 1870.

**Theorem** *1.3.1. (Cantor's uniqueness theorem). If the series (2.3.1) converges everywhere to zero, then $a_n = b_n = 0$ for all $n \in \mathbb{N}$.*

Let us briefly discuss the proof of Theorem 1.3.1. The first who systematically studied everywhere convergent trigonometric series was Riemann, in his habilitation thesis (1854). He had the idea to introduce the function

$$F(x) = \frac{a_0}{4}x^2 + \sum_{n=1}^{\infty} \frac{a_n \cos nx + b_n \sin nx}{n^2} \qquad (1.3.2)$$

obtained by formally integrating an everywhere convergent series (1.3.1) twice. Riemann assumed that the coefficients an, bn are bounded, in which case the series (1.3.2) converges uniformly and hence $F(x)$ is a continuous function on $\mathbb{R}$ (note that $F(x)$ is not periodic if $a_0 \neq 0$). He then proved that the Schwartz second derivative

$$D^2F(x) = \frac{F(x+h) - 2F(x) + F(x-h)}{h^2} \qquad (1.3.3)$$

exists, and is equal to (1.3.1).Cantor proved that the coefficients $a_n, b_n$ are tending to zero (and in particular, they are bounded). If we now assume that (1.3.1) converges everywhere to zero, then $D^2F(x) = 0$. It is then possible to prove that $F(x)$ is linear,which quite easily implies that $a_n = b_n = 0$. For more details see [6,chapt.I].

Let us consider the uniqueness problem for a trigonometric expansion which converges almost everywhere. That is, suppose a function $f(x)$ admits a trigonometric expansion such that (1.3.1) holds for almost every $x$. Is the expansion unique? Equivalently, suppose that (1.3.1) converges to zero for almost every $x$, does it follow that $a_n = b_n = 0$ for all $n$?

Lebesgue developed his theory of measure and integration in the years 1902-1906. In the following years it became common to consider sets of measure zero "negligible".

**Theorem** *2.3.2.(Cantor-Lebesgue).If $a_n \cos nx + b_n \sin nx \to 0$ for all $x$ in some set $E$ of one-dimensional positive measure, then* $a_n, b_n \to 0$.

**Proof**. By Egorov's theorem we may assume that $u_n(x) = a_n \cos nx + b_n \sin nx$ tends to $0$ uniformly on some set $E$ of positive measure. Consider the equations $u_n(x) = a_n \cos nx + b_n \sin nx$ and $u_n(y) = a_n \cos ny + b_n \sin ny$ as a linear system with unknowns $a_n, b_n$. The determinant of this system is $\sin n(y-x)$. Since $E$ has positive measure, the set $\check{E} = \{y - x | x, y \in E\}$ contains some interval $(-\delta, \delta)$ (see [41], Lemma 3.37, p. 46), therefore for any sufficiently large $n$ there exist $x, y \in E$ such that $y - x = \frac{\pi}{2n}$. For such $x, y$ we have $\sin n(y-x) = 1$, hence the above system determines $a_n, b_n$ uniquely, $a_n = u_n(x) \sin ny - u_n(y) \sin nx; b_n = u_n(y) \cos nx - u_n(x) \cos ny$. Therefore

$|a_n|, |b_n| \leq \lim_{n\to\infty} \left( 2 \sup_{x \in E} |u_n(x)| \right) = 0,$ and so Theorem 2.3.2 is proved.

**Theorem** *2.3.3.(Menshov).There exists a non-zero series (2.3.1) which converges to zero for almost every $x$.*

**Lemma** *(Menshov).There exists continuous function $F(x)$ such that:*
*(1) $F(x) \not\equiv const$ on $[0, 2\pi]$;*
*(2) $F(x) = c$ for all $x$ in some set $P$ of Lebesgue measure zero;*
*(3) the equality*

$$\lim_{n\to\infty} \int F(\alpha) \cos n(\alpha - x) d\alpha = 0 \qquad (2.3.4)$$

is satisfied uniformly on $[0, 2\pi]$.

**Proof**. (a) We define a set $P$ in the following way.From the interval $[0, 2\pi]$ we remove a central open interval such that there remain *two closed intervals* of equal length $\pi$. From each of these two intervals we remove again a central interval such that there remain $4$ closed intervals of length $\frac{2\pi}{3}$. Continuing this process, on the $k$-th step there remain $2^k$ closed intervals of length $\frac{2\pi}{k+1}$.
(b) Supose that $\delta_i, 1 \leq i \leq 2^k - 1$ is any one of the intervals which was deleted from the interval $[0, 2\pi]$ in the $k$-th step of the above procedure.

# I.4.The uniqueness problem of the trigonometric expansion in dimentions $N \geq 1$. Cantor-Lebesque theorem in dimentions $N \geq 2$.

Let $\mathbb{T}^N = [0,1)^N \subset \mathbb{R}^N$ be the $N$ dimensional torus.Let $\{\mathbf{f}_n(x)\}_{n \in \mathbb{N}}$ be a real or complex valued system of functions that are in

$$L_2(\mathbb{T}^N) \triangleq \left\{ f(x) | f : \mathbb{T}^N \to \mathbb{C}; \int_{\mathbf{T}^N} | f(x) |^2 d^N x < \infty \right\}. \tag{1.4.1}$$

The inner products $\langle \bullet, \bullet \rangle : L_2(\mathbb{T}^N) \times L_2(\mathbb{T}^N) \to \mathbb{C}$ in $L_2(\mathbb{T}^N)$ is

$$\langle \mathbf{f}_n, \mathbf{f}_m \rangle = \int \mathbf{f}_n(x), \overline{\mathbf{f}_m(x)} d^N x \tag{1.4.2}$$

where the bar denotes complex conjugate.If satisfy

$$\begin{gathered} \langle \mathbf{f}_n, \mathbf{f}_m \rangle = 0 \text{ if } n \neq m, \\ \langle \mathbf{f}_n, \mathbf{f}_m \rangle = 1 \text{ if } n = m, \\ n, m \in \mathbb{N} \end{gathered} \tag{1.4.3}$$

we call the system $\{\mathbf{f}_n(x)\}_{n \in \mathbb{N}}$ orthonormal (**ON**).Given an **ON** system and a function $f(x)$ on $\mathbb{T}^N$ it is often possible to represent $f(x) \in L_2(\mathbb{T}^N)$ as an infite linear combination of the elements of the system.

**Definition** *1.4.1.If the linear combination,* $\sum_{n \in \mathbb{N}} a_n \mathbf{f}_n(x)$ *be everywhere pointwise convergent to the value* $f(x)$, *i.e.*

$$\forall x \in \mathbb{T}^N : f(x) = \lim_{n \to \infty} \sum_{n \in \mathbb{N}} a_n \mathbf{f}_n(x), \tag{1.4.4}$$

3.2.2.If the linear combination, $\sum_{n \in \mathbb{N}} a_n \mathbf{f}_n(x)$ be o.e. pointwise convergent to the value $f(x)$, i.e.

$$\begin{gathered} \forall x \in \mathbb{T}^N \backslash E \\ \mu(E) = 0 : \\ f(x) = \lim_{n \to \infty} \sum_{n \in \mathbb{N}} a_n \mathbf{f}_n(x), \end{gathered} \tag{1.4.5}$$

In 1971 R.Cooke proved Cantor-Lebesque theorem in two dimentions [31]:

# Chapter II.Analysis on ${}^{\#_{w_{\{1\}}}}\mathbb{R}_{w_{\{1\}}}$.

## II.1.Paraordered fields.

### II.1.1.Designations

Remind that $\alpha^{[n]}$ stands for $\alpha^{[n-1]} \wedge (\alpha^{[n-1]})^{[0]}$, where $\alpha^{[0]} \triangleq \alpha \wedge \neg_w \alpha, 1 \leq n < \omega$.

**Designations 2.1.1.** In this section we will be write for short $x =_{w_{[n]}} y$ instead $(x =_w y)^{[n]}, n = 1,2,\ldots$; and we will write for short $x <_{w_{[n]}} y$ instead $(x <_w y)^{[n]}, n = 1,2,\ldots$

**Remark 2.1.1**. Thus we will be write

$$x =_{w_{[1]}} y \quad \text{instead} \quad (x =_w y) \wedge \neg_w (x =_w y) \tag{2.1.1}$$

etc. and we will be write

$$x <_{w_{[1]}} y \quad \text{instead} \quad (x <_w y) \wedge \neg_w (x <_w y) \tag{2.1.2}$$

etc. and we will be write

$$x \leq_{w_{[1]}} y \quad \text{instead} \quad \left(x <_{w_{[1]}} y\right) \vee (x =_{w_{[1]}} y). \tag{2.1.3}$$

**Remark 2.0.2**. In this section, we will be distinguish:

(1) the relations:

(i) strong (consistent) equality denoted by $(\bullet =_s \bullet)$,

(ii) weak equality denoted by $(\bullet =_w \bullet)$,

(iii) weak (inconsistent) equalities denoted by $(\bullet =_{w_{[1]}} \bullet),\ldots,(\bullet =_{w_{[n]}} \bullet),\ldots, n = 1,2,\ldots$

(2) (i) strong (consistent) inequality denoted by $(\bullet <_s \bullet)$,

(ii) weak inequality denoted by $(\bullet <_w \bullet)$,

(iii) weak (inconsistent) inequalities denoted by $(\bullet <_{w_{[1]}} \bullet),\ldots,(\bullet <_{w_{[n]}} \bullet),\ldots, n = 1,2,\ldots$ .

(iv) weak (inconsistent) inequalities denoted by $(\bullet <_{w_{\{1\}}} \bullet),\ldots,(\bullet <_{w_{\{n\}}} \bullet),\ldots, n = 1,2,\ldots$ .

**Designations 2.1.2**. (**I**) We will be write for short:

(i) $(\bullet \neq_s^s \bullet)$ instead $\neg_s(\bullet =_s \bullet)$,
(ii) $(\bullet \neq_s^w \bullet)$ instead $\neg_w(\bullet =_s \bullet)$,
(iii) $(\bullet \neq_w^s \bullet)$ instead $\neg_s(\bullet =_w \bullet)$,
(iv) $(\bullet \neq_w^w \bullet)$ instead $\neg_w(\bullet =_w \bullet)$,
(**II**) We will be write for short:
(i) $x =_{w_{\{0\}}} y$ instead $\left[(x =_s y) \vee (x =_w y) \wedge \left[\neg_w(x =_{w_{[1]}} y)\right]\right]$,
(ii) $x =_{w_{\{1\}}} y$ instead $\left[(x =_s y) \vee (x =_w y) \vee (x =_{w_{[1]}} y) \wedge \left[\neg_w(x =_{w_{[2]}} y)\right]\right]$,
(iii) $x =_{w_{\{n\}}} y$ instead $\left[(x =_s y) \vee (x =_w y) \vee \ldots \vee (x =_{w_{[n]}} y) \wedge \left[\neg_w(x =_{w_{[n+1]}} y)\right]\right]$,
$n = 1,2,\ldots$
(iv) $x =_{w_{\{\omega\}}} y$ instead $(x =_s y) \vee (x =_{w_{\{0\}}} y) \vee \bigvee_{0<n<\omega}(x =_{w_{[n]}} y)$,
**Remark 2.1.3.**(i) Note that in general case $\neg_s(x <_w y) \nvdash y <_w x$, i.e. in general case

$$x <_w y \not\Leftrightarrow_s \neg_s(y <_w x). \tag{2.1.4}$$

We often will be write for short: $x \nless_w^s y$ instead $\neg_s(x <_w y)$.
(ii) For any $x$ and $y$ such that $\neg_s(x <_w y) \vdash y <_w x$ we will be write for short:

$$x \prec_w^s y \tag{2.1.5}$$

instead $x <_w y$, i.e.we will be write $x \prec_w^s y$ iff

$$\neg_s(x \prec_w y) \Leftrightarrow_s y \prec_w x. \tag{2.1.6}$$

We often will be write for short: $x \nprec_w^s y$ instead $\neg_s(x \prec_w^s y)$.
(iii) Note that in general case $\neg_w(x <_w y) \nvdash y <_w x$, i.e. in general case the statement $\neg_w(x <_w y)$ does not imply provability of the statement $y <_w x$ and therefore in general case

$$\neg_w(x <_w y) \not\Leftrightarrow_s y <_w x. \tag{2.1.7}$$

We often will be write for short: $x \nless_w^w y$ instead $\neg_w(x <_w y)$
(iv) For any $x$ and $y$ such that $\neg_w(x <_w y) \vdash y <_w x$ we will be write for short:

$$x \prec_w^w y \tag{2.1.8}$$

instead $x <_w y$, i.e.we will be write $x \prec_w^w y$ iff

$$\neg_w(x \prec_w y) \Leftrightarrow_s y \prec_w x. \tag{2.1.9}$$

We often will be write for short: $x \nprec_w^w y$ instead $\neg_w(x \prec_w^w y)$.
(v) $(x <_w y) \wedge \neg_w(x <_w y) \not\Leftrightarrow_s \neg_w[_s(y <_w x) \wedge \neg_w(y <_w x)]$ or $x \leq_{w_{[1]}} y \not\Leftrightarrow_s \neg_w\left(x \leq_{w_{[1]}} y\right)$ in general case, i.e.
**Designations 2.1.3**. (**I**) We will be write for short:
(ii) $x =_{w_{(1)}} y$ instead $[(x =_w y) \vee (x =_{w_{[1]}} y)]$ etc.
**Designations 2.1.4**. We will be write for short:
(i) $x <_{w_{\{0\}}} y$ instead $\left[(x <_s y) \vee (x <_w y) \wedge \left[\neg_w(x <_{w_{[1]}} y)\right]\right]$,
(ii) $x <_{w_{\{1\}}} y$ instead $\left[(x <_s y) \vee (x <_w y) \vee (x <_{w_{[1]}} y) \wedge \neg_w(x <_{w_{[2]}} y)\right]$,
(iii) $x <_{w_{\{n\}}} y$ instead $\left[(x <_s y) \vee (x <_w y) \vee \ldots \vee (x <_{w_{[n]}} y) \wedge \left[\neg_w(x <_{w_{[n+1]}} y)\right]\right]$,
$n = 1,2,\ldots$ .
**Designations 2.1.5**. We will be write for short:
(i) $x \leq_{w_{\{0\}}} y$ instead $[(x \leq_s y) \vee (x <_w y)] \wedge \neg_s(x \leq_{w_{[1]}} y)$,
(ii) $x \leq_{w_{\{1\}}} y$ instead $[(x \leq_s y) \vee (x \leq_w y) \vee x <_{w_{[1]}} y] \wedge \neg_s(x \leq_{w_{[2]}} y)$,

(iii) $x \leq_{w_{\{n\}}} y$ instead $[(x \leq_s y) \vee (x \leq_w y) \vee \ldots \vee (x \leq_{w_{[n]}} y)] \wedge \neg_s(x \leq_{w_{[n+1]}} y)$, $n = 1,2,\ldots$ .

**Designations 2.1.6**. We will be write for short:

(i) $x <_{w_{\{0\}}} y$ instead $[(x <_s y) \vee (x <_w y)] \wedge \neg_s(x <_{w_{[1]}} y)$,

(ii) $x <_{w_{\{1\}}} y$ instead $[(x <_s y) \vee (x <_w y) \vee x <_{w_{[1]}} y] \wedge \neg_s(x <_{w_{[2]}} y)$,

(iii) $x <_{w_{\{n\}}} y$ instead $[(x <_s y) \vee (x <_w y) \vee \ldots \vee (x <_{w_{[n]}} y)] \wedge \neg_s(x <_{w_{[n+1]}} y)$, $n = 1,2,\ldots$ .

**Designations 2.1.7**. We will be write for short:

(i) $x \leq_{w_{\{0\}}} y$ instead $[(x \leq_s y) \vee (x <_w y)] \wedge \neg_s(x \leq_{w_{[1]}} y)$,

(ii) $x \leq_{w_{\{1\}}} y$ instead $[(x \leq_s y) \vee (x \leq_w y) \vee x <_{w_{[1]}} y] \wedge \neg_s(x \leq_{w_{[2]}} y)$,

(iii) $x \leq_{w_{\{n\}}} y$ instead $[(x \leq_s y) \vee (x \leq_w y) \vee \ldots \vee (x \leq_{w_{[n]}} y)] \wedge \neg_s(x \leq_{w_{[n+1]}} y)$, $n = 1,2,\ldots$ .

**Remark 2.1.4**.(i) Note that in general case $\neg_s(x <_{w_{\{0\}}} y) \nvdash y <_{w_{\{0\}}} x$, i.e. in general case

$$x <_{w_{\{0\}}} y \not\Leftrightarrow_s \neg_s(y <_{w_{\{0\}}} x). \qquad (2.0.10)$$

We often will be write for short: $x \not<^s_{w_{\{0\}}} y$ instead $\neg_s(x <_{w_{\{0\}}} y)$.

(ii) For any $x$ and $y$ such that $\neg_s(x <_{w_{\{0\}}} y) \vdash y <_{w_{\{0\}}} x$ we will be write for short:

$$x \prec^s_{w_{\{0\}}} y \qquad (2.0.11)$$

instead $x <_{w_{\{0\}}} y$, i.e.we will be write $x \prec^s_{w_{\{0\}}} y$ iff

$$\neg_s(x \prec_{w_{\{0\}}} y) \Leftrightarrow_s y \prec_{w_{\{0\}}} x. \qquad (2.0.12)$$

We often will be write for short: $x \nprec^s_{w_{\{0\}}} y$ instead $\neg_s\left(x \prec^s_{w_{\{0\}}} y\right)$.

(iii) Note that in general case $\neg_w(x <_{w_{\{0\}}} y) \nvdash y <_{w_{\{0\}}} x$, i.e. in general case the statement

$\neg_w(x <_{w_{\{0\}}} y)$ does not imply provability of the statement $y <_{w_{\{0\}}} x$ and therefore in general case

$$\neg_w(x <_{w_{\{0\}}} y) \not\Leftrightarrow_s y <_{w_{\{0\}}} x. \qquad (2.1.13)$$

We often will be write for short: $x \not<^w_{w_{\{0\}}} y$ instead $\neg_w(x <_{w_{\{0\}}} y)$

(iv) For any $x$ and $y$ such that $\neg_w(x <_{w_{\{0\}}} y) \vdash y <_{w_{\{0\}}} x$ we will be write for short:

$$x \prec^w_{w_{\{0\}}} y \qquad (2.1.14)$$

instead $x <_{w_{\{0\}}} y$, i.e.we will be write $x \prec^w_{w_{\{0\}}} y$ iff

$$\neg_w(x \prec_w y) \Leftrightarrow_s y \prec_w x. \qquad (2.1.15)$$

We often will be write for short: $x \nprec^w_{w_{\{0\}}} y$ instead $\neg_w\left(x \prec^w_{w_{\{0\}}} y\right)$.

**Remark 2.0.5**.(i) Note that in general case $\neg_s(x <_{w_{\{0\}}} y) \nvdash y <_{w_{\{0\}}} x$, i.e. in general

case

$$x <_{w_{\{0\}}} y \not\Leftrightarrow_s \neg_s(y <_{w_{\{0\}}} x). \qquad (2.0.10)$$

We often will be write for short: $x \not<^s_{w_{\{0\}}} y$ instead $\neg_s(x <_{w_{\{0\}}} y)$.

(ii) For any $x$ and $y$ such that $\neg_s(x <_{w_{\{0\}}} y) \vdash y <_{w_{\{0\}}} x$ we will be write for short:

$$x \prec^s_{w_{\{0\}}} y \qquad (2.1.11)$$

**Proposition 2.1.1**. (i) $x <_{w_{\{0\}}} y$ or $y <_{w_{\{0\}}} x$ but not $x <_{w_{\{0\}}} y$ and $y <_{w_{\{0\}}} x$

simultaneously,

(ii) $x <_{w_{\{1\}}} y$ or $y <_{w_{\{1\}}} x$ but not $x <_{w_{\{1\}}} y$ and $y <_{w_{\{1\}}} x$ simultaneously,

(iii) $x <_{w_{\{n\}}} y$ or $y <_{w_{\{n\}}} x$ but not $x <_{w_{\{n\}}} y$ and $y <_{w_{\{n\}}} x$ simultaneously, $n = 1, 2, \ldots$ .

Proof.Immetiately from definitions.

# II.1.2.Basics about paraordered fields.

In this section we will define the notion of paraordered field, which is simply a field in the algebraic sens together with a total order which has a compatible behavior with the operations of the field.

**Definition 2.1.1**.A $w$-consistent field ($w$-field) is a $w$-set $\mathbb{k}_w$ ($w$-set $\mathbb{k}_w$) together with two binary operations $+_w$ (addition), $\times_w$ (product) which satisfy the following axioms:

($\mathrm{F}_1^w$) (i) $\forall x,y(x,y \in_w \mathbb{k}_w \Rightarrow_s x +_w y \in_w \mathbb{k}_w)$, (ii) $\forall x,y(x,y \in_w \mathbb{k}_w \Rightarrow_s x \times_w y \in_w \mathbb{k}_w)$.

($\mathrm{F}_2^w$) (i) $\forall x,y,z[x,y,z \in_w \mathbb{k}_w \Rightarrow_s (x +_w y) +_w z =_w x +_w (y +_w z)]$,

(ii) $\forall x,y,z[x,y,z \in_w \mathbb{k}_w \Rightarrow_s (x \times_w y) \times_w z =_w x \times_w (y \times_w z)]$.

($\mathrm{F}_3^w$) (i) $\forall x,y[x,y \in_w \mathbb{k}_w \Rightarrow_s (x +_w y) =_w (y +_w x)]$,

(ii) $\forall x,y[x,y \in_w \mathbb{k}_w \Rightarrow_s (x \times_w y) =_w (y \times_w x)]$.

($\mathrm{F}_4^w$) There exists a $w$-unique $w$-element $0_w \in_w \mathbb{k}_w$ such that

$\forall x \in_w \mathbb{k}_w[x +_w 0_w] =_w x$.

($\mathrm{F}_5^w$) There exists a $w$-unique $w$-element $1_w \in_w \mathbb{k}_w$ such that

$\forall x \in_w \mathbb{k}_w[x \times_w 1_w] =_w x$.

($\mathrm{F}_6^w$) $\forall x(x \in_w \mathbb{k}_w)\exists y(y \in_w \mathbb{k}_w)[x +_w y =_w 0_w]$.

($\mathrm{F}_7^w$) $\forall x,y,z[x,y,z \in_w \mathbb{k}_w \Rightarrow_s (x +_w y) \times_w z =_w x \times_w z +_w y \times_w z]$.

**Definition 2.1.2** ($w$-ordered $w$-field). An $w$-ordered $w$-field is a $w$-field $\mathbb{k}_w$ such that a binary

predicate $<_w$ is defined on the set $w$-$\mathbb{k}_w$, such that $w$- satisfies the following axioms :

($\Delta_1^w$) $\forall x,y \in_w \mathbb{k}_w$ one and only one of the following holds :

(i) $x <_w y$, (ii) $x =_w y$, (iii) $x \neq_w^w y$, (iv) $x \nless_w^w y$, (v) $y <_w x$, (vi) $y \nless_w^w x$.

($\Delta_2^w$) (i) $x <_w y \Leftrightarrow_s \neg_s(y <_w x)$, (ii) $x <_w y \not\Rightarrow_w y <_w x$, (iii) $x <_w y \not\Rightarrow_s y <_w x$

($\Delta_3^w$) $\forall x,y,z(x,y,z \in_w \mathbb{k}_w)[(x <_w y) \wedge ({}_s y <_w z) \Rightarrow_s x <_w z]$.

($\Delta_4^w$) $\forall x,y,z(x,y,z \in_w \mathbb{k}_w)[x <_w y \Rightarrow_s x +_w z <_w y +_w z]$.

($\Delta_5^w$) $\forall x,y,z(x,y,z \in_w \mathbb{k}_w)[(x <_w y) \wedge (0_w <_w z) \Rightarrow_s x \times_w z <_w y \times_w z]$.

**Designation 2.1.1**.A $w$-field $(\mathbb{k}_w, +_w, \times_w)$ which is an $w$-ordered $w$-field for $<_w$ will be noted $(\mathbb{k}_w, +_w, \times_w, <_w)$.

**Definition 2.1.2**.We say that an element $x \in_w \mathbb{k}_w$ is a $w$-positive element if $x_w > 0_w$. We denote $\mathbb{k}_w^+$ the set of all $w$-positive elements.

**Remark 2.1.1**.

**Definition 2.1.3**.The following function $|\cdot|_w : \mathbb{k}_w \to \mathbb{k}_w^+ \cup_w \{0_w\}_w$ is called $w$-absolute value and can always be defined on any $w$-ordered $w$-field.

$$|x|_w =_w \begin{cases} w\text{-}\max(-_w x, x) & \Leftrightarrow_s x \neq^s_w 0_w \wedge x \neq^w_w 0_w \\ 0_w & \Leftrightarrow_s x =_w 0_w \end{cases} \quad (2.1.12)$$

**Proposition 2.1.1.**

**Definition 2.1.1.**A $w$-field is an $w_{[1]}$-inconsistent set order one $\Bbbk_w$ ($w$-set $\Bbbk_w$) together with
two binary operations $+_w$ (addition), $\times_w$ (product) which satisfy the following axioms:
($F^w_1$) (i) $\forall x,y(x,y \in^s_w \Bbbk_w \Rightarrow_s x +_w y \in^s_w \Bbbk_w)$, (ii) $\forall x,y(x,y \in^s_w \Bbbk_w \Rightarrow_s x \times_w y \in^s_w \Bbbk_w)$
(iii) $\forall x,y(x,y \in_{w_{[1]}} \Bbbk_w \Rightarrow_s x +_w y \in_{w_{[1]}} \Bbbk_w)$, (iv) $\forall x,y(x,y \in_w \Bbbk_{w_{[1]}} \Rightarrow_s x \times_w y \in_{w_{[1]}} \Bbbk_w)$
( $\Bbbk_w$ is closed at least in paraconsistent sense order one under addition and product)
($F^w_2$) (i) $\forall x,y,z[x,y,z \in_w \Bbbk_w \Rightarrow_s (x +_w y) +_w z =_{w_{[1]}} x +_w (y +_w z)]$,
(ii) $\forall x,y,z[x,y,z \in_w \Bbbk_w \Rightarrow_s (x +_w y) +_w z =_{w_{[1]}} x +_w (y +_w z)]$,
(the binary operations are associative in paraconsistent sense order one)
($F^w_3$) (i) $\forall x,y[x,y \in_w \Bbbk_w \Rightarrow_s (x +_w y) =_{w_{[1]}} (y +_w x)]$,
($F^w_3$) (ii) $\forall x,y[x,y \in_w \Bbbk_w \Rightarrow_s (x \times_w y) =_{w_{[1]}} (y \times_w x)]$
(the binary operations are commutative in paraconsistent sense order one)
($F^w_4$) There exists a $w$-unique $w$-element $0_w \in_w \Bbbk_w$ such that
(i) $\forall x \in_w \Bbbk_w$

# II.2.Limits continuity, and the derivative

Any consistent sequence $\{a_n\}_{n\in\mathbb{N}}$ is a map $a : \mathbb{N} \to \mathbb{R}$ and, as such, has an paraconsistent extension to a map ${}^{\#_w}a : {}^{\#_w}\mathbb{N} \to {}^{\#_w}\mathbb{R}$. For any $n \in {}^{\#_w}\mathbb{N}$ we write $a_n^{\#_w} = {}^{\#_w}a(n)$. We use $\{a_n\}_{n\in {}^{\#_w}\mathbb{N}}$ or $\{a_n^{\#_w}\}_{n\in {}^{\#_w}\mathbb{N}}$ to denote the extended paraconsistent $w$-sequence.

For any elements $a, a' \in {}^{\#_w}\mathbb{R}$ we shall write $a \approx_w a'$ to mean that the difference $a -_w a'$ is infinitesimal at least in inconsistent sense.

**PROPOSITION 2.2.1.**(i) $\lim_{n\to\infty} a_n = a$ iff $a_\omega^{\#_{w^*}} \approx_{w^*} a$ for all $\omega \in {}^{\#_{w^*}}\mathbb{N}\backslash\mathbb{N}$.

(ii) $\lim_{n\to\infty} a_n = a \Rightarrow a_\omega^{\#_w} \approx_w a$ for all $\omega \in {}^{\#_w}\mathbb{N}\backslash\mathbb{N}$.

**Remark 2.2.1.** (i) Here the left-hand side of the equivalence in statement (i) has its standard meaning inside $V(\mathbb{R})$. The right-hand side is a statement about the weacly consistent extended universe $V({}^{\#_{w^*}}\mathbb{R})$.

(ii) (i) Here the left-hand side of the implication in statement (ii) has its standard meaning
inside $V(\mathbb{R})$. The right-hand side is a statement about the paraconsistent extended universe $V({}^{\#_w}\mathbb{R})$.

**Proof**. (i) If $\mathrm{llim}_{n\to\infty} a_n = a$ , then given any $\varepsilon > 0$ there is some $n \in \mathbb{N}$ such that the following statement is true in $V(\mathbb{R})$ :

$$\forall m(m \in \mathbb{N})[m \geq n \Rightarrow |a - a_m|] < \varepsilon \qquad (2.1.1)$$

By $w^*$-transfer the statement

$$\forall m\Big(m \in {}^{\#_{w^*}}\mathbb{N}\Big)\Big[m \geq n \Rightarrow \Big|a - a_m^{\#_{w^*}}\Big|\Big] < \varepsilon \qquad (2.1.2)$$

is true in $V({}^{\#_{w^*}}\mathbb{R})$. If $\omega \in {}^{\#_{w^*}}\mathbb{N}\backslash\mathbb{N}$, then $\Big|a - a_\omega^{\#_{w^*}}\Big| < \varepsilon$ is true in $V({}^{\#_{w^*}}\mathbb{R})$. Since this is true for all standard $\varepsilon > 0$, it means that the difference $a_\omega^{\#_{w^*}} -_{w^*} a$, is $w^*$-infinitesimal, i.e., $a_\omega^{\#_{w^*}} \approx_{w^*} a$.

We present bellow two versions of the proof of the converse of the statement (i):

(ii)

**Definition 2.2.1.**Let $(x_n)_{n\in\mathbb{N}}$ be a sequence $(x_n)_{n\in\mathbb{N}} \subset \mathbb{R}$. A $w_{\{0\}}$-hypersequence $\{x_n\}_{n\in_{w_{\{0\}}}\mathbb{N}^{\#_{w_{\{0\}}}}}$ that is a mapping ${}^{\#_{w_{\{0\}}}}[(x_n)_{n\in\mathbb{N}}] : \quad \mathbb{N}^{\#_{w_{\{0\}}}} \to \mathbb{R}^{\#_{w_{\{0\}}}}$.

**Definition 2.2.2**. Let $(x_n)_{n\in\mathbb{N}}$ be a sequence $(x_n)_{n\in\mathbb{N}} \subset \mathbb{R}$. A $w_1$-hypersequence $\{x_n\}_{n\in_{w_1}\mathbb{N}^{\#_{w_1}}}$ that is a mapping ${}^{\#_{w_{\{1\}}}}[(x_n)_{n\in\mathbb{N}}] : \quad \mathbb{N}^{\#_{w_{\{1\}}}} \to \mathbb{R}^{\#_{w_{\{1\}}}}$.

**Definition 2.2.3.**A $w_{\{0\}}$-hypersequence $\{x_n\}_{n\in_{w_{\{0\}}}\mathbb{N}^{\#_{w_{\{0\}}}}}$ is:

(i) $w_{\{0\}}$-increasing (or non-$w_{\{0\}}$-decreasing) if $x_n \leq_{w_{\{0\}}} x_{n+1}$ for all $n \in_{w_{\{0\}}} \mathbb{N}^{\#_{w_{\{0\}}}}$;

(ii) $w_{\{0\}}$-decreasing (or non-$w_{\{0\}}$-increasing) if $x_{n+1} \leq_{w_{\{0\}}} x_n$ for all $n \in_{w_{\{0\}}} \mathbb{N}^{\#_{w_{\{0\}}}}$;

(iii) strictly $w_{\{0\}}$-increasing if $x_n <_{w_{\{0\}}} x_{n+1}$ for all $n \in_{w_{\{0\}}} \mathbb{N}^{\#_{w_{\{0\}}}}$;

(iv) strictly $w_{\{0\}}$-decreasing if $x_{n+1} <_{w_{\{0\}}} x_n$ for all $n \in_{w_{\{0\}}} \mathbb{N}^{\#_{w_{\{0\}}}}$.

**Definition 2.2.4**. A $w_1$-hypersequence $\{x_n\}_{n\in_{w_1}\mathbb{N}^{\#_{w_1}}}$ is:

(i) $w_{\{1\}}$-increasing (or non-$w_1$-decreasing) if $x_n \leq_{w_{\{1\}}} x_{n+1}$ for all $n \in_{w_{\{1\}}} \mathbb{N}^{\#_{w_1}}$;

(ii) $w_{\{1\}}$-decreasing (or non-$w_{\{1\}}$-increasing) if $x_{n+1} \leq_{w_{\{1\}}} x_n$ for all $n \in_{w_{\{1\}}} \mathbb{N}^{\#_{w_{\{1\}}}}$;

(iii) strictly $w_{\{1\}}$-increasing if $x_n <_{w_{\{1\}}} x_{n+1}$ for all $n \in_{w_{\{1\}}} \mathbb{N}^{\#_{w_1}}$;

(iv) strictly $w_{\{1\}}$-decreasing if $x_{n+1} <_{w_{\{1\}}} x_n$ for all $n \in_{w_{\{1\}}} \mathbb{N}^{\#_{w\{1\}}}$.

**Definition 2.1.5**. A $w_{\{0\}}$-hypersequence is:

(i) $w_{\{0\}}$-monotone if it is either $w_{\{0\}}$-increasing or $w_{\{0\}}$-decreasing;

(ii) strictly $w_{\{0\}}$-monotone if it is either strictly $w_{\{0\}}$-increasing or $w_{\{0\}}$-strictly $w_{\{0\}}$-decreasing.

**Definition 2.1.6.**A $w_{\{1\}}$-hypersequence is:

(i) $w_{\{1\}}$-monotone if it is either $w_{\{1\}}$-increasing or $w_{\{1\}}$-decreasing;

(ii) strictly $w_{\{1\}}$-monotone if it is either strictly $w_{\{1\}}$-increasing or strictly $w_{\{1\}}$-decreasing.

**Definition 2.1.7.**We call $x \in_{w_{\{0\}}} \mathbb{R}^{\#_{w_{\{0\}}}}$ the $w_{\{0\}}$-limit of the $w_{\{0\}}$-hypersequence $\{x_n\}_{n\in_{w_{\{0\}}}\mathbb{N}^{\#_{w_{\{0\}}}}}$ if the following condition holds:for each hyperreal number $\varepsilon \in_{w_{\{0\}}} \mathbb{R}^{\#_{w_{\{0\}}}}$, $\varepsilon \;_{w_{\{0\}}} > 0_{w_{\{0\}}}$, there exists a hypernatural number $N \in_{w_{\{0\}}} \mathbb{N}^{\#_{w_{\{0\}}}}$ such that, for every hypernatural number $n \;_{w_{\{0\}}} \geq N$, we have $|x_n -_{w_{\{0\}}} x|_{w_{\{0\}}} <_{w_{\{0\}}} \varepsilon$. The $w_{\{0\}}$-hypersequence $\{x_n\}_{n\in_{w_{\{0\}}}\mathbb{N}^{\#_{w_{\{0\}}}}}$ is said to $w_{\{0\}}$-converge to or tend to the $w_{\{0\}}$-limit $x$, written $x_n \to_{w_{\{0\}}} x$ or $w_{\{0\}}\text{-}\lim_{n\to\infty^{w_{\{0\}}}} x_n =_{w_{\{0\}}} x.$ Symbolically, this is:

$$\forall\varepsilon(\varepsilon_{w_{\{0\}}} > 0_{w_{\{0\}}}) \searrow \Big[\exists N\big(N \in_{w_{\{0\}}} \mathbb{N}^{\#_{w_{\{0\}}}}\big)\forall n\big(n \in_{w_{\{0\}}} \mathbb{N}^{\#_{w_{\{0\}}}}\big)[n_{w_{\{0\}}} \geq N \Rightarrow_s |x_n -_{w_{\{0\}}} x|_{w_{\{0\}}} <_{w_{\{0\}}} \varepsilon]\Big]. \quad (2.1.3)$$

**Remark 2.1.2**. For $w_{\{0\}}$-hypersequences $\{x_n\}_{n\in_{w_{\{0\}}}\mathbb{N}^{\#_{w_{\{0\}}}}} \subset_{w_{\{0\}}} \mathbb{R}^{\#_{w_{\{0\}}}}$ it is also convenient to define the notions $x_n \to_{w_{\{0\}}} \infty^{\#_{w_{\{0\}}}}$ and $x_n \to_{w_{\{0\}}} -_{w_{\{0\}}} \infty^{\#_{w_{\{0\}}}}$ as $n \to_{w_{\{0\}}} \infty^{\#_{w_{\{0\}}}}$.

**Definition 2.1.8.** If $\{x_n\}_{n\in_{w_{\{0\}}}\mathbb{N}^{\#_{w_{\{0\}}}}} \subset_{w_{\{0\}}} \mathbb{R}^{\#_{w_{\{0\}}}}$ then $x_n \to_{w_{\{0\}}} \infty^{\#_{w_{\{0\}}}}$ as $n \to_{w_{\{0\}}} \infty^{\#_{w_{\{0\}}}}$ if for every positive hyperreal number $M \in_{w_{\{0\}}} \mathbb{R}^{\#_{w_{\{0\}}}}$ there exists an hyperinteger $N \in_{w_{\{0\}}} \mathbb{N}^{\#_{w_{\{0\}}}}$ such that $n_{\ w_{\{0\}}} \geq N \Rightarrow_s x_{n\ w_{\{0\}}} \geq M$ ($x_n \leq_{w_{\{0\}}} -_{w_{\{0\}}} M$), we say $\{x_n\}_{n\in_{w_{\{0\}}}\mathbb{N}^{\#_{w_{\{0\}}}}}$ has $w_{\{0\}}$-limit $\infty^{\#_{w_{\{0\}}}}(-_{w_{\{0\}}}\infty^{\#_{w_{\{0\}}}})$ and write $w_{\{0\}}\text{-}\lim_{n\to\infty^{\#_{w_{\{0\}}}}} x_n =_{w_{\{0\}}} \infty^{\#_{w_{\{0\}}}}(-_{w_{\{0\}}}\infty^{\#_{w_{\{0\}}}})$.

**Definition 2.1.9.** We call $x \in_{w_{\{1\}}} \mathbb{R}^{\#_{w_{\{1\}}}}$ the $w_{\{1\}}$-limit of the $w_{\{1\}}$-hypersequence $\{x_n\}_{n\in_{w_{\{1\}}}\mathbb{N}^{\#_{w_{\{1\}}}}}$ if the following condition holds: for each hyperreal number $\varepsilon \in_{w_{\{1\}}} \mathbb{R}^{\#_{w_{\{1\}}}}$, $\varepsilon_{\ w_{\{1\}}} > 0_{w_{\{1\}}}$, there exists a hypernatural number $N \in_{w_{\{1\}}} \mathbb{N}^{\#_{w_{\{1\}}}}$ such that, for every hypernatural number $n_{\ w_{\{1\}}} \geq N$, we have $|x_n -_{w_{\{1\}}} x|_{w_{\{1\}}} <_{w_{\{1\}}} \varepsilon.$ The $w_{\{1\}}$-hypersequence $\{x_n\}_{n\in_{w_{\{1\}}}\mathbb{N}^{\#_{w_{\{1\}}}}}$ is said to $w_{\{1\}}$-converge to or tend to the $w_{\{1\}}$-limit $x$, written $x_n \to_{w_{\{1\}}} x$ or $\lim_{n\to\infty^{w_{\{1\}}}} x_n =_{w_{\{1\}}} x.$ Symbolically, this is:

$$\forall\varepsilon(\varepsilon_{w_{\{1\}}} > 0_{w_{\{1\}}}) \searrow \Big[\exists N\big(N \in_{w_{\{1\}}} \mathbb{N}^{\#_{w_{\{1\}}}}\big)\forall n\big(n \in_{w_{\{1\}}} \mathbb{N}^{\#_{w_{\{1\}}}}\big)[n_{w_{\{1\}}} \geq N \Rightarrow_s |x_n -_{w_{\{1\}}} x|_{w_{\{1\}}} <_{w_{\{1\}}} \varepsilon]\Big].$$

**Remark 2.1.3**. For $w_{\{1\}}$-hypersequences $\{x_n\}_{n\in_{w_{\{1\}}}\mathbb{N}^{\#_{w_{\{1\}}}}} \subset_{w_{\{1\}}} \mathbb{R}^{\#_{w_{\{1\}}}}$ it is also convenient to define the notions $x_n \to_{w_{\{1\}}} \infty^{\#_{w_{\{1\}}}}$ and $x_n \to_{w_{\{1\}}} -_{w_{\{1\}}} \infty^{\#_{w_{\{1\}}}}$ as $n \to_{w_{\{1\}}} \infty^{\#_{w_{\{1\}}}}$.

**Definition 2.1.10**. If $\{x_n\}_{n\in_{w_{\{1\}}}\mathbb{N}^{\#_{w_{\{1\}}}}} \subset_{w_{\{1\}}} \mathbb{R}^{\#_{w_{\{1\}}}}$ then $x_n \to_{w_{\{1\}}} \infty^{\#_{w_{\{1\}}}}\left(-\infty^{\#_{w_{\{1\}}}}\right)$ as $n \to_{w_{\{1\}}} \infty^{\#_{w_{\{1\}}}}$ if for every positive hyperreal number $M \in_{w_{\{1\}}} \mathbb{R}^{\#_{w_{\{1\}}}}$ there exists an hyperinteger $\mathbf{N} \in_{w_{\{1\}}} \mathbb{N}^{\#_{w_{\{1\}}}}$ such that $n_{\ w_{\{1\}}} \geq \mathbf{N} \Rightarrow_s x_{n\ w_{\{1\}}} \geq M$ ($x_n \leq_{w_{\{1\}}} -_{w_{\{1\}}} M$), we say $\{x_n\}_{n\in_{w_{\{1\}}}\mathbb{N}^{\#_{w_{\{1\}}}}}$ has $w_{\{0\}}$-limit $\infty^{\#_{w_{\{1\}}}}(-_{w_{\{1\}}}\infty^{\#_{w_{\{1\}}}})$ and write $w_{\{1\}}\text{-}\lim_{n\to\infty^{\#_{w_{\{1\}}}}} x_n =_{w_{\{1\}}} \infty^{\#_{w_{\{1\}}}}(-_{w_{\{1\}}}\infty^{\#_{w_{\{1\}}}})$.

**Theorem 2.2.1**. (i) Every $w_{\{0\}}$-internal $w_{\{0\}}$-hyperbounded $w_{\{0\}}$-monotone $w_{\{0\}}$-hypersequence in $\mathbb{R}^{\#_{w_{\{0\}}}}$ has a $w_{\{0\}}$-limit in $\mathbb{R}^{\#_{w_{\{0\}}}}$.
(ii) Every $w_{\{0\}}$-external $w_{\{0\}}$-hyperbounded strictly $w_{\{0\}}$-monotone $w_{\{0\}}$-hypersequence in $\mathbb{R}^{\#_{w_{\{0\}}}}$ has a $w_{\{0\}}$-limit in $\mathbb{R}^{\#_{w_{\{0\}}}}$.

Proof: Suppose $\{x_n\}_{n\in_{w_{\{0\}}}\mathbb{N}^{\#_{w_{\{0\}}}}} \subset_{w_{\{0\}}} \mathbb{R}^{\#_{w_{\{0\}}}}$ and $\{x_n\}_{n\in_{w_{\{0\}}}\mathbb{N}^{\#_{w_{\{0\}}}}}$ is $w_{\{0\}}$-increasing (if $\{x_n\}_{n\in_{w_{\{0\}}}\mathbb{N}^{\#_{w_{\{0\}}}}}$ is $w_{\{0\}}$-decreasing, the argument is analogous). Since the set $w_{\{0\}}\text{-}\bigcup_n \{x_n\}_{n\in_{w_{\{0\}}}\mathbb{N}^{\#_{w_{\{0\}}}}}$ of $w_{\{0\}}$-hyperreals is $w_{\{0\}}$-hyperbounded above, it has a least $w_{\{0\}}$-hyperupper bound in $\mathbb{R}^{\#w_{\{0\}}}$, $x$, say. We claim that $x_n \to_{w_{\{0\}}} x$ as $n \to_{w_{\{0\}}} \infty^{\#w_{\{0\}}}$. In order to see this, note that $x_n \leq_{w_{\{0\}}} x$ for all $n \in_{w_{\{0\}}} \mathbb{N}^{\#w_{\{0\}}}$; but if $\varepsilon >_{w_{\{0\}}} 0_{w_{\{0\}}}$ then $x_k \geq_{w_{\{0\}}} x -_{w_{\{0\}}} \varepsilon$ for some $k$, as otherwise $x -_{w_{\{0\}}} \varepsilon$ would be an upper bound. Choose such $k = k(\varepsilon)$. Since $x_k \geq_{w_{\{0\}}} x -_{w_{\{0\}}} \varepsilon$, then $x_n \geq_{w_{\{0\}}} x -_{w_{\{0\}}} \varepsilon$ for all $n \geq_{w_{\{0\}}} k$ as the sequence is increasing. Hence $x -_{w_{\{0\}}} \varepsilon <_{w_{\{0\}}} x_n <_{w_{\{0\}}} x$ for all $n \geq_{w_{\{0\}}} k$. Thus $|x -_{w_{\{0\}}} x_n|_{w_{\{0\}}} <_{w_{\{0\}}} \varepsilon$ for $n \geq_{w_{\{0\}}} k$, and so $x_n \to_{w_{\{0\}}} x$ since $\varepsilon >_{w_{\{0\}}} 0_{w_{\{0\}}}$ is arbitrary.

(ii) Let $\mathbb{R}_{\mathbf{d}}^{\#w_{\{0\}}}$ be Dedekind completion of $\mathbb{R}^{\#w_{\{0\}}}$. Suppose $\{x_n\}_{n \in_{w_{\{0\}}} \mathbb{N}^{\#w_{\{0\}}}} \subset_{w_{\{0\}}} \mathbb{R}^{\#w_{\{0\}}}$ and $\{x_n\}_{n \in_{w_{\{0\}}} \mathbb{N}^{\#w_{\{0\}}}}$ is strictly $w_{\{0\}}$-increasing (if $\{x_n\}_{n \in_{w_{\{0\}}} \mathbb{N}^{\#w_{\{0\}}}}$ is strictly $w_{\{0\}}$-decreasing, the argument is analogous). Since the set $w_{\{0\}}\text{-}\bigcup_n \{x_n\}_{n \in_{w_{\{0\}}} \mathbb{N}^{\#w_{\{0\}}}}$ of $w_{\{0\}}$-hyperreals is $w_{\{0\}}$-hyperbounded above, it has a least $w_{\{0\}}$-hyperupper bound in $\mathbb{R}_{\mathbf{d}}^{\#w_{\{0\}}}$, $x$, say. Assume that $x = {}^{\#}a, a \in \mathbb{R}^{\#w_{\{0\}}}$, where ${}^{\#}a$ is image $a$ in $\mathbb{R}_{\mathbf{d}}^{\#w_{\{0\}}}$, see [46]. In this case argument is the same as above. Assume now that $x \neq {}^{\#}a, \forall a\left[a \in \mathbb{R}^{\#w_{\{0\}}}\right]$, i.e. $x$ is absorbtion number in $\mathbb{R}_{\mathbf{d}}^{\#w_{\{0\}}}$. We claim that again $x_n \to_{w_{\{0\}}} x$ as $n \to_{w_{\{0\}}} \infty^{\#w_{\{0\}}}$. In order to see this, note that $x_n \leq_{w_{\{0\}}} x$ for all $n \in_{w_{\{0\}}} \mathbb{N}^{\#w_{\{0\}}}$; but if $\varepsilon >_{w_{\{0\}}} 0_{w_{\{0\}}}$ then $x_k \geq_{w_{\{0\}}} x -_{w_{\{0\}}} \varepsilon$ for some $k$, as otherwise $x -_{w_{\{0\}}} \varepsilon$ would be an upper bound. Choose such

**Theorem 2.1.2**.

**Theorem 2.2.3**.(**Comparison Test**)

(i) If $0_{w_{\{0\}}} \leq_{w_{\{0\}}} x_n \leq_{w_{\{0\}}} y_n$ for all $n \geq_{w_{\{0\}}} N \in_{w_{\{0\}}} \mathbb{N}^{\#w_{\{0\}}}$, and $y_n \to_{w_{\{0\}}} 0$ as $n \to_{w_{\{0\}}} \infty^{\#w_{\{0\}}}$, then $x_n \to_{w_{\{0\}}} 0$ as $n \to_{w_{\{0\}}} \infty^{\#w_{\{0\}}}$.

(ii) If $x_n \leq_{w_{\{0\}}} y_n$ for all $n \geq_{w_{\{0\}}} N \in_{w_{\{0\}}} \mathbb{N}^{\#w_{\{0\}}}$, $x_n \to_{w_{\{0\}}} x$ as $n \to_{w_{\{0\}}} \infty^{\#w_{\{0\}}}$ and $y_n \to_{w_{\{0\}}} y$ as $n \to_{w_{\{0\}}} \infty^{\#w_{\{0\}}}$, then $x \leq_{w_{\{0\}}} y$.

(iii) In particular, if $x_n \leq_{w_{\{0\}}} a$ for all $n \geq_{w_{\{0\}}} N \in_{w_{\{0\}}} \mathbb{N}^{\#w_{\{0\}}}$ and $x_n \to_{w_{\{0\}}} x$ as $n \to_{w_{\{0\}}} \infty^{\#w_{\{0\}}}$, then $x \leq_{w_{\{0\}}} a$.

.

# II.3.Cauchy $w_{\{0\}}$-hypersequences and $w_{\{1\}}$-hypersequences.

# II.3.1.Cauchy $w_{\{0\}}$-hypersequences.

**Defnition 2.3.1**. A $w_{\{0\}}$-metric space $(X; d_{w_{\{0\}}})$ is a set $X$ together with a distance function $d_{w_{\{0\}}} : X \times_{w_{\{0\}}} X \to \mathbb{R}^{\#w_{\{0\}}}$ such that for all $x, y, z \in_{w_{\{0\}}} X$ the following hold:

1.$d_{w_{\{0\}}}(x,y)\ {}_{w_{\{0\}}}\geq 0_{w_{\{0\}}}, d_{w_{\{0\}}}(x,y) =_{w_{\{0\}}} 0_{w_{\{0\}}} \Leftrightarrow_s x =_{w_{\{0\}}} y$ (positivity),
2. $d_{w_{\{0\}}}(x,y) =_{w_{\{0\}}} d_{w_{\{0\}}}(y,x)$ (symmetry),
3. $d_{w_{\{0\}}}(x,y) \leq_{w_{\{0\}}} d_{w_{\{0\}}}(x,z) +_{w_{\{0\}}} d_{w_{\{0\}}}(z,y)$ (triangle inequality).
We denote the corresponding metric space by $(X;d_{w_{\{0\}}})$, to indicate that
a metric space is determined by both the set $X$ and the metric $d_{w_{\{0\}}}$.
**Definition 2.3.2**. Let $\{x_n\}_{w_{\{0\}}} \triangleq \{x_n\}_{n\in_{w_{\{0\}}} \mathbb{N}^{\#_{w_{\{0\}}}}} \subset_{w_{\{0\}}} X$ where $(X;d_{w_{\{0\}}})$ is a $w_{\{0\}}$-metric space. Then $\{x_n\}_{w_{\{0\}}}$ is a Cauchy $w_{\{0\}}$-hypersequence if for every $\varepsilon\ {}_{w_{\{0\}}}\geq 0_{w_{\{0\}}}$ there exists an hyperinteger $\mathbf{N} \in_{w_{\{0\}}} \mathbb{N}^{\#_{w_{\{0\}}}}$ such that

$$m,n\ {}_{w_{\{0\}}}\geq \mathbf{N} \Rightarrow_s d_{w_{\{0\}}}(x_m,x_n) <_{w_{\{1\}}} \varepsilon \qquad (2.3.1)$$

We sometimes write this as $d_{w_{\{0\}}}(x_m,x_n) \to_{w_{\{0\}}} 0_{w_{\{0\}}}$ as $m,n \to_{w_{\{0\}}} \infty^{\#_{w_{\{0\}}}}$.

.

# II.3.2.Cauchy $w_{\{1\}}$-hypersequences.

**Defnition 2.3.1**. A $w_{\{1\}}$-metric space $(X;d_{w_{\{1\}}})$ is a set $X$ together with a distance function $d_{w_{\{1\}}} : X \times_{w_{\{1\}}} X \to \mathbb{R}^{\#_{w_{\{1\}}}}$ such that for all $x,y,z \in_{w_{\{1\}}} X$ the following hold:
1.$d_{w_{\{1\}}}(x,y)\ {}_{w_{\{1\}}}\geq 0_{w_{\{1\}}}, d_{w_{\{1\}}}(x,y) =_{w_{\{1\}}} 0_{w_{\{1\}}} \Leftrightarrow_s x =_{w_{\{1\}}} y$ ($w_{\{1\}}$-positivity),
2. $d_{w_{\{1\}}}(x,y) =_{w_{\{1\}}} d_{w_{\{1\}}}(y,x)$ ($w_{\{1\}}$-symmetry),
3. $d_{w_{\{1\}}}(x,y) \leq_{w_{\{1\}}} d_{w_{\{1\}}}(x,z) +_{w_{\{1\}}} d_{w_{\{1\}}}(z,y)$ (triangle $w_{\{1\}}$-inequality).

We denote the corresponding metric space by $(X;d_{w_{\{1\}}})$, to indicate that
a metric space is determined by both the set $X$ and the metric $d_{w_{\{1\}}}$.
**Definition 2.3.2**. Let $\{x_n\}_{w_{\{1\}}} \triangleq \{x_n\}_{n\in_{w_{\{1\}}} \mathbb{N}^{\#_{w_{\{1\}}}}} \subset_{w_{\{1\}}} X$ where $(X;d_{w_{\{1\}}})$ is a $w_{\{1\}}$-metric
space. Then $\{x_n\}_{w_{\{1\}}}$ is a Cauchy $w_{\{0\}}$-hypersequence if for every $\varepsilon\ {}_{w_{\{0\}}}\geq 0_{w_{\{0\}}}$ there exists
an hyperinteger $\mathbf{N} \in_{w_{\{0\}}} \mathbb{N}^{\#_{w_{\{0\}}}}$ such that

$$m,n\ {}_{w_{\{0\}}}\geq \mathbf{N} \Rightarrow_s d_{w_{\{0\}}}(x_m,x_n) <_{w_{\{1\}}} \varepsilon \qquad (2.3.1)$$

We sometimes write this as $d_{w_{\{0\}}}(x_m,x_n) \to_{w_{\{0\}}} 0_{w_{\{0\}}}$ as $m,n \to_{w_{\{0\}}} \infty^{\#_{w_{\{0\}}}}$.

.

# II.4.$w_{\{0\}}$-Limits and $w_{\{1\}}$-Limits of Functions

.

## II.4.1.$w_{\{0\}}$-Limits of Functions

**Definition 2.4.1**.Let $f : A \to Y, A \subset_{w_{\{0\}}} X,$ where $(X;d_{w_{\{0\}}})$ and $Y$ is a $w_{\{0\}}$-metric spaces,
and let $a$ be a $w_{\{0\}}$-limit point of $A$. Suppose

$$\left[\left(\{x_n\}_{w_{\{0\}}} \subset_{w_{\{0\}}} A\backslash_{w_{\{0\}}}\{a\}_{w_{\{0\}}}\right) \wedge (x_n \to_{w_{\{0\}}} a)\right] \Rightarrow_s f(x_n) \to_{w_{\{0\}}} b. \qquad (2.4.1)$$

Then we say $f$ has $w_{\{0\}}$-limit $b$ at $a$ and write

$$w_{\{0\}}\text{-}\lim_{x_n \to_{w_{\{0\}}} a,\, x \in_{w_{\{0\}}} A} f(x) =_{w_{\{0\}}} b \qquad (2.4.2)$$

or

$$w_{\{0\}}\text{-}\lim_{x_n \to_{w_{\{0\}}} a} f(x) =_{w_{\{0\}}} b, \qquad (2.4.3)$$

where in the last notation the intended domain $A$ is understood from the context.

**Definition 2.4.2**.

.

# II.4.2.$w_{\{1\}}$-Limits of $w_{\{1\}}$-Functions

# II.5.$w_{\{0\}}$-Continuity at a point

**Definition 2.5.1**.Let $f : A \to Y, A \subset_{w_{\{0\}}} X,$ where $(X; d_{w_{\{0\}}})$ and $Y$ is a $w_{\{0\}}$-metric spaces,

and let $a \in_{w_{\{0\}}} A$. Then $f$ is $w_{\{0\}}$-continuous at $a$ if a is an $w_{\{0\}}$-isolated point of $A$, or if $a$ is a $w_{\{0\}}$-limit point of $A$ and

$$w_{\{0\}}\text{-}\lim_{x_n \to_{w_{\{0\}}} a,\, x \in_{w_{\{0\}}} A} f(x) =_{w_{\{0\}}} f(a). \qquad (2.5.1)$$

**Definition 2.5.2**.If $f$ is $w_{\{0\}}$-continuous at every $a \in_{w_{\{0\}}} A$ then we say $f$ is $w_{\{0\}}$-continuous.

The set of all such $w_{\{0\}}$-continuous $w_{\{0\}}$-functions is denoted by $C_{w_{\{0\}}}(A, Y)$.

**Example 2.5.1**. Define

$$f(x) =_{w_{\{0\}}} \begin{cases} x \times \left(w_{\{0\}}\text{-}\sin x^{-_{w_{\{0\}}} 1_{w_{\{0\}}}}\right) & \text{if } \neg_s (x =_{w_{\{0\}}} 0_{w_{\{0\}}}) \\ & \text{if } x =_{w_{\{0\}}} 0_{w_{\{0\}}} \end{cases} \qquad (2.5.2)$$

$f$ is $w_{\{0\}}$-continuous everywhere on $\mathbb{R}^{\#_{w_{\{0\}}}}$.

.

# II.6.Uniform $w_{\{0\}}$-convergence of functions

## II.6.1.Uniform $w_{\{0\}}$-convergence of functions

**Definition 2.6.1**. Let $f, f_n : S \to Y$ for every $n \in_{w_{\{0\}}} \mathbb{N}^{\#_{w_{\{0\}}}}$, where $S$ is any set and $(Y; d_{w_{\{0\}}})$ is a $w_{\{0\}}$-metric space. If $f_n(x) \to_{w_{\{0\}}} f(x)$ for all $x \in_{w_{\{0\}}} S$ then $f_n(x) \to_{w_{\{0\}}} f(x)$ pointwise on $S$.

**Definition 2.6.2**. Let $f, f_n : S \to Y$ for every $n \in_{w_{\{0\}}} \mathbb{N}^{\#_{w_{\{0\}}}}$, where $S$ is any set and $(Y; d_{w_{\{0\}}})$ is a $w_{\{0\}}$-metric space.If for every $\varepsilon\ {}_{w_{\{0\}}} > 0_{w_{\{0\}}}$ there exists $N \in_{w_{\{0\}}} \mathbb{N}^{\#_{w_{\{0\}}}}$ such that $n\ {}_{w_{\{0\}}} \geq N \Rightarrow_s d_{w_{\{0\}}}(f_n(x), f(x)) <_{w_{\{0\}}} \varepsilon$ for all $x \in_{w_{\{0\}}} S$ then we say $f_n(x) \to_{w_{\{0\}}} f(x)$ uniformly on $S$ and write $f_n(x) \rightrightarrows_{w_{\{0\}}} f(x)$.

.

# II.6.2.Uniform $w_{\{1\}}$-convergence of functions

# II.7.Uniform $w_{\{0\}}$-convergence and $w_{\{0\}}$-continuity

.

# II.8.The $w_{\{0\}}$-derivative of $w_{\{0\}}$-internal function of one variable

**Definition 2.8.1**.For each $w_{\{0\}}$-function $f : \mathbb{R}^{\#_{w_{\{0\}}}} \to \mathbb{R}^{\#_{w_{\{0\}}}}$, we define its $w_{\{0\}}$-derived function $f'^{w_{\{0\}}} : \mathbb{R}^{\#_{w_{\{0\}}}} \to \mathbb{R}^{\#_{w_{\{0\}}}}$ by setting, for every point $p \in_{w_{\{0\}}} \mathbb{R}^{\#_{w_{\{0\}}}}$

$$f'^{w_{\{0\}}}(p) = w_{\{0\}}\text{-}\lim_{x \to_{w_{\{0\}}} p} \frac{f(x) -_{w_{\{0\}}} f(p)}{x -_{w_{\{0\}}} p} \tag{2.8.1}$$

if this $w_{\{0\}}$-limit exists.If the $w_{\{0\}}$-limit in (2.8.2) exists, we call it the $w_{\{0\}}$-derivative of $f$ at $p$. If, in addition, this $w_{\{0\}}$-limit is $w_{\{0\}}$-finite or $w_{\{0\}}$-hyper finite, we say that $f$ is $w_{\{0\}}$-differentiable at $p$. If this holds for each $p \in_{w_{\{0\}}} B \subseteq_{w_{\{0\}}} \mathbb{R}^{\#_{w_{\{0\}}}}$, we say that $f$ has a $w_{\{0\}}$-derivative (respectively, is $w_{\{0\}}$-differentiable) on $B$, and we call the function $f'^{w_{\{0\}}}$ the $w_{\{0\}}$-derivative of $f$ on $B$. If the limit in (2.8.2) is one sided (with $x \to_{w_{\{0\}}} p_-$ or $x \to_{w_{\{0\}}} p_+$), we call it a one-sided (left or right) $w_{\{0\}}$-derivative at $p$, denoted $f_-'^{w_{\{0\}}}$ or $f_+'^{w_{\{0\}}}$.

**Definition 2.8.2**.Given any $w_{\{0\}}$-internal function $f_{int} \triangleq f^{\#_{w_{\{0\}}}} : \mathbb{R}^{\#_{w_{\{0\}}}} \to \mathbb{R}^{\#_{w_{\{0\}}}}$, where $f : \mathbb{R} \to \mathbb{R}$, we define its $n$-th $w_{\{0\}}$-derived function (or $w_{\{0\}}$-derived function of order $n \in_{w_{\{0\}}} \mathbb{N}^{\#_{w_{\{0\}}}}$),denoted $f^{(n)_{w_{\{0\}}}} : \mathbb{R}^{\#_{w_{\{0\}}}} \to \mathbb{R}^{\#_{w_{\{0\}}}}$, by transfer:

$$f_{int}^{(0)_{w_{\{0\}}}} =_{w_{\{0\}}} f^{\#_{w_{\{0\}}}}, f_{int}^{(1)_{w_{\{0\}}}} =_{w_{\{0\}}} \left(f^{(1)}\right)^{\#_{w_{\{0\}}}}, \ldots, f_{int}^{(n)_{w_{\{0\}}}} =_{w_{\{0\}}} \left(f^{(n)}\right)^{\#_{w_{\{0\}}}}. \tag{2.8.2}$$

**Definition 2.8.3**.We say that $f_{int}$ has $n$ $w_{\{0\}}$-derivatives at a point $p$ iff the $w_{\{0\}}$-limits

# II.9.The $w_{\{0\}}$-integral.

## II.9.1.The $w_{\{0\}}$-Internal $w_{\{0\}}$-integral.

In this section we deal with $w_{\{0\}}$-internal function $f_{int} \triangleq f^{\#_{w_{\{0\}}}}$ where $f:\mathbb{R}\to\mathbb{R}$, defined on a $w_{\{0\}}$-finite interval $[a,b]^{\#_{w_{\{0\}}}}$. A $w_{\{0\}}$-hyper finite internal partition of $[a,b]^{\#_{w_{\{0\}}}}$ is a $w_{\{0\}}$-hyper finite set of $w_{\{0\}}$-subintervals

$$\left[x_{0_{w_{\{0\}}}}, x_{1_{w_{\{0\}}}}\right]^{\#_{w_{\{0\}}}}, \left[x_{1_{w_{\{0\}}}}, x_{2_{w_{\{0\}}}}\right]^{\#_{w_{\{0\}}}}, \ldots, \left[x_{N-_{w_{\{0\}}}1_{w_{\{0\}}}}, x_N\right]^{\#_{w_{\{0\}}}}, \qquad (2.9.1)$$

where $N \in_{w_{\{0\}}} \mathbb{N}^{\#_{w_{\{0\}}}}$ and where $\{x_n\}_{n=_{w_{\{0\}}}0_{w_{\{0\}}}}^{N}$ is an $w_{\{0\}}$-internal $w_{\{0\}}$-sequence such that

$$a =_{w_{\{0\}}} x_{0_{w_{\{0\}}}} <_{w_{\{0\}}} x_{1_{w_{\{0\}}}} <_{w_{\{0\}}} x_{2_{w_{\{0\}}}} <_{w_{\{0\}}} \ldots <_{w_{\{0\}}} x_N =_{w_{\{0\}}} b \qquad (2.9.2)$$

Thus, any $w_{\{0\}}$-internal set of $N +_{w_{\{0\}}} 1_{w_{\{0\}}}$ points satisfying (2.9.2) defines an $w_{\{0\}}$-internal partition $P_{w_{\{0\}}}(N)$ of $[a,b]^{\#_{w_{\{0\}}}}$, which we denote by

$$P_{w_{\{0\}}}(N) =_{w_{\{0\}}} \left\{x_{0_{w_{\{0\}}}}, x_{1_{w_{\{0\}}}}, \ldots, x_N\right\}_{w_{\{0\}}}. \qquad (2.9.3)$$

The points $x_{0_{w_{\{0\}}}}, x_{1_{w_{\{0\}}}}, \ldots, x_N$ are the partition points of $P_{w_{\{0\}}} \triangleq P_{w_{\{0\}}}(N)$. The largest of the lengths of the $w_{\{0\}}$-subintervals (3.1.1) is the norm of $P_{w_{\{0\}}}$, written as $\|P_{w_{\{0\}}}\|$; thus,

$$\|P_{w_{\{0\}}}\| =_{w_{\{0\}}} \max_{1<_{w_{\{0\}}} i <_{w_{\{0\}}} N} \left(x_i -_{w_{\{0\}}} x_{i-_{w_{\{0\}}}1_{w_{\{0\}}}}\right), \qquad (2.9.4)$$

where RHS of (2.9.4) is defined by $w_{\{0\}}$-transfer.

If $P_{w_{\{0\}}}$ and $P'_{w_{\{0\}}}$ are partitions of $[a,b]^{\#_{w_{\{0\}}}}$, then $P'_{w_{\{0\}}}$ is a refinement of $P_{w_{\{0\}}}$ if every partition point of $P_{w_{\{0\}}}$ is also a partition point of $P'_{w_{\{0\}}}$; that is, if $P'_{w_{\{0\}}}$ is obtained by inserting additional points between those of $P_{w_{\{0\}}}$.

**Definition 2.9.1.**If $w_{\{0\}}$-internal $f$ is defined on $[a,b]^{\#_{w_{\{0\}}}}$, then $w_{\{0\}}$-internal $w_{\{0\}}$-hyper finite $w_{\{0\}}$-sum

$$\sigma_{w_{\{0\}}}^{N} = w_{\{0\}}\text{-}\sum_{j=_{w_{\{0\}}}1_{w_{\{0\}}}}^{N} f(c_j) \times_{w_{\{0\}}} \left(x_j -_{w_{\{0\}}} x_{j-_{w_{\{0\}}}1_{w_{\{0\}}}}\right), \qquad (2.9.5)$$

where $x_{j-_{w_{\{0\}}}1_{w_{\{0\}}}} <_{w_{\{0\}}} c_j <_{w_{\{0\}}} x_j, 1_{w_{\{0\}}} <_{w_{\{0\}}} j <_{w_{\{0\}}} N$ and where $\{c_j\}_{j=_{w_{\{0\}}}1_{w_{\{0\}}}}^{N}$ is any $w_{\{0\}}$-internal $w_{\{0\}}$-hyper finite sequence, is a Riemann $w_{\{0\}}$-hyper finite $w_{\{0\}}$-sum of $f$ over
the partition $P_{w_{\{0\}}}(N)$. We will say more simply that $\sigma_{w_{\{0\}}}^{N}$ is a Riemann $w_{\{0\}}$-hyper finite $w_{\{0\}}$-sum of $f$ over $[a,b]^{\#_{w_{\{0\}}}}$.

## II.9.2.The $w_{\{1\}}$-External $w_{\{1\}}$-integral.

In this section we deal with $w_{\{0\}}$-internal function $f_{int} \triangleq f^{\#_{w_{\{0\}}}}$ where $f:\mathbb{R}\to\mathbb{R}$, defined on a $w_{\{0\}}$-finite interval $[a,b]^{\#_{w_{\{0\}}}}$. A $w_{\{0\}}$-hyper finite internal partition of $[a,b]^{\#_{w_{\{0\}}}}$ is a $w_{\{0\}}$-hyper finite set of $w_{\{0\}}$-subintervals

$$\left[x_{0_{w_{\{0\}}}}, x_{1_{w_{\{0\}}}}\right]^{\#_{w_{\{0\}}}}, \left[x_{1_{w_{\{0\}}}}, x_{2_{w_{\{0\}}}}\right]^{\#_{w_{\{0\}}}}, \ldots, \left[x_{N-_{w_{\{0\}}}1_{w_{\{0\}}}}, x_N\right]^{\#_{w_{\{0\}}}}, \qquad (2.9.1)$$

.

# Chapter III.Analysis on ${}^{\#_{w_{\{1\}}}}\mathbb{R}_{w_{\{1\}}}$.

## II.8.Internal and external series in $\mathbb{R}^{\#_{w_{\{n\}}}}$ of $w_{\{n\}}$-hyperreals.

## II.8.1.Internal series in $\mathbb{R}^{\#_{w_{\{0\}}}}$ of $w_{\{0\}}$-hyperreals.

**Definition 2.8.1**. We call $w_{\{0\}}$-series in $\mathbb{F}^{\#_{w_{\{0\}}}}$ (which means $\mathbb{R}^{\#_{w_{\{0\}}}}$ or $\mathbb{C}^{\#_{w_{\{0\}}}}$) any pair $(f^{\#_{w_{\{0\}}}}, g^{\#_{w_{\{0\}}}})$ of $w_{\{0\}}$-hypersequences, where $f^{\#_{w_{\{0\}}}} : \mathbb{N}^{\#_{w_{\{0\}}}} \to \mathbb{F}^{\#_{w_{\{0\}}}}$ is $w_{\{0\}}$-internal mapping which defines the general terms of the $w_{\{0\}}$-series, also noted $x_n =_{\#_{w_{\{0\}}}} f^{\#_{w_{\{0\}}}}(n)$ and

$g : \mathbb{N}^{\#_{w_{\{0\}}}} \to \mathbb{F}^{\#_{w_{\{0\}}}}$ represents the sequence of partial $w_{\{0\}}$-internal $w_{\{0\}}$-hyperfinite sums $s_n$, i.e.

$$s_n =_{\#_{w_{\{0\}}}} w_{\{0\}}\text{-}\sum_{i=0}^{n} x_n =_{\#_{w_{\{0\}}}} g^{\#_{w_{\{0\}}}}(n). \tag{2.8.1}$$

Instead of $(f^{\#_{w_{\{0\}}}}, g^{\#_{w_{\{0\}}}})$, the $w_{\{0\}}$-series in $\mathbb{F}^{\#_{w_{\{0\}}}}$ is frequently marked as an $w_{\{0\}}$-hyperinfinite sum

$$w_{\{0\}}\text{-}\sum_{n=0}^{\infty^{\#_{w_{\{0\}}}}} x_n \tag{2.8.2}$$

**Definition 2.8.2**. We say that the $w_{\{0\}}$-series $(f^{\#_{w_{\{0\}}}}, g^{\#_{w_{\{0\}}}})$ is $w_{\{0\}}$-convergent to $s \in_{w_{\{0\}}} \mathbb{R}^{\#_{w_{\{0\}}}}$, respectively $s$ is the sum of the $w_{\{0\}}$-series, iff the $w_{\{0\}}$-hypersequence $\{s_n\}_{n \in_{w_{\{0\}}} \mathbb{N}^{\#_{w_{\{0\}}}}}$, of partial sums, $w_{\{0\}}$-converges to $s$, and we note

$$w_{\{0\}}\text{-}\sum_{n=0}^{\infty^{\#_{w_{\{0\}}}}} x_n =_{w_{\{0\}}} w_{\{0\}}\text{-}\lim_{n \to \infty^{\#_{w_{\{0\}}}}} s_n =_{w_{\{0\}}} s. \tag{2.8.3}$$

**Theorem 2.8.1**.(The general Cauchy's criterion). The $w_{\{0\}}$-series (2.8.2) is $w_{\{0\}}$-convergent iff for any $\varepsilon >_{w_{\{0\}}} 0_{w_{\{0\}}}$ we can find $n_0(\varepsilon) \in_{w_{\{0\}}} \mathbb{N}^{\#_{w_{\{0\}}}}$ such that

$$\left| w_{\{0\}}\text{-}\sum_{n +_{w_{\{0\}}} 1_{w_{\{0\}}}}^{n +_{w_{\{0\}}} p} x_n \right| <_{w_{\{0\}}} \varepsilon \tag{2.8.4}$$

holds for all $n >_{w_{\{0\}}} n_0(\varepsilon)$ and arbitrary $p \in_{w_{\{0\}}} \mathbb{N}^{\#_{w_{\{0\}}}}$.

Proof.The assertion of the theorem reformulates in terms of $\varepsilon$ and $n_0(\varepsilon)$ the fact that a $w_{\{0\}}$-series (2.8.2) is $w_{\{0\}}$-convergent iff the hypersequence $\{s_n\}_{n \in_{w_{\{0\}}} \mathbb{N}^{\#_{w_{\{0\}}}}}$ of partial sums is $w_{\{0\}}$-fundamental. This is valid in both $\mathbb{R}^{\#_{w_{\{0\}}}}$ and $\mathbb{C}^{\#_{w_{\{0\}}}}$.

**Theorem 2.8.2**. If $w_{\{0\}}$-series (2.8.2) is convergent, then $x_n \to_{w_{\{0\}}} 0^{\#_{w_{\{0\}}}}$ as $n \to_{w_{\{0\}}} \infty^{\#_{w_{\{0\}}}}$.

Proof. Take $p = 1$ in the (2.8.4) above.

**Example 2.8.1.** In order to get the complete answer about the $w_{\{0\}}$-convergence of the $w_{\{0\}}$- hyperinfinite geometric series $w_{\{0\}}\text{-}\sum_{n=0}^{\infty^{\#w_{\{0\}}}} z^n$ we consider two cases:

(i) If $|z|_{w_{\{0\}}} <_{w_{\{0\}}} 1^{\#w_{\{0\}}}$, then $z_n \to_{w_{\{0\}}} 0^{\#w_{\{0\}}}$, and consequently $s_n \to_{w_{\{0\}}} s$, where

$$s_n =_{w_{\{0\}}} w_{\{0\}}\text{-}\sum_{i=_{w_{\{0\}}} 0^{\#w_{\{0\}}}}^{i=_{w_{\{0\}}} n-_{w_{\{0\}}} 1^{\#w_{\{0\}}}} z^n =_{w_{\{0\}}} \frac{1-_{w_{\{0\}}} z^n}{1-_{w_{\{0\}}} z} \to_{w_{\{0\}}} \frac{1}{1-_{w_{\{0\}}} z} =_{w_{\{0\}}} s. \qquad (2.8.5)$$

(ii) If $|z|_{w_{\{0\}}} \;_{w_{\{0\}}}\!\geq 1^{\#w_{\{0\}}}$, then the series is $w_{\{0\}}$-divergent because the general term is not

tending to zero $0^{\#w_{\{0\}}}$ (as the above Theorem 2.8.2 states).

**Theorem 2.8.3.** (The $1^{st}$ criterion of comparison) Let

## II.8.2. $w_{\{1\}}$-Internal and $w_{\{1\}}$-external series in $\mathbb{R}^{\#w_{\{1\}}}$ of $w_{\{1\}}$-hyperreals.

## II.9. $w_{\{n\}}$-Internal and $w_{\{n\}}$-External series of $w_{\{n\}}$-functions

## II.9.1. $w_{\{0\}}$-Internal series of $w_{\{0\}}$-functions

**Definition 2.9.1.** Let $D \subset_{w_{\{0\}}} \mathbb{R}^{\#w_{\{0\}}}$ be a fixed domain, and let $\mathcal{F}\left(D,\mathbb{R}^{\#w_{\{0\}}}\right) \triangleq \left(\mathbb{R}^{\#w_{\{0\}}}\right)^D$ be the set of all $w_{\{0\}}$-functions $f : D \to \mathbb{R}^{\#w_{\{0\}}}$. Any $w_{\{0\}}$-function $F : \mathbb{N}^{\#w_{\{0\}}} \to \mathcal{F}\left(D,\mathbb{R}^{\#w_{\{0\}}}\right)$ is called $w_{\{0\}}$-hypersequence of ($w_{\{0\}}$-hyperreal) $w_{\{0\}}$-functions. Most frequently it is marked by mentioning the terms $(f_n)_{n\in_{w_{\{0\}}}\mathbb{N}^{\#w_{\{0\}}}}$ or $(f_n)$, where $f_n =_{w_{\{0\}}} F(n)$, and $n$ is an arbitrary $w_{\{0\}}$-hypernatural number.

**Definition 2.9.2.** We say that a number $x \in_{w_{\{0\}}} D$ is a point of $w_{\{0\}}$-convergence of $(f_n)$ if the numerical sequence $(f_n(x))_{n\in_{w_{\{0\}}}\mathbb{N}^{\#w_{\{0\}}}}$ is $w_{\{0\}}$-convergent. The set of all such points forms the set (or domain) of $w_{\{0\}}$-convergence, denoted $w_{\{0\}}\text{-}D_c$. The resulting function, say $\varphi : w_{\{0\}}\text{-}D_c \to \mathbb{R}^{\#w_{\{0\}}}$, expressed at any $x \in_{w_{\{0\}}} D_c$ by

$$\varphi(x) =_{w_{\{0\}}} w_{\{0\}}\text{-}\lim_{n\to\infty^{\#w_{\{0\}}}} f_n(x) \qquad (2.9.1)$$

is called $w_{\{0\}}$-limit of the given $w_{\{0\}}$-sequences of $w_{\{0\}}$-functions. Alternatively we say that $\varphi$ is the (point-wise) $w_{\{0\}}$-limit of $(f_n)$, $(f_n)$ $p$-tends to $\varphi$ and we abreviate

$$\varphi \overset{p}{=}_{w_{\{0\}}} w_{\{0\}}\text{-}\lim_{n\to\infty^{\#w_{\{0\}}}} f_n. \qquad (2.9.2)$$

**Remark 2.9.1**. The notions of series of functions, partial sums, infinite sum, domain of convergence, etc., are similarly defined in $\mathcal{F}\left(D,\mathbb{R}^{\#_{w_{\{0\}}}}\right)$.

**Definition 2.9.3**.A functional $w_{\{0\}}$-series is a series

$$w_{\{0\}}\text{-}\sum_{n=0}^{\infty^{\#_{w_{\{0\}}}}} u_n(x) \tag{2.9.3}$$

where each term of the series $u_n(x)$ is a $w_{\{0\}}$-function on an interval $I$.

We can also de ne pointwise $w_{\{0\}}$-convergence for functional $w_{\{0\}}$-series:

**Definition 2.9.4**.The functional $w_{\{0\}}$-series (2.9.3) is pointwise $w_{\{0\}}$-convergent for each $x \in_{w_{\{0\}}} I$ if the $w_{\{0\}}$-limit

$$w_{\{0\}}\text{-}\sum_{n=0}^{\infty^{\#_{w_{\{0\}}}}} u_n(x) =_{w_{\{0\}}} w_{\{0\}}\text{-}\lim_{N\to\infty^{\#_{w_{\{0\}}}}} \sum_{n=0}^{N} u_n(x) \tag{2.9.4}$$

exists for each $x \in_{w_{\{0\}}} I$.

,

# II.2.THE INTEGRAL.

Let $f : I \to \mathbb{R}_+$ be a positive continuous function, where $I$ is some interval in $\mathbb{R}$. Let $[a,b] \subseteq I$ and let $\Delta x$ be a positive real. The Riemann sum is defined as

$$\sum_a^b f(x)\Delta x = \sum_{i=0}^{n-1} f(x_i)\Delta x + f(x_n)(b - n\Delta x), \qquad 2.2.1$$

where $n$ is the largest integer such that $a + n\Delta x$ and where

$$x_0 = a, x_1 = a + \Delta x, \ldots, x_n, = a + n\Delta x. \qquad 2.2.2$$

**Remark 2.2.1**.Note that it may happen that $n\Delta x < b < (n+1)\Delta x$. [Since $f$ is positive and continuous we have formed the Riemann sum as the sum of the rectangles over each subinterval with height equal to the value of $f(x)$ at the left end of the base of the rectangle.]The Riemann sum (2.2.1) for fixed $a,b$ is a function of $\Delta x$. By extension and transfer this function is also defined for positive $w^*$-infinitesimals $d^{\#_{w^*}}x$. We get a corresponding hyperfinite sum

$$\sum_a^b f(x)d^{\#_{w^*}}x, \tag{2.2.3}$$

where the number $n \in {}^{\#_{w^*}}\mathbb{N}\backslash\mathbb{N}$ in Eq.(2.2.1) is now an $w^*$-infinite number.

**Remark 2.2.2**.Note that the Riemann sum given by Eq.(2.2.3) is a finite $w^*$-hyperreal number; thus it has a $w^*$-standard part.

**Definition 2.2.1**.Let $[a,b] \subseteq I$ and let $d^{\#_{w^*}}x$ be a positive $w^*$-infinitesimal. The definiteinte integral of $f$ from $a$ to $b$ with respect to $d^{\#_{w^*}}x$ is the $w^*$-standard part of the Riemann sum,

$$\int_a^b f(x)dx \triangleq w^*\text{-}\mathbf{st}\left(\sum_a^b [{}^{\#_{w^*}}f(x)]d^{\#_{w^*}}x\right). \qquad (2.2.4)$$

**Remark 2.2.3**.Note that this definition depends upon the choice of infinitesimal dx. But it can be immediately proved that if $dx$ and $du$ are two positive $w^*$-infinitesimals, then

$$\int_a^b f(x)dx \triangleq \int_a^b f(u)du. \qquad (2.2.5)$$

Note that the $x$ in $f(x)$ and $u$ in f$(u)$ are dummy variables; the $d^{\#_{w^*}}x$ and the $d^{\#_{w^*}}u$ are not.

**Notation.2.2.1**.

**Remark 2.2**.It may be convenient to let the internal space by the $\#_{w^*}$-transform or by the $\#_w$-transform in general case of a classical standard measure space. For instance,if $(\mathbb{R},\mathcal{F},\mu)$ is the standard Lebesgue space on $\mathbb{R}$, our internal starting point could be the $w^*$-internal measure space $({}^{\#_{w^*}}\mathbb{R},{}^{\#_{w^*}}\mathcal{F},{}^{\#_{w^*}}\mu)$ and $w$-internal measure space $({}^{\#_w}\mathbb{R},{}^{\#_w}\mathcal{F},{}^{\#_w}\mu)$. Here ${}^{\#_{w^*}}\mu$ and ${}^{\#_w}\mu$ is finitely, hence finitely, hence hyperfinitely, additive on the $w^*$-internal algebra ${}^{\#_{w^*}}\mathcal{F}$and $w$-internal algebra ${}^{\#_w}\mathcal{F}$ correspondingly.Of course $\sigma$-additivity is lost in the transition. However, it is restored by passing to the associated Loeb space.By transfer we can write down integrals

$$\int_A {}^{\#_{w^*}}f(x)d^{\#_{w^*}}\mu(x) \qquad (2.2.5)$$
$$A \in_{w^*} {}^{\#_{w^*}}\mathcal{F}$$

and

$$\int_A {}^{\#_w}f(x)d^{\#_w}\mu(x) \qquad (2.2.5)$$
$$A \in_w {}^{\#_w}\mathcal{F},$$

which however must be handled with some care: no countable manipulations are allowed.

# Chapter III.

## III.1.Riemann’s non differentiable function.

According to Weierstrass [32],in a talk to the Royal Academy of Sciences in Berlin on 18 July 1872, Riemann introduced the function:

$$\Re(x) = \sum_{n=1}^{\infty} \frac{\sin(\pi n^2 x)}{n^2}. \tag{3.1.1}$$

in order to warn that continuous functions need not have a derivative.Not succeeding in verifying that $\Re(x)$ is nowhere differentiable, Weierstrass proved this property instead for the series

$$W(x) = \sum_{n=1}^{\infty} b^n \cos(a^n t), 0 < b < 1, 0 \leq a. \tag{3.1.2}$$

This appeared first in print in Du-Bois-Reymond [33]. According to Butzer and Stark [34], there are no other known sources which confirm Riemann's role in the story.Hardy [35,pp.322-323] proved that Riemann's function $\Re(x)$ is not differentiable in any irrational point $x \in \mathbb{R}$ and also $\Re(x)$ is not differentiable in a some class of rational point $x \in \mathbb{Q}$ .Gerver [36] succeeded in 1970 in showing that at every rational point $r = p/q$ with $p$ and $q$ both odd, $\Re(x)$ is differentiable, and has derivative equal to $-1/2$ at $r$. Furthermore he showed that at all other rational points the function is not differentiable. Other,shorter proofs were given by Smith [37], Quefelec [38], Mohr [39], Itatsu [40], Luther [41] and Holschneider and Tchamitchian [42]. For previous reviews on Riemann's function, see Neuenschwander [43] and Segal [44]; the literature list of [34] contains many further references abaut the

Riemann's function $\Re(x)$. In paper [45] Gerver introduced the function:

$$G_{3,\beta}(x) = \sum_{n=1}^{\infty} \frac{\exp(in^3 x)}{n^\beta}. \tag{3.1.3}$$

For reals $2 < \beta < 4$,in [45] directed analyze the behavior,near the points $y = \frac{p\pi}{q}$ of

(3.1.3).considered as a function of $x$ ,and expand this series into a constant term, a term on the order of quantity $z_1(x) = \left(x - \frac{p\pi}{q}\right)^{\frac{\beta-1}{3}}$, a term linear in $z_2(x) = \left(x - \frac{p\pi}{q}\right)$ a "chirp" term on the order of quantity $z_2(x) = \left(x - \frac{p\pi}{q}\right)^{\frac{2\beta-1}{4}}$, and an error term on the order $z_2(x) = \left(x - \frac{p\pi}{q}\right)^{\frac{\beta}{2}}$. At every such rational point, the left and right derivatives are either both finite (and equal) or both infinite, in contrast with the quadratic series, where the derivative is often finite on one side and infinite on the other. However, in the cubic series, again in contrast with the quadratic case, the chirp term generally has a different set of frequencies and amplitudes on the right and left sides. Finally, in [45] was shown that almost every irrational point can be closely approximated, in a suitable Diophantine sense, by rational points where the cubic series has an infinite derivative. This implies that when

$$\beta \leq \frac{\sqrt{97} - 1}{4} = 2.212\ldots, \tag{3.1.4}$$

both the real and imaginary parts of the cubic series are differentiable almost nowhere. At the same time it is necessary to note that in spite of a big progress obtained in the considered studies area, any general absence criterions of the finite almost everywhere derivate for absolutely convergent trigonometrical series was not obtained. In [22]-[23], using the methods of paralogical nonstandard analysis, was obtained the general criterion of the absence almost everywhere finite derivative for the following continuous function $\Re(x;\omega_1(n),\omega_2(n))$ :

$$\Re(x;\omega_1(n),\omega_2(n)) = \sum_{n=1}^{\infty} \frac{\exp(i \bullet x \bullet \omega_1(n))}{\omega_2(n)},$$
$$\omega_1 : \mathbb{N} \to \mathbb{N}, \omega_2 : \mathbb{N} \to \mathbb{N}, \qquad (3.1.5)$$
$$\sum_{n=1}^{\infty} \frac{1}{|\omega_2(n)|} < \infty.$$

It is shown in [22]-[23] that under condition

$$\sum_{n=1}^{\infty} \left( \frac{\omega_1(n)}{\omega_2(n)} \right)^2 = \infty \qquad (3.1.6)$$

function $\Re(x;\omega_1(n),\omega_2(n))$ does not have a finite derivate on a quantity of a positive measure. Particularly we shall reinforce the foregoing Gerver's result by showing that inequality (3.1.4) is possible to change by inequality $\beta \leq 4$, at least for a quantity of points of a positive measure.

.

# III.2.Non standard proof of the non-differentiability of the Riemann function $\Re(x)$.

Non-differentiable Riemann function $\Re(x)$ is defined by

$$\Re(x) = \sum_{n=1}^{\infty} \frac{\sin(n^2 x)}{n^2}, \qquad (3.2.1)$$

see subsection III.1.

**Theorem 3.2.1**. $\Re(x)$ is not a.e. differentiable on $[0,\pi]$.

Proof. See Remark 3.2.1 etc.

**Remark 3.2.1**.Remind that there exist imbeding

$$j_{\#w_{\{1\}}} : \mathbb{N} \hookrightarrow_{w_{\{1\}}} \mathbb{N}^{\#w_{\{1\}}} \qquad (3.2.2)$$

and there exist imbeding

$$j_{\#w_{\{1\}}} : \mathbb{R} \hookrightarrow_{w_{\{1\}}} \mathbb{R}^{\#w_{\{1\}}} \qquad (3.2.3)$$

such that

$$j_{\#w_{\{1\}}}(\mathbb{N}) =_{w_{\{1\}}} \mathbb{N}_{\#w_{\{1\}}} \subsetneqq_{\#w_{\{1\}}} \mathbb{N}^{\#w_{\{1\}}} \qquad (3.2.4)$$

and

$$j_{\#w_{\{1\}}}(\mathbb{R}) =_{w_{\{1\}}} \mathbb{R}_{\#w_{\{1\}}} \subsetneqq_{\#w_{\{1\}}} \mathbb{R}^{\#w_{\{1\}}} \quad (3.2.5)$$

correspondingly.

**Notation.3.2.1.**(i) We will use the following notation $j_{\#w_{\{1\}}}(n) \triangleq n_{\#w_{\{1\}}}, n \in \mathbb{N}$ and $j_{\#w_{\{1\}}}(x) \triangleq_{\#w_{\{1\}}} x_{\#w_{\{1\}}}, x \in \mathbb{R}$, $j_{\#w_{\{1\}}}(x \times y) \triangleq_{\#w_{\{1\}}} x_{\#w_{\{1\}}} \times_{\#w_{\{1\}}} y_{\#w_{\{1\}}}$,
$j_{\#w_{\{1\}}}(\infty) \triangleq_{\#w_{\{1\}}} \infty_{\#w_{\{1\}}}$, etc.

(**ii**) we often letter for short: simply $n$ instead $n_{\#w_{\{1\}}}$, simply $x$ instead $x_{\#w_{\{1\}}}$, simply $x \times y$ instead $x_{\#w_{\{1\}}} \times_{\#w_{\{1\}}} y_{\#w_{\{1\}}}$, etc.

(iii) We will use the following notation $f^{\#w_{\{1\}}}, \mu^{\#w_{\{1\}}}, \int^{\#w_{\{1\}}}, \ldots$ instead ${}^{\#w_{\{1\}}}f$, $\left({}^{\#w_{\{1\}}}\mu\right)$, $\left({}^{\#w_{\{1\}}}\int\right)$, $j_{\#w_{\{1\}}}(\mathbf{T}^1) =_{w_{\{1\}}} \mathbf{T}^1_{w_{\{1\}}}$ etc.

(iv) we let for short $j_{\#w_{\{1\}}}(\sin(\pi n^2 x)) =_{w_{\{1\}}} \sin_{\#w_{\{1\}}}\left(n^2_{\#w_{\{1\}}} x\right)$, where $x \in_{\#w_{\{1\}}} \mathbf{T}^1_{w_{\{1\}}}$, etc.

**Definition 3.2.1.**We define now a $w_{\{1\}}$-function $\Re_{\#w_{\{1\}}} : \mathbf{T}^1_{w_{\{1\}}} \to_{w_{\{1\}}} \mathbf{T}^1_{w_{\{1\}}}$ :

$$\Re_{\#w_{\{1\}}}\left(j_{\#w_{\{1\}}}(x)\right) =_{w_{\{1\}}} j_{\#w_{\{1\}}}(\Re((x))) =_{w_{\{1\}}}$$
$$j_{\#w_{\{1\}}}\left(\sum_{n=1}^{\infty} \frac{\sin(n^2 x)}{n^2}\right) =_{w_{\{1\}}} \sum_{n_{w_{\{1\}}} =_{w_{\{1\}}} 1_{w_{\{1\}}}}^{\infty_{w_{\{1\}}}} \left(\frac{\sin_{\#w_{\{1\}}}\left(n^2_{\#w_{\{1\}}} x\right)}{n^2_{\#w_{\{1\}}}}\right), \quad (3.2.6)$$
$$x \in_{\#w_{\{1\}}} \mathbf{T}^1_{w_{\{1\}}}$$

**Definition 3.2.2.**We define now a $w_{\{1\}}$-function $\Re^{\dagger}_{\#w_{\{1\}}} : \mathbf{T}^{1\#w_{\{1\}}} \to_{w_{\{1\}}} \mathbf{T}^{1\#w_{\{1\}}}$ :

$$\Re^{\dagger}_{\#w_{\{1\}}}(x) =_{w_{\{1\}}}$$
$$=_{w_{\{1\}}} \textbf{Ext-}w_{\{1\}}\text{-} \sum_{n =_{w_{\{1\}}} 1_{w_{\{1\}}}}^{M} \frac{c_n \times_{\#w_{\{1\}}} \sin^{\#w_{\{1\}}}\left(n^2 \times_{\#w_{\{1\}}} x\right)}{n^2}, \quad (3.2.7)$$

where $M \in_{w_{\{1\}}} \mathbb{N}^{\#w_{\{1\}}} \setminus_{w_{\{1\}}} \mathbb{N}_{\#w_{\{1\}}}$ and

$$c_n =_{w_{\{1\}}} \begin{cases} 1_{\#w_{\{1\}}} & \text{iff} \quad n \in_{w_{\{1\}}} \mathbb{N}_{\#w_{\{1\}}} \\ 0_{w_{\{1\}}} & \text{iff} \quad n \in_{w_{\{1\}}} \mathbb{N}^{\#w_{\{1\}}} \setminus_{w_{\{1\}}} \mathbb{N}_{\#w_{\{1\}}} \end{cases} \quad (3.2.8)$$

**Remark 3.2.3.**Note that for any $x \in_{\#w_{\{1\}}} \mathbf{T}^1_{w_{\{1\}}}$ :

$$\Re_{\#w_{\{1\}}}(x) =_{\#w_{\{1\}}} \Re^{\dagger}_{\#w_{\{1\}}}(x). \quad (3.2.9)$$

**Remark 3.2.4.**We assume now that a Riemann function $\Re(x)$ is differentiable almost everywhere in the sense of the Lebesgue measure $d\mu = d\mu_{\mathcal{L}}$, i.e., a.e. the derivative $\Re'(x)$ exists and finite, i.e., $\exists \xi(x)$ such that a.e. $\xi(x) < \infty$ and

$$\text{a.e.}: \ \Re'(x) < \xi(x) < \infty. \quad (3.2.10)$$

**Remark 3.2.4.**Therefore (i) from Eq.(3.2.9) by $w_{\{1\}}$-transfer it follows that a $w_{\{1\}}$-function $\Re^{\dagger}_{\#_{w_{\{1\}}}}(x)$ is $\#_{w_{\{1\}}}$-differentiable $\#_{w_{\{1\}}}$-almost everywhere on $\mathbf{T}^{1\#_{w_{\{1\}}}}$ in the sense of the $w_{\{1\}}$-transfered Lebesgue measure $d^{\#_{w_{\{1\}}}}\mu = d^{\#_{w_{\{1\}}}}\mu_{\mathcal{L}}$. (ii) By $w_{\{1\}}$-transfer from (3.2.10) we obtain

$$w_{\{1\}}\text{-a.e.: } \frac{d^{\#_{w_{\{1\}}}}}{d^{\#_{w_{\{1\}}}}x}\Re^{\dagger}_{\#_{w_{\{1\}}}}(x) <_{\#_{w_{\{1\}}}} \xi^{\#_{w_{\{1\}}}}(x), \tag{3.2.11}$$

where

$$w_{\{1\}}\text{-a.e.: } \xi^{\#_{w_{\{1\}}}}(x) <_{\#_{w_{\{1\}}}} \infty^{\#_{w_{\{1\}}}} \tag{3.2.12}$$

From Eqs.(3.2.7)-(3.2.8) by $w_{\{1\}}$-differentiation one obtains

$$\begin{gathered}\frac{d^{\#_{w_{\{1\}}}}}{d^{\#_{w_{\{1\}}}}x}\Re^{\dagger}_{\#_{w_{\{1\}}}}(x) =_{w_{\{1\}}} \\ =_{w_{\{1\}}} \frac{d^{\#_{w_{\{1\}}}}}{d^{\#_{w_{\{1\}}}}x}\left[\mathbf{Ext}\text{-}w_{\{1\}}\text{-}\sum_{n=_{w_{\{1\}}}1_{w_{\{1\}}}}^{M} \frac{c_n \times_{\#_{w_{\{1\}}}} \sin^{\#_{w_{\{1\}}}}\left(n^2 \times_{\#_{w_{\{1\}}}} x\right)}{n^2}\right] \\ =_{w_{\{1\}}} \mathbf{Ext}\text{-}w_{\{1\}}\text{-}\sum_{n=_{w_{\{1\}}}1_{w_{\{1\}}}}^{M} \frac{d^{\#_{w_{\{1\}}}}}{d^{\#_{w_{\{1\}}}}x}\left[\frac{c_n \times_{\#_{w_{\{1\}}}} \sin^{\#_{w_{\{1\}}}}\left(n^2 \times_{\#_{w_{\{1\}}}} x\right)}{n^2}\right] \\ =_{w_{\{1\}}w_{\{1\}}} \mathbf{Ext}\text{-}w_{\{1\}}\text{-}\sum_{n=_{w_{\{1\}}}1_{w_{\{1\}}}}^{M} c_n \times_{\#_{w_{\{1\}}}} \cos^{\#_{w_{\{1\}}}}\left(n^2 \times_{\#_{w_{\{1\}}}} x\right).\end{gathered} \tag{3.2.13}$$

Thus finally we obtain

$$\begin{gathered}\Phi(x) =_{w_{\{1\}}} \frac{d^{\#_{w_{\{1\}}}}}{d^{\#_{w_{\{1\}}}}x}\Re^{\dagger}_{\#_{w_{\{1\}}}}(x) =_{w_{\{1\}}} \\ =_{w_{\{1\}}} \mathbf{Ext}\text{-}w_{\{1\}}\text{-}\sum_{n=_{w_{\{1\}}}1_{w_{\{1\}}}}^{M} c_{n_{\#_{w_{\{1\}}}}} \times_{\#_{w_{\{1\}}}} \cos^{\#_{w_{\{1\}}}}\left(n^2 \times_{\#_{w_{\{1\}}}} x\right).\end{gathered} \tag{3.2.14}$$

**Remark 3.2.5**. Note that a $w_{\{1\}}$-function $\Phi(x)$ is not $w_{\{1\}}$-a.e. $w_{\{1\}}$-finite on $\mathbf{T}^{1\#_{w_{\{1\}}}}$, i.e.

$$\neg_s\left[w_{\{1\}}\text{-a.e. : } \Phi(x) <_{\#_{w_{\{1\}}}} \infty^{\#_{w_{\{1\}}}} \middle| x \in_{w_{\{1\}}} \mathbf{T}^{1\#_{w_{\{1\}}}}\right] \tag{3.2.15}$$

In order to proof (3.2.15) we calculate now the $w_{\{1\}}$-integral

$$\int_{[-\pi,\pi]^{w_{\{1\}}}} \Phi(x) \times_{w_{\{1\}}} \Phi(x) d^{\#_{w_{\{1\}}}}\mu(x).$$

From Eq.(3.2.14) one obtains

$$\Phi(x)\times_{w_{\{1\}}}\Phi(x)=_{w_{\{1\}}}$$

$$=_{w_{\{1\}}}\textbf{Ext-}w_{\{1\}}\text{-}\sum_{n_{\#w_{\{1\}}}=_{w_{\{1\}}}1_{w_{\{1\}}}}^{M}\sum_{m_{\#w_{\{1\}}}=_{w_{\{1\}}}1_{w_{\{1\}}}}^{M}c_{n_{\#w_{\{1\}}}}\times_{\#w_{\{1\}}}c_{m_{\#w_{\{1\}}}}\times_{w_{\{1\}}}\tag{3.2.16}$$
$$\times_{w_{\{1\}}}\cos^{\#_{w_{\{1\}}}}\left(n^2_{\#w_{\{1\}}}\times_{\#w_{\{1\}}}x\right)\times_{w_{\{1\}}}\cos^{\#_{w_{\{1\}}}}\left(m^2_{\#w_{\{1\}}}\times_{\#w_{\{1\}}}x\right).$$

From Eq.(3.2.16) by $w_{\{1\}}$-integratiion one obtains

$$\int_{[-\pi,\pi]^{w_{\{1\}}}}\Phi(x)\times_{w_{\{1\}}}\Phi(x)d^{\#_{w_{\{1\}}}}\mu(x)=_{w_{\{1\}}}$$
$$=_{w_{\{1\}}}\textbf{Ext-}w_{\{1\}}\text{-}\sum_{n_{\#w_{\{1\}}}=_{w_{\{1\}}}1_{w_{\{1\}}}}^{M}\sum_{m_{\#w_{\{1\}}}=_{w_{\{1\}}}1_{w_{\{1\}}}}^{M}c_n\times_{\#w_{\{1\}}}c_m\times_{w_{\{1\}}}\tag{3.2.17}$$
$$\times_{w_{\{1\}}}\int_{[0,\pi]^{w_{\{1\}}}}\cos^{\#_{w_{\{1\}}}}\left(n^2\times_{\#w_{\{1\}}}x\right)\times_{w_{\{1\}}}\cos^{\#_{w_{\{1\}}}}\left(m^2\times_{\#w_{\{1\}}}x\right)d^{\#_{w_{\{1\}}}}\mu(x).$$

where by $d^{\#_{w_{\{1\}}}}\mu(x)$ we denote $\#_{w_{\{1\}}}$- transfered standard Lebesgue measure $d\mu(x)$ on $[0,\pi]$.

Note that

$$\int_{[-\pi,\pi]^{w_{\{1\}}}}\cos^{\#_{w_{\{1\}}}}\left(n^2\times_{\#w_{\{1\}}}x\right)\times_{w_{\{1\}}}\cos^{\#_{w_{\{1\}}}}\left(n^2\times_{\#w_{\{1\}}}x\right)d^{\#_{w_{\{1\}}}}\mu(x)=_{w_{\{1\}}}\pi_{w_{\{1\}}}\tag{3.2.18}$$

and

$$\int_{[-\pi,\pi]^{w_{\{1\}}}}\cos^{\#_{w_{\{1\}}}}\left(n^2\times_{\#w_{\{1\}}}x\right)\times_{w_{\{1\}}}\cos^{\#_{w_{\{1\}}}}\left(m^2\times_{\#w_{\{1\}}}x\right)d^{\#_{w_{\{1\}}}}\mu(x)=_{w_{\{1\}}}0_{w_{\{1\}}}\tag{3.2.19}$$

iff $\neg_s\left(n_{\#w_{\{1\}}}=_{w_{\{1\}}}m_{\#w_{\{1\}}}\right)$. Then from Eq.(3.2.17) and Eqs.(3.2.18)-(3.2.19) one obtains

$$\int_{[-\pi,\pi]^{w_{\{1\}}}}\Phi(x)\times_{w_{\{1\}}}\Phi(x)d^{\#_{w_{\{1\}}}}\mu(x)=_{w_{\{1\}}}$$
$$\pi_{w_{\{1\}}}\times_{w_{\{1\}}}\left(\textbf{Ext-}w_{\{1\}}\text{-}\sum_{n=_{w_{\{1\}}}1_{w_{\{1\}}}}^{M}c_n\right)=_{w_{\{1\}}}\pi_{w_{\{1\}}}\times_{w_{\{1\}}}\Omega,\tag{3.2.20}$$

where

$$\Omega=_{w_{\{1\}}}\textbf{Ext-}w_{\{1\}}\text{-}\sum_{n=_{w_{\{1\}}}1_{w_{\{1\}}}}^{M}c_n\tag{3.2.21}$$

and therefore

$$\int_{[-\pi,\pi]^{w_{\{1\}}}}\Phi(x)\times_{w_{\{1\}}}\Phi(x)d^{\#_{w_{\{1\}}}}\mu(x)=_{w_{\{1\}}}\pi_{w_{\{1\}}}\times_{w_{\{1\}}}\Omega.\tag{3.2.22}$$

**Remark 3.2.6**. Note that obviously $\Omega\in_{w_{\{1\}}}\mathbb{R}^{\#_{w_{\{1\}}}}\backslash_{w_{\{1\}}}\mathbb{R}^{\#_{w_{\{1\}}}}_{\mathbf{fin}}$ and therefore (3.2.15) holds.

But (3.2.15) contradicts with (3.2.11). This contradiction finalized the proof.

# III.3.Non standard proof of the non-differentiability of the Generalized Riemann function $\Re(x;\omega_1(n),\omega_2(n))$.

**Theorem 3.3.1**. Let $\Re(x;\omega_1(n),\omega_2(n))$ be the continuous function

$$\Re(x;\omega_1(n),\omega_2(n)) = \sum_{n=1}^{\infty} \frac{\exp(i \bullet x \bullet \omega_1(n))}{\omega_2(n)}, \tag{3.3.1}$$

where $\omega_1 : \mathbb{N} \to \mathbb{N}$, $\omega_1 : \mathbb{N} \to \mathbb{N}$ and the following conditions holds:

(i) $\forall n \forall m[(\omega_1(n) = \omega_1(m)) \Leftrightarrow n = m]$,

(ii) $\sum_{n=1}^{\infty} \frac{1}{|\omega_2(n)|} < \infty$ and

(iii)

$$\sum_{n=1}^{\infty}\left(\frac{\omega_1(n)}{\omega_2(n)}\right)^2 = \infty. \tag{3.3.2}$$

Then a function $\Re(x;\omega_1(n),\omega_2(n))$ does not have a finite derivative on a set $\mathcal{F} \subseteq [-\pi,\pi]$ of a positive Lebesgue measure $\mu_{\mathcal{L}}(\mathcal{F}) > 0$.

**Proof**. Similarly to proof of the Theorem 3.3.1.

**Definition 3.3.1**.We define now a $w_{\{1\}}$-function $\Re_{\#w_{\{1\}}}(x;\omega_1(n),\omega_2(n)) : \mathbf{T}^1_{w_{\{1\}}} \to_{w_{\{1\}}} \mathbf{T}^1_{w_{\{1\}}}$ :

$$\begin{gathered}\Re_{\#w_{\{1\}}}(x;\omega_1(n),\omega_2(n)) =_{w_{\{1\}}} j_{\#w_{\{1\}}}(\Re(x;\omega_1(n),\omega_2(n))) =_{w_{\{1\}}} \\ j_{\#w_{\{1\}}}\left(\sum_{n=1}^{\infty}\frac{\sin(\omega_1(n)\times x)}{\omega_2(n)}\right) =_{w_{\{1\}}} \\ \sum_{n_{w_{\{1\}}}=_{w_{\{1\}}}1_{w_{\{1\}}}}^{\infty_{w_{\{1\}}}}\left(\frac{\sin_{\#w_{\{1\}}}\left(\omega_{1w_{\{1\}}}\left(n_{w_{\{1\}}}\right)\times_{w_{\{1\}}} x\right)}{\omega_{2w_{\{1\}}}\left(n_{w_{\{1\}}}\right)}\right), \\ \omega_{1w_{\{1\}}}\left(n_{w_{\{1\}}}\right) =_{w_{\{1\}}} j_{\#w_{\{1\}}}(\omega_1(n)), \\ \omega_{2w_{\{1\}}}\left(n_{w_{\{1\}}}\right) =_{w_{\{1\}}} j_{\#w_{\{1\}}}(\omega_2(n)).\end{gathered} \tag{3.3.3}$$

**Definition 3.3.2**.We define now a $w_{\{1\}}$-function $\Re^{\dagger}_{\#w_{\{1\}}} : \mathbf{T}^{1\#w_{\{1\}}} \to_{w_{\{1\}}} \mathbf{T}^{1\#w_{\{1\}}}$ :

$$\begin{gathered}\Re^{\dagger}_{\#w_{\{1\}}}(x;\omega_1(n),\omega_2(n)) =_{w_{\{1\}}} \\ =_{w_{\{1\}}} \textbf{Ext-}w_{\{1\}}\text{-}\sum_{n=_{w_{\{1\}}}1_{w_{\{1\}}}}^{M}\frac{c_n\times_{\#w_{\{1\}}}\sin^{\#w_{\{1\}}}\left(\omega_1^{w_{\{1\}}}(n)\times_{w_{\{1\}}} x\right)}{\omega_2^{w_{\{1\}}}(n)},\end{gathered} \tag{3.3.4}$$

where $M \in_{w_{\{1\}}} \mathbb{N}^{\#w_{\{1\}}} \backslash_{w_{\{1\}}} \mathbb{N}_{\#w_{\{1\}}}$ and

$$c_n =_{w_{\{1\}}} \begin{cases} 1_{\#w_{\{1\}}} & \text{iff} \quad n \in_{w_{\{1\}}} \mathbb{N}_{\#w_{\{1\}}} \\ 0_{w_{\{1\}}} & \text{iff} \quad n \in_{w_{\{1\}}} \mathbb{N}^{\#w_{\{1\}}} \backslash_{w_{\{1\}}} \mathbb{N}_{\#w_{\{1\}}} \end{cases} \quad (3.3.5)$$

**Remark 3.3.1.** Note that for any $x \in_{\#w_{\{1\}}} \mathbf{T}^1_{w_{\{1\}}}$ :

$$\Re_{\#w_{\{1\}}}(x;\omega_1(n),\omega_2(n)) =_{\#w_{\{1\}}} \Re^{\dagger}_{\#w_{\{1\}}}(x;\omega_1(n),\omega_2(n)). \quad (3.3.6)$$

**Remark 3.3.2.** We assume now that a Generalized Riemann function $\Re(x;\omega_1(n),\omega_2(n))$ is differentiable almost everywhere in the sense of the Lebesgue measure $d\mu = d\mu_{\mathcal{L}}$, i.e., a.e. on $\mathbf{T}^1$ the derivative $\Re'(x;\omega_1(n),\omega_2(n))$ exists and finite, i.e., $\exists\xi(x)$ such that a.e. $\xi(x) < \infty$ and

$$\text{a.e.:}\ \Re'(x;\omega_1(n),\omega_2(n)) < \xi(x) < \infty. \quad (3.3.7)$$

**Remark 3.3.3.** Therefore (i) from Eq.(3.3.6) by $w_{\{1\}}$-transfer it follows that a $w_{\{1\}}$-function $\Re^{\dagger}_{\#w_{\{1\}}}(x;\omega_1(n),\omega_2(n))$ is $\#_{w_{\{1\}}}$-differentiable $\#_{w_{\{1\}}}$-almost everywhere on $\mathbf{T}^{1\#w_{\{1\}}}$ in the sense of the $w_{\{1\}}$-transfered Lebesgue measure $d^{\#w_{\{1\}}}\mu = d^{\#w_{\{1\}}}\mu_{\mathcal{L}}$. (ii) By $w_{\{1\}}$-transfer from (3.3.7) we obtain

$$w_{\{1\}}\text{-a.e.:}\ \frac{d^{\#w_{\{1\}}}}{d^{\#w_{\{1\}}}x}\Re^{\dagger}_{\#w_{\{1\}}}(x;\omega_1(n),\omega_2(n)) <_{\#w_{\{1\}}} \xi^{\#w_{\{1\}}}(x), \quad (3.3.8)$$

where

$$w_{\{1\}}\text{-a.e.:}\ \xi^{\#w_{\{1\}}}(x) <_{\#w_{\{1\}}} \infty^{\#w_{\{1\}}} \quad (3.3.9)$$

From Eqs.(3.3.4)-(3.3.5) by $w_{\{1\}}$-differentiation one obtains

$$\begin{aligned}
&\frac{d^{\#w_{\{1\}}}}{d^{\#w_{\{1\}}}x}\Re^{\dagger}_{\#w_{\{1\}}}(x;\omega_1(n),\omega_2(n)) =_{w_{\{1\}}} \\
&=_{w_{\{1\}}} \frac{d^{\#w_{\{1\}}}}{d^{\#w_{\{1\}}}x}\left[\mathbf{Ext}\text{-}w_{\{1\}}\text{-}\sum_{n=_{w_{\{1\}}}1_{w_{\{1\}}}}^{M} \frac{c_n \times_{\#w_{\{1\}}} \sin^{\#w_{\{1\}}}\left(\omega_1^{w_{\{1\}}}(n) \times_{w_{\{1\}}} x\right)}{\omega_2^{\#w_{\{1\}}}(n)}\right] \\
&=_{w_{\{1\}}} \mathbf{Ext}\text{-}w_{\{1\}}\text{-}\sum_{n=_{w_{\{1\}}}1_{w_{\{1\}}}}^{M} \frac{d^{\#w_{\{1\}}}}{d^{\#w_{\{1\}}}x}\left[\frac{c_n \times_{\#w_{\{1\}}} \sin^{\#w_{\{1\}}}\left(\omega_1^{w_{\{1\}}}(n) \times_{w_{\{1\}}} x\right)}{\omega_2^{\#w_{\{1\}}}(n)}\right] \quad (3.3.10) \\
&=_{w_{\{1\}}w_{\{1\}}} \mathbf{Ext}\text{-}w_{\{1\}}\text{-}\sum_{n=_{w_{\{1\}}}1_{w_{\{1\}}}}^{M} c_n \times_{\#w_{\{1\}}} \frac{\omega_1^{w_{\{1\}}}(n)}{\omega_2^{\#w_{\{1\}}}(n)} \times_{\#w_{\{1\}}} \\
&\times_{\#w_{\{1\}}} \cos^{\#w_{\{1\}}}\left(\omega_1^{w_{\{1\}}}(n) \times_{w_{\{1\}}} x\right).
\end{aligned}$$

Thus finally we obtain

$$\Phi(x;\omega_1(n),\omega_2(n)) =_{w_{\{1\}}} \frac{d^{\#_{w_{\{1\}}}}}{d^{\#_{w_{\{1\}}}}x} \Re^{\dagger}_{\#_{w_{\{1\}}}}(x;\omega_1(n),\omega_2(n)) =_{w_{\{1\}}}$$

$$=_{w_{\{1\}}} \mathbf{Ext}\text{-}w_{\{1\}}\text{-}\sum_{n=_{w_{\{1\}}}1_{w_{\{1\}}}}^{M} c_n \times_{\#_{w_{\{1\}}}} \frac{\omega_1^{w_{\{1\}}}(n)}{\omega_2^{\#_{w_{\{1\}}}}(n)} \times_{\#_{w_{\{1\}}}} \tag{3.3.11}$$

$$\times_{\#_{w_{\{1\}}}} \cos^{\#_{w_{\{1\}}}}\left(\omega_1^{w_{\{1\}}}(n) \times_{w_{\{1\}}} x\right).$$

**Remark 3.3.4**. Note that a $w_{\{1\}}$-function $\Phi(x;\omega_1(n),\omega_2(n))$ is not $w_{\{1\}}$-a.e. $w_{\{1\}}$-finite on $\mathbf{T}^{1\#_{w_{\{1\}}}}$, i.e.

$$\neg_s\left[w_{\{1\}}\text{-a.e. : } \Phi(x;\omega_1(n),\omega_2(n)) <_{\#_{w_{\{1\}}}} \infty^{\#_{w_{\{1\}}}} \middle| x \in_{w_{\{1\}}} \mathbf{T}^{1\#_{w_{\{1\}}}}\right] \tag{3.3.12}$$

In order to proof (3.2.12) we calculate now the $w_{\{1\}}$-integral

$$\int_{[-\pi,\pi]^{w_{\{1\}}}} \left[\Phi(x;\omega_1(n),\omega_2(n)) \times_{w_{\{1\}}} \Phi(x;\omega_1(n),\omega_2(n))\right] d^{\#_{w_{\{1\}}}}\mu(x). \tag{3.3.13}$$

From Eq.(3.2.11) one obtains

$$\Phi(x;\omega_1(n),\omega_2(n)) \times_{w_{\{1\}}} \Phi(x;\omega_1(n),\omega_2(n)) =_{w_{\{1\}}}$$

$$=_{w_{\{1\}}} \mathbf{Ext}\text{-}w_{\{1\}}\text{-}\sum_{n_{\#_{w_{\{1\}}}} =_{w_{\{1\}}} 1_{w_{\{1\}}}}^{M} \sum_{m_{\#_{w_{\{1\}}}} =_{w_{\{1\}}} 1_{w_{\{1\}}}}^{M} c_n \times_{\#_{w_{\{1\}}}} c_m \times_{w_{\{1\}}}$$

$$\times_{w_{\{1\}}} \frac{\omega_1^{w_{\{1\}}}(n)}{\omega_2^{\#_{w_{\{1\}}}}(n)} \times_{w_{\{1\}}} \frac{\omega_1^{w_{\{1\}}}(m)}{\omega_2^{\#_{w_{\{1\}}}}(m)} \times_{w_{\{1\}}} \tag{3.2.16}$$

$$\times_{w_{\{1\}}} \cos^{\#_{w_{\{1\}}}}\left(\omega_1^{w_{\{1\}}}(n) \times_{w_{\{1\}}} x\right) \times_{w_{\{1\}}} \cos^{\#_{w_{\{1\}}}}\left(\omega_1^{w_{\{1\}}}(m) \times_{w_{\{1\}}} x\right).$$

# IV.1.Non standard proof of the Carleson's theorem.

Let us consider Fourier series in space $\mathcal{L}_2(\mathbf{T}^1)$

$$\sum_{n=0}^{\infty} c_n \exp(inx), \tag{4.1.1}$$

where $\mathbf{T}^1 = [-\pi,\pi]$, such that

$$\sum_{n=0}^{\infty} |c_n|^2 < \infty. \tag{4.1.2}$$

**Remark 4.1.1.**Note that in this section we will be consider more general trigonometric series such that

$$\sum_{k=1}^{\infty} c_k \exp(ixn_k(x)), \qquad (4.1.3)$$

where $\mathbf{T}^1 = [-\pi,\pi]$, $n_k \to \infty$ if $k \to \infty$ and

$$\sum_{k=1}^{\infty} |c_k|^2 < \infty, \qquad (4.1.4)$$

or

$$\sum_{k=1}^{\infty} c_{n_k(x)} \exp(ixn_k(x)), \qquad (4.1.5)$$

where $\mathbf{T}^1 = [-\pi,\pi]$, $n_k \to \infty$ if $k \to \infty$ and

$$\sum_{k=1}^{\infty} |c_{n_k(x_k)}|^2 < \infty \qquad (4.1.6)$$

$\forall k (x_k \in \mathbf{T}^1)$

(**I**) Now we go to prove that under the condition (4.1.6) the following statement holds: for any sequence $p_k$ a.e. on $\mathbf{T}^1$

$$\lim_{p_k \to \infty, k \to \infty} \left| \sum_{k=1}^{p_k} c_{n_k(x)} \exp(ixn_k(x)) \right|^2 < \infty. \qquad (4.1.7)$$

(1) In contrary with (4.1.7) we assume now that : a.e. on $\mathbf{T}^1$

$$\left| \sum_{k=1}^{\infty} c_{n_k(x)} \exp(ixn_k(x)) \right|^2 = \infty. \qquad (4.1.8)$$

Let $x \in \mathbf{T}^1$be a real number and there exists a sequence $\{n_k(x)\}_{k \in \mathbb{N}}$ such that

$$\lim_{m_q \to \infty, q \to \infty} \left| \sum_{k=1}^{m_q} c_{n_k(x)} \exp(ixn_k(x)) \right|^2 = \infty. \qquad (4.1.9)$$

**Remark 4.1.2**.Let $x \in \mathbf{T}^1$be a real number. Note that a sequence $\{n_k(x)\}_{k \in \mathbb{N}}$ mentioned above in Eq.(4.1.7) in general case is not unique and there exists infinite set of the sequences $\{n_k^l(x)\}_{k \in \mathbb{N}}, l = 1,2,\ldots$such that for any $l \in \mathbb{N}$

$$\lim_{m_q \to \infty, q \to \infty} \left| \sum_{k=1}^{m_q} c_{n_k^l(x)} \exp(ixn_k^l(x)) \right|^2 = \infty. \qquad (4.1.10)$$

**Remark 4.1.3**.Note that any sequence $\{n_k^l(x)\}_{k \in \mathbb{N}}, l = 1,2,\ldots$ mentioned above in Eq.(4.1.8) depend on number $x \in \mathbf{T}^1$ and we will be denoted such sequences $\{n_k^l\}_{k \in \mathbb{N}}$ by $\{n_k^l(x)\}_{k \in \mathbb{N}}$ or by $\{n_k(x)\}_{k \in \mathbb{N}}, \{m_k(x)\}_{k \in \mathbb{N}}, \{r_k(x)\}_{k \in \mathbb{N}}$ etc.

**Remark 4.1.4**.Note that from (4.1.7) it follows that : a.e. on $\mathbf{T}^1$

$$\lim_{m_q \to \infty} \left| \sum_{k=0}^{m_q} c_{n_k(x)} \exp(ixn_k(x)) \right|^2 = \infty, \qquad (4.1.11)$$

where $m_q \to \infty, q \to \infty$.From (4.1.11) by $\#_{w_{\{1\}}}$- transfer it follows that : $\#_{w_{\{1\}}}$- a.e. on $\mathbf{T}^1_{\#_{w_{\{1\}}}} \triangleq (\mathbf{T}^1)^{\#_{w_{\{1\}}}}$

$$\left| \#_{w_{\{1\}}}\textbf{-Ext}\sum_{k=_{w_{\{1\}}}0_{w_{\{1\}}}}^{M} c^{\#_{w_{\{1\}}}}_{n_k^{\#_{w_{\{1\}}}}(x)} \exp\left(ix \times_{\#_{w_{\{1\}}}} n_k^{\#_{w_{\{1\}}}}(x)\right)\right|^2 =_{\#_{w_{\{1\}}}}$$

$$\#_{w_{\{1\}}}\textbf{-Ext}\sum_{k_1=_{w_{\{1\}}}0_{w_{\{1\}}}}^{M} \#_{w_{\{1\}}}\textbf{-Ext}\sum_{k_2=_{w_{\{1\}}}0_{w_{\{1\}}}}^{M} c^{\#_{w_{\{1\}}}}_{n_{k_1}^{\#_{w_{\{1\}}}}(x)} \overline{c^{\#_{w_{\{1\}}}}_{n_{k_2}^{\#_{w_{\{1\}}}}(x)}} \times_{\#_{w_{\{1\}}}} \tag{4.1.12}$$

$$\exp\left(ix \times_{\#_{w_{\{1\}}}} n_{k_1}^{\#_{w_{\{1\}}}}(x)\right) \times_{\#_{w_{\{1\}}}} \overline{\exp\left(ix \times_{\#_{w_{\{1\}}}} n_{k_2}^{\#_{w_{\{1\}}}}(x)\right)} =_{\#_{w_{\{1\}}}}$$

$$=_{\#_{w_{\{1\}}}} N_M(x),$$

where $M \in_{w_{\{1\}}} \mathbb{N}^{\#_{w_{\{1\}}}} \backslash \mathbb{N}_{w_{\{1\}}}$ and

$$N_M(x) \in_{w_{\{1\}}} \mathbb{N}^{\#_{w_{\{1\}}}} \backslash \mathbb{N}_{\#_{w_{\{1\}}}}, \tag{4.1.13}$$

where $\#_{w_{\{1\}}}$-sequence

$$\left\{ n_k^{\#_{w_{\{1\}}}}(x) \right\}_{k \in_{\#_{w_{\{1\}}}} \mathbb{N}^{\#_{w_{\{1\}}}}} \tag{4.1.14}$$

is obtained by using $\#_{w_{\{1\}}}$- transfer from the standard sequence $\{n_k(x)\}_{k\in\mathbb{N}}$, i.e.

$$\left\{ n_k^{\#_{w_{\{1\}}}}(x) \right\}_{k \in_{w_{\{1\}}} \mathbb{N}^{\#_{w_{\{1\}}}}} =_{w_{\{1\}}} (\{n_k(x)\}_{k\in\mathbb{N}})^{\#_{w_{\{1\}}}}. \tag{4.1.15}$$

**Remark 4.1.5**. We introduce now $w_{\{1\}}$-inconsistent hyperintegers $\mathbf{n}_k^{\#}$ corresponding to trigonometric series (4.1.5) by the following way

$$\bigwedge_{x \in_{w_{\{1\}}} \mathbf{T}_{\#}^1} \left[ \mathbf{n}_k^{\#} =^{w}_{w_{\{1\}}} n_k^{\#_{w_{\{1\}}}}(x) \right]. \tag{4.1.16}$$

Note that for any $w_{\{1\}}$-inconsistent hyperintegers $\mathbf{n}_k^{\#}$ and $\mathbf{m}_k^{\#}$ the following property holds

$$\mathbf{n}_k^{\#} =_{w_{\{1\}}} \mathbf{m}_k^{\#} \Leftrightarrow \forall x \left[ \left( \mathbf{n}_k^{\#} =^{w}_{w_{\{1\}}} n_k^{\#_{w_{\{1\}}}}(x) \right) \Leftrightarrow \left( \mathbf{m}_k^{\#} =^{w}_{w_{\{1\}}} n_k^{\#_{w_{\{1\}}}}(x) \right) \right]. \tag{4.1.17}$$

**Notation 4.1.1**. We often abbreviate for short

$$\mathbf{n}_{k,x}^{\#} =^{w}_{\#_{w_{\{1\}}}} n_k^{\#_{w_{\{1\}}}}(x), \tag{4.1.18}$$

where $x \in_{w_{\{1\}}} \mathbf{T}_{\#}^1$, instead (4.1.16).

**Definition 4.1.1**. For any $w_{\{1\}}$-inconsistent hyperinteger $\mathbf{n}_k^{\#}$ we define a $w_{\{1\}}$-set $\mathbf{Val}(\mathbf{n}_k^{\#})$
by

$$\forall x \left[ n_k^{\#_{w_{\{1\}}}}(x) \in_{w_{\{1\}}} \mathbf{Val}(\mathbf{n}_k^{\#}) \Leftrightarrow \mathbf{n}_{k,x}^{\#} =^{w}_{\#_{w_{\{1\}}}} n_k^{\#_{w_{\{1\}}}}(x) \right]. \tag{4.1.19}$$

Note that for any $w_{\{1\}}$-inconsistent hyperintegers $\mathbf{n}_k^{\#}$ and $\mathbf{m}_k^{\#}$ the following property holds

$$\mathbf{n}_k^{\#} =_{w_{\{1\}}} \mathbf{m}_k^{\#} \Leftrightarrow \mathbf{Val}(\mathbf{n}_k^{\#}) =_{w_{\{1\}}} \mathbf{Val}(\mathbf{m}_k^{\#}). \tag{4.1.20}$$

**Definition 4.1.2**.For any $w_{\{1\}}$-inconsistent hyperintegers $\mathbf{n}_k^{\#}$ and $\mathbf{m}_k^{\#}$ we define now the relation $\mathbf{n}_k^{\#} \subseteqq_{w_{\{1\}}} \mathbf{m}_k^{\#}$ :

$$\mathbf{n}_k^{\#} \subseteqq_{w_{\{1\}}} \mathbf{m}_k^{\#} \Leftrightarrow \mathbf{Val}(\mathbf{n}_k^{\#}) \subseteqq_{w_{\{1\}}} \mathbf{Val}(\mathbf{m}_k^{\#}). \tag{4.1.21}$$

**Remark 4.1.6**.(i) The vector $w_{\{1\}}$-addition $\mathbf{n}_k^{\#} +_{\#_{w_{\{1\}}}} \mathbf{m}_k^{\#}$ of $w_{\{1\}}$-inconsistent hyperintegers
$\mathbf{n}_k^{\#}$ and $\mathbf{m}_k^{\#}$ is defined by

$$\bigwedge_{x \in_{w_{\{1\}}} \mathbf{T}_{\#}^1} \left[ \mathbf{n}_k^{\#} +_{w_{\{1\}}} \mathbf{m}_k^{\#} =_{w_{\{1\}}}^{w} n_k^{\#_{w_{\{1\}}}}(x) +_{w_{\{1\}}} m_k^{\#_{w_{\{1\}}}}(x) \right] \quad (4.1.22)$$

or

$$\forall (x \in_{w_{\{1\}}} \mathbf{T}_{\#}^1) \left[ \mathbf{n}_{k,x}^{\#} +_{w_{\{1\}}} \mathbf{m}_{k,x}^{\#} =_{w_{\{1\}}}^{w} n_k^{\#_{w_{\{1\}}}}(x) +_{w_{\{1\}}} m_k^{\#_{w_{\{1\}}}}(x) \right] \quad (4.1.23)$$

(ii) The vector $w_{\{1\}}$-multiplication on scalar $\alpha$ takes any scalar $\alpha \in_{\#_{w_{\{1\}}}} \mathbb{R}^{\#_{w_{\{1\}}}}$ or $\alpha \in_{w_{\{1\}}} \mathbb{C}^{\#_{w_{\{1\}}}}$ and any $w_{\{1\}}$-inconsistent hyperinteger $\mathbf{n}_k^{\#}$ and gives $w_{\{1\}}$-inconsistent hyperreal number or nonstandard complex number defined by

$$\bigwedge_{x \in_{w_{\{1\}}} \mathbf{T}_{\#}^1} \left[ \alpha \times_{w_{\{1\}}} \mathbf{n}_k^{\#} =_{\#_{w_{\{1\}}}}^{w} \alpha \times_{w_{\{1\}}} n_k^{\#_{w_{\{1\}}}}(x) \right]. \quad (4.1.24)$$

or

$$\forall x (x \in_{w_{\{1\}}} \mathbf{T}_{\#}^1) \left[ \mathbf{n}_{k,x}^{\#} \times_{w_{\{1\}}} \mathbf{m}_{k,x}^{\#} =_{w_{\{1\}}}^{w} n_k^{\#_{w_{\{1\}}}}(x) \times_{w_{\{1\}}} m_k^{\#_{w_{\{1\}}}}(x) \right] \quad (4.1.25)$$

(iv) Note that the following properties holds:
(a) $w_{\{1\}}$-associativity of vector $w_{\{1\}}$-addition:
$(\mathbf{n}_k^{\#} +_{w_{\{1\}}} \mathbf{m}_k^{\#}) +_{w_{\{1\}}} \mathbf{k}_k^{\#} =_{w_{\{1\}}} \mathbf{n}_k^{\#} +_{w_{\{1\}}} (\mathbf{m}_k^{\#} +_{w_{\{1\}}} \mathbf{k}_k^{\#}),$
(b) $w_{\{1\}}$-associativity of vector $w_{\{1\}}$-multiplication:
$\left(\mathbf{n}_k^{\#} \times_{\#_{w_{\{1\}}}} \mathbf{m}_k^{\#}\right) \times_{w_{\{1\}}} \mathbf{m}_k^{\#} =_{\#_{w_{\{1\}}}} \mathbf{n}_k^{\#} \times_{\#_{w_{\{1\}}}} \left(\mathbf{m}_k^{\#} \times_{\#_{w_{\{1\}}}} \mathbf{k}_k^{\#}\right),$
(c) $w_{\{1\}}$-commutativity of vector $w_{\{1\}}$-addition:
$\mathbf{n}_k^{\#} +_{w_{\{1\}}} \mathbf{m}_k^{\#} =_{w_{\{1\}}} \mathbf{m}_k^{\#} +_{w_{\{1\}}} \mathbf{n}_k^{\#},$
(d) $w_{\{1\}}$-commutativity of vector $w_{\{1\}}$-multiplication:
$\mathbf{n}_k^{\#} \times_{w_{\{1\}}} \mathbf{m}_k^{\#} =_{w_{\{1\}}} \mathbf{m}_k^{\#} \times_{w_{\{1\}}} \mathbf{n}_k^{\#},$
(c) inverse elements of vector $w_{\{1\}}$-addition: $\mathbf{n}_k^{\#} +_{w_{\{1\}}} (-_{w_{\{1\}}} \mathbf{n}_k^{\#}) =_{w_{\{1\}}} \mathbf{0}_{w_{\{1\}}}$
(d) compatibility of vector $w_{\{1\}}$-multiplication on scalars $\alpha, \beta \in_{w_{\{1\}}} \mathbb{C}^{\#_{w_{\{1\}}}}$ with multiplication
in field $\mathbb{C}^{\#_{w_{\{1\}}}}$:
$(\alpha \times_{w_{\{1\}}} \mathbf{n}_k^{\#}) \times_{w_{\{1\}}} (\beta \times_{w_{\{1\}}} \mathbf{m}_k^{\#}) =_{w_{\{1\}}} (\alpha \times_{w_{\{1\}}} \beta) \times_{w_{\{1\}}} (\mathbf{n} \times_{w_{\{1\}}} \mathbf{m}_k^{\#}).$

**Definition 4.1.6**. Let $\mathfrak{I}^{\#}$ be a $w_{\{1\}}$-set of the all $w_{\{1\}}$-inconsistent hyperintegers $\mathbf{n}_k^{\#}$ with binary operations $+_{w_{\{1\}}}, \times_{w_{\{1\}}}$, etc. defined above.The tuple $\{\mathfrak{I}^{\#}, +_{w_{\{1\}}}, \times_{w_{\{1\}}}\}$ is an inconsistent $\mathbb{C}^{\#_{w_{\{1\}}}}$-algebra and we will be denoted this algebra by $\mathfrak{I}^{\#}$.

**Remark 4.1.7**. We introduce now $w_{\{1\}}$-inconsistent complex nonstandard numbers $\mathbf{c}_k^{\#}$ corresponding to trigonometric series (4.1.5) by the following way

$$\bigwedge_{x \in_{w_{\{1\}}} \mathbf{T}_{\#}^1} \left[ \mathbf{c}_{\mathbf{n}_k^{\#}}^{\#} =_{w_{\{1\}}}^{w} c_{n_k^{\#_{w_{\{1\}}}}(x)}^{\#_{w_{\{1\}}}} \right]. \quad (4.1.26)$$

**Notation 4.1.2**. We often abbreviate for short $\mathbf{c}_k^{\#}$ instead notation $\mathbf{c}_{\mathbf{n}_k^{\#}}^{\#}$, i.e.

$$\mathbf{c}_k^{\#} \triangleq \mathbf{c}_{\mathbf{n}_k^{\#}}^{\#}. \quad (4.1.27)$$

**Remark 4.1.8**. Note that for any $w_{\{1\}}$-inconsistent numbers $\mathbf{c}_{k_1}^{\#}$ and $\mathbf{c}_{k_2}^{\#}$ the following

property holds

$$\mathbf{c}^{\#}_{k_1} =_{w_{\{1\}}} \mathbf{c}^{\#}_{k_2} \Leftrightarrow \forall x\left[\left(\mathbf{c}^{\#}_{k_1} =^{w}_{w_{\{1\}}} n^{\#_{w_{\{1\}}}}_{k_1}(x)\right) \Leftrightarrow \left(\mathbf{c}^{\#}_{k_2} =^{w}_{w_{\{1\}}} n^{\#_{w_{\{1\}}}}_{k_2}(x)\right)\right]. \quad (4.1.28)$$

**Notation 4.1.3**. We often abbreviate for short

$$\mathbf{c}^{\#}_{k,x} =^{w}_{\#_{w_{\{1\}}}} c^{\#_{w_{\{1\}}}}_{k}(x), \quad (4.1.29)$$

where $x \in_{w_{\{1\}}} \mathbf{T}^1_{\#}$, instead (4.1.28).

**Definition 4.1.7**. For any $w_{\{1\}}$-inconsistent number $\mathbf{c}^{\#}_{k}$ we define a $w_{\{1\}}$-set $\mathbf{Val}(\mathbf{c}^{\#}_{k})$ by

$$\forall x\left[c^{\#_{w_{\{1\}}}}_{k}(x) \in_{w_{\{1\}}} \mathbf{Val}(\mathbf{c}^{\#}_{k}) \Leftrightarrow \mathbf{c}^{\#}_{k,x} =^{w}_{\#_{w_{\{1\}}}} c^{\#_{w_{\{1\}}}}_{k}(x)\right]. \quad (4.1.30)$$

Note that for any $w_{\{1\}}$-inconsistent numbers $\mathbf{c}^{\#}_{k_1}$ and $\mathbf{c}^{\#}_{k_2}$ the following property holds

$$\mathbf{c}^{\#}_{k_1} =_{w_{\{1\}}} \mathbf{c}^{\#}_{k_2} \Leftrightarrow \mathbf{Val}(\mathbf{c}^{\#}_{k_1}) =_{w_{\{1\}}} \mathbf{Val}(\mathbf{c}^{\#}_{k_2}). \quad (4.1.31)$$

**Remark 4.1.9**.(i) The vector $w_{\{1\}}$-addition $\mathbf{c}^{\#}_{k_1} +_{\#_{w_{\{1\}}}} \mathbf{c}^{\#}_{k_2}$ of $w_{\{1\}}$-inconsistent numbers $\mathbf{c}^{\#}_{k_1}$ and $\mathbf{c}^{\#}_{k_2}$ is defined by

$$\bigwedge_{x \in_{w_{\{1\}}} \mathbf{T}^1_{\#}} \left[\mathbf{c}^{\#}_{k_1} +_{\#_{w_{\{1\}}}} \mathbf{c}^{\#}_{k_2} =^{w}_{w_{\{1\}}} c^{\#_{w_{\{1\}}}}_{k_1}(x) +_{\#_{w_{\{1\}}}} c^{\#_{w_{\{1\}}}}_{k_2}(x)\right] \quad (4.1.32)$$

or

$$\forall x(x \in_{w_{\{1\}}} \mathbf{T}^1_{\#})\left[\mathbf{c}^{\#}_{k_1,x} +_{\#_{w_{\{1\}}}} \mathbf{c}^{\#}_{k_2,x} =^{w}_{w_{\{1\}}} c^{\#_{w_{\{1\}}}}_{k_1}(x) +_{\#_{w_{\{1\}}}} c^{\#_{w_{\{1\}}}}_{k_2}(x)\right] \quad (4.1.33)$$

(ii) The mixed $w_{\{1\}}$-addition $\mathbf{c}^{\#}_{k_1} \oplus_{\#_{w_{\{1\}}}} \mathbf{c}^{\#}_{k_2}$ of $w_{\{1\}}$-inconsistent numbers $\mathbf{c}^{\#}_{k_1}$ and $\mathbf{c}^{\#}_{k_2}$ is defined by

$$\bigwedge_{x,y \in_{w_{\{1\}}} \mathbf{T}^1_{\#}} \left[\mathbf{c}^{\#}_{k_1} \oplus_{\#_{w_{\{1\}}}} \mathbf{c}^{\#}_{k_2} =^{w}_{w_{\{1\}}} c^{\#_{w_{\{1\}}}}_{k_1}(x) +_{\#_{w_{\{1\}}}} c^{\#_{w_{\{1\}}}}_{k_2}(y)\right] \quad (4.1.34)$$

or

$$\forall x \forall y(x,y \in_{w_{\{1\}}} \mathbf{T}^1_{\#})\left[\mathbf{c}^{\#}_{k_1} \oplus_{\#_{w_{\{1\}}}} \mathbf{c}^{\#}_{k_2} =^{w}_{w_{\{1\}}} c^{\#_{w_{\{1\}}}}_{k_1}(x) +_{\#_{w_{\{1\}}}} c^{\#_{w_{\{1\}}}}_{k_2}(y)\right] \quad (4.1.35)$$

(iii) Note that the following properties holds:

(a) $w_{\{1\}}$-associativity of vector $w_{\{1\}}$-addition:

$\left(\mathbf{c}^{\#}_{k_1} +_{\#_{w_{\{1\}}}} \mathbf{c}^{\#}_{k_2}\right) +_{\#_{w_{\{1\}}}} \mathbf{c}^{\#}_{k_3} =_{\#_{w_{\{1\}}}} \mathbf{c}^{\#}_{k_1} +_{\#_{w_{\{1\}}}} \left(\mathbf{c}^{\#}_{k_2} +_{\#_{w_{\{1\}}}} \mathbf{c}^{\#}_{k_3}\right),$

(b) $w_{\{1\}}$-associativity of vector $w_{\{1\}}$-multiplication:

$\left(\mathbf{c}^{\#}_{k_1} \times_{\#_{w_{\{1\}}}} \mathbf{c}^{\#}_{k_2}\right) \times_{\#_{w_{\{1\}}}} \mathbf{c}^{\#}_{k,3} =_{\#_{w_{\{1\}}}} \mathbf{c}^{\#}_{k_1} \times_{\#_{w_{\{1\}}}} \left(\mathbf{c}^{\#}_{k_2} \times_{\#_{w_{\{1\}}}} \mathbf{c}^{\#}_{k,3}\right),$

(c) $w_{\{1\}}$-commutativity of vector $w_{\{1\}}$-addition:

$\mathbf{c}^{\#}_{k_1} +_{\#_{w_{\{1\}}}} \mathbf{c}^{\#}_{k_2} =_{\#_{w_{\{1\}}}} \mathbf{c}^{\#}_{k_2} +_{\#_{w_{\{1\}}}} \mathbf{c}^{\#}_{k_1},$

(d) $w_{\{1\}}$-commutativity of vector $w_{\{1\}}$-multiplication:

$\mathbf{c}^{\#}_{k_1} \times_{\#_{w_{\{1\}}}} \mathbf{c}^{\#}_{k_2} =_{\#_{w_{\{1\}}}} \mathbf{c}^{\#}_{k_2} \times_{\#_{w_{\{1\}}}} \mathbf{c}^{\#}_{k_1},$

(e) $w_{\{1\}}$-associativity of mixed $w_{\{1\}}$-addition:

$\left(\mathbf{c}^{\#}_{k_1} \oplus_{\#_{w_{\{1\}}}} \mathbf{c}^{\#}_{k_2}\right) \oplus_{\#_{w_{\{1\}}}} \mathbf{c}^{\#}_{k_3} =_{\#_{w_{\{1\}}}} \mathbf{c}^{\#}_{k_1} \oplus_{\#_{w_{\{1\}}}} \left(\mathbf{c}^{\#}_{k_2} \oplus_{\#_{w_{\{1\}}}} \mathbf{c}^{\#}_{k_3}\right),$

(f) $w_{\{1\}}$-commutativity of mixed $w_{\{1\}}$-addition:

$\mathbf{c}^{\#}_{k_1} \oplus_{\#_{w_{\{1\}}}} \mathbf{c}^{\#}_{k_2} =_{\#_{w_{\{1\}}}} \mathbf{c}^{\#}_{k_1} \oplus_{\#_{w_{\{1\}}}} \mathbf{c}^{\#}_{k_2}.$

**Definition 4.1.3**.Let $\{z_k(x)\}_{k \in \mathbb{N}}$ be any sequence of functions $\{z_k(x)\}_{k \in \mathbb{N}}$ such that $z_k : \mathbf{T}^1 \to \mathbb{C}$.Assume that a $\#_{w_{\{1\}}}$-sequence

$$\left\{z_k^{\#_{w_{\{1\}}}}(x)\right\}_{k\in_{\#_{w_{\{1\}}}}\mathbb{N}^{\#_{w_{\{1\}}}}} \quad (4.1.36)$$

is obtained by using $\#_{w_{\{1\}}}$- transfer from the standard sequence $\{z_k(x)\}_{k\in\mathbb{N}}$, i.e.

$$\left\{z_k^{\#_{w_{\{1\}}}}(x)\right\}_{k\in_{w_{\{1\}}}\mathbb{N}^{\#_{w_{\{1\}}}}} =_{w_{\{1\}}} \left(\{z_k(x)\}_{k\in\mathbb{N}}\right)^{\#_{w_{\{1\}}}}. \quad (4.1.37)$$

We introduce now $w_{\{1\}}$-inconsistent nonstandard complex numbers $\mathbf{z}_k^{\#}$ corresponding to sequence (4.1.37) by the following way

$$\bigwedge_{x\in_{w_{\{1\}}}\mathbf{T}_{\#}^{1}}\left[\mathbf{z}_k^{\#} =_{w_{\{1\}}}^{w} z_k^{\#_{w_{\{1\}}}}(x)\right]. \quad (4.1.38)$$

**Remark 4.1.8**. Note that for any $w_{\{1\}}$-inconsistent numbers $\mathbf{c}_{k_1}^{\#}$ and $\mathbf{c}_{k_2}^{\#}$ the following property holds

$$\mathbf{z}_{k_1}^{\#} =_{w_{\{1\}}} \mathbf{z}_{k_2}^{\#} \Leftrightarrow \forall x\left[\left(\mathbf{z}_{k_1}^{\#} =_{w_{\{1\}}}^{w} z_{k_1}^{\#_{w_{\{1\}}}}(x)\right) \Leftrightarrow \left(\mathbf{z}_{k_2}^{\#} =_{w_{\{1\}}}^{w} z_{k_2}^{\#_{w_{\{1\}}}}(x)\right)\right]. \quad (4.1.39)$$

**Notation 4.1.4**. We often abbreviate for short

$$\mathbf{z}_{k,x}^{\#} =_{\#_{w_{\{1\}}}}^{w} z_k^{\#_{w_{\{1\}}}}(x), \quad (4.1.40)$$

where $x \in_{w_{\{1\}}} \mathbf{T}_{\#}^{1}$, instead (4.1.38).

**Definition 4.1.4**. For any $w_{\{1\}}$-inconsistent number $\mathbf{z}_k^{\#}$ we define a $w_{\{1\}}$-set $\mathbf{Val}(\mathbf{z}_k^{\#})$ by

$$\forall x\left[z_k^{\#_{w_{\{1\}}}}(x) \in_{w_{\{1\}}} \mathbf{Val}(\mathbf{z}_k^{\#}) \Leftrightarrow \mathbf{z}_{k,x}^{\#} =_{\#_{w_{\{1\}}}}^{w} z_k^{\#_{w_{\{1\}}}}(x)\right]. \quad (4.1.41)$$

Note that for any $w_{\{1\}}$-inconsistent numbers $\mathbf{z}_{k_1}^{\#}$ and $\mathbf{z}_{k_2}^{\#}$ the following property holds

$$\mathbf{z}_{k_1}^{\#} =_{w_{\{1\}}} \mathbf{z}_{k_2}^{\#} \Leftrightarrow \mathbf{Val}(\mathbf{z}_{k_1}^{\#}) =_{w_{\{1\}}} \mathbf{Val}(\mathbf{z}_{k_2}^{\#}). \quad (4.1.42)$$

**Definition 4.1.5**.(i) The vector $w_{\{1\}}$-addition $\mathbf{z}_{k_1}^{\#} +_{\#_{w_{\{1\}}}} \mathbf{z}_{k_2}^{\#}$ of $w_{\{1\}}$-inconsistent numbers $\mathbf{z}_{k_1}^{\#}$ and $\mathbf{c}_{k_2}^{\#}$ is defined by

$$\bigwedge_{x\in_{w_{\{1\}}}\mathbf{T}_{\#}^{1}}\left[\mathbf{z}_{k_1}^{\#} +_{\#_{w_{\{1\}}}} \mathbf{z}_{k_2}^{\#} =_{w_{\{1\}}}^{w} z_{k_1}^{\#_{w_{\{1\}}}}(x) +_{\#_{w_{\{1\}}}} z_{k_2}^{\#_{w_{\{1\}}}}(x)\right] \quad (4.1.43)$$

or

$$\forall x(x \in_{w_{\{1\}}} \mathbf{T}_{\#}^{1})\left[\mathbf{z}_{k_1,x}^{\#} +_{\#_{w_{\{1\}}}} \mathbf{z}_{k_2,x}^{\#} =_{w_{\{1\}}}^{w} z_{k_1}^{\#_{w_{\{1\}}}}(x) +_{\#_{w_{\{1\}}}} z_{k_2}^{\#_{w_{\{1\}}}}(x)\right] \quad (4.1.44)$$

(ii) The mixed $w_{\{1\}}$-addition $\mathbf{z}_{k_1}^{\#} \oplus_{\#_{w_{\{1\}}}} \mathbf{z}_{k_2}^{\#}$ of $w_{\{1\}}$-inconsistent numbers $\mathbf{z}_{k_1}^{\#}$ and $\mathbf{z}_{k_2}^{\#}$ is defined by

$$\bigwedge_{x,y\in_{w_{\{1\}}}\mathbf{T}_{\#}^{1}}\left[\mathbf{z}_{k_1}^{\#} \oplus_{\#_{w_{\{1\}}}} \mathbf{z}_{k_2}^{\#} =_{w_{\{1\}}}^{w} z_{k_1}^{\#_{w_{\{1\}}}}(x) +_{\#_{w_{\{1\}}}} z_{k_2}^{\#_{w_{\{1\}}}}(y)\right] \quad (4.1.45)$$

or

$$\forall x\forall y(x,y \in_{w_{\{1\}}} \mathbf{T}_{\#}^{1})\left[\mathbf{z}_{k_1}^{\#} \oplus_{\#_{w_{\{1\}}}}^{\#} \mathbf{z}_{k_2}^{\#} =_{w_{\{1\}}}^{w} z_{k_1}^{\#_{w_{\{1\}}}}(x) +_{\#_{w_{\{1\}}}} z_{k_2}^{\#_{w_{\{1\}}}}(y)\right] \quad (4.1.46)$$

(iii) Note that the following properties holds:

(a) $w_{\{1\}}$-associativity of vector $w_{\{1\}}$-addition:

$\left(\mathbf{z}^{\#}_{k_1} +_{\#_{w_{\{1\}}}} \mathbf{z}^{\#}_{k_2}\right) +_{\#_{w_{\{1\}}}} \mathbf{z}^{\#}_{k_3} =_{\#_{w_{\{1\}}}} \mathbf{z}^{\#}_{k_1} +_{\#_{w_{\{1\}}}} \left(\mathbf{z}^{\#}_{k_2} +_{\#_{w_{\{1\}}}} \mathbf{z}^{\#}_{k_3}\right),$

(b) $w_{\{1\}}$-associativity of vector $w_{\{1\}}$-multiplication:

$\left(\mathbf{z}^{\#}_{k_1} \times_{w_{\{1\}}} \mathbf{z}^{\#}_{k_2}\right) \times_{w_{\{1\}}} \mathbf{z}^{\#}_{k,3} =_{w_{\{1\}}} \mathbf{z}^{\#}_{k_1} \times_{\#_{w_{\{1\}}}} \left(\mathbf{z}^{\#}_{k_2} \times_{w_{\{1\}}} \mathbf{z}^{\#}_{k,3}\right),$

(c) $w_{\{1\}}$-commutativity of vector $w_{\{1\}}$-addition:

$\mathbf{z}^{\#}_{k_1} +_{w_{\{1\}}} \mathbf{z}^{\#}_{k_2} =_{w_{\{1\}}} \mathbf{z}^{\#}_{k_2} +_{w_{\{1\}}} \mathbf{z}^{\#}_{k_1},$

(d) $w_{\{1\}}$-commutativity of vector $w_{\{1\}}$-multiplication:

$\mathbf{z}^{\#}_{k_1} \times_{w_{\{1\}}} \mathbf{z}^{\#}_{k_2} =_{w_{\{1\}}} \mathbf{z}^{\#}_{k_2} \times_{w_{\{1\}}} \mathbf{z}^{\#}_{k_1},$

(e) $w_{\{1\}}$-associativity of mixed $w_{\{1\}}$-addition:

$\left(\mathbf{z}^{\#}_{k_1} \oplus_{\#_{w_{\{1\}}}} \mathbf{z}^{\#}_{k_2}\right) \oplus_{\#_{w_{\{1\}}}} \mathbf{z}^{\#}_{k_3} =_{\#_{w_{\{1\}}}} \mathbf{z}^{\#}_{k_1} \oplus_{\#_{w_{\{1\}}}} \left(\mathbf{z}^{\#}_{k_2} \oplus_{\#_{w_{\{1\}}}} \mathbf{z}^{\#}_{k_3}\right),$

(f) $w_{\{1\}}$-commutativity of mixed $w_{\{1\}}$-addition:

$\mathbf{z}^{\#}_{k_1} \oplus_{w_{\{1\}}} \mathbf{z}^{\#}_{k_2} =_{w_{\{1\}}} \mathbf{z}^{\#}_{k_1} \oplus_{w_{\{1\}}} \mathbf{z}^{\#}_{k_2}.$

**Definition 4.1.7.**

**Definition 4.1.8.**of Let $\{\mathbf{z}^{\#}_k\}_{k \in_{w_{\{1\}}} \mathbb{N}^{\#_{w_{\{1\}}}}}$ be a $w_{\{1\}}$-sequence of $w_{\{1\}}$-inconsistent numbers

$\mathbf{z}^{\#}_k, k \in_{w_{\{1\}}} \mathbb{N}^{\#_{w_{\{1\}}}}, \mathbf{m} \in_{w_{\{1\}} w_{\{1\}}} \mathbb{N}^{\#_{w_{\{1\}}}}$. External vector $w_{\{1\}}$-summation of the sequence $\{\mathbf{z}^{\#}_k\}_{k =_{w_{\{1\}}} 0_{w_{\{1\}}}}$ may be defined recursively by using external induction principle as follows:

if $\mathbf{m}$ is any $w_{\{1\}}$-hyperinteger, then the recursion schemata reads

$$\#_{w_{\{1\}}}\text{-}\mathbf{Ext} \sum_{k =_{w_{\{1\}}} 0_{w_{\{1\}}}}^{0_{w_{\{1\}}}} \mathbf{z}^{\#}_k =_{w_{\{1\}}} \mathbf{z}^{\#}_{0_{w_{\{1\}}}},$$

$$\#_{w_{\{1\}}}\text{-}\mathbf{Ext} \sum_{k =_{w_{\{1\}}} 0_{w_{\{1\}}}}^{\mathbf{m}} \mathbf{z}^{\#}_k =_{w_{\{1\}}} \#_{w_{\{1\}}}\text{-}\mathbf{Ext} \sum_{k =_{w_{\{1\}}} 0_{w_{\{1\}}}}^{\mathbf{m} -_{w_{\{1\}}} 1_{w_{\{1\}}}} \mathbf{z}^{\#}_k +_{w_{\{1\}}} \mathbf{z}^{\#}_{\mathbf{m}}. \qquad (4.1.47)$$

**Propozition 4.1.1.**

**Definition 4.1.7.**

**Definition 4.1.8.**

$$\#_{w_{\{1\}}}\text{-} \oplus_{w_{\{1\}}} \text{-}\mathbf{Ext} \sum_{k =_{w_{\{1\}}} 0_{w_{\{1\}}}}^{M} \mathbf{c}^{\#}_{k,i} \quad (4.1.)$$

$$\#_{w_{\{1\}}} - \mathbf{Ext} \bigoplus_{k =_{w_{\{1\}}} 0_{w_{\{1\}}}}^{M} \mathbf{c}^{\#}_{k,i} \quad (4.1.)$$

**Definition 4.1.7.**

**Definition 4.1.8.** We define now a function $\mathbf{Exp}(u, \mathbf{n}^{\#}_k) : \mathbf{T}^1_{\#} \times_{w_{\{1\}}} \mathfrak{I}^{\#} \to \mathfrak{R}^{\#}$ by

$$\exp(iu \times_{w_{\{1\}}} \mathbf{n}^{\#}_k) =^{w}_{w_{\{1\}}} \exp(iu \times_{w_{\{1\}}} \mathbf{n}^{\#}_{k,x}) \qquad (4.1.)$$

where $u,x \in_{w_{\{1\}}} \mathbf{T}^1_{\#}, \mathbf{n}_k^{\#} \in_{\#_{w_{\{1\}}}} \Im^{\#}$ and $k \in_{w_{\{1\}}} \mathbb{N}^{\#_{w_{\{1\}}}}$.

**Remark 3.3.5**. We introduce now a $w_{\{1\}}$-function $\Psi_M(x)$ by the following way

$$\Psi_M(x) =_{\#_{w_{\{1\}}}} \#_{w_{\{1\}}}\textbf{-Ext}\sum_{k=_{w_{\{1\}}} 0_{w_{\{1\}}}}^{M} c_{\mathbf{n}_k^{\#}} \times \exp\left(ix \times_{\#_{w_{\{1\}}}} \mathbf{n}_k^{\#}\right). \quad (3.3.16)$$

From Eq.(3.3.8) and Eq.(3.3.10) we obtain

$$\Psi(x) =_{\#_{w_{\{1\}}}} \#_{w_{\{1\}}}\textbf{-Int}\sum_{k=_{w_{\{1\}}} 0_{w_{\{1\}}}}^{M} c_{\mathbf{n}_k^{\#}} \exp\left(ix \times_{\#_{w_{\{1\}}}} \mathbf{n}_k^{\#}\right) =_{\#_{w_{\{2\}}}} N_M(x). \quad (3.3.12)$$

By $\#_{w_{\{1\}}}$- integration From Eq.(3.3.11) we obtain

$$\Delta(M) =_{\#_{w_{\{1\}}}} \#_{w_{\{1\}}}\textbf{-Ext}\int_{\mathbf{T}^1_{\#_{w_{\{1\}}}}} \Psi_M^2(x) d^{\#_{w_{\{1\}}}}\mu(x) =_{\#_{w_{\{1\}}}}$$

$$\#_{w_{\{1\}}}\textbf{-Ext}\int_{\mathbf{T}^1_{\#_{w_{\{1\}}}}} \left[\#_{w_{\{1\}}}\textbf{-Int}\sum_{k=_{w_{\{1\}}} 0_{w_{\{1\}}}}^{M} c_{\mathbf{n}_k^{\#}} \exp\left(ix \times_{\#_{w_{\{1\}}}} \mathbf{n}_k^{\#}\right)\right]^2 d^{\#_{w_{\{1\}}}}\mu(x) =_{\#_{w_{\{1\}}}} \quad (3.3.13)$$

$$=_{\#_{w_{\{1\}}}} \#_{w_{\{1\}}}\textbf{-Int}\sum_{k=_{w_{\{0\}}} 0_{w_{\{0\}}}}^{M} \left|c_{\mathbf{n}_k^{\#}}\right|^2,$$

where by $d^{\#_{w_{\{1\}}}}\mu(x)$ we denote $\#_{w_{\{1\}}}$- transfered standard Lebesgue measure $d\mu(x)$ on $\mathbf{T}^1$. From (3.3.2) by $\#_{w_{\{1\}}}$- transfer it follows that $M \in_{w_{\{1\}}} \mathbb{N}^{\#_{w_{\{1\}}}} \backslash \mathbb{N}$

$$\Delta(M) =_{\#_{w_{\{1\}}}} \#_{w_{\{1\}}}\textbf{-Ext}\sum_{k=_{w_{\{1\}}} 0_{w_{\{1\}}}}^{M} \left|c_{\mathbf{n}_k^{\#}}\right|^2 \in_{w_{\{1\}}} \mathbb{R}_{\mathbf{fin}}^{\#_{w_{\{1\}}}}, \quad (3.3.14)$$

i.e. the quantity $\Delta(M)$ alwais is $\#_{w_{\{1\}}}$- finite i.e.

$$\neg_s\left[\Delta(M) \in_{w_{\{1\}}} \mathbb{R}^{\#_{w_{\{1\}}}} \backslash_{w_{\{1\}}} \mathbb{R}_{\mathbf{fin}}^{\#_{w_{\{1\}}}}\right]. \quad (3.3.15)$$

From RHS of the Eq.(3.3.12) By $\#_{w_{\{1\}}}$- integration we obtain

$$\#_{w_{\{1\}}}\textbf{-Ext}\int_{\mathbf{T}^1_{\#_{w_{\{1\}}}}} N_M^2(x) d^{\#_{w_{\{1\}}}}\mu(x) =_{\#_{w_{\{2\}}}}$$

$$=_{\#_{w_{\{2\}}}} \int_{\mathbf{T}^1_{\#_{w_{\{1\}}}}} \left[\#_{w_{\{1\}}}\textbf{-Ext}\sum_{k=_{w_{\{1\}}} 0_{w_{\{1\}}}}^{M} c_{\mathbf{n}_k^{\#}} \exp\left(ix \times_{\#_{w_{\{1\}}}} \mathbf{n}_k^{\#}\right)\right]^2 d^{\#_{w_{\{1\}}}}\mu(x) =_{\#_{w_{\{1\}}}} \quad (3.3.16)$$

$$=_{\#_{w_{\{1\}}}} \#_{w_{\{1\}}}\textbf{-Ext}\sum_{k=_{w_{\{1\}}} 0_{w_{\{1\}}}}^{M} \left|c_{\mathbf{n}_k^{\#}}\right|^2 =_{\#_{w_{\{1\}}}} \Delta(M)$$

But on other hand from Eq.(3.3.8) By $\#_{w_{\{1\}}}$- integration we obtain

$$\int_{\mathbf{T}^1_{\#_{w_{\{1\}}}}}^{\#_{w_{\{1\}}}} N_M(x) d^{\#_{w_{\{1\}}}}\mu(x) =_{\#_{w_{\{1\}}}}$$

$$=_{\#_{w_{\{1\}}}} \int_{\mathbf{T}^1_{\#_{w_{\{1\}}}}}^{\#_{w_{\{1\}}}} \left[ \#_{w_{\{1\}}}\text{-}\mathbf{Int} \sum_{k=_{w_{\{0\}}} 0_w}^{M} c_{n_k(x)}^{\#_{w_{\{1\}}}} \exp\left( ix \times_{\#_{w_{\{1\}}}} n_k^{\#_{w_{\{1\}}}}(x) \right) \right] d^{\#_{w_{\{1\}}}}\mu(x) \quad (3.3.17)$$

$$\in_{w_{\{1\}}} \mathbb{R}^{\#_{w_{\{1\}}}} \backslash_{w_{\{1\}}} \mathbb{R}_{\mathbf{fin}}^{\#_{w_{\{1\}}}}.$$

Obviously by (3.3.15),Eq.(3.3.16) and Eq.(3.3.17) one obtains a contradiction.

(**II**) Now we go to prove that: a.e. on $\mathbf{T}^1$

$$\exists \lim_{n\to\infty} \sum_{n=0}^{n} c_n \exp(inx) < \infty. \quad (3.3.17)$$

It follows from (**I**) that

$$\exists \lim_{k\to\infty} \sum_{n_k=0}^{n_k} c_{n_k} \exp(in_k x) \Rightarrow \lim_{k\to\infty} \sum_{n_k=0}^{n_k} c_{n_k} \exp(in_k x) < \infty, \quad (3.3.18)$$

where $n_k \to \infty$ iff $k \to \infty$.
We assume now that: a.e. on $\mathbf{T}^1$

$$\neg\exists \lim_{n\to\infty} \sum_{n=0}^{n} c_n \exp(inx). \quad (3.3.19)$$

Let $z \in \mathbf{T}^1$be a real number such that

$$\neg\exists \lim_{n\to\infty} \sum_{n=0}^{n} c_n \exp(inz). \quad (3.3.20)$$

Notice that (3.3.20) meant that there exists countable sequence $\{n_k(z)\}_{k\in\mathbb{N}}$ such that

$$\neg\exists \lim_{k\to\infty} \sum_{n_k=0}^{n_k(z)} c_{n_k(z)} \exp(izn_k(z)). \quad (3.3.21)$$

Notice that (3.3.21) meant that there exists $\varepsilon(z) > 0$ and $N \in \mathbb{N}$ such that

$$\sum_{n_k(z)=n_{k_1}(z)}^{n_k(z)=n_{k_2}(z)} c_{n_k(z)} \exp(izn_k(z)) \geq \varepsilon(z), \quad (3.3.22)$$

where $n_{k_2}(z) > n_{k_1}(z) \geq N$.
From (3.3.22) by $\#_{w_{\{1\}}}$- transfer it follows that : $\#_{w_{\{1\}}}$- a.e. on $\mathbf{T}^1_{\#_{w_{\{1\}}}} \triangleq (\mathbf{T}^1)^{\#_{w_{\{1\}}}}$

$$\#_{w_{\{1\}}}\text{-}\mathbf{Ext} \sum_{n_{k_1}^{\#_{w_{\{1\}}}}(x)}^{n_{k_2}^{\#_{w_{\{1\}}}}(x)} c^{\#_{w_{\{1\}}}}_{n_k^{\#_{w_{\{1\}}}}(x)} \exp\left( ix \times_{\#_{w_{\{1\}}}} n_k^{\#_{w_{\{1\}}}}(x) \right)_{\#_{w_{\{1\}}}} \geq \varepsilon^{\#_{w_{\{1\}}}}(x), \quad (3.3.23)$$

where $n_{k_1}^{\#_{w_{\{1\}}}}(x) \in_{w_{\{1\}}} \mathbb{N}^{\#_{w_{\{1\}}}} \backslash \mathbb{N}$ and $n_{k_2}^{\#_{w_{\{1\}}}}(x) \in_{w_{\{1\}}} \mathbb{N}^{\#_{w_{\{1\}}}} \backslash \mathbb{N}$ and where a sequence

$$\left\{ n_k^{\#_{w_{\{1\}}}}(x) \right\}_{k \in_{\#_{w_{\{1\}}}} \mathbb{N}^{\#_{w_{\{1\}}}}}$$

is obtained by using $\#_{w_{\{1\}}}$- transfer from sequence $\{n_k(x)\}_{k\in\mathbb{N}}$,i.e.

$$\left\{ n_k^{\#_{w_{\{1\}}}}(x) \right\}_{k \in_{w_{\{1\}}} \mathbb{N}^{\#_{w_{\{1\}}}}} =_{\#_{w_{\{1\}}}} (\{n_k(x)\}_{k\in\mathbb{N}})^{\#_{w_{\{1\}}}}. \quad (3.3.24)$$

**Remark 3.3.5**. We introduce now $w_{\{1\}}$-inconsistent numbers by the following way

$$\bigwedge_{x \in_{w_{\{1\}}} \mathbf{T}^1_{\#}} \left[ \mathbf{n}^{\#}_k =_{\#_{w_{\{1\}}}} n_k^{\#_{w_{\{1\}}}}(x) \right]. \quad (3.3.25)$$

**Remark 3.3.6**. We introduce now a $w_{\{1\}}$-function $\Psi_M(x)$ by the following way

$$\Psi_{\mathbf{n}^{\#}_{k_1},\mathbf{n}^{\#}_{k_2}}(x) =_{\#_{w_{\{1\}}}} \#_{w_{\{1\}}}\textbf{-Ext-}\sum_{\mathbf{n}^{\#}_{k_1}}^{\mathbf{n}^{\#}_{k_2}} c^{\#_{w_{\{1\}}}}_{\mathbf{n}^{\#}_k} c_{\mathbf{n}^{\#}_k} \exp\left(ix \times_{\#_{w_{\{1\}}}} \mathbf{n}^{\#}_k\right). \quad (3.3.26)$$

From Eq.(3.3.23) and Eq.(3.3.25) we obtain

$$\Psi(x) =_{\#_{w_{\{1\}}}} \#_{w_{\{1\}}}\textbf{-Ext-}\sum_{\mathbf{n}^{\#}_{k_1}}^{\mathbf{n}^{\#}_{k_2}} c_{\mathbf{n}^{\#}_k} \exp\left(ix \times_{\#_{w_{\{1\}}}} \mathbf{n}^{\#}_k\right) =_{\#_{w_{\{2\}}}} \varepsilon^{\#_{w_{\{1\}}}}(x). \quad 3.3.12)$$

By $\#_{w_{\{1\}}}$- integration From Eq.(3.3.26) we obtain

$$\Delta(M) =_{\#_{w_{\{1\}}}} \#_{w_{\{1\}}}\textbf{-Ext-}\int_{\mathbf{T}^1_{\#_{w_{\{1\}}}}} \Psi^2_M(x) d^{\#_{w_{\{1\}}}}\mu(x) =_{\#_{w_{\{1\}}}}$$

$$\int_{\mathbf{T}^1_{\#_{w_{\{1\}}}}}^{\#_{w_{\{1\}}}} \left[ \#_{w_{\{1\}}}\textbf{-Int} \sum_{k=_{w_{\{1\}}} 0_{w_{\{1\}}}}^{M} c_{\mathbf{n}^{\#}_k} \exp\left(ix \times_{\#_{w_{\{1\}}}} \mathbf{n}^{\#}_k\right) \right]^2 d^{\#_{w_{\{1\}}}}\mu(x) =_{\#_{w_{\{1\}}}} \quad (3.3.13)$$

$$=_{\#_{w_{\{1\}}}} \#_{w_{\{1\}}}\textbf{-Int} \sum_{k=_{w_{\{0\}}} 0_{w_{\{0\}}}}^{M} \left| c_{\mathbf{n}^{\#}_k} \right|^2,$$

where by $d^{\#_{w_{\{1\}}}}\mu(x)$ we denote $\#_{w_{\{1\}}}$- transfered standard Lebesgue measure $d\mu(x)$ on $\mathbf{T}^1$. From (3.3.2) by $\#_{w_{\{1\}}}$- transfer it follows that $M \in_{w_{\{1\}}} \mathbb{N}^{\#_{w_{\{1\}}}} \backslash \mathbb{N}$

$$\Delta(M) =_{\#_{w_{\{1\}}}} \#_{w_{\{1\}}}\textbf{-Ext} \sum_{k=_{w_{\{1\}}} 0_{w_{\{1\}}}}^{M} \left| c_{\mathbf{n}^{\#}_k} \right|^2 \in_{w_{\{1\}}} \mathbb{R}^{\#_{w_{\{1\}}}}_{\mathbf{fin}}, \quad (3.3.14)$$

i.e. the quantity $\Delta(M)$ alwais is $\#_{w_{\{1\}}}$- finite i.e.

$$\neg_s \left[ \Delta(M) \in_{w_{\{1\}}} \mathbb{R}^{\#_{w_{\{1\}}}} \backslash_{w_{\{1\}}} \mathbb{R}^{\#_{w_{\{1\}}}}_{\mathbf{fin}} \right]. \quad (3.3.15)$$

From RHS of the Eq.(3.3.12) By $\#_{w_{\{1\}}}$- integration we obtain

$$\int_{\mathbf{T}^1_{\#_{w_{\{1\}}}}}^{\#_{w_{\{1\}}}} N^2_M(x) d^{\#_{w_{\{1\}}}}\mu(x) =_{\#_{w_{\{2\}}}}$$

$$=_{\#_{w_{\{2\}}}} \int_{\mathbf{T}^1_{\#_{w_{\{1\}}}}}^{\#_{w_{\{1\}}}} \left[ \#_{w_{\{1\}}}\textbf{-Int} \sum_{k=_{w_{\{1\}}} 0_{w_{\{1\}}}}^{M} c_{\mathbf{n}^{\#}_k} \exp\left(ix \times_{\#_{w_{\{1\}}}} \mathbf{n}^{\#}_k\right) \right]^2 d^{\#_{w_{\{1\}}}}\mu(x) =_{\#_{w_{\{1\}}}} \quad (3.3.16)$$

$$=_{\#_{w_{\{1\}}}} \#_{w_{\{1\}}}\textbf{-Int} \sum_{k=_{w_{\{1\}}} 0_{w_{\{1\}}}}^{M} \left| c_{\mathbf{n}^{\#}_k} \right|^2 =_{\#_{w_{\{1\}}}} \Delta(M)$$

But on other hand from Eq.(3.3.8) By $\#_{w_{\{1\}}}$- integration we obtain

$$\int_{\mathbf{T}^1_{\#_{w_{\{1\}}}}} {}^{\#_{w_{\{1\}}}} N_M(x) d^{\#_{w_{\{1\}}}}\mu(x) =_{\#_{w_{\{1\}}}}$$

$$=_{\#_{w_{\{1\}}}} \int_{\mathbf{T}^1_{\#_{w_{\{1\}}}}} {}^{\#_{w_{\{1\}}}} \left[ \#_{w_{\{1\}}}\text{-}\mathbf{Int} \sum_{k=_{w_{\{0\}}} 0_w}^{M} c^{\#_{w_{\{1\}}}}_{n_k(x)} \exp\left( ix \times_{\#_{w_{\{1\}}}} n_k^{\#_{w_{\{1\}}}}(x) \right) \right] d^{\#_{w_{\{1\}}}}\mu(x) \qquad (3.3.17)$$

$$\in_{w_{\{1\}}} \mathbb{R}^{\#_{w_{\{1\}}}} \backslash_{w_{\{1\}}} \mathbb{R}^{\#_{w_{\{1\}}}}_{\mathbf{fin}}.$$

Obviously by (3.3.15),Eq.(3.3.16) and Eq.(3.3.17) one obtains a contradiction.

# Apendix 1.Paraconsistent Nonstandard Arithmetic

**Designations 1.1**. We will be write for short:

(i) $x =_{w_{\{0\}}} y$ instead $\left[ (x =_s y) \vee \left[ (x =_w y) \wedge \neg_s (x =_{w_{[1]}} y) \right] \right]$,

(ii) $x =_{w_{\{1\}}} y$ instead $\left[ (x =_s y) \vee (x =_{w_{\{0\}}} y) \vee \left[ (x =_{w_{[1]}} y) \wedge \neg_s (x =_{w_{[2]}} y) \right] \right]$,

(iii) $x =_{w_{\{n\}}} y$ instead $\left[ (x =_s y) \vee (x =_{w_{\{0\}}} y) \vee \ldots \vee \left[ (x =_{w_{[n]}} y) \wedge \neg_s (x =_{w_{[n+1]}} y) \right] \right]$, $n = 1, 2, \ldots$ .

**Designations 1.2**. We will be write for short:

(i) $x \neq^s_{w_{\{0\}}} y$ instead $\neg_s (x =_{w_{\{0\}}} y)$, i.e. instead

$\neg_s \left[ \left[ (x =_s y) \vee \left[ (x =_w y) \wedge \neg_s (x =_{w_{[1]}} y) \right] \right] \right]$,

(ii) $x =_{w_{\{1\}}} y$ instead $\left[ (x =_s y) \vee (x =_{w_{\{0\}}} y) \vee \left[ (x =_{w_{[1]}} y) \wedge \neg_s (x =_{w_{[2]}} y) \right] \right]$,

(iii) $x =_{w_{\{n\}}} y$ instead $\left[ (x =_s y) \vee (x =_{w_{\{0\}}} y) \vee \ldots \vee \left[ (x =_{w_{[n]}} y) \wedge \neg_s (x =_{w_{[n+1]}} y) \right] \right]$, $n = 1, 2, \ldots$ .

**The Theory $\mathbf{PA}_s$**

**The Theory $\mathbf{PA}_{w_{\{0\}}}$**

Let $\mathbb{N}_{w_{\{0\}}}$ be a set containing an $w_{\{0\}}$-element $0_{w_{\{0\}}}$, and let $\mathbf{S}_{w_{\{0\}}} : \mathbb{N}_{w_{\{0\}}} \to \mathbb{N}_{w_{\{0\}}}$ be a $w_{\{0\}}$-function satisfying the following postulates:

$\mathbf{PA}_{w_{\{0\}}} 0$ : $0_{w_{\{0\}}} =_{w_{\{0\}}} 0_{w_{\{0\}}}$,

$\mathbf{PA}_{w_{\{0\}}} 1$ : $\mathbf{S}_{w_{\{0\}}}(x) \neq^s_{w_{\{0\}}} 0_{w_{\{0\}}}$, for all $x \in_{w_{\{0\}}} \mathbb{N}_{w_{\{0\}}}$.

$\mathbf{PA}_{w_{\{0\}}} 2$ : $\forall x, y \in_{w_{\{0\}}} \mathbb{N}_{w_{\{0\}}} \{ [\mathbf{S}_{w_{\{0\}}}(x) =_{w_{\{0\}}} \mathbf{S}_{w_{\{0\}}}(y)] \Rightarrow_s x =_{w_{\{0\}}} y \}$,

$\mathbf{PA}_{w_{\{0\}}} 3$ : Let $A$ be any $w_{\{0\}}$-subset of $\mathbb{N}_{w_{\{0\}}}$ which contains $0_{w_{\{0\}}}$ and which is closed under $\mathbf{S}_{w_{\{0\}}}$ i.e. $\mathbf{S}_{w_{\{0\}}}(x) \in_{w_{\{0\}}} A$ for all $x \in_{w_{\{0\}}} A$. Then $A =_{w_{\{0\}}} \mathbb{N}_{w_{\{0\}}}$.

$\mathbf{PA}_{w_{\{0\}}} 4$. $\forall y_1 \ldots \forall y_k [A(0_{w_{\{0\}}}) \wedge \forall x (A(x) \Rightarrow_s A(\mathbf{S}_{w_{\{0\}}}(x)))]$

where $A$ is any formula whose free variables are among $x, y_1, y_k$.

**The Theory $\mathbf{PA}_{w_{\{1\}}}$**

Let $\mathbb{N}_{w_{\{1\}}}$ be a $w_{\{1\}}$-set containing an $w_{\{1\}}$-element $0_{w_{\{1\}}}$, and let $\mathbf{S}_{w_{\{1\}}} : \mathbb{N}_{w_{\{1\}}} \to \mathbb{N}_{w_{\{1\}}}$ be a

$w_{\{1\}}$-function satisfying the following postulates:

$\mathbf{PA}_{w_{\{1\}}}0$ : $0_{w_{\{0\}}} =_{w_{\{0\}}} 0_{w_{\{0\}}}$,

$\mathbf{PA}_{w_{\{1\}}}1$ : $\mathbf{S}_{w_{\{0\}}}(x) \neq^{s}_{w_{\{0\}}} 0_{w_{\{0\}}}$, for all $x \in_{w_{\{0\}}} \mathbb{N}^{\#}_{w_{\{1\}}}$.

$\mathbf{PA}_{w_{\{1\}}}2$ : If $\mathbf{S}_{w_{\{0\}}}(x) =_{w_{\{0\}}} \mathbf{S}_{w_{\{0\}}}(y)$ then $x =_{w_{\{0\}}} y$, for all $x,y \in_{w_{\{0\}}} \mathbb{N}_{w_{\{1\}}}$.

$\mathbf{PA}_{w_{\{1\}}}4$.

$\mathbf{PA}_{w_{\{1\}}}5$. There exists $w_{\{1\}}$-subset $\mathbb{N}_{w_{\{1\}}} \subset_{w_{\{1\}}} \mathbb{N}^{\#}_{w_{\{1\}}}$ such that the following statement holds:

any $w_{\{1\}}$-subset $X \subset_{w_{\{1\}}} \mathbb{N}_{w_{\{1\}}}$ has a strong $w_{\{1\}}$-complement $\mathbb{N}_{w_{\{1\}}} \backslash^{s}_{w_{\{1\}}} X$ in $\mathbb{N}_{w_{\{1\}}}$.

**Definition 1.1**.The condition that $X$ has a strong $<_{w_{\{1\}}}$ -least element reads

$$\exists x\left(x \in_{w_{\{1\}}} X\right)\left[\forall y \in_{w_{\{1\}}} X \neg_{\mathbf{s}}\left(y <_{w_{\{1\}}} x\right)\right]. \qquad (1.2)$$

**Definition 1.2**.

**Remark 1.1**.

**Theorem 1.1**. $\mathbb{N}_{w_{\{1\}}}$ is a strong well-$w_{\{1\}}$-ordered $w_{\{1\}}$-set.

**Proof**.We will prove by using strong (or complete) induction.

Let $X$ be a $w_{\{1\}}$-nonempty $w_{\{1\}}$-subset of $\mathbb{N}_{w_{\{1\}}}$. Suppose $X$ does not have a $<_{w_{\{1\}}}$ -least element. Then consider the set $\mathbb{N}_{w_{\{1\}}} \backslash^{s}_{w_{\{1\}}} X$.

Case 1) $\mathbb{N}_{w_{\{1\}}} \backslash^{s}_{w_{\{1\}}} X =_{w_{\{1\}}} \varnothing_{w_{\{1\}}}$. Then $X =_{w_{\{1\}}} \mathbb{N}_{w_{\{1\}}}$ and so $0_{w_{\{1\}}}$ is a strong $<_{w_{\{1\}}}$ -least element. Contradiction.

Case 2) $\neg_{\mathbf{s}}\left(\mathbb{N}_{w_{\{1\}}} \backslash^{s}_{w_{\{1\}}} X =_{w_{\{1\}}} \varnothing_{w_{\{1\}}}\right)$. There exists an $n \in_{w_{\{1\}}} \mathbb{N}_{w_{\{1\}}} \backslash^{s}_{w_{\{1\}}} X$ such that for all $k <_{w_{\{1\}}} n$; $k \in_{w_{\{1\}}} \mathbb{N}_{w_{\{1\}}} \backslash^{s}_{w_{\{1\}}} X$.

(Note that $n$ necessarily exists because $0_{w_{\{1\}}} \in_{w_{\{1\}}} \mathbb{N}_{w_{\{1\}}} \backslash^{s}_{w_{\{1\}}} X$, else $0_{w_{\{1\}}} \in_{w_{\{1\}}} X$ and would

be a $<_{w_{\{1\}}}$ - least $w_{\{1\}}$-element of $X$.)

Since we have supposed that N  X does not have a least element, thus n =2 X.

Using strong induction, we see that for all k < n; k 2 N  X and n 2 N  X. We can conclude

n 2 N  X for all n 2 N. Thus N  X = N implies X = ?.

This is a contradiction to X being a nonempty subset of N.

**Notation 1.1**. We often abbreviate for short $x \prec_{w_{\{1\}}} y$ instead

$$\neg_{w}\left(y <_{w_{\{1\}}} x\right). \qquad (1.1)$$

**Definition 1.1**.The condition that $X$ has a weak $<_{w_{\{1\}}}$ -least $w_{\{1\}}$-element reads

$$\exists x\left(x \in_{w_{\{1\}}} X\right)\left[\forall y \in_{w_{\{1\}}} X \neg_{w}\left(y <_{w_{\{1\}}} x\right)\right]. \qquad (1.2)$$

or

$$\exists x\left(x \in_{w_{\{1\}}} X\right)\left[\forall y \in_{w_{\{1\}}} X\left(x \prec_{w_{\{1\}}} y\right)\right]. \quad (1.3)$$

**Definition 1.1**.Assume that the condition that $X$ has a weak $<_{w_{\{1\}}}$ -least element is satisfied and let $x$ be a weak $<_{w_{\{1\}}}$ -least $w_{\{1\}}$-element of the $X$. We will say that $x$ is inconsistent if the following statement is true

$\mathbf{PA}_{w_{\{0\}}}3.$ **The a weak $w_{\{1\}}$-well-ordering principle**:

(i) every non-empty $w_{\{1\}}$-set of natural numbers $w_{\{1\}}$-contains a weak $<_{w_{\{1\}}}$ -least element

or in the following equivalent form

(ii) every non-empty $w_{\{1\}}$-set of natural numbers contains a $\prec_{w_{\{1\}}}$ -least $w_{\{1\}}$-element

**Remark 1.1**.We remind that

$$y <_{w_{\{1\}}} x \Leftrightarrow \neg_{\mathbf{s}}\left(x <_{w_{\{1\}}} y\right). \quad (1.4)$$

**Theorem 1.1**.

Proof. Assume that (i) and (ii) are both true statements.

Let $X_P$ be the $w_{\{1\}}$-set of all natural numbers for which $P(y)$ is false, i.e.

$$\forall y \in_{w_{\{1\}}} X_P \Leftrightarrow \neg_{\mathbf{s}} P(y). \quad (1.5)$$

If $X_P$ is $w_{\{1\}}$-empty set then we are done, so assume that $X_P$ is not $w_{\{1\}}$-empty.Then,by the weak well $w_{\{1\}}$-ordering principle, $X_P$ has a weak $<_{w_{\{1\}}}$ -least $w_{\{1\}}$-member let's say $x$, i.e.

$$\exists x\left(x \in_{w_{\{1\}}} X_P\right)\left[\forall y \in_{w_{\{1\}}} X_P\left(x \prec_{w_{\{1\}}} y\right)\right]. \quad (1.6)$$

Since $x$ is the weak $<_{w_{\{1\}}}$ -least $w_{\{1\}}$-member of $X_P$ it follows that $P\left(x -_{w_{\{1\}}} 1_{w_{\{1\}}}\right)$ is true. But this means, by (ii) above, that $P(x)$ is true. We have a contradiction and so our assumption that $\neg_{\mathbf{s}}\left(X_P =_{w_{\{1\}}} \varnothing_{w_{\{1\}}}\right)$ must be wrong.